\documentclass[OPT,biber]{nowfnt_arxiv} 

\usepackage[utf8]{inputenc}

\newcommand{\ie}{\emph{i.e.}}

\usepackage{caption}
\usepackage{subcaption}
\usepackage{typearea}

\def\jobname{ms} 
\usepackage{psfrag}
\usepackage{auto-pst-pdf} 
\usepackage{dsfont}

\usepackage{subfiles} 
\graphicspath{{figures/}{../figures/}} 
\usepackage{amssymb,amsmath}
\usepackage{algorithmic,algorithm}

\usepackage{pifont} 
\usepackage{booktabs} 
\usepackage{tikz}
\usepackage{appendix}
\usepackage{pgfkeys}
\usetikzlibrary{shapes,arrows}
\usepackage{pgfplots}
\pgfplotsset{compat=1.17}
\newcommand{\togglecodeurl}{\url{https://github.com/AdrienTaylor/AccelerationMonograph}}

\title{Acceleration Methods}

\maintitleauthorlist{
Alexandre d'Aspremont \\
CNRS \& Ecole Normale Sup\'erieure, Paris\\
aspremon@ens.fr
\and
Damien Scieur \\
Samsung SAIT AI Lab \& Mila, Montreal\\
damien.scieur@gmail.com
\and
Adrien Taylor \\
INRIA \& Ecole Normale Sup\'erieure, Paris\\
adrien.taylor@inria.fr
}

\issuesetup
{%
 copyrightowner={A.~Heezemans and M.~Casey},
 volume        = xx,
 issue         = xx,
 pubyear       = 2020,
 isbn          = xxx-x-xxxxx-xxx-x,
 eisbn         = xxx-x-xxxxx-xxx-x,
 doi           = 10.1561/XXXXXXXXX,
 firstpage     = 1, 
 lastpage      = 18
 }

\author[1]{d'Aspremont,Alexandre}
\author[2]{Scieur,Damien}
\author[3]{Taylor,Adrien}

\affil[1]{CNRS \& Ecole Normale Sup\'erieure, Paris; aspremon@ens.fr}
\affil[2]{Samsung SAIT AI Lab \& Mila, Montreal; damien.scieur@gmail.com}
\affil[3]{INRIA \& Ecole Normale Sup\'erieure, Paris; adrien.taylor@inria.fr}

\DeclareCiteCommand{\citeyear}
    {}
    {\bibhyperref{\printdate}}
    {\multicitedelim}
    {}

\newcommand{\BEAS}{\begin{eqnarray*}}
\newcommand{\EEAS}{\end{eqnarray*}}
\newcommand{\BEA}{\begin{eqnarray}}
\newcommand{\EEA}{\end{eqnarray}}
\newcommand{\BEQ}{\begin{equation}}
\newcommand{\EEQ}{\end{equation}}
\newcommand{\BIT}{\begin{itemize}}
\newcommand{\EIT}{\end{itemize}}
\newcommand{\BNUM}{\begin{enumerate}}
\newcommand{\ENUM}{\end{enumerate}}

\newcommand{\BA}{\begin{array}}
\newcommand{\EA}{\end{array}}

\newcommand{\ones}{\mathbf 1}

\newcommand{\reals}{{\mathbb R}}

\newcommand{\symm}{{\mbox{\bf S}}}
\newcommand{\Span}{\mbox{\textrm{span}}}

\newcommand{\Null}{{\textrm{Null}}}

\newcommand{\idm}{\mathbf{I}}
\newcommand{\HH}{\textbf{H}}
\newcommand{\GG}{\textbf{G}}

\newcommand{\EE}{\textbf{E}}

\newcommand{\Co}{{\mathop {\bf Co}}}

\newcommand{\Ext}{{\mathop {\bf Ext}}}

\newcommand{\argmin}{\mathop{\rm argmin}}

\newcommand{\prox}{{\mathrm{prox}}}

\newcommand{\dom}{\mathop{\bf dom}}

\newcommand{\cond}{q}
\newcommand{\defeq}{\triangleq}

\PassOptionsToPackage{hyphens}{url}\usepackage{hyperref}

\begin{document}

\maketitle

\clearpage 

\thispagestyle{empty}
~
\vfill

\begin{center}
    \textbf{Typos and Errata} \\
\end{center}

~\\
~\\
Please do not hesitate to contact us for any typo or error that you can find.
A list of identified typos and errata can be found at\\ \href{https://accelerationmethods.github.io/AccelerationMethodsWebsite/}{https://accelerationmethods.github.io/AccelerationMethodsWebsite/}. 
    
\vfill
~

\clearpage 

\tableofcontents
\clearpage 

\begin{abstract}

This monograph covers some recent advances in a range of acceleration techniques frequently used in convex optimization. We first use quadratic optimization problems to introduce two key families of methods, namely momentum and nested optimization schemes. They coincide in the quadratic case to form the \textit{Chebyshev method}.

We discuss momentum methods in detail, starting with the seminal work of~\citet{Nest83} and structure convergence proofs using a few master templates, such as that for \emph{optimized gradient methods}, which provide the key benefit of showing how momentum methods optimize convergence guarantees. We further cover proximal acceleration, at the heart of the \emph{Catalyst} and \emph{Accelerated Hybrid Proximal Extragradient} frameworks, using similar algorithmic patterns.

Common acceleration techniques rely directly on the knowledge of some of the regularity parameters in the problem at hand. We conclude by discussing \emph{restart} schemes, a set of simple techniques for reaching nearly optimal convergence rates while adapting to unobserved regularity parameters. 
\end{abstract}

\chapter{Introduction}\label{c-intro}

Optimization methods are a core component of the modern numerical toolkit. In many cases, iterative algorithms for solving convex optimization problems have reached a level of efficiency and reliability comparable to that of advanced linear algebra routines. This is largely true for medium scale-problems where interior point methods reign supreme, but less so for large-scale problems where the complexity of first-order methods is not as well understood and efficiency remains a concern. 

The situation has improved markedly in recent years, driven in particular by the emergence of a number of applications in statistics, machine learning, and signal processing. Building on Nesterov's path-breaking algorithm from the 80's, several accelerated methods and numerical schemes have been developed that both improve the efficiency of optimization algorithms and refine their complexity bounds. Our objective in this monograph is to cover these recent developments using a few master templates. 

The methods described in this manuscript can be arranged in roughly two categories. The first, stemming from the work of \citet{Nest83}, produces variants of the gradient method with accelerated worst-case convergence rates that are provably optimal under classical regularity assumptions. The second uses outer iteration (a.k.a. nested) schemes to speed up convergence. In this second setting, accelerated schemes run both an inner loop and an outer loop, with the inner iterations being solved by classical optimization methods, and the outer loop containing the acceleration mechanism.

\paragraph{Direct acceleration techniques.} Ever since the original algorithm by \citet{Nest83}, the acceleration phenomenon was regarded as somewhat of a mystery. While accelerated gradient methods can be seen as iteratively building a model for the function and using it to guide gradient computations, the argument is essentially algebraic and is simply an effective exploitation of regularity assumptions. This approach of collecting inequalities induced by regularity assumptions and cleverly chaining them to prove convergence was also used in e.g.,~\citep{Beck09}, to produce an optimal proximal gradient method. There too, however, the proof yielded little evidence as to why the method is actually faster.

Fortunately, we are now better equipped to push the proof mechanisms much further. Recent advances in the programmatic design of optimization algorithms allow us to design and analyze algorithms by following a more principled approach. In particular, the \emph{performance estimation approach}, pioneered by~\citet{Dror14}, can be used to design optimal methods from scratch, selecting algorithmic parameters to optimize worst-case performance guarantees~\citep{Dror14,kim2016optimized}. Primal dual optimality conditions on the design problem then provide a blueprint for the accelerated algorithm structure and for its convergence proof.

Using this framework, acceleration is no longer a mystery: it is the main objective in the design of the algorithm. We recover the usual ``soup of regularity inequalities'' that forms the template of classical convergence proofs, but the optimality conditions of the design problem explicitly produce a method that optimizes the convergence guarantee. In this monograph, we cover accelerated first-order methods using  this systematic template and describe a number of convergence proofs for classical variants of the accelerated gradient method, such as those of Nesterov (\citeyear{Nest83,Nest03a}), \citet{Beck09,tseng2008accelerated} as well as more recent ones~\citep{kim2016optimized}.

\paragraph{Nested acceleration schemes.} The second category of acceleration techniques that we cover in this monograph is composed of outer iteration schemes, in which classical optimization algorithms are used as a black-box in the inner loop and acceleration is produced by an argument in the outer loop. We describe three acceleration results of this type.

The first scheme is based on nonlinear acceleration techniques. Based on arguments dating back to \citep{Aitk27,Wynn56,Ande87}, these techniques use a weighted average of iterates to extrapolate a better candidate solution than the last iterate. We begin by describing the Chebyshev method for solving quadratic problems, which interestingly qualifies both as a gradient method and as an outer iteration scheme. It takes its name from the use of Chebyshev polynomial coefficients to approximately minimize the gradient at the extrapolated solution. The argument can be extended to non-quadratic optimization problems provided the extrapolation procedure is regularized.

The second scheme, due to \citep{guler1992new,monteiro2013accelerated,lin2015universal} relies on a conceptual accelerated proximal point algorithm, and uses classical iterative methods to approximate the proximal point in an inner loop. In particular, this framework produces accelerated gradient methods (in the same sense as Nesterov's acceleration) when the approximate proximal points are computed using linearly converging gradient-based optimization methods, taking advantage of the fact that the inner problems are always strongly convex.

Finally, we describe restart schemes. These techniques exploit regularity properties called H\"olderian error bounds, which extend strong convexity properties near the optimum and hold almost generically, to improve the convergence rates of most first-order methods. The parameters of the H\"olderian error bounds are usually unknown, but the restart schemes are robust: that is, they are adaptive to the H\"olderian parameters and their empirical performance is excellent on problems with reasonable precision targets.

\paragraph{Content and organization.}
We present a few convergence acceleration techniques that are particularly relevant in the context of (first-order) convex optimization. Our summary includes our own points of view on the topic and is focused on techniques that have received substantial attention since the early 2000's, although some of the underlying ideas are much older. We do not pretend to be exhaustive, and we are aware that valuable references might not appear below.

The sections can be read nearly independently. However, we believe the insights of some sections can benefit the understanding of others. In particular, Chebyshev acceleration (Section~\ref{c-Cheb}) and nonlinear acceleration (Section~\ref{c-RNA}) are clearly complementary readings. Similarly, Chebyshev acceleration (Section~\ref{c-Cheb}) and Nesterov acceleration (Section~\ref{c-Nest}), Nesterov acceleration (Section~\ref{c-Nest}) and proximal acceleration (Section~\ref{c-prox}), as well as Nesterov acceleration (Section~\ref{c-Nest}) and restart schemes (Section~\ref{c-restart}) certainly belong together.

\paragraph{Prerequisites and complementary readings.} This monograph is not meant to be a general-purpose manuscript on convex optimization, for which we refer the reader to the now classical references~\citep{boyd2004convexopt,bonnans2006numerical,nocedal2006numerical}. Other directly related references are provided in the text.

We assume the reader to have a working knowledge of base linear algebra and convex analysis (such as of subdifferentials), as we do not detail the corresponding technical details while building on them. Classical references on the latter include~\citep{Rock70,rockafellar2009variational,hiriart2013convex}.  


\chapter{Chebyshev Acceleration}\label{c-Cheb}

While ``Chebyshev polynomials are everywhere dense in numerical analysis,'' we would like to argue here that Chebyshev polynomials also provide one of the most direct and intuitive explanations for acceleration arguments in first-order methods. That is, one can form linear combinations of past gradients for optimizing a worst-case guarantee on the distance to an optimal solution. In quadratic optimization, these linear combinations emerge from a Chebyshev minimization problem, whose solution can also be computed iteratively, thereby yielding an algorithm called the {\em Chebyshev method} \citep{Nemic84}. The Chebyshev method traces its roots to at least \citet{Flan50}, who credit Tuckey and Grosch. Its recurrence matches asymptotically the one of the heavy-ball method and is detailed below.

\enlargethispage{\baselineskip}
\section{Introduction}\label{s:cheb-intro}
In this section, we demonstrate basic acceleration results on quadratic minimization problems. In such problems, optimal points are the solutions of a linear system, and the basic gradient method can be seen as a simple iterative solver for this linear system. In this context, acceleration methods can be obtained using a classical argument involving Chebyshev polynomials.

Analyzing this simple scenario is useful in two ways. First, recursive formulations of the Chebyshev argument yield a basic algorithmic template for designing accelerated methods and provide first approach to their structures, such as the presence of a momentum term. Second, the arguments are robust to perturbations of the quadratic function $f$ and hence apply in more generic contexts. This property enables acceleration in a wider range of applications, which we cover later in the Section~\ref{c-RNA} and Section~\ref{c-Nest}.

For now, consider the following unconstrained quadratic minimization problem 
\BEQ\label{eq:quad-prob}
\mbox{minimize} ~ \left\{f(x) \triangleq \frac{1}{2} \langle x;\HH x\rangle - \langle b;x\rangle \right\}
\EEQ
in the variable $x\in\reals^d$, where $\HH\in\symm_d$ (the set of symmetric matrices of size $d\times d$) is the Hessian of $f$. We further assume that $f$ is both smooth and strongly convex, i.e., that there exist some $L>\mu>0$ such that $\mu \idm \preceq \HH \preceq L \idm$. The reasoning of this section readily extends to the case where $\mu$ is the smallest {\em nonzero} eigenvalue of $\HH$. We start by analyzing the convergence of the fixed step gradient method (Algorithm~\ref{alg:cheb-grad}) for solving~\eqref{eq:quad-prob}.

\begin{algorithm}[h]
  \caption{Gradient method}
  \label{alg:cheb-grad}
  \begin{algorithmic}[1]
    \REQUIRE
      A differentiable convex function $f$, initial point $x_0$, step size $\gamma >0$, budget $N$.
    \FOR{$k=1,\ldots,N$}
      \STATE $x_{k}= x_{k-1} - \gamma \nabla f(x_{k-1}) $
    \ENDFOR
    \ENSURE Approximate solution $x_{N}$.
  \end{algorithmic}
\end{algorithm}

For problem~\eqref{eq:quad-prob}, the iteration reads
\[
    x_{k+1} = (\idm - \gamma \HH) x_{k} + \gamma b,
\]
and calling $x_\star$ the optimum of problem~\eqref{eq:quad-prob} (satisfying  $\HH x_\star=b$) yields
\BEQ\label{eq:error-bnd}
    x_{k+1}-x_\star = (\idm - \gamma \HH) (x_{k} - x_\star).
\EEQ
This means the iterates of gradient descent $x_{k}-x_\star$ can be computed from $x_{0}-x_\star$ via $x_{k}-x_\star= P_k^{\text{Grad}}(\HH)(x_0-x_\star)$ using the matrix polynomial 
\BEQ
    P_k^{\text{Grad}}(\HH) = (\idm-\gamma \HH)^k. \label{eq:polynomial_gradient_descent}
\EEQ
Suppose we set the step size $\gamma$ to ensure
\[
\|\idm - \gamma \HH\|_2 < 1,
\]
where $\|\cdot\|_2$ stands for the operator $\ell_2$ norm. Then,~\eqref{eq:error-bnd} controls the convergence with
\BEQ\label{eq:grad-bnd}
\|x_{k}-x_\star\|_2 \leq  \|\idm - \gamma \HH\|_2^k ~\|x_{0} - x_\star\|_2, \quad \mbox{for $k\geq 0$}.
\EEQ
Because the matrix $\HH$ is symmetric and hence diagonalizable in an orthogonal basis, given $\gamma>0$, we obtain
\BEAS
\|\idm - \gamma \HH\|_2 & \leq & \max_{\mu \idm \preceq H \preceq L \idm} \| \idm - \gamma \HH \|_2 \\
& \leq & \max_{\mu \leq \lambda \leq L} \left|1-\gamma \lambda\right| \\
& \leq & \max_{\mu \leq \lambda \leq L} \max \left\{\gamma \lambda-1\; ; \; 1-\gamma \lambda \right\}\\
& \leq & \max \left\{\gamma L-1\; ; \; 1-\gamma \mu \right\}.
\EEAS
To get the best possible worst-case convergence rate, we now minimize this quantity in $\gamma$ by solving
\BEQ\label{eq:rate_convergence_gradient_descent}
\min_{\gamma} ~\max\left\{\gamma L-1, 1-\gamma \mu \right\} = \frac{L-\mu}{L+\mu}.
\EEQ
The optimal step size is obtained when both terms in the max are equal, reaching:
\BEQ\label{eq:step}
    \gamma=\frac{2}{L+\mu}.
\EEQ
Denoting by $\kappa\triangleq \tfrac{L}{\mu} \geq 1$ the {\em condition number} of the function $f$, the bound in~\eqref{eq:grad-bnd} finally becomes
\BEQ\label{eq:error-bnd-kappa}
\|x_{k}-x_\star\|_2 \leq  \left(\frac{\kappa-1}{\kappa+1}\right)^k \|x_{0} - x_\star\|_2, \quad \mbox{for $k\geq 0$},
\EEQ
which is a worst-case guarantee for the gradient method when minimizing smooth strongly convex quadratic functions.

\section{Optimal Methods and Minimax Polynomials}

\begin{figure}
    \centering
    \includegraphics[width=0.6\textwidth]{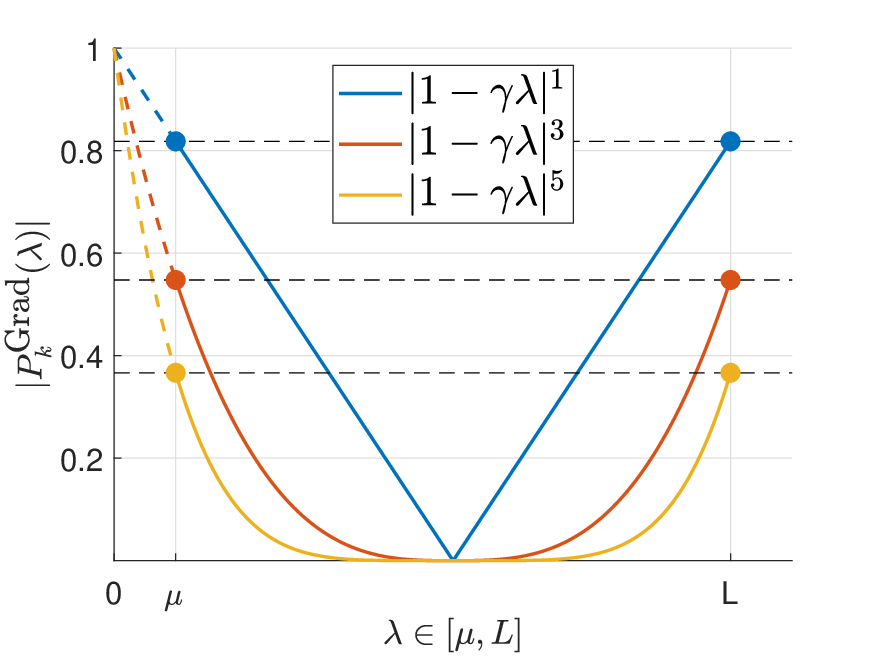}
    \caption{We plot $|P_k^{\text{Grad}}(\lambda)|$ (for the optimal $\gamma$ in~\eqref{eq:step}) for $k\in\{1,3,5\}$, $\mu = 1$, $L=10$. Note that the polynomials satisfy $P_k^{\text{Grad}}(0)=1$. The rate is equal to the largest value of $|P_k^{\text{Grad}}(\lambda)|$ on the interval, which is achieved at the boundaries (where $\lambda$ is either equal to $\mu$ or to $L$).}
    \label{fig:gradient_polynomial}
\end{figure}

In Equation~\eqref{eq:grad-bnd} above, we saw that the worst-case convergence rate of the gradient method on quadratic functions can be controlled by the spectral norm of a matrix polynomial. Figure~\ref{fig:gradient_polynomial} plots the polynomial $P_k^{\text{Grad}}$ for several degrees $k$. We can extend this reasoning further to produce methods with accelerated worst-case convergence guarantees.

\subsection{First-Order Methods and Matrix Polynomials}
The bounds derived above for gradient descent can be extended to a broader class of first-order methods for quadratic optimization. We consider first-order algorithms in which each iterate belongs to the span of previous gradients, i.e.
\begin{equation}\label{eq:first-order-cheby}
    x_{k+1} \in x_0 + \Span\left\{ \nabla f(x_0),\, \nabla f(x_1),\,\ldots,\, \nabla f(x_{k}) \right\},
\end{equation}
and show that the iterates can be written using matrix polynomials as in~\eqref{eq:polynomial_gradient_descent} above.

\newpage\begin{proposition}\label{prop:first_order_to_polynomial}
  Let $x_0\in\mathbb{R}^d$ and $f$ be a quadratic function defined as in~\eqref{eq:quad-prob} with $\mu \idm \preceq \HH \preceq L \idm$ for some $L>\mu>0$. The sequence $\{x_k\}_{k=0,1,\ldots}$ satisfies
    \BEQ\label{eq:span-fo}
        x_{k+1} \in x_0 + \Span\left\{ \nabla f(x_0),\, \nabla f(x_1),\,\ldots,\, \nabla f(x_{k}) \right\},
    \EEQ
    for all $k=0,1,\ldots$, if and only if the errors $\{x_k-x_\star\}_{k=0,1,\ldots}$ can be written as
    \BEQ\label{eq:pol-fo}
        x_{k}-x_\star = P_{k}(\HH)(x_0-x_\star),
    \EEQ
    for all $k=0,1\ldots$, for some sequence of polynomials $\{P_{k}\}_{k=0,1,\ldots}$ with $P_k$ of degree at most ${k}$ and $P_{k}(0)=1$.
\end{proposition}
\enlargethispage{\baselineskip}
\begin{proof}
Since $\nabla f(x)$ is the gradient of a quadratic function, it reads
\[
    \nabla f(x) = \HH x-b = \HH(x-x_\star)
\]
for any $x_\star$  satisfying $\HH x_\star = b$, where $\HH$ is symmetric. We have
\begin{align*}
    x_0-x_\star &= 1 \cdot (x_0-x_\star) \\&= P_0(\HH ) (x_0-x_\star).
\end{align*}
We now show recursively that $x_k-x_\star = P_k(\HH)(x_0-x_\star)$, where $P_k$ is a residual polynomial of degree at most $k$. Our assumption about the iterates~\eqref{eq:span-fo} implies that, for some sequence of coefficients $ \{\alpha_i^{(k+1)}\}_{i=0,\ldots,k}$,
\begin{align*}
    x_{k+1} - x_\star & = x_0-x_\star + \sum_{i=0}^{k}\alpha_i^{(k+1)} \nabla f(x_i).
\end{align*}
Assuming recursively that~\eqref{eq:pol-fo} holds for all indices $i \leq k$,
\begin{align*}
    x_{k+1} - x_\star & = x_0-x_\star + \sum_{i=0}^{k}\alpha_i^{(k+1)}\HH P_i(\HH ) (x_0-x_\star)\\
    & = \left(\idm+\HH  \sum_{i=0}^{k}\alpha_i^{(k+1)} P_i(\HH )\right) (x_0-x_\star).
\end{align*}
Then, by writing $P_{k+1}(x) = 1+x\sum_{i=0}^{k} \alpha_i^{(k+1)} P_i(x)$, we have 
\[
    x_{k+1}-x_\star = P_{k+1}(\HH )(x_0-x_\star)
\]
with $P_{k+1}(0) = 1$ and $\deg(P_{k+1}) \leq k+1$. Since the proof is a sequence of equalities, the equivalence readily follows.
\end{proof}

Given a class $\mathcal{M}$ of problem matrices $\HH$, Proposition~\ref{prop:first_order_to_polynomial} provides a way to design algorithms. Indeed, we can extract a first-order method from a sequence of polynomials $\{P_k\}_{k=0,\ldots, N}$. We can therefore use tools from approximation theory to find optimal polynomials and extract corresponding methods from them. Given a matrix class $\mathcal{M}$, this involves minimizing the worst-case convergence bound over $\HH \in \mathcal{M}$ by solving
\begin{equation}\label{eq:optimal_poly_pb_class}
    P_k^* = \argmin_{\substack{P \in\mathcal{P}_k,\\ P(0)=1}} ~\max_{\HH \in \mathcal{M}} \|P(\HH )\|_2
\end{equation}
where ${P}_k$ is the set of polynomials of degree at most $k$. The polynomial $P^*_k$ is an optimal polynomial for $\mathcal{M}$ and yields a (worst-case) optimal algorithm for the class~$\mathcal{M}$. In terms of the notation used in the proof of Proposition~\ref{prop:first_order_to_polynomial}, we are by construction looking for coefficients $\{\alpha_{i}^{(j)}\}_i$ depending only on the problem class $\mathcal{M}$, but \textit{not} on a specific instance of $\HH$. 

\section{The Chebyshev Method}\label{s:cheb-cheb}

In the case where $\mathcal{M}$ is the set of positive definite matrices with a bounded spectrum, namely
\[
    \mathcal{M} = \{ \HH  \in\symm_d: 0 \prec \mu \idm \preceq \HH  \preceq L \idm \},
\]
the optimal polynomial can be found by solving
\BEQ\label{eq:problem_chebyshev}
    P_k^* = \argmin_{\substack{P \in\mathcal{P}_k,\\ P(0)=1}} ~ \max_{\lambda \in [\mu,\, L]} |P(\lambda)| 
\EEQ
Polynomials that solve~\eqref{eq:problem_chebyshev} are derived from \textit{Chebyshev polynomials of the first kind} in approximation theory and can be formed explicitly to produce an optimal algorithm called the {\em Chebyshev method}. This section describes this method and provides its corresponding worst-case convergence guarantees.

\subsection{Shifted Chebyshev Polynomials}

We now explicitly introduce the Chebyshev polynomials. A more complete treatment of these polynomials is available in, e.g.,~\citet{mason2002chebyshev}. Chebyshev polynomials of the first kind are defined recursively as follows
\begin{equation}\label{eq:cheby_rec}
\begin{aligned}
    \mathcal{T}_0(x) & = 1, \\
    \mathcal{T}_1(x) & = x, \\
    \mathcal{T}_k(x) & = 2x \mathcal{T}_{k-1}(x) - \mathcal{T}_{k-2}(x), \quad\mbox{for $k \geq 2$}.
\end{aligned}
\end{equation}

There exists a compact explicit solution for Chebyshev polynomials that involves trigonometric functions:
\BEQ\label{eq:cheby_explicit_trigono}
    \mathcal{T}_k(x) = 
    \begin{cases}
        \cos(k\,\text{acos}(x)) & x\in[-1,1],\\
        \text{cosh}(k\,\text{acosh}(x)) & x>1,\\
        (-1)^k\text{cosh}(k\,\text{acosh}(-x)) & x<1.\\
    \end{cases}
\EEQ
It is possible to show that Chebyshev polynomials satisfy the minimax property
\[
    \frac{1}{2^{k-1}} \mathcal{T}_k = \argmin_{\substack{\deg(P)\leq k \\ \text{$P$ is monic}}} \max_{[-1,1]} |P(x)|,
\]
where a monic polynomial is a polynomial whose coefficient associated with the highest power is equal to one. From this minimax definition, that defines the minimal polynomial over $[-1,1]$, we can transform it by \textit{shifting} it to the interval $[\mu,L]$, then \textit{rescaling} it to obtain a polynomial such that $P(0)=1$. More precisely, using a simple linear mapping  from $[\mu,\, L]$ to $[-1,\,1]$,
\[
    x\rightarrow t^{[\mu,L]}(x) = \frac{2x-(L+\mu)}{L-\mu},
\]
we obtain \textit{shifted Chebyshev polynomials}:
\BEQ
    C_k^{[\mu,L]}(x) = \frac{\mathcal{T}_k\big(t^{[\mu,L]}(x)\big)}{\mathcal{T}_k\big(t^{[\mu,L]}(0)\big)}. \label{eq:shifted_cheby_definition}
\EEQ
where we have enforced the normalization constraint $C^{[\mu,L]}_k(0) = 1$. Under these transformations, the shifted Chebyshev polynomials keep some minimax property, and can be shown to be solutions to~\eqref{eq:problem_chebyshev}.

\begin{figure}
    \centering
    \includegraphics[width=0.6\textwidth]{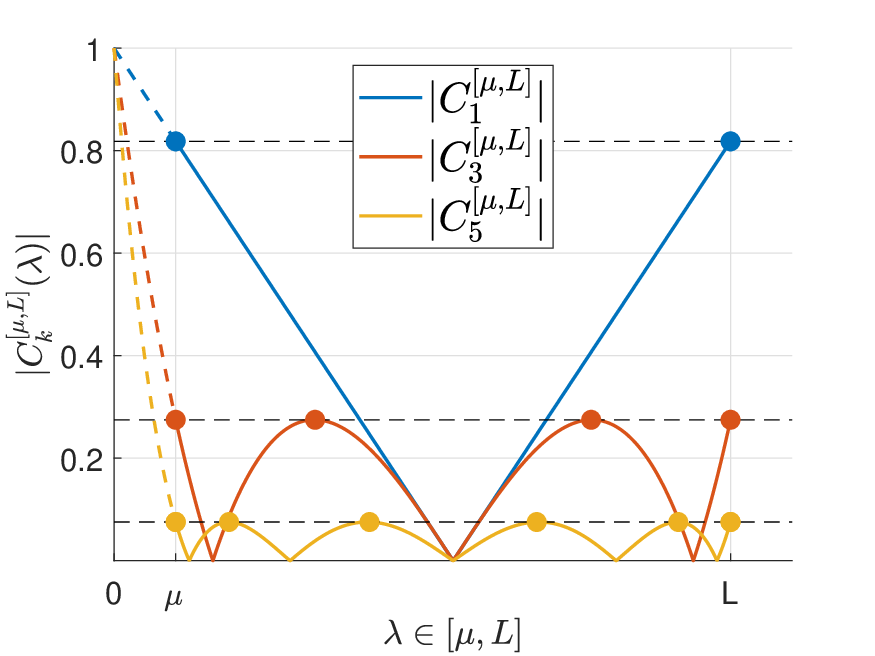}
    \caption{We plot the absolute value of $C_1^{[\mu,L]}(x)$,  $C_3^{[\mu,L]}(x)$ and $C_5^{[\mu,L]}(x)$ for $\lambda\in [\mu, L]$, where $\mu = 1$ and $L=10$. Note that the polynomials satisfy $C_k^{[\mu,L]}(0)=1$. The maximum value of the image of $[\mu,L]$ by $C^{[\mu,L]}_k$ decreases rapidly as $k$ grows, implying an accelerated rate of convergence. \label{fig:cheb}}
\end{figure}

More formally, \citet{golub1961chebyshev} characterize the $\ell_\infty$ optimality on the interval $[\mu,L]$ using an \textit{equi-oscillation} argument (see, e.g., \citep{suli2003introduction}); i.e., they show that the solution of~\eqref{eq:problem_chebyshev} is $P_k^* = C^{[\mu,L]}_k$. The equioscillation property of the shifted Chebyshev polynomial is clearly visible on Figure \ref{fig:cheb}, where the polynomial hits its maximum value $k+1$ times on the interval $[\mu,L]$.

\subsection{Chebyshev Algorithm} 

The following recursion follows from~\eqref{eq:cheby_rec} together with~\eqref{eq:shifted_cheby_definition} and a few simplifications:
\BEA\label{eq:shift-cheb}
    C_0^{[\mu,L]}(x) & = & 1,\nonumber\\
    C_1^{[\mu,L]}(x) & = & 1 - \frac{2}{L+\mu} x,\\
    C_k^{[\mu,L]}(x) & = & \frac{2\delta_{k}}{L-\mu}\left(L+\mu - 2x\right)C_{k-1}^{[\mu,L]}(x)\nonumber\\
    && +\left(1-\frac{2\delta_{k}(L+\mu)}{L-\mu}\right)C_{k-2}^{[\mu,L]}(x), \quad \mbox{for $k\geq 2$,}\nonumber
\EEA
where $\delta_1 = \frac{L-\mu}{L+\mu}$ and 
\[
    \delta_{k} = -\frac{\mathcal{T}_{k-1}\big(t^{[\mu,L]}(0))}{\mathcal{T}_k\big(t^{[\mu,L]}(0))} = \frac{1}{2\frac{L+\mu}{L-\mu}-\delta_{k-1}}, \quad \mbox{for $k\geq 2$.}
\]
The sequence $\delta_k$ ensures $C_k^{[\mu,L]}(0)=1$. We present this recursion for simplicity, but one should note that it might present some numerical stability issues in practice. There exist numerically more stable first-order methods based on Chebyshev polynomials, see, e.g.,~\citep[Algorithm 1]{gutknecht2002chebyshev}. We plot  $C_1^{[\mu,L]}(x)$, $C_3^{[\mu,L]}(x)$ and $C_5^{[\mu,L]}(x)$ in Figure~\ref{fig:cheb} for illustration.

\paragraph{Computational details for the shifted Chebyshev.} Let us quickly detail how to arrive to~\eqref{eq:shift-cheb} using~\eqref{eq:cheby_rec} and~\eqref{eq:shifted_cheby_definition}. We first expand $\mathcal{T}_k\big(t^{[\mu,L]}(x)\big)$ using~\eqref{eq:cheby_rec},
\begin{align*}
    \mathcal{T}_k\big(t^{[\mu,L]}(x)\big) & = 2t^{[\mu,L]}(x) \mathcal{T}_{k-1}\big(t^{[\mu,L]}(x)\big) - \mathcal{T}_{k-2}\big(t^{[\mu,L]}(x)\big),\\
    & = \frac{2(2x-L-\mu)}{L-\mu}\mathcal{T}_{k-1}\big(t^{[\mu,L]}(x)\big) - \mathcal{T}_{k-2}\big(t^{[\mu,L]}(x)\big).
\end{align*}
Since $C_k^{[\mu,L]}(x) = \frac{\mathcal{T}_k\big(t^{[\mu,L]}(x)\big)}{\mathcal{T}_k\big(t^{[\mu,L]}(0)\big)}$, we substitute $\mathcal{T}_k\big(t^{[\mu,L]}(x)\big)$ by $\mathcal{T}_k\big(t^{[\mu,L]}(0)) \cdot C_k^{[\mu,L]}(x)$ in the equation above and obtain
\begin{align*}
    \mathcal{T}_k\big(t^{[\mu,L]}(0)) \cdot C_k^{[\mu,L]}(x) =& \frac{2(2x-L-\mu)}{L-\mu}\mathcal{T}_{k-1}\big(t^{[\mu,L]}(0)) \cdot C_{k-1}^{[\mu,L]}(x) \\
    &  - \mathcal{T}_{k-2}\big(t^{[\mu,L]}(0))\cdot C_{k-2}^{[\mu,L]}(x),\\
    C_k^{[\mu,L]}(x) =& \frac{2(2x-L-\mu)}{L-\mu}\frac{\mathcal{T}_{k-1}\big(t^{[\mu,L]}(0))}{\mathcal{T}_k\big(t^{[\mu,L]}(0))} \cdot C_{k-1}^{[\mu,L]}(x) \\ 
    &  - \frac{\mathcal{T}_{k-2}\big(t^{[\mu,L]}(0))}{\mathcal{T}_k\big(t^{[\mu,L]}(0))} \cdot C_{k-2}^{[\mu,L]}(x),\\
    C_k^{[\mu,L]}(x) =& -\frac{2(2x-L-\mu)}{L-\mu} \delta_k C_{k-1}^{[\mu,L]}(x) \\
    & - \delta_{k-1}\delta_k  C_{k-2}^{[\mu,L]}(x).
\end{align*}
For obtaining a simple recursion on $\delta_k$, note that $C_k^{[\mu,L]}(0)=1$ for all $k\geq 0$ by construction. It follows that
\[
    C_k^{[\mu,L]}(0) = \frac{2(L+\mu)}{L-\mu} \delta_k \underbrace{C_{k-1}^{[\mu,L]}(0)}_{=1} - \delta_{k-1}\delta_k  \underbrace{C_{k-2}^{[\mu,L]}(0)}_{=1} = 1,
\]
and thereby
\[
    \delta_{k-1}\delta_k  =  1-2\delta_k\frac{L+\mu}{L-\mu}  \quad \text{and} \quad \delta_k = \frac{1}{2\frac{L+\mu}{L-\mu}-\delta_{k-1}}.
\]

\subsection{Chebyshev and Polyak's Heavy-Ball Methods}\label{s:HeavyBall}
We now present the resulting algorithm, called \textit{Chebyshev semi-iterative method} \citep{Golu61}. We define iterates using $C_k^{[\mu,L]}(x)$ as follows,
\[
    x_{k}-x_\star = C_k^{[\mu,L]}(\HH )(x_0-x_\star).
\]
The recursion in~\eqref{eq:shift-cheb} then yields
\BEAS
    x_{k}-x_\star & = &  \frac{2\delta_{k}}{L-\mu}\left((L+\mu)\idm - 2\HH \right)(x_{k-1}-x_\star) \\
    & & +\left(1-2\delta_{k}\frac{L+\mu}{L-\mu}\right)(x_{k-2}-x_\star).
\EEAS
Since the gradient of the function reads
\[
    \nabla f(x_k) = \HH (x_{k}-x_\star),
\]
we can simplify away $x_\star$ to get the following recursion:
\[
    x_{k} = \frac{2\delta_{k}}{L-\mu}\left((L+\mu)x_{k-1} - 2\nabla f(x_{k-1})\right) +\left(1-2\delta_{k}\frac{L+\mu}{L-\mu}\right)x_{k-2},
\]
which describes iterates of the Chebyshev method. We summarize it as Algorithm~\ref{alg:cheb}. By construction, the Chebyshev method is a worst-case optimal first-order method for minimizing quadratics whose spectrum lies in $[\mu, L]$. Surprisingly, its iteration structure is simple and somewhat intuitive: it involves a gradient step with variable step size ${4\delta_{k}}/({L-\mu})$, combined with a variable momentum term. 

\begin{algorithm} 
  \caption{Chebyshev's method} \label{alg:cheb}
    \begin{algorithmic}[1]
  \REQUIRE An $L$-smooth $\mu$-strongly convex quadratic $f$, initial point $x_0$ and budget~$N$.
  \STATE Set $\delta_1 = \frac{L-\mu}{L+\mu},\;\; x_1 = x_0 - \frac{2}{L+\mu}\nabla f(x_0)$.
  \FOR{$k=2,\ldots,N$}
    \STATE Set $\delta_{k} = \frac{1}{2\frac{L+\mu}{L-\mu}-\delta_{k-1}},$
    \STATE $ x_{k} = x_{k-1} - \frac{4\delta_{k}}{L-\mu}\nabla f(x_{k-1}) +\left(1-2\delta_{k}\frac{L+\mu}{L-\mu}\right)(x_{k-2}-x_{k-1}).$
    \ENDFOR
    \ENSURE Approximate solution $x_{N}$.
    \end{algorithmic}
\end{algorithm}

Perhaps more surprisingly, the Chebyshev method has a \emph{stationary regime} that is even simpler. Indeed, when $k\rightarrow \infty$, the coefficients of the recursion from Algorithm~\ref{alg:cheb} converge to the ones of {\em Polyak's heavy-ball method},
\[
    x_{k} = x_{k-1} - \frac{4}{(\sqrt{L}+\sqrt{\mu})^2} \nabla f(x_{k-1}) + \frac{(\sqrt{L}-\sqrt{\mu})^2}{(\sqrt{L}+\sqrt{\mu})^2}(x_{k-1}-x_{k-2}).
\]
To see this, it suffices to compute the limit of $\delta_k$, written as $\delta_{\infty}$, by solving
\[
    \delta_{\infty} = \frac{1}{2\frac{L+\mu}{L-\mu}-\delta_{\infty}},
\]
reaching
\begin{equation}
    \delta_{\infty} = \frac{\sqrt{L}-\sqrt{\mu}}{\sqrt{L}+\sqrt{\mu}}. \label{eq:delta_infty}
\end{equation}
We obtain Polyak's heavy-ball method by replacing $\delta_k$ with $\delta_{\infty}$ in Algorithm~\ref{alg:cheb}.

\subsection{Worst-case Convergence Bounds}\label{s:cheby_wc}
The \textit{shifted} Chebyshev polynomials are solutions to~\eqref{eq:problem_chebyshev}. Therefore, using the same trick as for gradient descent, we can obtain the following worst-case bound for the Chebyshev method
\BEQ\label{eq:bound_using_cheb_poly}
    \|x_k-x_\star\|_2 \leq \|C_k^{[\mu,L]}(\HH )(x_0-x_\star)\|_2 \leq \|x_0-x_\star\|_2\max_{x\in[\mu,L]} |C_k^{[\mu,L]}(x)|. 
\EEQ
The maximum value is determined by evaluating the polynomial at one of the extremities of the interval~\citep[Chapter 2]{mason2002chebyshev} (see also Figure \ref{fig:cheb}), i.e, 
\[
    \max_{x\in[\mu,L]} |C_k^{[\mu,L]}(x)| = C_k^{[\mu,L]}(L). 
\]
Using \eqref{eq:shifted_cheby_definition} and \eqref{eq:cheby_explicit_trigono} successively,
\[
    |C_k^{[\mu,L]}(L)| = \frac{1}{|\mathcal{T}_k\big(t^{[\mu,L]}(0)\big)|} = \frac{1}{\text{cosh}\left( k \; \text{acosh}\left( \frac{L+\mu}{L-\mu} \right)\right)}.
\]
We obtain the worst-case convergence guarantee of the Chebyshev method, as stated in the following theorem.

\begin{theorem}\label{thm:cheb_rate_convergence}
  Let $x_0\in\mathbb{R}^d$ and $f$ be a quadratic function defined as in~\eqref{eq:quad-prob} with $\mu \idm \preceq \HH \preceq L \idm$ for some $L>\mu>0$. For any $N\in\mathbb{N}$, the iterates of the Chebyshev method (Algorithm~\ref{alg:cheb}) satisfy
    \BEQ \label{eq:rate_convergence_cheby}
        \| x_N-x_\star \|_2 \leq \frac{2}{\xi^N + \xi^{-N}}\|x_0-x_\star\|_2 \quad \text{where }\; \xi = \frac{\sqrt{\frac{L}{\mu}}+1}{\sqrt{\frac{L}{\mu}}-1}.
    \EEQ
\end{theorem}
\begin{proof}
    We bound \eqref{eq:bound_using_cheb_poly} as follows,
    \BEAS
        \|x_{N}-x_\star\|_2 & \leq & \|x_0-x_\star\|_2\max_{x\in[\mu,L]} |C_N^{[\mu,L]}(x)| \\
        & = & \|x_0-x_\star\|_2 \frac{1}{\text{cosh}\left( N \; \text{acosh}\left( \frac{L+\mu}{L-\mu} \right)\right)}.
    \EEAS
    First, we evaluate the $\text{acosh}$ term,
    \BEAS
        \text{acosh}\left( \frac{L+\mu}{L-\mu} \right) & = & \text{ln}\left(\frac{L+\mu}{L-\mu}+\sqrt{\left(\frac{L+\mu}{L-\mu}\right)^2-1}\right), \\
        & = & \text{ln}\left( \xi \right),\quad \text{where }\;\xi = \frac{\sqrt{\frac{L}{\mu}}+1}{\sqrt{\frac{L}{\mu}}-1}.
    \EEAS
    After plugging this result into the $\text{cosh}$, we get
    \[
        \frac{1}{\text{cosh}\left( N \; \text{ln}\left( \xi \right)\right)} = \frac{2}{e^{N \; \text{ln}\left( \xi \right)} + e^{-N \; \text{ln}\left( \xi \right)}}= \frac{2}{\xi^N + \xi^{-N}},
    \]
    thereby reaching the desired result.
\end{proof}

It may be difficult to compare the convergence rate of the Chebyshev method with that of gradient descent, due to its more complex expression. However, by neglecting the denominator term $\xi^{-N}$, we obtain the following upper bound:
\[
    \|x_{N}-x_\star\|_2 \leq 2 \left(\frac{\sqrt{\kappa}-1}{\sqrt{\kappa}+1}\right)^N\|x_0-x_\star\|_2.
\]

Note that the convergence rate of Polyak's heavy-ball method matches (up to a multiplicative factor) that of Chebyshev's method asymptotically, which is better than that of gradient descent in \eqref{eq:error-bnd-kappa}, which reads
\[
    \|x_{N}-x_\star\|_2 \leq \left(\frac{{\kappa}-1}{{\kappa}+1}\right)^N\|x_0-x_\star\|_2.
\]

We summarize this result in the following corollary, which compares the number of iterations required to reach a target accuracy $\epsilon$.

\begin{corollary}\label{cor:iter_comp}
    Let $x_0\in\mathbb{R}^d$ and $f$ be a quadratic function defined as in~\eqref{eq:quad-prob} with $\mu \idm \preceq \HH \preceq L \idm$ for some $L>\mu>0$. After respectively
    \begin{itemize}
        \item $N \geq \frac{L}{2\mu}~\log\left( \frac{\|x_0-x_\star\|_2}{\epsilon} \right) $ iterations of gradient descent (Algorithm~\ref{alg:cheb-grad} with~\eqref{eq:step}), or
        \item $N \geq \sqrt{\frac{L}{2\mu}}~\log\left( \frac{\|x_0-x_\star\|_2}{\epsilon} \right) $ iterations of Chebyshev's method (Algorithm~\ref{alg:cheb}),
    \end{itemize}
    we have that
    \[
        \|x_N-x_\star\|_2 \leq \epsilon.
    \]
\end{corollary}
\begin{proof}
    For gradient descent, a sufficient condition on the number of iterations required to reach an accuracy of $\epsilon$ reads
    \[
        \|x_N-x_\star\|_2 \leq  \left(\frac{{\kappa}-1}{{\kappa}+1}\right)^N \|x_0-x_\star\|_2 \leq \epsilon.
    \]
    Taking the log on both sides, we get
    \[
        N \geq  \frac{ \log\left( \frac{\|x_0-x_\star\|_2}{\epsilon} \right)}{\log \left(\frac{{\kappa}+1}{{\kappa}-1}\right)}.
    \]
    Using the bound $\log\left(\frac{\frac{1}{x}+1}{\frac{1}{x}-1}\right) > 2x$, the above condition can be simplified to the following stronger condition on $N$
    \[
        N \geq \frac{\kappa}{2}\log\left( \frac{\|x_0-x_\star\|_2}{\epsilon} \right).
    \]
    This gives the desired result for gradient descent. With the same approach, we also get the result for the Chebyshev algorithm.
\end{proof}

This corollary shows that the Chebyshev method can be~$\sqrt{\kappa}$ \textit{faster} than gradient descent. This translate to a speedup factor of $100$ in problems with a (reasonable) condition number of $10^4$, which is very significant. 

\paragraph{Worst-case optimality of Chebyshev's method.} When the dimension $d$ of the ambient space is sufficiently large, and without further assumptions on the spectrum of~$\HH$, the worst-case guarantee on Chebyshev's method is essentially unimprovable. Informally, given a budget $N$, a problem class $\mathcal{M}$ and some $R>0$, the best worst-case guarantee on the distance to optimality $\|x_N-x_\star\|_2$ that can be achieved by a first-order method is given by
\begin{equation}\label{eq:wc_cg_dist}
\max_{\substack{\HH\in\mathcal{M},\,x_\star,x_0\in\mathbb{R}^d\\ \|x_0-x_\star\|_2 \leq R}} \;\; \min_{\substack{P \in\mathcal{P}_N,\\ P(0)=1}} \|P(\HH)(x_0-x_\star)\|_2,
\end{equation}
which corresponds to the worst-case performance of the best performing method on any problem of the class. More precisely, for any first-order method satisfying $x_k \in x_0 +\Span\left\{ \nabla f(x_0),\, \nabla f(x_1),\,\ldots,\, \nabla f(x_{k-1}) \right\}$ for all $k\geq 1$ (the ``span assumption'') applied on the quadratic problem~\eqref{eq:quad-prob}, it holds that:
\begin{align*}
         \Span&\left\{ \nabla f(x_0),\, \ldots,\, \nabla f(x_{k-1}) \right\}\\&\subseteq \Span\left\{ H(x_0-x_\star),\, \ldots,\, H^{k}(x_{0}-x_\star) \right\},
\end{align*}
and therefore, the conjugate gradient-like method 
\begin{equation}\label{eq:CG_dist}
\begin{aligned}
         x_{k}=\argmin_{x} &\,\|x-x_\star\|_2\\ \text{s.t. }& x\in x_0 + \Span\{ H(x_0-x_\star),\, \ldots,\, H^{k}(x_{0}-x_\star) \},
\end{aligned}
\end{equation}
is \emph{instance-optimal}: it achieves the best worst-case performance on any problem instance. It follows that the worst-case performance of any first-order method satisfying the span assumption can only be worse than that of~\eqref{eq:CG_dist}. The worst-case performance of~\eqref{eq:CG_dist} being given by~\eqref{eq:wc_cg_dist}, we have that for any initialization $x_0\in\mathbb{R}^d$ any first-order method satisfying the span assumption in Equation~\eqref{eq:first-order-cheby}, there exists at least one problem on which
\begin{equation}\label{eq:wc_lb}
    \|x_N-x_\star\|_2 \geq \max_{\substack{\HH\in\mathcal{M},\,x_\star\in\mathbb{R}^d\\ \|x_0-x_\star\|_2 \leq R}} \;\; \min_{\substack{P \in\mathcal{P}_N,\\ P(0)=1}} \|P(\HH)(x_0-x_\star)\|_2,
\end{equation}
where $x_N\in\mathbb{R}^d$ is the output of the first-order method under consideration, and $x_\star$ is the optimal point of the problem. 

In other words, the $\max$ term in~\eqref{eq:wc_lb} searches for the ``most difficult'' quadratic function, while the $\min$ term represents the \textit{best} first-order method for a specific quadratic function. Of course, the method~\eqref{eq:CG_dist} is much more powerful than the Chebyshev method, since it is optimal \textit{for any specific function}. However, it is possible to show that despite being more powerful, this optimal algorithm has the same worst-case performance as that of Chebyshev's method when the dimension $d$ of the ambient space is large enough. The lower bound result is summarized by the next theorem.

\begin{theorem}\citep[Proposition~12.3.2]{nemirovskinotes1995}\label{thm:quad_lb}
    Let $N,d\in\mathbb{N}$ such that $d\geq N+1$, and let $x_0\in\mathbb{R}^d$.
    There exists $\HH \in\symm_d : \mu \idm \preceq \HH  \preceq L \idm$ and $x_\star \in\mathbb{R}^d$ such that any sequence $\{x_k\}_{k=0, \ldots, N}$ generated by any first-order method satisfying \eqref{eq:first-order-cheby}, and initiated at $x_0$ for minimizing the quadratic function $f$ in the form~\eqref{eq:quad-prob} satisfies
\begin{align*}
        \|x_N-x_\star\|_2 \geq \|x_0-x_\star\|_2 \min_{\substack{P \in\mathcal{P}_N,\\ P(0)=1}} \max_{\lambda \in [\mu,\, L]} |P(\lambda)|= \frac{2 }{\xi^N + \xi^{-N}}\|x_0-x_\star\|_2,
\end{align*}
with $ \xi = \frac{\sqrt{\frac{L}{\mu}}+1}{\sqrt{\frac{L}{\mu}}-1}$.
\end{theorem}
More details on this topic can be found in~\citep[Section 12.3]{nemirovskinotes1995}.

\section{Notes and References}\label{s:cheby_notesandrefs}

The Chebyshev method presented in this section is \textit{worst-case} optimal for the class of quadratic functions with Hessians $\HH$ satisfying $\mu \idm\preceq \HH \preceq L\idm$. More detailed discussions and developments on the topic of Chebyshev polynomials, for quadratic minimization, are provided in~\citep{Nemi83,nemirovsky1992information,Nest03a}, as well as in the lecture notes~\citep[Chapter 10]{nemirovskinotes1995}. Those references include the treatment of the case where the smallest eigenvalue is $\mu=0$. Finally, one should note that the optimal convergence bounds achieved by the Chebyshev method requires knowledge of the problem class parameters, $\mu$ and $L$, which might or might not be an issue, depending on the problem at hand. 

Probably the most celebrated method for unconstrained quadratic optimization problems is the conjugate gradient (CG) method. Its origin is usually attributed to \citet{stiefel1952methods,straeter1971extension}. As for the setup of this section, it turns out that CG methods are \emph{instance-optimal}, in the sense that they are the best performing first-order methods on every particular problem instance in the range of unconstrained quadratic minimization problems (in particular, the CG variant presented in~\eqref{eq:CG_dist} achieves the lower bound from Theorem~\ref{thm:quad_lb}). The classical CG produces iterates $\{x_{k}\}_{k\geq 0}$ such that
\begin{equation*}
\begin{aligned}
x_{k+1}\in\argmin_{x} &\, f(x) \\\text{s.t. }& x\in x_0 + \Span\{H(x_0-x_\star),\, \ldots,\, H^{k}(x_{0}-x_\star) \},
\end{aligned}
\end{equation*}
which admits efficient formulation; see, e.g.,~\citep{nocedal2006numerical}. Another variant of CG is often referred to as~MINRES, which produces iterates $\{x_{k}\}_{k\geq0}$ in the form
\begin{equation*}
\begin{aligned}
x_{k+1}\in\argmin_{x} &\, \|\nabla f(x)\|_2 \\\text{s.t. }& x\in x_0 + \Span\{H(x_0-x_\star),\, \ldots,\, H^{k}(x_{0}-x_\star) \}.
\end{aligned}
\end{equation*}
Its generalization GMRES~\citep{saad1986gmres} is popular for solving linear systems of the form $\HH x =b$ when $\HH$ is not required to be either symmetric or invertible.

Yet another alternative for dealing with quadratic minimization is to resort on \emph{Anderson}-type acceleration schemes. As for conjugate gradient methods, those schemes do not readily extend beyond quadratic minimization with the same nice theoretical guarantees. This is the topic of the next section.

Beyond quadratic optimization, properties of Chebyshev polynomials is the focus of~\citep{mason2002chebyshev}. The use of Chebyshev polynomials in the context of solving linear systems is covered at length in~\citep{fischerpolynomial}. In particular, the theory of~\citep{fischerpolynomial} can be instantiated for the convex quadratic minimization in average-case analyses, where Chebyshev polynomials (along with their heavy-ball limits) also naturally appear~\citep{pedregosa2020acceleration,lacotte2020optimal,scieur2020universal}.

\chapter{Nonlinear Acceleration}\label{c-RNA}

In this section, we see that the main argument used in the Chebyshev method can be adapted beyond quadratic problems. The extension that we present here, called nonlinear acceleration, follows a pattern that is known in numerical analysis as {\em vector extrapolation methods}: it seeks to accelerate the convergence of sequences by extrapolation using nonlinear averages. Different such strategies are known under various names, starting with Aitken's $\Delta^2$ \citep{Aitk27}, Wynn's epsilon algorithm \citep{Wynn56}, and Anderson acceleration \citep{anderson1965iterative}; a survey of these techniques can be found in \citep{Sidi86}. The vector extrapolation techniques, generic by nature, can be applied to optimization, as explained in what follows.

\section{Introduction}

This section focuses on the convex minimization problem:
\BEQ\label{eq:prob-nacc}
    \mbox{minimize} ~ f(x)
\EEQ
in the variable $x\in\mathbb{R}^d$. We assume $f$ to be twice continuously differentiable in a neighborhood of its minimizer~$x_\star$, and denote by $f_{\star} = f(x_\star)$ the minimum of $f$. 

We aim at adapting some of the ideas behind Chebyshev's acceleration (see Section~\ref{c-Cheb}) to a broader class of convex minimization problems beyond quadratic minimization. These adaptations stems from a local quadratic approximation of the objective:
\BEQ\label{eq:local-qp}
f(x) = f_\star + \frac12 \, \langle x-x_\star; \HH (x-x_\star)\rangle + o(\|x-x_\star\|_2^2),
\EEQ
where $\HH=\nabla^2f(x_\star)\in\symm_d$ (the set of symmetric $d\times d$ matrices) is the Hessian of $f$ at~$x_\star$, which we assume to satisfy $\mu \idm \preceq \HH \preceq L \idm$ for some $0<\mu<L$. Of course, neglecting the second-order term in~\eqref{eq:local-qp} allows recovering a quadratic minimization problem for which one could apply Chebyshev's method as is.  

We recall that the Chebyshev method is the first-order method associated with the \textit{best worst-case polynomial}. In short, the $k$th iteration of Chebyshev's method consists in combining previous gradients for minimizing a worst-case convergence bound over all $\mu$-strongly convex $L$-smooth quadratic problems in $\mathbb{R}^d$ (in the form~\eqref{eq:quad-prob} and with $d\geq k+1$):
\begin{align}\label{eq:opt_wc_method}
\begin{aligned}
\alpha_k &=\argmin_{\{\alpha^{(i)}\}_{i}}\, 
    \max_{\substack{\mu \idm \preceq \HH \preceq L \idm\\x_0,\ldots,x_k,\,x_\star,\,b\in\mathbb{R}^{k+1}}} \, &&\frac{\|x_k-x_\star\|_2^2}{\|x_0-x_\star\|^2_2}\\ 
    && \text{s.t. }& x_k=x_0-\sum_{i=0}^{k-1}\alpha^{(i)}\nabla f(x_i),\\
    & && \quad\HH x_\star=b,\\
    x_k &=x_0-\sum_{i=0}^{k-1}\alpha_{k}^{(i)}\nabla f(x_i),
\end{aligned}
\end{align}
so that $\alpha_k$ does not depend on the particular problem instance $(\HH,x_\star)$ and on the initialization $x_0$, but only on the problem class described by $\mu$ and $L$ (we note that in Section~\ref{c-Cheb},~\eqref{eq:opt_wc_method} was expressed in terms of optimizing a polynomial~\eqref{eq:problem_chebyshev}). A natural alternative to Chebyshev's method consists in choosing those weights {\em adaptively}. That is, depending on the particular instance of the problem at hand. For doing so, we have to choose another way to measure performance (because minimizing $\|x_k-x_\star\|_2$ would require knowledge of $x_\star$); one such possibility is to rely on function values or gradient norms. One could then rely on conjugate gradient-type methods which are very attractive for unconstrained quadratic minimization (see, e.g., discussions in Section~\ref{s:cheby_notesandrefs}). In this section, we consider the case where a first-order optimization method provided us with a sequence of pairs $\{(x_i, \nabla f(x_i)\}_{i=0,\ldots,k}$ satisfying (for $i=1,\ldots,k$)
\BEQ\label{eq:first_order_algorithm}
    x_{i} \in x_0 + \Span\left\{ \nabla f(x_0),\, \nabla f(x_1),\,\ldots,\, \nabla f(x_{i-1}) \right\},
\EEQ
and we study methods producing approximations of $x_\star$ as linear combinations $\sum_{i=0}^kc_ix_i$ of the previous iterates $\{x_i\}_{i=0,\ldots,k}$. For choosing the corresponding weights, we minimize the norm of the gradient at the approximated point. In unconstrained convex quadratic minimization problems, this approach is closely related to the so-called MINRES~\citep{paige1975solution} and GMRES \citep{saad1986gmres} methods (conjugate gradient-type methods minimizing gradient norms; see discussions by~\citet{walker2011anderson} and Section \ref{s:cheby_notesandrefs}). That is, when $f$ is quadratic, we choose the weights $\{c_i\}_{i=0,\ldots,k}$ by solving
\BEQ\label{eq:best_nonlinear}
    c_\star = \argmin_{c} \left\{ \left\|\nabla f \left({\textstyle \sum_{i=0}^kc_ix_i}\right) \right\|_2^2:\, \ones^T c = 1\right\}.
\EEQ
Whereas the new approximation is a linear combination of previous iterates $\{x_i\}_{i=0,\ldots,k}$, the coefficients $\{c_i\}_{i=0,\ldots,k}$ depend \textit{nonlinearly} on both $\nabla f$ and on $\{x_i\}_{i=0,\ldots,k}$. This technique is known under a few different names including those of \emph{Anderson acceleration} and \emph{minimal polynomial extrapolation} (see discussions and references in Section \ref{s:relatedw} for more details). In this section, we refer to all these methods as ``nonlinear acceleration'' techniques.

\section{Nonlinear Acceleration for Quadratic Minimization}

In this section, we present the simplest form of nonlinear acceleration, which is often referred to as the \textit{offline} nonlinear acceleration mechanism. We start with the main arguments underlying the technique, and present a few variants later in this section.

The core idea of the mechanism is to use a sequence of iterates $\{x_i\}_{i=0,\ldots,k}$ provided by a first-order method for solving~\eqref{eq:prob-nacc}. On this basis, we generate a new approximation of a solution to~\eqref{eq:prob-nacc} as a linear combination of past iterates, in the form $x_{\text{extr}}=\sum_{i=0}^kc_ix_i$. The point $x_{\text{extr}}$ is commonly referred to as the \textit{extrapolation} and can be chosen in different ways. In classical nonlinear acceleration mechanisms, it is chosen for making $\|\nabla f(x_{\text{extr}})\|_2$ small, as in~\eqref{eq:best_nonlinear}.  In general, solving~\eqref{eq:best_nonlinear} is just as costly as solving \eqref{eq:prob-nacc}, but the mechanism turns out to have an efficient formulation when minimizing quadratic functions of the form
\BEQ \label{eq:quadratic_function_nonlinear_acceleration}
    f(x) = \frac12 \langle x-x_\star ; \HH(x-x_\star) \rangle + f_\star.
\EEQ 
In this case, it is possible to find an explicit formula for~\eqref{eq:best_nonlinear}. Indeed, the gradient of $f$ is then a linear function, and because the coefficients $\{c_i\}_{i=0,\ldots,k}$ sum to one, the gradient of the linear combination is equal to a linear combination of gradients:
\BEQ \label{eq:equality_grad_combination}
    \nabla f \left(\sum_{i=0}^kc_ix_i\right) =\HH\left(\sum_{i=0}^kc_ix_i - x_\star\right) = \sum_{i=0}^kc_i \HH\left(x_i - x_\star\right) = \sum_{i=0}^kc_i \nabla f(x_i) .
\EEQ
It follows that~\eqref{eq:best_nonlinear} reduces to a simple quadratic program involving gradients of the past iterates. It can be formulated as
\BEQ
    c_\star = \argmin_{c} \left\{ \left\| {\textstyle \sum_{i=0}^kc_i} \nabla f \left(x_i\right) \right\|_2^2:\, c^T\ones = 1\right\}. \label{eq:nonlinear_acc_subproblem}
\EEQ
For convenience, we use the following more compact form in the sequel
\BEQ
    c_\star = \argmin_{c^T\ones = 1} \|\GG c\|_2^2, \label{eq:nonlinear_acc_subproblem_matrix}
\EEQ
where $\GG = [\nabla f(x_0), \ldots, \nabla f(x_k)]$ is the matrix formed by concatenating past gradients. This quadratic subproblem requires solving a small linear system of $(k+1)$ equations. When $\GG^T \GG$ is invertible, an explicit solution is provided by 
\[
    c_\star = \frac{z}{\sum_{i=0}^{k}z_i}, \quad \text{where }z = (\GG^T\GG)^{-1}\ones.
\]
This mechanism is summarized in Algorithm~\ref{alg:nonlinear_acceleration}.

\begin{algorithm}[!ht] 
  \caption{Nonlinear acceleration (offline version)}
  \label{alg:nonlinear_acceleration}
  \begin{algorithmic}[1]
    \REQUIRE Sequence of pairs $\{(x_i,\,\nabla f(x_i))\}_{i=0,\ldots,k}$.
      \STATE Form the matrix $\GG = [\nabla f(x_0),\,\ldots,\,\nabla f(x_k)]$, and compute $\GG^T\GG$.
      \STATE Solve the linear system $(\GG^T\GG)z = \ones$, and set $c = \frac{z}{z^T\ones}$.
      \STATE Form the extrapolated point $x_{\text{extr}} = \sum_{i=0}^{k} c_ix_i$.
    \ENSURE Approximate solution $x_{\text{extr}}$.
  \end{algorithmic}
\end{algorithm}

\subsection{Worst-case Convergence Bounds}
In this section, we quantify the accuracy of nonlinear acceleration. In particular, we show that it is instance-optimal and achieves the same worst-case convergence rate as that Chebyshev's method (see Theorem~\ref{thm:cheb_rate_convergence} and Theorem~\ref{thm:quad_lb}) in the worst-case, as soon as the sequence $\{x_i\}_{i=0,\ldots,k}$ is generated by a reasonable first-order method. Before going into the analysis, we introduce a few technical ingredients specifying what is a \emph{reasonable} first-order method. In short, we require that $\{x_i\}_{i=0,\ldots,k}$ is obtained from a ``nondegenerate'' first-order method. That is, we assume that the method uses $\nabla f(x_i)$ non-trivially for generating $x_{i+1}$ (for all $i=0,\ldots,k-1$).

\begin{definition}[Nondegenerate first-order method] \label{def:nondegenerate_fom} 
    Let $x_0\in\mathbb{R}^d$ be an initial point, and $f:\mathbb{R}^d\rightarrow\mathbb{R}$ be a continuously differentiable convex function. A first-order method generates sequences of iterates $\{x_i\}_{i=0,1,\ldots}$ such that for all $i=1,2,\ldots$
    \[ x_i\in x_0+\Span\{\nabla f(x_0),\ldots,\nabla f(x_{i-1})\}.\]
    A first-order method is \emph{nondegenerate} if for all continuously differentiable convex function $f:\mathbb{R}^d\rightarrow\mathbb{R}$, all $i=0,1,\ldots$, and all $x_0\in\mathbb{R}^d$ there exists some $\{\alpha_j^{(i)}\}_{j=0,\ldots,i}\subset \mathbb{R}$ with $\alpha^{(i)}_i \neq 0$ such that     
    \begin{equation}\label{eq:nondegenerate_fom}
        x_{i+1} = x_0 + \sum_{j=0}^i \alpha^{(i)}_j \nabla f(x_j).
\end{equation}
\end{definition}

The proposition below shows that when the sequence~$\{x_i\}_{i=0,1,\ldots}$ is generated by a nondegenerate first-order method, the gradient of any iterate $x_i$ can be written using a polynomial of degree \textit{exactly} $i$.

\begin{proposition} \label{prop:nondegenerate_polynomials}
    Let $x_0\in\mathbb{R}^d$ be an initial point, $f:\mathbb{R}^d\rightarrow \mathbb{R}$ be a quadratic function in the form~\eqref{eq:quadratic_function_nonlinear_acceleration} with $\mu\idm\preceq\HH\preceq L\idm$, and let $\{x_i\}_{i=0,1,\ldots}$ be generated by a nondegenerate first-order method. Then, for all $i=0,1,\ldots$, there exists a polynomial $P_i$ of degree \textit{exactly} $i$ such that $P_i(0) = 1$ and
    \[
        \nabla f(x_i) = P_i(\HH)\nabla f(x_0).
    \]
\end{proposition}
\begin{proof}
    We proceed by induction. First, we have $P_0 = 1$ ($P_0$ has degree $0$ and $P_0(0)=1$) and hence
    \[
        \nabla f(x_0) = P_0(\HH)\nabla f(x_0),
    \]
    thereby trivially reaching the desired conclusion for $i=0$.
    
    We proceed with the induction hypothesis, assuming that the desired result holds for~$x_i$. By Definition~\ref{def:nondegenerate_fom} (nondegenerate first-order method), and because $f$ is a quadratic function \eqref{eq:quadratic_function_nonlinear_acceleration}, we have
    \begin{align*}
        \nabla f(x_{i+1}) = &\nabla f\left(x_0 + \sum_{j=0}^i \alpha^{(i)}_j \nabla f(x_j)\right)\\
        =&\HH\left(x_0-x_{\star} + \sum_{j=0}^i \alpha^{(i)}_j \nabla f(x_j)\right)\\
        =&\HH(x_0-x_{\star}) + \sum_{j=0}^i \alpha^{(i)}_j \HH\nabla f(x_j)\\
        =&\nabla f(x_0) + \sum_{j=0}^i \alpha^{(i)}_j \HH\nabla f(x_j).
    \end{align*}
    Thanks to the induction hypothesis, we have $\nabla f(x_j) = P_j(\HH)\nabla f(x_0)$ with $P_j(0)=1$ and $\deg(P_j) = j$ for $j\leq i$. For showing that $P_{i+1}$ satisfies the desired claim, we start by expressing $ \nabla f(x_{i+1})$ in terms of a polynomial:
    \[
         \nabla f(x_{i+1}) = \underbrace{\left(P_0(\HH) + \sum_{j=0}^i \alpha^{(i)}_j\HH P_j(\HH)\right)}_{=P_{i+1}(\HH)}  \nabla f(x_0).
    \]
    It is relatively straightforward to verify $P_{i+1}(0)=1$ using the previous expression:
    \[
        P_{i+1}(0) = P_0(0) + \sum_{j=0}^i \alpha^{(i)}_j \cdot 0 \cdot P_j(0) = P_0(0) = 1.
    \]
    Finally, a minor reorganization of the expression of $P_{i+1}$ allows writing
    \[
        P_{i+1}(\HH) = \underbrace{P_0(\HH) + \sum_{j=0}^{i-1} \alpha^{(i)}_j\HH P_j(\HH)}_{\text{degree} \leq i} + \alpha_{i}^{(i)} \HH P_{i}(\HH).
    \]
    Nondegeneracy of the first-order method implies that there exists some $\alpha_{i}^{(i)} \neq 0$ such that the previous expression holds. Finally it follows from $\deg(P_{i}) = i$ (induction hypothesis) that $\deg(\HH P_{i}(\HH))=i+1$, thereby reaching the desired claim.
\end{proof}

Equipped with previous technical ingredients, one can show that Algorithm~\ref{alg:nonlinear_acceleration} is ``instance-optimal'' when applied to a nondegenerate first-order method. This means that the nonlinear acceleration algorithm finds the \textit{best polynomial} given a \textit{specific} quadratic function $f$---in opposition to the Chebyshev method that finds the best polynomial for a \textit{class} of functions (Theorem~\ref{thm:cheb_rate_convergence}). In other words, nonlinear acceleration \textit{adaptively} looks for the best combination of previous iterates given the information stored in the previous gradients, while Chebyshev's method uses the same worst-case optimal polynomial in all cases. Moreover, Algorithm~\ref{alg:nonlinear_acceleration} does \textit{not} require knowledge of the smoothness or strong convexity parameters. 

\begin{theorem}[Instance-optimality of nonlinear acceleration] \label{thm:rate_nonlinear_acc}
Let $x_0\in\mathbb{R}^d$ be an initial point, $f$ be the quadratic function from \eqref{eq:quadratic_function_nonlinear_acceleration} with $\mu\idm\preceq\HH\preceq L\idm$, and $\{x_i\}_{i=0,\ldots,k}$ be generated by a nondegenerate first-order method initiated from $x_0$. For any $k\geq 0$, it holds that
    \BEQ \label{eq:rate_convergence_quad_nonlinear_accel}
    \begin{aligned}
    	\| \nabla f(x_{\text{extr}}) \|_2 &= \;\; \min_{c^T\ones=1} \left\| \nabla f\left( {\textstyle \sum_{i=0}^k} c_ix_i\right) \right\|_2 \\&= \displaystyle \min_{\substack{P\in\mathcal{P}_{k},\\ P(0)=1}} \| P(\HH)\nabla f(x_0) \|_2,
    \end{aligned}
    \EEQ
    where $x_{\text{extr}}$ is obtained from Algorithm~\ref{alg:nonlinear_acceleration} $\{ (x_i,\nabla f(x_i)) \}_{i=0,\ldots,k}$, and $\mathcal{P}_k$ is the set of polynomials of degree at most $k$.
\end{theorem}

\begin{proof}
    Since the coefficients $c_i$ sum to one, we have the following equalities:
    \BEAS
        \nabla f(x_{\text{extr}}) & = & \nabla f\left( \sum_{i=0}^{k}c_ix_i\right)= \HH\left( \sum_{i=0}^{k}c_i(x_i-x_{\star})\right) \\
         & = & \left( \sum_{i=0}^{k}c_i\HH(x_i-x_{\star})\right) =  \sum_{i=0}^{k}c_i\nabla f(x_i).
    \EEAS
    Therefore,
    \[
        \|\nabla f(x_{\text{extr}})\|_2 = \left\|\sum_{i=0}^{k}c_i\nabla f(x_i)\right\|_2.
    \]
    We now use the definition of a first-order method, which yields iterates $x_{i}$ such that
    \[
        x_{i} \in x_0 + \Span\{ \nabla f(x_0),\,\ldots,\, \nabla f(x_{i-1}) \}.
    \]
    From Section~\ref{s:cheb-cheb}, Proposition \ref{prop:first_order_to_polynomial}, it follows that $x_i$ can be written as
    \[
        x_{i} = x_{\star} + P_i(\HH)(x_0-x_{\star}), \quad P_i \in \mathcal{P}_i, \quad P_i(0)=1.
    \]
    It also holds that its gradient can be written using the same polynomial:
    \BEA
        \nabla f(x_{i}) & = & \HH(x_i-x_{\star}) = \HH P_i(\HH)(x_0-x_{\star}) = P_i(\HH)(\HH(x_0-x_{\star})) \nonumber \\
        & = & P_i(\HH)\nabla f(x_0). \label{eq:equality_grad_poly}
    \EEA
     By substituting this expression in the objective of \eqref{eq:nonlinear_acc_subproblem}, we obtain
    \BEQ \label{eq:temp_rna_to_poly}
        \min_{c^T\ones = 1}\left\| \sum_{i=0}^k c_i \nabla f(x_i) \right\|_2^2 = \min_{c^T\ones = 1} \left\| \left(\sum_{i=0}^k c_i P_i(\HH)\right)\nabla f(x_0) \right\|_2^2.
    \EEQ
    Since the iterates $\{x_i\}_{i=0,\ldots,k}$ are generated by a nondegenerate first-order method, all polynomials $P_i$ have differents degrees. Therefore, the polynomials $\{P_i\}_{i=0,\ldots,k}$ are linearly independent, and hence $\{P_i\}_{i=0,\ldots, k}$ is a basis for the space~$\mathcal{P}_k$. Finally, because $P_i(0) = 1$ and $c^T\ones=1$, we can rephrase the objective~\eqref{eq:temp_rna_to_poly} as
    \[
        \min_{c^T\ones = 1} \left\| \left(\sum_{i=0}^k c_i P_i(\HH)\right)\nabla f(x_0) \right\|_2 = \min_{\substack{P\in\mathcal{P}_k,\\P(0)=1}} \left\| P(\HH)\nabla f(x_0) \right\|_2.
    \]
    The convergence rate is thus given by
    \BEAS
        \|\nabla f(x_{\text{extr}})\|_2 & =  & \left\| \sum_{i=0}^{k}c_i\nabla f(x_i)\right\|_2 = \min_{\substack{P\in\mathcal{P}_k,\\P(0)=1}} \left\| P(\HH)\nabla f(x_0) \right\|_2.
    \EEAS
\end{proof}

This theorem shows that the worst-case convergence rate is essentially controlled by the optimal value of a minimization problem. The minimum in \eqref{eq:rate_convergence_quad_nonlinear_accel} can be bounded using a Chebyshev argument similar to the main argument used in Section~\ref{c-Cheb}.

\begin{corollary}\label{cor:wc_rate_na}
    Let $x_0\in\mathbb{R}^d$ be an initial point, $f$ be the quadratic function from \eqref{eq:quadratic_function_nonlinear_acceleration} with $\mu\idm\preceq\HH\preceq L\idm$, and $\{x_i\}_{i=0,\ldots,k}$ be generated by a nondegenerate first-order method initiated from $x_0$. Then, for any $k\geq 0$, it holds that
    \[
        \| \nabla f(x_{\text{extr}}) \|_2 \leq \dfrac{2}{\xi^k + \xi^{-k}}\|\nabla f(x_0)\|_2, \quad \xi = \frac{\sqrt{\frac{L}{\mu}}+1}{\sqrt{\frac{L}{\mu}}-1}.
    \]
    where $x_{\text{extr}}$ is the output of Algorithm~\ref{alg:nonlinear_acceleration} applied to~$\{ (x_i,\,\nabla f(x_i)) \}_{i=0,\ldots,k}$.
\end{corollary}
\begin{proof}
    We use the shifted Chebyshev polynomial from Theorem \ref{thm:cheb_rate_convergence} as a feasible solution of the minimization problem \eqref{eq:rate_convergence_quad_nonlinear_accel}.
\end{proof}
As for the Chebyshev method, it is possible to show that the convergence rate cannot be improved as it matches that of the corresponding lower bound, which can be obtained by adapting Theorem~\ref{thm:quad_lb} to gradient norms; see e.g.,~\citep[Proposition~12.3.2]{nemirovskinotes1995}.

\begin{remark}[Finite-time convergence]\label{rem:instance-opt}
When~$\{(x_i, \,\nabla f(x_i))\}_{i=0,\ldots,k}$ is generated by a nondegenerate first-order method and when $k$ is large enough, nonlinear acceleration eventually converges \textit{exactly} to the minimizer of the quadratic function $f$ (this follows easily from analogies with conjugate gradient-type methods; see, e.g.,Section \ref{s:cheby_notesandrefs}). More formally, for all $k\geq d$ it holds that
\BEQ
    x_{\text{extr}} = x_{\star},
\EEQ
where $x_{\text{extr}}$ is the output of Algorithm~\ref{alg:nonlinear_acceleration} applied to~$\{(x_i, \,\nabla f(x_i))\}_{i=0,\ldots,k}$. This is a natural consequence of the fact that $x_{\text{extr}}$ is the best point in $ x_0+\Span\{\nabla f(x_0),\ldots,\nabla f(x_{k-1})\}$, as provided by~\eqref{eq:nonlinear_acc_subproblem}. 
\end{remark}

\subsection{Computational Complexity}\label{rem:comp_complexity_offline}

 The computational complexity of Algorithm~\ref{alg:nonlinear_acceleration} is $O(dk^2 + k^3)$, where $k$ is the length of the sequence of gradients and $d$ is the dimension of the ambient space. The first term originates from the matrix-matrix multiplication $\GG^T\GG$ in step~1, and the second one comes from solving a $(k+1)\times(k+1)$ matrix in step~2. The length $k$ being often much smaller than $d$, the resulting complexity is typically~$O(dk^2)$.

\paragraph{Low-rank updates.} When using the nonlinear acceleration method in parallel with a first-order method generating a growing sequence of iterates $\{x_i\}_{i=0,\ldots,k}$, an extrapolation step can be computed each time a new iterate is produced. we can reduce the per iteration complexity up to $O(dk)$ by computing the matrix $\textbf{G}^T\textbf{G}$ and the coefficients $c_\star$ using low-rank updates \citep{sidi1991efficient}.

\paragraph{Limited-memory.} Because the iteration complexity of nonlinear acceleration grows with the length of the input sequence, it might become costly to compute an extrapolation. It is therefore common to use nonlinear acceleration only with the last few iterates produced by the first-order method. More formally, if we impose a maximum memory of $m$ pairs, we compute the extrapolation as follow when $k\geq m$. More formally, if we impose a maximum memory of $m$ pairs, one can use Algorithm~\ref{alg:nonlinear_acceleration} with input $\{(x_{k-m+i},\,\nabla f(x_{k-m+i}))\}_{i=1,\ldots,m}$.

\subsection{Online Nonlinear Acceleration}

So far, we have seen a post-processing procedure that generates an extrapolated point $x_{\text{extr}}$ from a sequence of pairs $\{x_i,\nabla f(x_i)\}_{i=0,\ldots,k}$. If this sequence is generated by a nondegenerate first-order method, then Corollary~\ref{cor:wc_rate_na} shows that the gradient of the extrapolated point $x_{\text{extr}}$ converges to zero at an optimal worst-case convergence rate, without any hyper-parameters. Perhaps surprisingly, Algorithm~\ref{alg:nonlinear_acceleration} is itself not a nondegenerate first-order method, and can therefore not be used recursively as is.

In what follows, we introduce a \emph{mixing} parameter. This parameter transforms the nonlinear acceleration mechanism to a nondegenerate first-order method without hurting its worst-case performance. This enables using nonlinear acceleration recursively for generating the whole sequence $\{x_i\}_{i=1,2,\ldots}$. This technique is often referred to as \emph{online nonlinear acceleration}.

\paragraph{Mixing parameter.} The idea underlying the mixing parameter is fairly simple: instead of combining previous iterates, we combine \textit{gradient steps} as follows:
\begin{equation} \label{eq:extrapolation_with_mixing}
    x_{\text{extr}}^{\text{mixing}} = \sum_{i=0}^k c_i \Big(x_i-h\nabla f(x_i)\Big),
\end{equation}
as provided by Algorithm \ref{alg:nonlinear_acceleration_mixing}. Furthermore, the use of an appropriate step size $h$ can even slightly improve the worst-case convergence speed of Algorithm \ref{alg:nonlinear_acceleration}.
\begin{algorithm}[!ht] 
  \caption{Nonlinear acceleration (with mixing)}
  \label{alg:nonlinear_acceleration_mixing}
  \begin{algorithmic}[1]
    \REQUIRE Sequence of pairs $\{(x_i,\,\nabla f(x_i))\}_{i=0,\ldots,k}$, mixing parameter $h$.
      \STATE Form the matrix $\GG = [\nabla f(x_0),\,\ldots,\,\nabla f(x_k)]$, and compute $\GG^T\GG$.
      \STATE Solve the linear system $(\GG^T\GG)z = \ones$, and compute $c = \frac{z}{z^T\ones}$.
      \STATE Form the extrapolated point $x_{\text{extr}} = \sum_{i=0}^{k} c_i\big(x_i-h\nabla f(x_i)\big)$.
    \ENSURE Approximate solution $x_{\text{extr}}$.
  \end{algorithmic}
\end{algorithm}
    
Intuitively, this \textit{mixing} between iterates and gradients emulates a gradient step on the extrapolated point from Algorithm \ref{alg:nonlinear_acceleration}, that we call here $x_{\text{extr}}^{\text{offline}} = \sum_{i=0}^k c_i x_i$,
\begin{align*}
    x_{\text{extr}}^{\text{mixing}} =& \sum_{i=0}^k c_i x_i-h \sum_{i=0}^k c_i \nabla f(x_i) \\
    =&  \sum_{i=0}^k c_i x_i-h \nabla f\left(\sum_{i=0}^k c_i x_i\right) \\
    =& x_{\text{extr}}^{\text{offline}} - h\nabla f(x_{\text{extr}}^{\text{offline}}),
\end{align*}
where the second equality comes from the fact that $f$ is a quadratic function, and that $\sum_{i=0}^k c_i = 1$.

This mixing parameter requires tuning one hyper-parameter $h$, which can be chosen in various ways. Proposition \ref{prop:performance_mixing_nonlinear} shows that the mixing strategy slightly improves the performance of nonlinear acceleration if $h$ is set properly.

\begin{proposition}\label{prop:performance_mixing_nonlinear}
Let $x_0\in\mathbb{R}^d$ be an initial point, $f$ be the quadratic function from \eqref{eq:quadratic_function_nonlinear_acceleration} with $\mu\idm\preceq\HH\preceq L\idm$, and $\{x_i\}_{i=0,\ldots,k}$ be generated by a nondegenerate first-order method initiated from $x_0$. Then, for any $k\geq 0$, it holds that
    \[
        \| \nabla f(x_{\text{extr}}^{\text{mixing}}) \|_2 \leq C_h\|\nabla f(x_{\text{extr}}^{\text{offline}})\|_2,
    \]
    where $x_{\text{extr}}^{\text{offline}}$ is obtained from Algorithm~\ref{alg:nonlinear_acceleration} applied to~$\{ (x_i,\,\nabla f(x_i)) \}_{i=0,\ldots, k}$, $x_{\text{extr}}^{\text{mixing}}$ is the extrapolation with mixing from \eqref{eq:extrapolation_with_mixing}, and 
    \[
        C_h = \max\left\{ 1-\mu h \;;\; L h-1 \right\}.
    \]
    Moreover, the factor $C_h$ is guaranteed to be smaller than one if $h\in(0,\,\frac{2}{L})$, and takes its minimal value at $h=\frac{2}{L+\mu}$.
\end{proposition}
\begin{proof}
    We start with the identity
    \[
        x_{\text{extr}}^{\text{mixing}} = x_{\text{extr}}^{\text{offline}} - h\nabla f(x_{\text{extr}}^{\text{offline}}).
    \]
    Because $f$ is the quadratic function \eqref{eq:quadratic_function_nonlinear_acceleration}, the gradient of $x_{\text{extr}}^{\text{mixing}}$ reads
    \begin{align*}
        \nabla f(x_{\text{extr}}^{\text{mixing}}) & = \HH(x_{\text{extr}}^{\text{mixing}}-x_\star)\\
        & = \HH(x_{\text{extr}}^{\text{offline}} - h\nabla f(x_{\text{extr}}^{\text{offline}})-x_\star) \\
        & = \nabla f(x_{\text{extr}}^{\text{offline}}) - \HH h \nabla f(x_{\text{extr}}^{\text{offline}})\\
        & = (I - \HH h) \nabla f(x_{\text{extr}}^{\text{offline}}).
    \end{align*}
    Therefore, we have the bound
    \[
        \|\nabla f(x_{\text{extr}}^{\text{mixing}})\|_2 \leq \| I - \HH h\|_2 \| \nabla f(x_{\text{extr}}^{\text{offline}})\|_2,
    \]
    and the desired result follows from $\mu \idm \preceq \HH \preceq L\idm$.
\end{proof}

\paragraph{Online nonlinear acceleration method.} As previously underlined, the mixing parameter transforms the nonlinear acceleration method into a nondegenerate first-order method. We can therefore use it recursively. The online variant of the nonlinear acceleration technique, with limited memory, is provided in Algorithm~\ref{alg:online_nonlinear_acceleration}, using Algorithm~\ref{alg:nonlinear_acceleration_mixing} as a subroutine. One should note that when $m=\infty$ (no memory restriction), the worst-case performance of offline version of nonlinear acceleration with mixing (Algorithm~\ref{alg:nonlinear_acceleration_mixing}) is also valid for its online variant (Algorithm~\ref{alg:online_nonlinear_acceleration}). It follows from Proposition~\ref{prop:performance_mixing_nonlinear} that the worst-case performance of Algorithm~\ref{alg:nonlinear_acceleration_mixing} and Algorithm~\ref{alg:online_nonlinear_acceleration} is no worse than that of offline version of nonlinear acceleration (Algorithm~\ref{alg:nonlinear_acceleration}), provided by Corollary~\ref{cor:wc_rate_na}.

\begin{algorithm}[!ht] 
  \caption{Nonlinear acceleration (online version, limited memory)}
  \label{alg:online_nonlinear_acceleration}
  \begin{algorithmic}[1]
    \REQUIRE A differentiable function~$f$, initial point~$x_0$, mixing parameter~$h$, maximum memory parameter $m$ (optional, $m=\infty$ by default).
    \STATE \textbf{Initialize} Empty sequence $\mathcal{S}$ of pairs iterate/gradient.
      \FOR{$k=0,\ldots$}
        \STATE Compute $\nabla f(x_k)$; append the pair $\mathcal{S} \leftarrow \mathcal{S} \cup \{(x_k,\,\nabla f(x_k))\}$.
        \IF{$k\geq m$}
            \STATE Discard the oldest pair from $\mathcal{S}$.
        \ENDIF
        \STATE Compute the extrapolation $x_{k+1} = \text{[Algorithm \ref{alg:nonlinear_acceleration_mixing}]}(\mathcal{S},h)$.
      \ENDFOR
    \ENSURE Approximate solution $x_{k+1}$.
  \end{algorithmic}
\end{algorithm}

In the next section, we see that nonlinear acceleration technique might suffer from serious instability issues when applied beyond quadratic minimization. Perhaps luckily, a simple regularization technique allows stabilizing the procedure beyond quadratics.

\section{Regularized Nonlinear Acceleration Beyond Quadratics}
Nonlinear acceleration suffers some serious drawbacks when used outside the restricted setting of quadratic functions. In fact, Algorithm~\ref{alg:nonlinear_acceleration} and Algorithm~\ref{alg:nonlinear_acceleration_mixing} are numerically highly unstable. This problem originates from the conditioning of the matrix $\GG^T\GG$, used to compute the coefficients $\{c_i\}_{i=0,\ldots,k}$. To illustrate this statement, assume we run Algorithm \ref{alg:nonlinear_acceleration} with a noisy sequence of gradients $\{(x_{i},\,\nabla f(x_i)+e_i)\}_{i=0,\ldots,k}$ where the sequence $\{e_i\}_{i=0,\ldots,k}$ is such that $\|e_i\|_2\leq \epsilon$ for some $\epsilon>0$. \citet[Proposition 3.1]{scieur2016regularized} show that the relative distance between $\tilde c$ (the coefficients computed using the noisy sequence above) and its noise-free version $c$ satisfies
\BEQ
    \frac{\| c-\tilde c \|_2}{\|c\|_2} = O \left(\|\EE\|_2\;\; \| ((\GG+\EE)^T(\GG+\EE))^{-1} \|_2\right). \label{eq:bound_perturbation}
\EEQ
where $\EE \triangleq [e_0,\,\ldots,\,e_k]$ is the noise matrix. In this bound, the perturbation impacts the solution proportionally to its norm and to the conditioning of $\GG+\EE$.

Unfortunately, even for small perturbations, the condition number of $\GG+\EE$ and the norm of the vector $c$ are usually \textit{huge}. In fact, $\GG$ has a \textit{Krylov} matrix structure, which is notoriously poorly conditioned \citep{Tyrt94}. Thereby, even a very small perturbation $\EE$ might have a significant impact on performance. For illustrating this, let us briefly illustrate the link between $\GG$ and Krylov matrices: consider using gradient descent with step size $\gamma$ on a quadratic function; the iterates follow the rule
\[
    x_{k+1}-x_\star = (1-\gamma \HH)^k(x_0-x_\star) \quad \Leftrightarrow \quad \nabla f(x_{k+1}) = (1-\gamma \HH)^k\nabla f(x_0).
\]
Thereby, a matrix $\GG$ formed by these expressions has the form
\[
    \GG = \Big[ \nabla f(x_0),\, (\idm-\gamma \HH)\nabla f(x_0),\, (\idm-\gamma \HH)^2\nabla f(x_0),\, \ldots \Big],
\]
which shows that $\GG$ is in fact a \textit{Krylov} matrix---by definition, a Krylov matrix $K$ associated with a matrix $A$ and vector $v$ is defined as $K = [v,\, Av, \, A^2v, \ldots]$.

\begin{figure}
     \centering
     \begin{subfigure}[b]{0.45\textwidth}
         \centering
         \includegraphics[width=\textwidth]{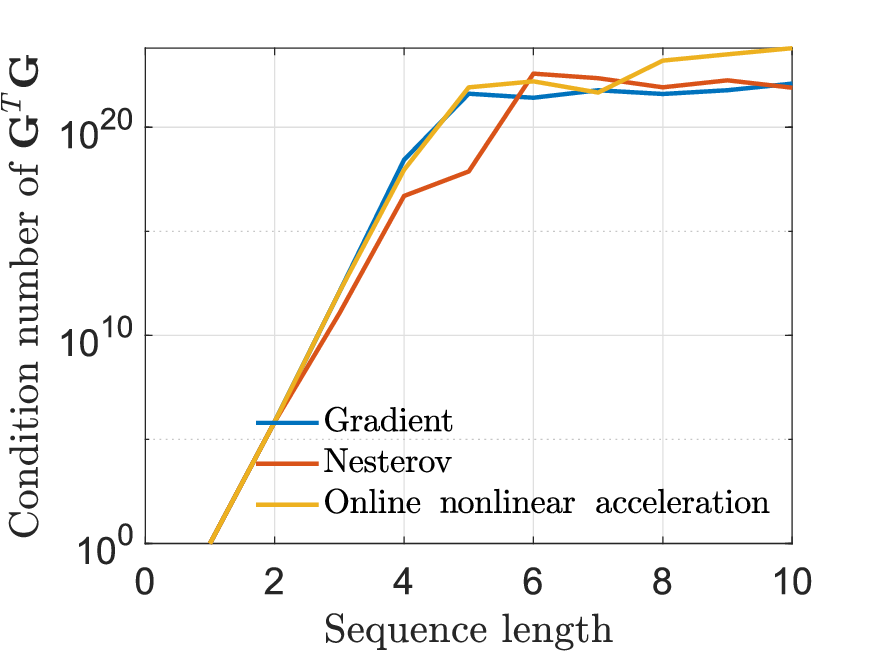}
         \caption{}
         \label{fig:sensitivity_nonlinear_acceleration_cong_gg}
     \end{subfigure}
     \hfill
     \begin{subfigure}[b]{0.45\textwidth}
         \centering
         \includegraphics[width=\textwidth]{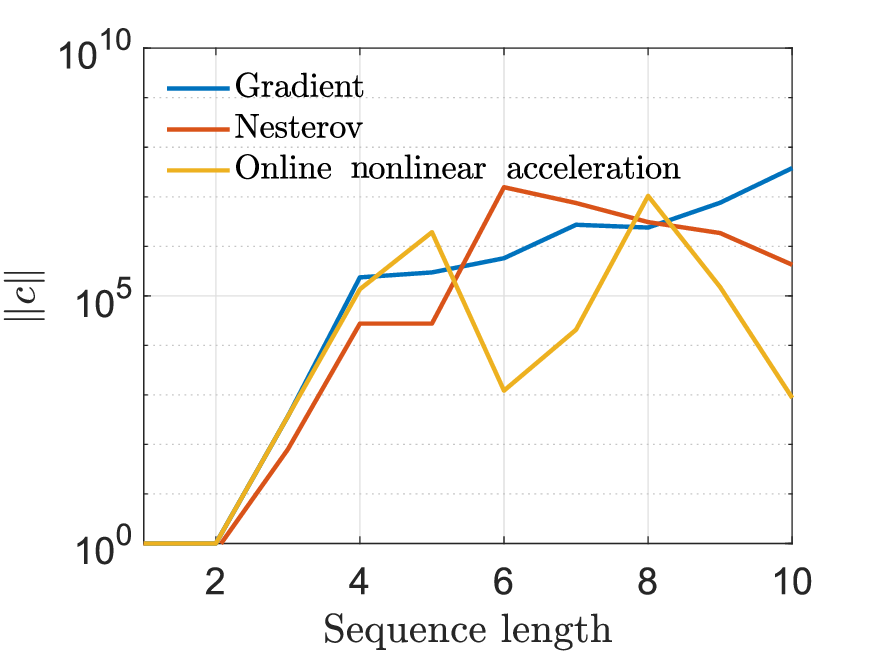}
         \caption{}
         \label{fig:sensitivity_nonlinear_acceleration_norm_c}
     \end{subfigure}
        \caption{Illustration of the sensitivity of nonlinear acceleration when applying nonlinear acceleration to gradient descent, to Nesterov's method (see Section \ref{c-Nest}), and using its online variant (Algorithm~\ref{alg:online_nonlinear_acceleration}) to minimize some random quadratic function. Figure \ref{fig:sensitivity_nonlinear_acceleration_cong_gg}: the condition number of the matrix $\GG^T\GG$, which grows exponentially with its size (the plateau on the right is caused by numerical errors). Figure \ref{fig:sensitivity_nonlinear_acceleration_norm_c}: the norm of the vector of coefficients $c$.}
        \label{fig:sensitivity_nonlinear_acceleration}
\end{figure}

In Figure~\ref{fig:sensitivity_nonlinear_acceleration}, we show the norm of $c$ and the condition number of the matrix $\GG^T\GG$ when it is formed from iterates of gradient descent, accelerated gradient descent (see Section~\ref{c-Nest}), and nonlinear acceleration (in the online setting, see Algorithm~\ref{alg:online_nonlinear_acceleration}) for minimizing some randomly generated quadratic function. Figure \ref{fig:sensitivity_nonlinear_acceleration_cong_gg} shows that even after 3 iterations, the system can already be considered singular (i.e., the condition number exceeds $10^{16}$).

For stabilizing the method, it is common to regularize the linear system. The resulting algorithm is often referred to as regularized nonlinear acceleration (RNA)~\citep{scieur2016regularized}. The following section is devoted to some theoretical properties of this method.

\subsection{Regularized Nonlinear Acceleration} 

Regularized nonlinear acceleration (RNA) consists of using Algorithm~\ref{alg:nonlinear_acceleration} with a regularization, thereby rendering the method less sensitive to noise. In short, the base operation underlying RNA is to solve
\BEQ \label{eq:rna_equation}
    \argmin_{c^T\ones = 1} \frac{\| \GG c \|_2^2}{\|\GG\|_2^2} + \lambda  \| c-c_{\text{ref}} \|_2^2
\EEQ
in the variable $c\in\mathbb{R}^k$, where $c_{\text{ref}}$ is some reference vector. The effect of regularization is therefore to force $c$ to be close to $c_{\text{ref}}$. Of course, it makes more sense to pick $c_{\text{ref}}$ summing to one. A common choice is to pick $c_{\text{ref}} = \ones/k$, which would enforce the procedure to be ``not too far'' from a simple averaging of the iterates. Another possibility is $c_{\text{ref}} = [\textbf{0}_{k-1}, 1]$, which puts more weight on the last iterate. The division by $\|\GG\|_2^2$ is for scaling purposes, as it makes $\lambda$ dimensionless. The resulting method is slightly more complicated than its previous version without regularization and is provided in Algorithm~\ref{algo:rna_general}. When $c_{\text{ref}} =\ones/k$, the procedure simplifies to Algorithm~\ref{algo:rna_simplified}.

\begin{algorithm}[ht] 
  \caption{Regularized nonlinear acceleration (RNA)}
  \label{algo:rna_general}
  \begin{algorithmic}[1]
      \REQUIRE Sequence of pairs $\{(x_i,\,\nabla f(x_i))\}_{i=0,\ldots,k}$, mixing parameter $h$, regularization term $\lambda>0$,  reference vector $c_{\text{ref}}$.
      \STATE Form $\GG = [\nabla f(x_0),\,\ldots,\,\nabla f(x_k)]$, compute $\mathcal{G} = \frac{\GG^T\GG}{\|\GG^T\GG\|_2}$.
      \STATE Solve the linear system $(\mathcal{G}+\lambda \idm)w = \lambda c_{\text{ref}}$.
      \STATE Solve the linear system $(\mathcal{G}+\lambda \idm)z = \ones$.
      \STATE Compute the coefficients $ c = w +z\frac{(1-w^T\ones)}{z^T\ones} $
      \STATE Form the extrapolated point $x_{\text{extr}} = \sum_{i=0}^{k} c_i \left(x_i - h\nabla f(x_i)\right)$.
    \ENSURE Approximate solution $x_{\text{extr}}$.
  \end{algorithmic}
\end{algorithm}

\begin{algorithm}[ht] 
  \caption{Regularized nonlinear acceleration (with $c_{\text{ref}} = \ones/k$)}
  \label{algo:rna_simplified}
  \begin{algorithmic}[1]
      \REQUIRE Sequence of pairs $\{(x_i,\,\nabla f(x_i))\}_{i=0,\ldots,k}$, mixing parameter $h$, regularization term $\lambda>0$.
      \STATE Form $\GG = [\nabla f(x_0),\,\ldots,\,\nabla f(x_k)]$ and compute $\mathcal{G} = \frac{\GG^T\GG}{\|\GG^T\GG\|_2}$.
      \STATE Solve the linear system $(\mathcal{G}+\lambda \idm)z = \ones/k$.
      \STATE Compute the coefficients $ c =\frac{z}{z^T\ones}$.
      \STATE Form the extrapolated point $x_{\text{extr}} = \sum_{i=0}^{k} c_i \left(x_i - h\nabla f(x_i)\right)$.
    \ENSURE Approximate solution $x_{\text{extr}}$.
  \end{algorithmic}
\end{algorithm}

\paragraph{Online regularized nonlinear acceleration.} As in the quadratic case, Algorithm~\ref{algo:rna_general} and Algorithm~\ref{algo:rna_simplified} could be used as subroutines in the online nonlinear acceleration method (Algorithm~\ref{alg:online_nonlinear_acceleration}), thereby forming the regularized version of the online acceleration algorithm.

\subsection{Perturbed Linear Gradients}

In this section, we consider the problem of minimizing a twice continuously differentiable convex function $f$, as in~\eqref{eq:prob-nacc}, beyond quadratic problems. For doing that, we introduce \textit{perturbed linear gradients}. As before, we consider iterates $\{x_i\}_{i=0,1,\ldots}$ originating from a first-order method satisfying
\BEAS
     x_{k+1} & \in & x_0 + \Span\{ \nabla f(x_0),\,\nabla f(x_1),\,\ldots,\nabla f(x_k) \}.
\EEAS
However, $\nabla f(x)$ is now no longer the gradient of a quadratic function. Instead, gradients of $f$ can be written as a sum of the gradients of a quadratic function with a perturbation term $e(x)$, as follows:
\BEQ
    \nabla f(x) = \HH(x-x_\star) + e(x), \quad \text{with $\HH$ : } \textbf{0} \prec \mu \idm \preceq \HH \preceq L\idm. \label{eq:model_gradient}
\EEQ
Indeed, it follows from twice continuous differentiability of $f$ that 
\[
    f(x)=f_\star + \underbrace{\langle \nabla f(x_\star)}_{=0} ; x-x_\star\rangle + \frac12 \langle x-x_\star; \nabla^2 f(x_\star)(x-x_\star)\rangle + O(\|x-x_\star\|_2^3).
\]
Therefore, $f(x)$ can be approximated by the quadratic function~\eqref{eq:model_gradient} with $\HH=\nabla^2f(x_\star)$. Similarly, its gradient reads
\[
    \nabla f(x) = \underbrace{\nabla f(x_\star)}_{=0} + \nabla^2 f(x_\star) (x-x_\star) + e(x) = \HH(x-x_\star) + e(x),
\]
where $e(x)$ is the first-order Taylor remainder of the gradient. Thus, minimizing a non-quadratic function is equivalent to minimizing a perturbed quadratic one with a second-order error on its gradient:
\begin{equation} \label{eq:error_order_twice_diff_functions}
    e(x_k) =  \nabla f(x_k) - \HH (x_k-x_\star) \quad \left(\Rightarrow \;\; \|e(x_k)\|_2= O (\|x_k-x_\star\|_2^2) \right),
\end{equation}
where $\HH = \nabla^2 f(x_\star)$.

\subsection{Convergence Bound}

Using a perturbation argument, it is possible to derive a convergence guarantee for RNA. We state here a simplified version of \citep[Theorem 3.2]{scieur2018online}, which describes how regularization balances acceleration and stability in Algorithm~\ref{algo:rna_simplified}. We discuss convergence rates in greater detail in what follows.

\begin{theorem} \label{thm:convergence_rna}
Let $x_0\in\mathbb{R}^d$, $f$ be a twice continuously differentiable function with minimum $x_\star$ and whose gradient $\nabla f$ follows~\eqref{eq:model_gradient} with $\|e(x_i)\|_2 \leq \epsilon$. Let $\{x_i\}_{i=0,\ldots,k}$ be generated by a nondegenerate first-order method initiated from $x_0$. Then, for any $k\geq 0$, it holds that
\begin{align*}
        &\|\HH(x_{\text{extr}}-x_{\star})\|_2 \\
        &\leq \|\idm-h\HH\|_2 \left( \underbrace{V_k^{[\mu,L]}(\lambda)\|\HH(x_0-x_{\star})\|_2}_{\textbf{acceleration}} + \underbrace{O\left( \sqrt{1+\frac{1}{\lambda}} \epsilon\right)}_{\textbf{stability}} \right),
\end{align*}
    where $x_{\text{extr}}$ is the output of Algorithm \ref{algo:rna_simplified} applied to~$\{ (x_i,\,\nabla f(x_i)) \}_{i=0,\ldots, k}$ with parameters $h,\,\lambda>0$, and $V_k^{[\mu,L]}(\lambda)$ is a constant that corresponds to the maximum value on the interval $[\mu, L]$ of the \textit{regularized Chebyshev polynomial}, i.e,
    \begin{equation} \label{regularized_cheby}
    \begin{aligned}
        V_k^{[\mu,L]}(\lambda) =& \max_{x \in [\mu,L]} |C_k^{[\mu,L],\lambda}(x)|,\\
        \text{with }& C_k^{[\mu,L],\lambda} =\argmin_{\substack{P\in\mathcal{P}_k,\\ P(0)=1}} \;\; \max_{x \in [\mu,L]} P^2(x) + \lambda \| \GG^T\GG \|_2 \|P\|_2^2,
	\end{aligned}
    \end{equation}
    where $\|P\|_2$ is the norm of the vector of coefficients of the polynomial~$P$.
\end{theorem}
This theorem states that regularization helps stabilizing the algorithm while slowing down the convergence rate. The regularized Chebyshev polynomial is somehow a mid-point between the classical shifted Chebyshev polynomial $C_k^{[\mu,L]}$ (from \eqref{eq:shifted_cheby_definition}) and the polynomial whose coefficients are defined by $\ones/k$ (the polynomial that averages the iterates $\{x_i\}_{i=0,\ldots,k}$). By construction, its maximum value is always larger than that of the Chebyshev polynomial, but the norm of its coefficients is smaller. Unfortunately, there is as yet no known explicit expression of the regularized Chebyshev polynomial. To the best of our knowledge, its value can nevertheless be computed numerically~\citep{Barr20}.

This mid-point between Chebyshev coefficients and the simple averaging of iterates is also natural in the context of noisy iterations. When the noise is negligible, a small regularization parameter combines the iterates $\{x_i\}_{i=0,\ldots, k}$ using nearly the classical Chebyshev weights. When the noise is more substantial, a larger regularization parameter brings the vector of coefficients $c$ closer to the average $\ones/k$, thereby improving the ``stability term'' while rendering the ``acceleration'' less effective.

\subsection{Asymptotic Convergence Rate} 

We briefly discuss the behavior of RNA when the initial point $x_0$ approaches the solution $x_{\star}$. In particular, the next proposition shows that if the perturbation magnitude $\epsilon$ decreases faster than $\| \HH(x_0-x_{\star}) \|_2$, the parameter $\lambda$ can be adjusted to ensure an asymptotic convergence rate comparable to that of the Chebyshev method on quadratic problems (see Section~\ref{c-Cheb}). 

Informally, the theorem exploits the fact that as $x_0$ approaches $x_\star$, $f$ gets closer to its quadratic approximation around $x_\star$. Thereby, an appropriate tuning of RNA allows matching (asymptotically) the convergence rate of nonlinear acceleration on quadratics (see Theorem~\ref{thm:rate_nonlinear_acc}).

\begin{proposition} \label{prop:asymptotic_rna}
    Let $x_0\in\mathbb{R}^d$, $f$ be a twice continuously differentiable function with minimum $x_\star$, whose gradient $\nabla f$ follows~\eqref{eq:model_gradient} with $\|e(x_i)\|_2 \leq \epsilon$. Let $\{x_i\}_{i=0,\ldots,k}$ be generated by a nondegenerate first-order method initiated from $x_0$. If we have
    \[
       \epsilon = O(\|x_0-x_\star\|_2^\alpha), \quad \alpha > 1,
    \] 
    and if we set $\lambda \propto \|x_0-x_\star\|_2^{s}$ (proportional to $\|x_0-x_\star\|_2^{s}$), where $0<s<2(\alpha-1)$, then it holds that
    \[
        \lim_{x_0\rightarrow x_\star} \frac{\|\HH(x_{\text{extr}}-x_{\star})\|_2}{\|\HH(x_0-x_{\star})\|_2} \leq \|\idm-h\HH\|_2  \frac{2}{\xi^k + \xi^{-k}} \quad \text{where }\; \xi = \frac{\sqrt{\frac{L}{\mu}}+1}{\sqrt{\frac{L}{\mu}}-1},
    \]
    where $x_{\text{extr}}$ is the output of Algorithm \ref{algo:rna_simplified} applied to the sequence  $\{ (x_i,\nabla f(x_i)) \}_{i=0,\ldots, k}$ with parameters $h,\,\lambda>0$.
\end{proposition}
\begin{proof}
    To simplify the notation, set $ R \triangleq \|\HH(x_0-x_{\star})\|_2$. We start from the result of Theorem \ref{thm:convergence_rna} and divide both sides by $R$:
    \[
        \frac{\|\HH(x_{\text{extr}}-x_{\star})\|_2}{R} \leq \|\idm-h\HH\|_2 \left( V_k^{[\mu,L]}(\lambda) + O\left( \sqrt{1+\frac{1}{\lambda}} \frac{\epsilon}{R}\right) \right).
    \]
    Since $\lambda\propto R^s$ and $\epsilon = O(R^\alpha)$,
    \begin{align*}
        &\frac{\|\HH(x_{\text{extr}}-x_{\star})\|_2}{R} \\
        &\quad\leq \|\idm-h\HH\|_2 \left( V_k^{[\mu,L]}(\lambda) + \sqrt{O(R^{2(\alpha-1)}) +  O(R^{2(\alpha-1)-s})}\right).
    \end{align*}
    When $x_0\rightarrow x_\star$, we have $R\rightarrow 0$ and
    \begin{align*}
        R^{2(\alpha-1)} &\rightarrow 0 & \text{(since $\alpha > 1$),}\\
        R^{2(\alpha-1)-s} &\rightarrow 0 & \text{(since $s < 2(\alpha-1)$)}.
    \end{align*}
    Finally, in \eqref{regularized_cheby} the regularization parameter $\lambda \|\GG^T\GG\|_2 = O(R^2) \rightarrow 0$. Since the (non-regularized) shifted Chebyshev polynomial $C_k^{[\mu,L]}=C_k^{[\mu,L],0}$ is a feasible solution of \eqref{regularized_cheby}, we have the following bounds:
    \begin{align*}
        \left(V_k^{[\mu,L]}(\lambda)\right)^2 & = \max_{x \in [\mu,L]} |C_k^{[\mu,L],\lambda}(x)|^2,\\
        & \leq \left\{\max_{x \in [\mu,L]} |C_k^{[\mu,L],\lambda}(x)|^2\right\} + \lambda \| \GG^T\GG \|_2 \| C_k^{[\mu,L],\lambda}\|_2^2,\\
        & = \min_{\substack{P\in\mathcal{P}_k,\\ P(0)=1}} \;\; \max_{x \in [\mu,L]} P^2(x) + \lambda \| \GG^T\GG \|_2 \|P\|_2^2,\\
        & \leq \left\{\max_{x \in [\mu,L]} |C_k^{[\mu,L]}(x)|^2\right\} + \lambda \| \GG^T\GG \|_2 \| C_k^{[\mu,L]}\|_2^2.
    \end{align*}
    As $\lambda \rightarrow 0$ we have that the upper bound on $(V_k^{[\mu,L]}(\lambda)$ converges to the maximum value of the regular (shifted) Chebyshev polynomial, thereby reaching the desired claim.
\end{proof}

In short, the previous theorem states that the asymptotic convergence rate matches the rate of Chebyshev's method as soon as $\lambda \propto \|x_0-x_\star\|_2^s$ is decreasing (condition $s>0$), but not too quickly compared to the perturbation magnitude $\epsilon=O(\|x_0-x_\star\|_2^\alpha)$ (condition $s<2(\alpha-1)$), which is achievable only when $\alpha>1$. This condition is met, for instance, when accelerating twice continuously differentiable functions with gradient descent: the error decreases as $O(\|x_0-x_{\star}\|_2^2)$, see \eqref{eq:error_order_twice_diff_functions}, and therefore $\alpha = 2>1$.

\section{Extensions}

The previous sections presented the nonlinear acceleration mechanism for unconstrained convex quadratic minimization. It also contained an analysis of its regularized version when applied beyond quadratics. In this section, we briefly cover two natural extensions: (i)~the application of nonlinear acceleration to iterates that are corrupted by a stochastic noise, and (ii)~the application of nonlinear acceleration to constrained/composite convex optimization problems when a projection/proximal operator is used.

\paragraph{Stochastic gradients.}
In the common situation where the first-order method under consideration only has access to stochastic estimates $\tilde \nabla f(x)$ of the gradient (satisfying $\mathbb{E}[\tilde \nabla f(x_i)] = \nabla f(x_i)$) one can adapt the perturbation model
\[
    \tilde \nabla f(x_i) = \HH(x_i-x_{\star}) + e_i,
\]
to $e_i$ being the sum of a stochastic noise with a Taylor remainder. This is typically the case when applying RNA to stochastic gradient descent (SGD) and related methods. In this case, Theorem~\ref{thm:convergence_rna} holds in expectation under standard assumptions~\citep{Scie17}, such as a bounded variance of $e_i$. However, the asymptotic convergence result from Proposition~\ref{prop:asymptotic_rna} may not be achieved. Indeed, in this setting, Proposition~\ref{prop:asymptotic_rna} also holds in expectation under the condition
\[
    \mathbb{E}[\|e(x_i)\|_2] = O(\|x_0-x_\star\|_2^\alpha), \quad \alpha > 1.
\]
Unfortunately, an asymptotic acceleration is not always possible. For instance, when trying to accelerate the fixed step SGD, we have $\epsilon = O(1)$ (i.e., $\alpha = 0$) and the asymptotic guarantee does not apply. This is probably not a surprise as this SGD does not converge to the optimum, hence there is no apparent reason for any sequence extrapolation technique to work at all. Fortunately, Algorithm~\ref{algo:rna_general} does usually work for ``variance reduced'' first-order methods~\citep{Scie17}, such as SAG~\citep{schmidt2017minimizing}, SAGA~\citep{defazio2014saga}, or SVRG \citep{johnson2013accelerating}.

\paragraph{Nonlinear acceleration of the proximal gradient method.} It is common to apply first-order methods to composite convex minimization problems of the form:
\BEQ \label{eq:composite_optim}
    \min_{x\in\mathbb{R}^d} \{ F(x) \defeq f(x) + h(x)\},
\EEQ
where $f$ is a smooth strongly convex function  (this class of functions is used intensively in Section~\ref{c-Nest}; see Definition~\ref{def:smoothstrconvex}) and $h$ is a closed, proper, and convex function (i.e., $h$ has a closed, non-empty, and convex epigraph) and whose \emph{proximal operator} is available:
\BEQ \label{eq:proximal_operator}
    \prox_{\gamma h}(x) \triangleq \argmin_{z} \left\{\gamma h(z) + \frac{1}{2}\|x-z\|_2^2\right\}
\EEQ
for some step size  $\gamma > 0$. Problem~\eqref{eq:composite_optim} can then be approached iteratively via the proximal gradient method:
\BEQ
    x_{k+1} = \prox_{\gamma h}( x_k - \gamma\nabla f(x_k)). \label{eq:proximal_gradient}
\EEQ
We omit most of the details on proximal algorithms; see Section~\ref{c-Nest} and Section~\ref{c-prox} for more details and references. For instance, when $h$ is the indicator function of a non-empty closed convex set $\mathcal{C}$, the proximal operator corresponds to an orthogonal projection onto $\mathcal{C}$ and the proximal gradient method reduces to the projected gradient method.

Unfortunately, a naive use of nonlinear acceleration on the iterates $\{x_i\}_{i=0,\ldots,k}$ does not immediately work in this context, for several reasons. In particular, it is not possible to ensure that the extrapolated point $x_{\text{extr}}$ belongs to $\dom(h)$ (or to the set $\mathcal{C}$ when $h$ is an indicator function for $\mathcal{C}$). Moreover, due to the use of the proximal operator, the iterates $\{x_{i}\}_{i=0,\ldots,k}$ do not necessarily satisfy the span assumption~\eqref{eq:first_order_algorithm}.

Recently, \citet{mai2019anderson} adapted the Anderson Acceleration method to handle a large class of constrained and non-smooth composite problems. The main idea is as follows: instead of accelerating the sequence $\{x_i\}_{i=0,\ldots,k}$ generated by the proximal gradient method~\eqref{eq:proximal_gradient}, we accelerate an alternate sequence $\{z_i\}_{i=0,\ldots,k}$ which satisfies
\[
    z_{i+1} = \prox_{\gamma h} (z_{i}) - \gamma \nabla f\left( \prox_{\gamma h} (z_{i}) \right), \quad z_1 = x_0-\gamma \nabla f(x_0).
\]
This sequence $\{z_i\}_{i=0,\ldots,k}$ corresponds to the sequence generated by~\eqref{eq:proximal_gradient} with the ordering of the gradient and proximal steps being swapped.

This trick allows obtaining convergence bounds for nonlinear acceleration in the proximal setup under very few changes in the algorithm. In particular, \citet{mai2019anderson} show that using Algorithm~\ref{algo:rna_prox} in the presence of a proximal operator does not change the convergence analysis---using  Clarke's generalized Jacobian~\citep{clarke1990optimization}, semi-smoothness~\citep{mifflin1977semismooth,qi1993nonsmooth} and assuming that the function $h$ is twice \emph{epi-differentiable} and that $h$ is twice-differentiable around the solution $x_\star$. We refer the reader to~\citep[Section 13]{rockafellar2009variational} for a comprehensive treatment of epi-differentiability.

\begin{algorithm}[ht] 
  \caption{Online regularized nonlinear acceleration with a proximal operator}
  \label{algo:rna_prox}
  \begin{algorithmic}[1]
    \REQUIRE Differentiable function $f$, closed proper convex function $h$ with proximal operator available, initial point $x_0$, step size $\gamma$, regularization term $\lambda$ and reference vector $c_{\text{ref}}$.
    \STATE \textbf{Initialize} $z_1 = x_0-\gamma \nabla f(x_0)$, $x_1 = \prox_{\gamma h}(z_1)$, empty sequence $S$ of pairs iterate/gradient.
    \FOR{$k=1\ldots$}  
        \STATE Compute $g_k = \frac{\gamma\nabla f(x_k)+z_{k}-x_k}{\gamma}$, then append the pair $S\leftarrow S \cup \{(z_k,\, g_k)\}$.
        \STATE Compute the extrapolation $z_{k+1} = \text{[Algorithm~\ref{algo:rna_general}]}(\mathcal{S}, \gamma, \lambda,c_{\text{ref}})$.
        \STATE $x_{k+1} = \prox_{\gamma h}(z_{k+1})$.
        \ENDFOR
    \ENSURE Approximate solution $x_{k+1}$.
  \end{algorithmic}
\end{algorithm}

\section{Globalization Strategies and Speeding-up Heuristics}

As for many standard optimization methods, such as quasi-Newton methods, RNA only has local convergence guarantees beyond quadratics. Therefore, it is common to embed the mechanism with some globalization strategies, a.k.a. safeguards. Those strategies ensure not to deteriorate too much the performance of the initial first-order method in cases where RNA is used beyond its guaranteed range of applications. Those strategies can also be seen as speeding-up heuristics.

\paragraph{Descent condition.} It is in general not guaranteed that the extrapolated point $x_{\text{extr}}$ is better than any iterate of the sequence $\{x_i\}_{i=0,\ldots,k}$ produced by the original method. This situation might for example occur when extrapolating with a bad mixing or regularization parameter, or simply when the error terms are too large. One classical way of limiting the impact of such problems is by checking some descent condition. For instance, one might consider ``accepting'' $x_{\text{extr}}$ only if it is better than previous iterates $\{x_i\}_{i=0,\ldots,k}$:
\[
    f\left(x_{\text{extr}}\right) < \min_{i\in\{ 0,\ldots, k\}} f(x_i),
\]
and to discard it otherwise.

\paragraph{Line-search.} Nonlinear acceleration requires the selection of a mixing parameter, which might be difficult to tune in practice. One common trick is to choose it via a line-search strategy. That is, defining:
\[
    x_{\text{extr}}(h) = \sum_{i=0}^k \left(c_i x_i - h\nabla f(x_i)\right),
\]
one can choose $h$ by approximately solving $\argmin_{h} f\left(  x_{\text{extr}}(h) \right)$.

\section{Notes and References}\label{s:relatedw}
Nonlinear acceleration techniques have been studied extensively during recent decades, and excellent reviews can be found in \citep{smith1987extrapolation,jbilou1991some,brezinski1991extrapolation,jbilou1995analysis,jbilou2000vector,brezinski2001convergence,brezinski2019genesis}. The first usage of an acceleration technique for fixed point iteration can be traced back to \citep{gekeler1972solution,brezinski1971algorithme,brezinski1970application}.

There are numerous independent works leading to methods similar to those described here. The most classical, and probably the most similar, is \textit{Anderson acceleration} \citep{anderson1965iterative}, which corresponds exactly to the online mode of nonlinear acceleration (without regularization). Despite it being an old algorithm, there has been a recent uptake of interest in the convergence analysis \citep{walker2011anderson, toth2015convergence} of Anderson acceleration thanks to its good empirical performance, and strong connection with quasi-Newton methods \citep{fang2009two}.

Other versions of nonlinear acceleration use different arguments but behave similarly. For instance, minimal polynomial extrapolation (MPE), which uses the properties of the minimal polynomial of a matrix \citep{cabay1976polynomial}; reduced rank extrapolation (RRE); and the Mesina method \citep{mevsina1977convergence,eddy1979extrapolating} are also variants of Anderson acceleration. The properties and equivalences of these approaches have been studied extensively during the past decades \citep{sidi1988extrapolation,ford1988recursive,sidi1991efficient,jbilou1991some,sidi1998upper,sidi2008vector,sidi2017minimal,sidi2017vector,brezinski2018shanks,brezinski2020shanks}. Unfortunately, these methods do not extend well to nonlinear functions, especially due to conditioning problems \citep{sidi1986convergence,sidi1988convergence,scieur2016regularized}. Recent works have nevertheless proven the convergence of such methods, provided that good conditioning of the linear system \citep{sidi2019convergence} can be ensured.

There are also other classes of nonlinear acceleration algorithms, based on existing algorithms, for accelerating the convergence of scalar sequences \citep{brezinski1975generalisations}. For instance, the topological epsilon vector algorithm (TEA) extends the idea of the scalar $
\varepsilon$-algorithm of \citep{Wynn56} to vectors.

\chapter{Nesterov Acceleration}\label{c-Nest} 

\enlargethispage{-\baselineskip}
This section presents a systematic interpretation of the acceleration of the gradient method stemming from Nesterov's original work~\citep{Nest83}. The early parts of the section are devoted to the gradient method and the ``optimized gradient method,'' due to~\citet{Dror14} and~\citet{kim2016optimized}. The motivations and ideas underlying the latter are intuitive and very similar to those behind the introduction of Chebyshev methods for optimizing quadratic functions (see Section~\ref{c-Cheb}). Furthermore, the optimized gradient method has a relatively simple format and proof and can be used as an inspiration for developing numerous variants with wider ranges of applications, including Nesterov's early accelerated gradient methods~\citep{Nest83,Nest13} and the fast iterative shrinkage-thresholding algorithm~\citep[FISTA]{Beck09}. Although some parts of this section are more technical, we believe all the ideas can be reasonably well understood even when skipping, or skimming through the algebraic proofs. The section and the proofs are organized so that each time an additional ingredient (strong convexity, constraints, etc.) is included, its inclusion only requires a few additional ingredients compared to the previous (simpler) proofs of the base versions of the method. 

We start with the theory and interpretation of acceleration in a simple setting: smooth unconstrained convex minimization in a Euclidean space. All subsequent developments follow from the same template, namely a linear combination of regularity inequalities, with additional ingredients being added one by one. The next part is devoted to methods that take advantage of strong convexity by using the same ideas and algorithmic structures. On the way, we provide a few different (equivalent) templates for the algorithms, since in more advanced settings, those templates do not generalize in the same way. We then recap and discuss a few practical extensions for handling constrained problems, nonsmooth regularization terms, unknown problem parameters/line-searches, and non-Euclidean geometries. Finally, we briefly discuss a popular ordinary differential equation (ODE)-based interpretation of Nesterov's method. Techniques for obtaining the worst-case analyses presented throughout this text are presented in Appendix~\ref{a-WC_FO}, and notebooks for simpler reproduction of the proofs are provided in Section~\ref{s:notes_ref_chapt_nest}.

\section{Introduction}
In the first part of this section, we consider smooth unconstrained convex minimization problems. This type of problems is a direct extension of unconstrained convex quadratic minimization problems where the quadratic function has eigenvalues bounded above by some constant. More precisely, we consider the simple unconstrained differentiable convex minimization problem
\begin{equation}
f_\star=\min_{x\in\mathbb{R}^d} f(x),\label{eq:smoothcvx}
\end{equation}
where $f$ is convex with an $L$-Lipschitz gradient (we call such functions convex and $L$-smooth, see Definition~\ref{def:smoothstrconvex} below), and we assume throughout that there exists a minimizer~$x_\star$. The goal of the methods presented below is to find a candidate solution $x$ satisfying $f(x)-f_\star\leq \epsilon$ for some $\epsilon>0$. Depending on the target application, other quality measures, such as guarantees on $\| \nabla f(x)\|_2$ or $\|x-x_\star\|_2$, might be preferred. We refer to Section~\ref{s:notes_ref_chapt_nest} ``Changing the performance measure'', for discussions on this topic.

We start with the analysis of gradient descent and then show that its iteration complexity can be significantly improved using an acceleration technique proposed by~\citet{Nest83}. 

After presenting the theory for the smooth convex case we see how it goes in the smooth strongly convex one. This class of problems extends to that of unconstrained convex quadratic minimization problems where the quadratic function has eigenvalues respectively bounded above and below by some constants~$L$ and~$\mu$.
\begin{definition}\label{def:smoothstrconvex} Let $0\leq\mu< L<+\infty$. A continuously differentiable function $f:\mathbb{R}^d\rightarrow\mathbb{R}$ is $L$-smooth and $\mu$-strongly convex (denoted $f\in\mathcal{F}_{\mu,L}$) if and only if 
\begin{itemize}
    \item ($L$-smoothness) for all $x,y\in\mathbb{R}^d$, it holds that
    \begin{equation}\label{eq:smooth}
    f(x)\leq f(y)+\langle \nabla f(y);x-y\rangle+\frac{L}{2}\|x-y\|^2_2,
    \end{equation}
    \item ($\mu$-strong convexity)  for all $x,y\in\mathbb{R}^d$, it holds that
    \begin{equation}\label{eq:strcvx}
    f(x)\geq f(y)+\langle \nabla f(y);x-y\rangle+\frac{\mu}{2}\|x-y\|^2_2.
    \end{equation}
\end{itemize} Furthermore, we denote by $\cond=\tfrac{\mu}{L}$ the inverse condition number (that is,  $\cond=\tfrac{1}{\kappa}$, with $\kappa$ is the usual condition number as used e.g., in Section~\ref{c-Cheb}) of functions in the class $\mathcal{F}_{\mu,L}$.
\end{definition}\paragraph{Notation.}We use the notation $\mathcal{F}_{0,L}$ for the set of smooth convex functions. By extension, we use $\mathcal{F}_{0,\infty}$ for the set of (possibly non-differentiable) proper closed convex functions (i.e., convex functions whose epigraphs are non-empty closed convex sets). Finally, we denote by $\partial f(x)$ the subdifferential of $f$ at $x\in\mathbb{R}^d$ and by $g_f(x)\in\partial f(x)$ a particular subgradient of $f$ at~$x$.

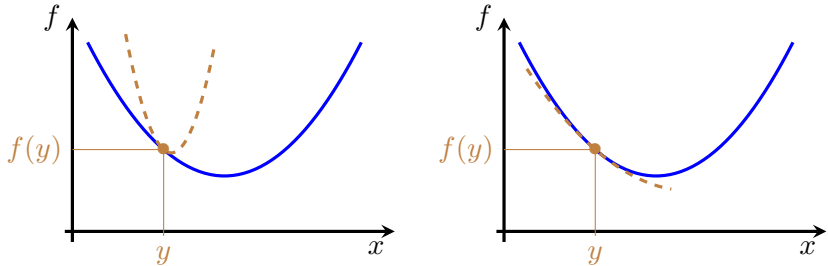
\begin{figure}[!ht]
    \centering
        \begin{tabular}{cc}
            \begin{tikzpicture}[yscale=0.35,xscale=0.2]
                \draw[very thick, black, -stealth] (-10.5,-0.1) -- (11.2,-0.1);
		        \draw (11.2,-0.1) node[below left]{$x$};
		        \draw[very thick, black, -stealth] (-10,-0.5) -- (-10,8);
		        \draw (-10.1,8) node[left]{$f$};
		        \draw(0,0) [very thick, color=blue, domain=-3*3:3*3, samples=50] plot({\x}, {2+(\x/4)^2});
		        \draw(0,0) [dashed,very thick, color=brown, domain=-6.5:-0.7, samples=50] plot({\x}, {2+1-1/2*(\x+4)+1/2*(\x+4)^2});
		        \node [color=brown] (poi) at (-4,3) {$\bullet$};
		        \node [color=brown, dashed] (poix) at (-4,-1) {$y$};
		        \node [color=brown, dashed] (poiy) at (-12.5,3) {$f(y)$};
		        \draw [color=brown] (poi.center) -- (poiy.east);
		        \draw [color=brown] (poi.center) -- (poix.north);
            \end{tikzpicture}
             &
            \begin{tikzpicture}[yscale=0.35,xscale=0.2]
                \draw[very thick, black, -stealth] (-10.5,-0.1) -- (11.2,-0.1);
		        \draw (11.2,-0.1) node[below left]{$x$};
		        \draw[very thick, black, -stealth] (-10,-0.5) -- (-10,8);
		        \draw (-10.1,8) node[left]{$f$};
		        \draw(0,0) [very thick, color=blue, domain=-3*3:3*3, samples=50] plot({\x}, {2+(\x/4)^2});
		        \draw(0,0) [dashed,very thick, color=brown, domain=-8.5:1, samples=50] plot({\x}, {2+1-1/2*(\x+4)+((\x+4)/5)^2});
		        \node [color=brown] (poi) at (-4,3) {$\bullet$};
		        \node [color=brown, dashed] (poix) at (-4,-1) {$y$};
		        \node [color=brown, dashed] (poiy) at (-12.5,3) {$f(y)$};
		        \draw [color=brown] (poi.center) -- (poiy.east);
		        \draw [color=brown] (poi.center) -- (poix.north);
            \end{tikzpicture}
        \end{tabular}
    \caption{Let $f(\cdot)$ (blue) be a differentiable function. (Left)~Smoothness: $f(\cdot)$ (blue) is $L$-smooth if and only if it is upper bounded by $ f(y)+\langle \nabla f(y);.-y\rangle+\tfrac{L}{2}\|.-y\|^2_2$ (dashed, brown) for all $y$. (Right)~Strong convexity: $f(\cdot)$ (blue) is $\mu$-strongly convex if and only if it is lower bounded by $ f(y)+\langle \nabla f(y);.-y\rangle+\tfrac{\mu}{2}\|.-y\|^2_2$ (dashed, brown) for all $y$.}
    \label{fig:smooth_strcvx}
\end{figure}

\paragraph{Smooth strongly convex functions.}
Figure~\ref{fig:smooth_strcvx} provides an illustration of the global quadratic {\em upper approximation}  (with curvature $L$) on $f(\cdot)$ due to smoothness and of the global quadratic {\em lower approximation} (with curvature $\mu$) on $f(\cdot)$ due to strong convexity.

A number of inequalities can be written to characterize functions in $\mathcal{F}_{\mu,L}$: see, for example,~\citet[Theorem 2.1.5]{Nest03a}. When analyzing methods for minimizing functions in this class, it is crucial to have the \emph{right} inequalities at our disposal, as worst-case analyses essentially boil down to appropriately combining such inequalities. We provide the most important inequalities along with their interpretations and proofs in Appendix~\ref{a-inequalities}. In this section, we only use three. First, we use the quadratic upper and lower bounds arising from the definition of smooth strongly convex functions, that is,~\eqref{eq:smooth} and~\eqref{eq:strcvx}. For some analyses however, we need an additional inequality, provided by the following theorem. This inequality is often referred to as an \emph{interpolation} (or \emph{extension}) inequality. Its proof is relatively simple: it only consists of requiring all quadratic lower bounds from~\eqref{eq:strcvx} to be below all quadratic upper bounds from~\eqref{eq:smooth} (details in Appendix~\ref{s:ineq_eucl}). It can be shown that worst-case analyses of all first-order methods for minimizing smooth strongly convex functions can be performed using \emph{only} this inequality for some specific values of $x$ and $y$ (details in Appendix~\ref{a-WC_FO}).

\clearpage
\begin{theorem} Let $f$ be a continuously differentiable function. $f$ is $L$-smooth and $\mu$-strongly convex (possibly with $\mu=0$) if and only if for all $x,y\in\mathbb{R}^d$, it holds that

\begin{equation}\label{eq:interpolation_strcvx}
\begin{aligned}
f(x)\geq f(y)&+\langle \nabla f(y);x-y\rangle+\frac{1}{2L}\|\nabla f(y)-\nabla f(x)\|^2_2\\&+\frac{\mu L}{2(L-\mu)}\|x-y-\tfrac1L(\nabla f(x)-\nabla f(y))\|^2_2.
\end{aligned}
\end{equation}
\end{theorem}
As discussed later in Section~\ref{sec:LO_OGM}, this inequality has some flaws. Therefore, we only use~\eqref{eq:smooth} and~\eqref{eq:strcvx} whenever possible.

Before continuing to the next section, we mention that both smoothness and strong convexity are strong assumptions. More generic assumptions are discussed in Section~\ref{c-restart} to obtain improved rates under weaker assumptions.

\section{Gradient Method and Potential Functions}
In this section, we analyze gradient descent using the concept of {\em potential functions}. The resulting proofs are technically simple, although they might not seem to provide any direct intuition about the method at hand. We use the same ideas to analyze a few improvements on gradient descent before providing interpretations underlying this mechanism.

\subsection{Gradient Descent}
The simplest and probably most natural method for minimizing differentiable functions is gradient descent. It is often attributed to~\citet{cauchy1847methode} and consists of iterating
\[ x_{k+1}=x_k-\gamma_k \nabla f(x_k),\]
where $\gamma_k$ is some step size. There are many different techniques for picking $\gamma_k$, the simplest of which is to set $\gamma_k=1/L$, assuming $L$ is known---otherwise, line-search techniques are typically used; see Section~\ref{sec:nest_practice}. Our present objective is to bound the number of iterations required by gradient descent to obtain an approximate minimizer $x_k$ of $f$ that satisfies $f(x_k)-f_\star\leq \epsilon$. 

\subsection{A Simple Proof Mechanism: Potential Functions}\label{sec:GM_pot}
Potential (or Lyapunov/energy) functions are classical tools for proving convergence rates in the first-order literature, and a nice recent review of this topic is given by~\citet{bansal2019potential}. For gradient descent, the idea consists in recursively using a simple inequality (proof below),
\[
(k+1)(f(x_{k+1})-f_\star)+\frac{L}2\|x_{k+1}-x_\star\|^2_2\leq k(f(x_{k})-f_\star)+\frac{L}2\|x_{k}-x_\star\|^2_2,
\]
that is valid for all $f\in\mathcal{F}_{0,L}$ and all $x_k\in\mathbb{R}^d$ when $x_{k+1}=x_{k}-\tfrac1L \nabla f(x_{k})$. In this context, we refer to 
\[\phi_k\defeq k(f(x_{k})-f_\star)+\frac{L}2\|x_{k}-x_\star\|^2_2\]
as a potential and use $\phi_{k+1}\leq\phi_k$ as the building block for the worst-case analysis. Once such a {\em potential inequality} $\phi_{k+1}\leq\phi_k$ is established, a worst-case guarantee can easily be deduced through a recursive argument, yielding
\begin{equation}\label{eq:chained_pot}
N(f(x_N)-f_\star)\leq \phi_N\leq \phi_{N-1}\leq\hdots\leq \phi_0=\frac{L}{2}\|x_0-x_\star\|^2_2,
\end{equation}
and hence, $f(x_N)-f_\star\leq \tfrac{L}{2N}\|x_0-x_\star\|^2_2$. We also conclude that the worst-case accuracy of gradient descent is $O(N^{-1})$ or equivalently, that its iteration complexity is $O(\epsilon^{-1})$. Therefore, the main inequality to be proved for this worst-case analysis to work is the \emph{potential} inequality $\phi_{k+1}\leq \phi_k$. In other words, the analysis of $N$ iterations of gradient descent is reduced to the analysis of a single iteration, using an appropriate potential. This kind of approach was already used for example by~\citet{Nest83}, and many different variants of the potential function can be used to prove convergence of gradient descent and related methods in similar ways.

\begin{theorem}\label{thm:pot_GM}
Let $f$ be an $L$-smooth convex function, $x_\star\in\mathrm{argmin}_x\, f(x)$, and $k\in\mathbb{N}$. For any $A_k\geq 0$ and $x_k\in\mathbb{R}^d$, it holds that
\begin{equation*}
\begin{aligned}
A_{k+1}(f(x_{k+1})-f_\star)&+\frac{L}{2}\|x_{k+1}-x_\star\|^2_2\\
&\leq A_{k}(f(x_{k})-f_\star)+\frac{L}{2}\|x_{k}-x_\star\|^2_2,
\end{aligned}
\end{equation*}
with $x_{k+1}=x_k-\tfrac1L \nabla f(x_k)$ and $A_{k+1}=1+A_k$.
\end{theorem}
\begin{proof} The proof consists of performing a weighted sum of the following inequalities:
\begin{itemize}
    \item convexity of $f$ between $x_k$ and $x_\star$, with weight $\lambda_1=A_{k+1}-A_k$:
    \[ 0\geq f(x_k)-f_\star+\langle \nabla f(x_k); x_\star-x_k\rangle,\]
    \item smoothness of $f$ between $x_k$ and $x_{k+1}$ with weight $\lambda_2=A_{k+1}$:
    \[ 0\geq f(x_{k+1})-\left(f(x_k)+\langle \nabla f(x_k);x_{k+1}-x_k\rangle+\frac{L}{2}\lVert x_k-x_{k+1}\rVert^2_2\right).\]
    The last inequality is often referred to as the \emph{descent lemma} since substituting $x_{k+1}$ allows to obtain $f(x_{k+1})\leq f(x_k)-\frac1{2L}\|\nabla f(x_k)\|^2_2$.
\end{itemize}
The weighted sum forms a valid inequality:
\begin{equation*}
\begin{aligned}
0\geq &\lambda_1[f(x_k)-f_\star+\langle \nabla f(x_k); x_\star-x_k\rangle]\\
&+\lambda_2[f(x_{k+1})-(f(x_k)+\langle \nabla f(x_k);x_{k+1}-x_k\rangle +\frac{L}{2}\lVert x_k-x_{k+1}\rVert^2_2)].
\end{aligned}
\end{equation*}
Using $x_{k+1}=x_k-\tfrac1L \nabla f(x_k)$, this inequality can be rewritten (by completing the squares or simply extending both expressions and verifying that they match on a term-by-term basis) as follows:
\begin{equation*}
\begin{aligned}
0\geq &(A_k+1)(f(x_{k+1})-f_\star)+\frac{L}{2}\lVert x_{k+1}-x_\star\rVert^2_2\\&- A_k(f(x_{k})-f_\star)-\frac{L}{2}\lVert x_{k}-x_\star\rVert^2_2+\frac{A_{k+1}-1}{2L}\lVert \nabla f(x_k)\rVert^2_2\\&- (A_{k+1}-A_k-1)\langle \nabla f(x_k);x_k-x_\star\rangle,
\end{aligned}
\end{equation*}
which can be reorganized and simplified to
\begin{equation*}
\begin{aligned}
(A_k&+1)(f(x_{k+1})-f_\star)+\frac{L}{2}\lVert x_{k+1}-x_\star\rVert^2_2\\\leq& A_k(f(x_{k})-f_\star)+\frac{L}{2}\lVert x_{k}-x_\star\rVert^2_2-\frac{A_{k+1}-1}{2L}\lVert \nabla f(x_k)\rVert^2_2\\&\quad +(A_{k+1}-A_k-1)\langle \nabla f(x_k);x_k-x_\star\rangle \\\leq& A_k(f(x_{k})-f_\star)+\frac{L}{2}\lVert x_{k}-x_\star\rVert^2_2,
\end{aligned}
\end{equation*}
where the last inequality follows from picking $A_{k+1}=A_k+1$ and neglecting the last residual term $-\tfrac{A_k}{2L}\|\nabla f(x_k)\|^2_2$ (which is nonpositive) on the right-hand side.
\end{proof}

A convergence rate for gradient descent can be obtained directly as a consequence of Theorem~\ref{thm:pot_GM}, following the reasoning of~\eqref{eq:chained_pot}, and the worst-case guarantee corresponds to $f(x_N)-f_\star=O(A_N^{-1})=O(N^{-1})$. We detail this in the next corollary.

\begin{corollary}\label{cor:pot_gm}
Let $f$ be an $L$-smooth convex function, and $x_\star\in\mathrm{argmin}_x\, f(x)$. For any $N\in\mathbb{N}$, the iterates of gradient descent with step size $\gamma_0=\gamma_1=\hdots=\gamma_N=\tfrac1L$ satisfy
\[f(x_N)-f_\star \leq \frac{L\|x_0-x_\star\|^2_2}{2N}.\]
\end{corollary}
\begin{proof}
Following the reasoning of~\eqref{eq:chained_pot}, we recursively use Theorem~\ref{thm:pot_GM}, starting with $A_0=0$. That is, we define 
\[\phi_k\defeq A_k(f(x_k)-f_\star)+\frac{L}{2}\|x_k-x_\star\|^2_2\]
and recursively use the inequality $\phi_{k+1}\leq\phi_k$ from Theorem~\ref{thm:pot_GM}, with $A_{k+1}=A_k+1$ and $A_0=0$; hence, $A_k=k$. We thus obtain 
\[A_N (f(x_N)-f_\star)\leq\phi_N\leq\hdots\leq \phi_0=\frac{L}{2}\|x_0-x_\star\|^2_2,\]
resulting in the desired statement \[f(x_N)-f_\star\leq \frac{L}{2A_N}\|x_0-x_\star\|^2_2=\frac{L}{2N}\|x_0-x_\star\|^2_2.\qedhere\]
\end{proof}

\subsection{How Conservative is this Worst-case Guarantee?}
Before moving to other methods, we show that the worst-case rate $O(N^{-1})$ of gradient descent is attained on very simple problems, motivating the search for alternate methods with better guarantees. This rate is observed on, e.g., all functions that are nearly linear over large regions. One such common function is the Huber loss (with $x_\star=0$, arbitrarily):
\[ f(x)=\left\{\begin{array}{ll}
    a_\tau |x| - b_\tau \quad & \text{if } |x|\geq \tau, \\
    \frac{L}{2}x^2 & \text{otherwise,}
\end{array} \right. \]
with $a_\tau=L\tau$ and $b_\tau=-\frac{L}{2}\tau^2$ to ensure its continuity and differentiability. On this function, as long as the iterates of gradient descent satisfy $|x_k|\geq \tau$, they behave as if the function were linear, and the gradient is constant. It is therefore relatively easy to explicitly compute all iterates. In particular, by picking $\tau=\frac{|x_0|}{2N+1}$, we get $f(x_N)-f_\star=\frac{L\|x_0-x_\star\|^2_2}{2(2N+1)}$ and reach the $O(N^{-1})$ worst-case bound; see \citet[Theorem 3.2]{Dror14}. Therefore, it appears that the worst-case bound from Corollary~\ref{cor:pot_gm} for gradient descent can only be improved in terms of the constants, but the rate itself is the best possible one for this simple method; see, for example,~\citep{Dror14,drori2014contributions} for the corresponding tight expressions.

In the next section, we show that similar reasoning based on potential functions produces methods with improved worst-case convergence rate $O(N^{-2})$, compared to the $O(N^{-1})$ of vanilla gradient descent.

\section{Optimized Gradient Method}\label{s:OGM}
Given that the complexity bound for gradient descent cannot be improved, it is reasonable to look for alternate, hopefully better, methods. In this section, we show that accelerated methods can be designed by optimizing their worst-case performance. To do so, we start with a reasonably broad family of candidate first-order methods described by
\begin{equation}\label{eq:ogm_start}
\begin{aligned}
y_1 &= y_0 - h_{1,0} \nabla f(y_0),\\
y_2 &= y_1 - h_{2,0} \nabla f(y_0) - h_{2,1} \nabla f(y_1),\\
y_3 &= y_2 - h_{3,0} \nabla f(y_0) - h_{3,1} \nabla f(y_1) - h_{3,2} \nabla f(y_2),\\
&\vdots\\
y_N &= y_{N-1} - \sum_{i=0}^{N-1} h_{N,i} \nabla f(y_i).
\end{aligned}
\end{equation}
Of course, methods in this form are impractical since they require keeping track of all previous gradients. Neglecting this potential problem for now, one possibility for choosing the step size $\{h_{i,j}\}$ is to solve a minimax problem:
\begin{equation}\label{eq:pep_opt_h}
 \min_{\{h_{i,j}\}}\,\, \max_{f\in\mathcal{F}_{0,L}} \left\{ \frac{f(y_N)-f_\star}{\|y_0-x_\star\|^2_2} \, : \, y_N \text{ obtained from~\eqref{eq:ogm_start} and $y_0$}\right\}.    
\end{equation}
In other words, we are looking for the best possible worst-case ratio among methods of the form~\eqref{eq:ogm_start}. 
Of course, different target notions of accuracy could be considered instead of $\left.(f(y_N)-f_\star)\middle/\|y_0-x_\star\|^2_2\right.$, but we proceed with this notion for now.

It turns out that~\eqref{eq:pep_opt_h} has a clean solution, obtained by~\citet{kim2016optimized}, based on clever reformulations and relaxations of~\eqref{eq:pep_opt_h} developed by~\citet{Dror14} (some details are provided in Section~\ref{s:notes_ref_chapt_nest}). Furthermore, this method has ``factorized'' forms that do not require keeping track of previous gradients. The optimized gradient method (OGM) is parameterized by a sequence $\{\theta_{k,N}\}_k$ that is constructed recursively starting from $\theta_{-1,N}=0$ (or equivalently $\theta_{0,N}=1$), using
\begin{equation}\label{eq:thetas}
\begin{aligned}
\theta_{k+1,N}= \left\{\begin{array}{ll}
        \frac{1+\sqrt{4\theta_{k,N}^2+1}}{2}  \, & \text{if } k\leq N-2  \\
        \frac{1+\sqrt{8\theta_{k,N}^2+1}}{2}  \, & \text{if } k=N-1.
      \end{array}\right.
\end{aligned}
\end{equation}
We also mention that optimized gradient methods can be stated in various equivalent formats, we provide two variants in Algorithm~\ref{alg:OGM_1} and Algorithm~\ref{alg:OGM_2} (a rigorous equivalence statement is provided in Appendix~\ref{s:eq_OGM}). While the shape of Algorithm~\ref{alg:OGM_2} is more common in accelerated methods, the equivalent formulation provided in Algorithm~\ref{alg:OGM_1} allows for slightly more direct proofs.
\begin{algorithm}[!ht]
  \caption{Optimized gradient method (OGM), form I}
  \label{alg:OGM_1}
  \begin{algorithmic}[1]
    \REQUIRE
     $L$-smooth convex function $f$, initial point $x_0$, and budget $N$.
    \STATE \textbf{Initialize} $z_0=y_0=x_0$ and $\theta_{-1,N}=0$.
    \FOR{$k=0,\ldots,N-1$}
      \STATE $\theta_{k,N}=\tfrac{1+\sqrt{4\theta_{k-1,N}^2+1}}{2}$
      \STATE $y_{k}=\left(1-\tfrac{1}{\theta_{k,N}}\right)x_{k}+\tfrac1{\theta_{k,N}}z_{k}$
      \STATE $x_{k+1}=y_k-\tfrac1L \nabla f(y_k)$
      \STATE $z_{k+1}=x_0-\tfrac2L\sum_{i=0}^k \theta_{i,N} \nabla f(y_k)$
    \ENDFOR
    \ENSURE Approximate solution $y_{N}=\left(1-\tfrac{1}{\theta_{N,N}}\right)x_{N}+\tfrac1{\theta_{N,N}}z_{N}$ with $\theta_{N,N}=\tfrac{1+\sqrt{8\theta_{N-1,N}^2+1}}{2}$.
  \end{algorithmic}
\end{algorithm}
\begin{algorithm}[!ht]
  \caption{Optimized gradient method (OGM), form II}
  \label{alg:OGM_2}
  \begin{algorithmic}[1]
    \REQUIRE
     $L$-smooth convex function $f$, initial point $x_0$, and budget $N$.
    \STATE \textbf{Initialize} $z_0=y_0=x_0$ and $\theta_{0,N}=1$.
    \FOR{$k=0,\ldots,N-1$}
      \STATE $\theta_{k+1,N}= \left\{\begin{array}{ll}
        \frac{1+\sqrt{4\theta_{k,N}^2+1}}{2}  \, & \text{if } k\leq N-2  \\
        \frac{1+\sqrt{8\theta_{k,N}^2+1}}{2}  \, & \text{if } k=N-1.
      \end{array}\right.$
      \STATE $x_{k+1}=y_k-\tfrac1L \nabla f(y_k)$
      \STATE $y_{k+1}=x_{k+1}+\tfrac{\theta_{k,N}-1}{\theta_{k+1,N}}(x_{k+1}-x_k)+\frac{\theta_{k,N}}{\theta_{k+1,N}}(x_{k+1}-y_k)$.
    \ENDFOR
    \ENSURE Approximate solution $y_{N}$.
  \end{algorithmic}
\end{algorithm}

Direct approaches to~\eqref{eq:pep_opt_h} are rather technical---see details in~\citep{Dror14,kim2016optimized}. However, showing that the OGM is indeed optimal on the class of smooth convex functions can be accomplished indirectly by providing an upper bound on its worst-case complexity guarantees and by showing that no first-order method can have a better worst-case guarantee on this class of problems. We detail a fully explicit worst-case guarantee for OGM in the next section. It consists in showing that
\begin{equation}\label{eq:pot_ogm_preview}
\phi_k\defeq 2\theta_{k-1,N}^2\left(f(y_{k-1})-f_\star-\frac{1}{2L}\lVert \nabla f(y_{k-1})\rVert^2_2\right)+\frac{L}{2}\|z_{k}-x_\star\|^2_2
\end{equation}
is a potential function for the optimized gradient method when \linebreak $k<N$ (Theorem~\ref{thm:pot_OGM}, below). For $k=N$, we need a minor adjustment (Lemma~\ref{thm:pot_OGM_final}, below) to obtain a bound on $f(y_N)-f_\star$ and not in terms of  $f(y_{N})-f_\star-\tfrac{1}{2L}\lVert \nabla f(y_{N})\rVert^2_2$, which appears in the potential.

As in the case of gradient descent, the proof relies on potential functions. Following the recursive argument from~\eqref{eq:chained_pot}, the convergence guarantee is driven by the convergence speed of $\theta_{k,N}^{-2}$ towards $0$. We note that when $k<N-1$,
\begin{equation}\label{eq:theta_evol}
\theta_{k+1,N}=\frac{1+\sqrt{4\theta_{k,N}^2+1}}{2}\geq \frac{1+2\theta_{k,N}}{2}=\theta_{k,N}+\frac12,  
\end{equation}
and therefore, $\theta_{k,N}\geq \tfrac{k}{2}+1$. We also directly obtain
\begin{equation}\label{eq:theta_evol2} \theta_{N,N}=\frac{1+\sqrt{8\theta_{N-1,N}^2+1}}{2}\geq \frac{1+\sqrt{2(N+1)^2+1}}{2}\geq \frac{N+1}{\sqrt{2}},
\end{equation}
and hence, $\theta_{N,N}^{-2}=O(N^{-2})$. Before providing the proof, we mention that it heavily relies on inequality~\eqref{eq:interpolation_strcvx} with $\mu=0$. This inequality is key for formulating~\eqref{eq:pep_opt_h} in a tractable way.

\subsection{A Potential for the Optimized Gradient Method}
The main point now is to prove that~\eqref{eq:pot_ogm_preview} is indeed a potential for the optimized gradient method. We emphasize again that our main motivation for proving this is to show that the OGM provides a good template algorithm for acceleration (i.e., a method involving two or three sequences) and that the corresponding potential functions can also be used as a template for the analysis of more advanced methods.

Note that the potential structure does not seem immediately intuitive: it was actually found using computer-assisted proof design techniques; see Section~\ref{s:notes_ref_chapt_nest} ``{On obtaining the proofs in this section}'' and Appendix~\ref{a-WC_FO} for further references. In particular, the following theorem can be found in~\citet[Theorem 11]{taylor19bach}.

\begin{theorem}\label{thm:pot_OGM}
Let $f$ be an $L$-smooth convex function, $x_\star\in\mathrm{argmin}_x\, f(x)$, and $N\in\mathbb{N}$. For any $k\in\mathbb{N}$ with $0\leq k \leq N-1$ and any $y_{k-1},z_k\in\mathbb{R}^d$, it holds that
\begin{equation*}
\begin{aligned}
2\theta_{k,N}^2&\left(f(y_{k})-f_\star-\frac{1}{2L}\lVert \nabla f(y_{k})\rVert^2_2\right)+\frac{L}{2}\|z_{k+1}-x_\star\|^2_2\\
&\leq 2\theta_{k-1,N}^2\left(f(y_{k-1})-f_\star-\frac1{2L}\lVert \nabla f(y_{k-1})\rVert^2_2\right)+\frac{L}{2}\|z_{k}-x_\star\|^2_2,
\end{aligned}
\end{equation*}
when $y_{k}$ and $z_{k+1}$ are obtained from Algorithm~\ref{alg:OGM_1}.
\end{theorem}
\begin{proof}
Recall that the algorithm can be written as
\begin{equation*}
\begin{aligned}
y_{k}&=\left(1-\frac{1}{\theta_{k,N}}\right)\left(y_{k-1}-\frac1L \nabla f(y_{k-1})\right)+\frac{1}{\theta_{k,N}}z_{k}\\
z_{k+1}&= z_{k}-\frac{2\theta_{k,N}}{L}\nabla f(y_{k}).
\end{aligned}
\end{equation*}
The proof consists of performing a weighted sum of the following inequalities.
\begin{itemize}
    \item Smoothness and convexity of $f$ between $y_{k-1}$ and $y_k$ with weight $\lambda_1={2\theta_{k-1,N}^2}$:
    \begin{align*}
    0\geq & f(y_k)-f(y_{k-1})+\langle \nabla f(y_k); y_{k-1}-y_k\rangle \\
    & \quad +\frac1{2L}\lVert \nabla f(y_k)-\nabla f(y_{k-1})\rVert^2_2.
    \end{align*}
    \item Smoothness and convexity of $f$ between $x_\star$ and $y_k$ with weight $\lambda_2={2\theta_{k,N}}$:
    \[ 0\geq f(y_k)-f_\star+\langle \nabla f(y_k); x_\star-y_k\rangle +\frac1{2L}\lVert \nabla f(y_k)\rVert^2_2. \]
\end{itemize}
Since the weights are nonnegative, the weighted sum produces a valid inequality:
\begin{equation}\label{eq:OGM_tot_WS}
\begin{aligned}
0\geq&  \lambda_1 \bigg[f(y_k)-f(y_{k-1})+\langle \nabla f(y_k); y_{k-1}-y_k\rangle +\frac1{2L}\lVert \nabla f(y_k)\\
& \quad -\nabla f(y_{k-1})\rVert^2_2\bigg]\\
& +\lambda_2  \left[f(y_k)-f_\star+\langle \nabla f(y_k); x_\star-y_k\rangle +\frac1{2L}\lVert \nabla f(y_k)\rVert^2_2\right],
\end{aligned}
\end{equation}
which (either by completing the squares or simply by extending both expressions and verifying that they match on a term-by-term basis) can be reformulated as 
\begin{equation*}
\begin{aligned}
0\geq&\lambda_1 \bigg[f(y_k)-f(y_{k-1})+\langle \nabla f(y_k); y_{k-1}-y_k\rangle \\
& \quad +\frac1{2L}\lVert \nabla f(y_k)-\nabla f(y_{k-1})\rVert^2_2\bigg]\\
& +\lambda_2  \left[f(y_k)-f_\star+\langle \nabla f(y_k); x_\star-y_k\rangle +\frac1{2L}\lVert \nabla f(y_k)\rVert^2_2\right]\\
=& 2\theta_{k,N}^2\left(f(y_{k})-f_\star-\frac{1}{2L}\lVert \nabla f(y_{k})\rVert^2_2\right)+\frac{L}{2}\|z_{k+1}-x_\star\|^2_2\\
&-2\theta_{k-1,N}^2\left(f(y_{k-1})-f_\star-\frac1{2L}\lVert \nabla f(y_{k-1})\rVert^2_2\right)-\frac{L}{2}\|z_{k}-x_\star\|^2_2\\
\end{aligned}
\end{equation*}

\begin{equation*}
\begin{aligned}
{\color{white}0\geq}&+\frac{2}{\theta_{k,N}}\left(\theta_{k-1,N}^2+\theta_{k,N}-\theta_{k,N}^2\right)\langle\nabla f(y_k); y_{k-1}-\tfrac{1}{L}\nabla f(y_{k-1})-z_k\rangle\\
&+2 \left(\theta_{k-1,N}^2+\theta_{k,N}-\theta_{k,N}^2\right) \left(f(y_k)-f_\star+\frac{1}{2 L}\|\nabla f(y_k)\|^2_2\right).
\end{aligned}
\end{equation*}
The desired conclusion follows from picking $\theta_{k,N}\geq\theta_{k-1,N}$ satisfying \[\theta_{k-1,N}^2+\theta_{k,N}-\theta_{k,N}^2=0,\] and hence the choice~\eqref{eq:thetas}, thus reaching
\begin{equation*}
\begin{aligned}
2\theta_{k,N}^2&\left(f(y_{k})-f_\star-\frac{1}{2L}\lVert \nabla f(y_{k})\rVert^2_2\right)+\frac{L}{2}\|z_{k+1}-x_\star\|^2_2\\&\leq
2\theta_{k-1,N}^2\left(f(y_{k-1})-f_\star-\frac1{2L}\lVert \nabla f(y_{k-1})\rVert^2_2\right)+\frac{L}{2}\|z_{k}-x_\star\|^2_2.\qedhere
\end{aligned}
\end{equation*}
\end{proof}

A final technical fix is required now. To show that the optimized gradient method is an optimal solution to~\eqref{eq:pep_opt_h}, we need an upper bound on the function values, rather than on the function values minus a squared gradient norm. This discrepancy is handled by the following technical lemma.
\begin{lemma}\label{thm:pot_OGM_final}
Let $f$ be an $L$-smooth convex function, $x_\star\in\mathrm{argmin}_x\, f(x)$, and $N\in\mathbb{N}$. For any $y_{N-1},z_{N}\in\mathbb{R}^d$, it holds that
\begin{equation*}
\begin{aligned}
\theta_{N,N}^2&\left(f(y_N)-f_\star\right)+\frac{L}{2}\|z_{N}-\tfrac{\theta_{N,N}}{L}\nabla f(y_N)-x_\star\|^2_2\\
&\leq 2\theta_{N-1,N}^2\left(f(y_{N-1})-f_\star-\frac1{2L}\lVert \nabla f(y_{N-1})\rVert^2_2\right)+\frac{L}{2}\|z_{N}-x_\star\|^2_2,
\end{aligned}
\end{equation*}
where $y_{N}$ is obtained from Algorithm~\ref{alg:OGM_1}.
\end{lemma}
\begin{proof}
The proof consists of performing a weighted sum of the following inequalities.
\begin{itemize}
    \item Smoothness and convexity of $f$ between $y_{N-1}$ and $y_N$ with weight $\lambda_1={2\theta_{N-1,N}^2}$:
    \begin{align*}
    0\geq& f(y_N)-f(y_{N-1})+\langle \nabla f(y_N); y_{N-1}-y_N\rangle \\
    &+\frac1{2L}\lVert \nabla f(y_N)-\nabla f(y_{N-1})\rVert^2_2.
    \end{align*}
    \item Smoothness and convexity of $f$ between $x_\star$ and $y_N$ with weight $\lambda_2=\theta_{N,N}$:
    \[ 0\geq f(y_N)-f_\star+\langle \nabla f(y_N); x_\star-y_N\rangle +\frac1{2L}\lVert \nabla f(y_N)\rVert^2_2. \]
\end{itemize}
Since the weights are nonnegative, the weighted sum produces a valid inequality:
\begin{equation*}
\begin{aligned}
0\geq & \lambda_1 \bigg[f(y_N)-f(y_{N-1})+\langle \nabla f(y_N); y_{N-1}-y_N\rangle \\
& \quad +\frac1{2L}\lVert \nabla f(y_N)-\nabla f(y_{N-1})\rVert^2_2\bigg]\\
&+ \lambda_2 \left[f(y_N)-f_\star+\langle \nabla f(y_N); x_\star-y_N\rangle +\frac1{2L}\lVert \nabla f(y_N)\rVert^2_2\right],
\end{aligned}
\end{equation*}
which can be reformulated as
\begin{equation*}
\begin{aligned}
0\geq &\theta_{N,N}^2 \left(f(y_N)-f_\star\right)+\frac{L}{2}\|z_{N}-\tfrac{\theta_{N,N}}{L}\nabla f(y_N)-x_\star\|^2_2\\
&-2\theta_{N-1,N}^2\left(f(y_{N-1})-f_\star-\frac1{2L}\lVert \nabla f(y_{N-1})\rVert^2_2\right)-\frac{L}{2}\|z_{N}-x_\star\|^2_2\\
&+\frac{1}{\theta_{N,N}} \left(2 \theta_{N-1,N}^2-\theta_{N,N}^2+\theta_{N,N}\right) \\
& \quad\quad \times\langle \nabla f(y_N);y_{N-1}-\tfrac{1}{L}\nabla f(y_{N-1})-z_N\rangle\\
&+\left(2 \theta_{N-1,N}^2-\theta_{N,N}^2+\theta_{N,N}\right) \left(f(y_N)-f_\star+\frac{1}{2 L}\|\nabla f(y_N)\|^2_2\right).
\end{aligned}
\end{equation*}
The conclusion follows from choosing $\theta_{N,N}\geq\theta_{N-1,N}$ such that \[2 \theta_{N-1,N}^2-\theta_{N,N}^2+\theta_{N,N}=0,\] thereby reaching the desired inequality:
\begin{equation*}
\begin{aligned}
\theta_{N,N}^2 &\left(f(y_N)-f_\star\right)+\frac{L}{2}\|z_{N}-\tfrac{\theta_{N,N}}{L}\nabla f(y_N)-x_\star\|^2_2\\
&\leq 2\theta_{N-1,N}^2\left(f(y_{N-1})-f_\star-\frac1{2L}\lVert \nabla f(y_{N-1})\rVert^2_2\right)+\frac{L}{2}\|z_{N}-x_\star\|^2_2.\qedhere
\end{aligned}
\end{equation*}
\end{proof}

By combining Theorem~\ref{thm:pot_OGM} and the technical Lemma~\ref{thm:pot_OGM_final}, we get the final worst-case performance bound of the OGM on function values, detailed in the corollary below.

\begin{corollary}\label{cor:OGM_bound}
Let $f$ be an $L$-smooth convex function, and $x_\star\in\mathrm{argmin}_x\, f(x)$. For any  $N\in\mathbb{N}$ and $x_0\in\mathbb{R}^d$, the output of the optimized gradient method (OGM, Algorithm~\ref{alg:OGM_1} or Algorithm~\ref{alg:OGM_2}) satisfies
\[ f(y_N)-f_\star \leq \frac{L\lVert x_0-x_\star\rVert^2_2}{2\theta_{N,N}^2}\leq \frac{L\|x_0-x_\star\|^2_2}{(N+1)^2}.\]
\end{corollary}
\begin{proof} Defining, for $k\in\{1,\hdots,N\}$
\[ \phi_k\defeq 2\theta_{k-1,N}^2\left(f(y_{k-1})-f_\star-\frac{1}{2L}\lVert \nabla f(y_{k-1})\rVert^2_2\right)+\frac{L}{2}\|z_{k}-x_\star\|^2_2,\]
and \[ \phi_{N+1}\defeq \theta_{N,N}^2\left(f(y_N)-f_\star\right)+\frac{L}{2}\|z_{N}-\tfrac{\theta_{N,N}}{L}\nabla f(y_N)-x_\star\|^2_2,\]
we reach the desired statement: \[\theta_{N,N}^2(f(y_N)-f_\star)\leq\phi_{N+1}\leq \phi_N\leq\hdots\leq \phi_0=\frac{L}2\|x_0-x_\star\|^2_2,\]
using Theorem~\ref{thm:pot_OGM} and technical Lemma~\ref{thm:pot_OGM_final}. We obtain the last bound by using $\theta_{N,N}\geq (N+1)/\sqrt{2}$; see~\eqref{eq:theta_evol2}.
\end{proof}
In the following section, we mostly use potential functions, relying directly on the function value $f(x_k)$ instead of $f(y_k)$ for practical reasons discussed below. Note that using the \emph{descent lemma} (i.e., the inequality $f(x_{k+1})\leq f(y_k)-\tfrac1{2L}\|\nabla f(y_k)\|^2_2$) directly on the potential function allows us to obtain a bound on $f(x_{k+1})$ for the OGM. This result can be found in~\citep[Theorem 3.1]{kim2017convergence} without the ``potential function'' mechanism.

\subsection{Optimality of Optimized Gradient Methods}

A nice, commonly used guide for designing {optimal} methods consists of constructing problems that are {difficult} for all methods within a certain class. This strategy results in \emph{lower complexity bounds}, and it is often deployed via the concept of \emph{minimax risk} (of a class of problems and a class of methods)---see, e.g.,~\citet{guzman2015lower}---which corresponds to the worst-case performance of the best method within the prescribed class. In this section, we briefly discuss such results in the context of smooth convex minimization, on the particular class of \emph{black-box} first-order methods. The term \emph{black-box} is used to emphasize that the method has no prior knowledge of $f$ (beyond the class of functions to which $f$ belongs, so methods are allowed to use $L$) and that it can only obtain information about $f$ by a evaluating its gradient/function value through an \emph{oracle}. 

Of particular interest to us,~\citet{drori2017exact} established that the worst-case performance achieved by the optimized gradient method (see Corollary~\ref{cor:OGM_bound}) on the class of smooth convex functions cannot in general be improved by any black-box first-order method.
\begin{theorem}\citep[Theorem 3]{drori2017exact}\label{thm:smooth_LB} Let $L>0$, $d,N\in\mathbb{N}$ with $d\geq N+1$. For any black-box first-order method that performs at most $N$ calls to the first-order oracle $(f(\cdot),\nabla f(\cdot))$, there exists a function $f\in\mathcal{F}_{0,L}(\mathbb{R}^d)$ and $x_0\in\mathbb{R}^d$ such that
\[ f(x_N)-f(x_\star)\geq \frac{L\lVert x_0-x_\star\rVert^2_2}{2\theta_{N,N}^2},\]
with $x_\star\in\mathrm{argmin}_x\, f(x)$, $x_N$ is the output of the method under consideration, and $x_0$ its input.
\end{theorem}
In the previous sections, we showed that $\theta_{N,N}^{2}\geq \tfrac{(N+1)^2}{2}$. It is also relatively easy to establish that \[\theta_{k+1,N}\leq \frac{1+\sqrt{(2\theta_{k,N}+1)^2}}{2}=1+\theta_{k,N},\]thereby obtaining $\theta_{k,N}\leq k+1$ (because $\theta_{0,N}=1$) as well as
\[\theta_{N,N}=\frac{1+\sqrt{8\theta_{N-1,N}^2+1}}{2}\leq\frac{1+\sqrt{8N^2+1}}{2}\leq \sqrt{2}N+1.\] We conclude (through Theorem~\ref{thm:smooth_LB}) that the lower bound has the form \[ f(x_N)-f(x_\star)\geq \frac{L\lVert x_0-x_\star\rVert^2_2}{2\theta_{N,N}^2}\geq \frac{L\|x_0-x_\star\|^2_2}{(\sqrt{2}N+1)^2}=\Omega(N^{-2}).\]

While Drori's approach to obtaining this lower bound is rather technical (a slightly simplified and weaker version of this result can be found in~\citep[Corollary 5]{drori2021exact}), there are simpler approaches that allow us to show that the rate $\Omega(N^{-2})$ (that is, neglecting the tight constants) cannot in general be beaten in black-box smooth convex minimization.
For one such example, we refer to~\citep[Theorem~2.1.6]{Nest03a}. In a closely related line of work,~\citep{nemirovsky1991optimality} established similar \emph{exact} bounds in the context of solving linear systems of equations and for minimizing convex quadratic functions (see also Section~\ref{s:cheby_wc}). For convex quadratic problems whose Hessian has bounded eigenvalues between $0$ and $L$, these lower bounds are attained by the Chebyshev (see Section~\ref{c-Cheb}) and by conjugate gradient methods~\citep{nemirovsky1991optimality,nemirovsky1992information}. 

Perhaps surprisingly, the conjugate gradient method also achieves the lower complexity bound of smooth convex minimization provided by Theorem~\ref{thm:smooth_LB}. Furthermore, the proof follows essentially the same structure as that for the OGM. In particular, it relies on the same potential function (see Appendix~\ref{s:conj_methods}).

\subsection{Optimized Gradient Method: Summary}\label{sec:LO_OGM}
Before going further, we quickly summarize what we have learned from the optimized gradient method. First of all, the optimized gradient method can be seen as a counterpart of the Chebyshev method for minimizing quadratics, applied to smooth convex minimization. It is an \emph{optimal} method in the sense that it has the smallest possible worst-case ratio $\tfrac{f(y_N)-f_\star}{\|y_0-x_\star\|^2_2}$ over the class $f\in\mathcal{F}_{0,L}$ among all black-box first-order methods, given a fixed computational budget of $N$ gradient evaluations. Furthermore, although this method has a few drawbacks (we mention a few below), it can be seen as a \emph{template} for designing other accelerated methods using the same algorithmic and proof structures. We extensively use variants of this template below. In other words, most variants of \emph{accelerated gradient methods} rely on the same two (or three) sequence structures, and on similar potential functions. Such variants usually rely on slight variations in the choice of the parameters used throughout the iterative process, typically involving less aggressive step size strategies (i.e., smaller values for $\{h_{i,j}\}$ in~\eqref{eq:ogm_start}).

Second, the OGM is not a very practical method as such: it is fined-tuned for unconstrained smooth convex minimization and does not readily extend to other situations, such as situations involving constraints, for which~\eqref{eq:interpolation_strcvx} does not hold in general; see the discussions in~\citep{drori2018properties} and Remark~\ref{rem:restrictedset}.

On the other hand, we see in what follows that it is relatively easy to design other methods that follow the same template and achieve the same $O(N^{-2})$ rate, while resolving the issues of the OGM listed above. Such methods use slightly less aggressive step size strategies, at the cost of being slightly suboptimal for~\eqref{eq:pep_opt_h}, i.e., they have slightly worse worst-case guarantees. In this vein, we start by discussing the original accelerated gradient method, proposed by~\citet{Nest83}.

\section{Nesterov's Acceleration}
Motivated by the format of the optimized gradient method, we detail a potential-based proof for Nesterov's method. We then quickly review the concept of \emph{estimate sequences} and show that they provide an interpretation of potential functions as increasingly good models of the function to be minimized. Finally, we extend these results to strongly convex minimization.

\subsection{Nesterov's Method, from Potential Functions}
In this section, we follow the algorithmic template provided by the optimized gradient method. In this spirit, we start by discussing the first accelerated method in its simplest form (Algorithm~\ref{alg:FGM_1}) as well as its potential function, originally proposed by~\citet{Nest83}, but the presentation here is different. 

Our goal is to derive the simplest algebraic proof for this scheme. We follow the algorithmic template of the optimized gradient method (which is further motivated in Section~\ref{s:ogm_sc}). Once a potential is chosen, the proofs are quite straightforward as simple combinations of inequalities and basic algebra. Our choice of potential function is not immediately obvious but allows for simple extensions afterwards. Other choices are possible, for example, incorporating $f(y_k)$ as (in the OGM) or additional terms such as $\|\nabla f(x_k)\|^2_2$. We pick a potential function similar to that used for gradient descent and that of the optimized gradient method, which is written
\[ 
\phi_k \defeq A_k (f(x_k)-f_\star)+\frac{L}{2}\|z_k-x_\star\|^2_2,
\]
where one iteration of the algorithm has the following form, reminiscent of the OGM:
\begin{equation}\label{eq:FGM_1}
\begin{aligned}
y_{k}&=x_k+\tau_k (z_k-x_k)\\
x_{k+1}&=y_{k}-\alpha_k \nabla f(y_{k})\\
z_{k+1}&=z_k-\gamma_k \nabla f(y_{k}).
\end{aligned}
\end{equation}
Our goal is to select algorithmic parameters $\{(\tau_k,\alpha_k,\gamma_k)\}_k$ so as to greedily make $A_{k+1}$ as large as possible as a function of $A_k$ since the convergence rate of the method is controlled by the inverse of the growth rate of $A_k$, i.e., $f(x_N)-f_\star=O(A_N^{-1})$. 

In practice, we can pick $A_{k+1}=A_k+\frac12(1+\sqrt{4 A_k+1})$ by choosing $\tau_k=1-A_k/A_{k+1}$, $\alpha_k=\frac1L$, and $\gamma_k=(A_{k+1}-A_k)/{L}$ (see Algorithm~\ref{alg:FGM_1}), and the proof is then quite compact.

\begin{algorithm}[!ht]
  \caption{Nesterov's method, form I}
  \label{alg:FGM_1}
  \begin{algorithmic}[1]
    \REQUIRE An $L$-smooth convex function $f$ and initial point $x_0$.
    \STATE \textbf{Initialize} $z_0=x_0$ and $A_0=0$.
    \FOR{$k=0,\ldots$}
      \STATE $a_{k}=\frac12(1+\sqrt{4 A_k+1})$
      \STATE $A_{k+1}=A_k+a_k$
      \STATE $y_{k}= x_k+(1-\frac{A_k}{A_{k+1}}) (z_k-x_k)$
      \STATE $x_{k+1}=y_{k}-\frac1L \nabla f(y_{k})$
      \STATE $z_{k+1}=z_k- \frac{A_{k+1}-A_k}{L}\nabla f(y_{k})$
    \ENDFOR
    \ENSURE Approximate solution $x_{k+1}$.
  \end{algorithmic}
\end{algorithm}

Before continuing to the proof of the potential inequality, we show that $A_k^{-1}=O(k^{-2})$. Indeed, we have,
\begin{equation}\label{eq:conv_Ak_fgm}
\begin{aligned}
A_{k}&=A_{k-1}+\frac{1+\sqrt{4A_{k-1}+1}}{2}\geq A_{k-1}+\frac12+\sqrt{A_{k-1}}\geq \left(\sqrt{A_{k-1}}+\frac12\right)^2\geq \frac{k^2}{4},
\end{aligned}
\end{equation}
where the last inequality follows from a recursive application of the previous one, along with $A_0=0$. 

\begin{theorem}\label{thm:Nest_first}
Let $f$ be an $L$-smooth convex function, $x_\star\in\mathrm{argmin}_x\, f(x)$, and $k\in\mathbb{N}$. For any $x_k,z_k\in\mathbb{R}^d$ and $A_k\geq 0$, the iterates of Algorithm~\ref{alg:FGM_1} satisfy
\[ A_{k+1}(f(x_{k+1})-f_\star)+\frac{L}{2}\|z_{k+1}-x_\star\|^2_2\leq A_k (f(x_k)-f_\star)+\frac{L}2 \|z_k-x_\star\|^2_2,\]
with $A_{k+1}=A_{k}+\tfrac{1+\sqrt{4A_{k}+1}}{2}$.
\end{theorem}
\begin{proof}
The proof consists of a weighted sum of the following inequalities.
\begin{itemize}
    \item Convexity of $f$ between $x_\star$ and $y_{k}$ with weight $\lambda_1=A_{k+1}-A_k$:
    \[ f_\star \geq f(y_{k})+\langle \nabla f(y_{k});x_\star-y_k\rangle.\]
    \item Convexity of $f$ between $x_k$ and $y_{k}$ with weight $\lambda_2=A_k$:
    \[ f(x_k) \geq f(y_{k})+\langle \nabla f(y_{k});x_k-y_{k}\rangle.\]
    \item Smoothness of $f$ between $y_{k}$ and $x_{k+1}$ (a.k.a., \emph{descent lemma}) with weight $\lambda_3=A_{k+1}$:
    \[ f(y_{k}) + \langle \nabla f(y_{k}); x_{k+1}-y_{k}\rangle +\frac{L}{2}\|x_{k+1}-y_{k}\|^2_2 \geq f(x_{k+1}).\]
\end{itemize}
We therefore arrive at the following valid inequality
\begin{equation*}
\begin{aligned}
0\geq& \lambda_1 [ f(y_{k})-f_\star +\langle \nabla f(y_{k});x_\star-y_{k}\rangle]\\&+ \lambda_2[f(y_{k})-f(x_k)+\langle \nabla f(y_{k});x_k-y_{k}\rangle]\\&+\lambda_3[f(x_{k+1})-f(y_{k}) - \langle \nabla f(y_{k}); x_{k+1}-y_{k}\rangle -\frac{L}{2}\|x_{k+1}-y_{k}\|^2_2].
\end{aligned}
\end{equation*}
For the sake of simplicity, we do not substitute $A_{k+1}$ by its expression until the last stage of the reformulation. Substituting $y_k$, $x_{k+1}$, and $z_{k+1}$ by their expressions in~\eqref{eq:FGM_1} along with $\tau_k=1-A_k/A_{k+1}$, $\alpha_k=\frac1L$, and $\gamma_k=\frac{A_{k+1}-A_k}{L}$, basic algebra shows that the previous inequality can be reorganized as 

\begin{equation*}
\begin{aligned}
0\geq &A_{k+1}(f(x_{k+1})-f_\star)+\frac{L}{2}\|z_{k+1}-x_\star\|^2_2\\
& -A_k (f(x_k)-f_\star)-\frac{L}2 \|z_k-x_\star\|^2_2\\
& +\frac{A_{k+1}-(A_k-A_{k+1})^2}{2 L}\|\nabla f(y_{k})\|^2_2.
\end{aligned}
\end{equation*}
The claim follows from selecting $A_{k+1}\geq A_k$ such that $A_{k+1}-(A_k-A_{k+1})^2=0$, thereby reaching
\begin{equation*}
\begin{aligned}
&A_{k+1}(f(x_{k+1})-f_\star)+\frac{L}{2}\|z_{k+1}-x_\star\|^2_2 \leq A_k (f(x_k)-f_\star)+\frac{L}2 \|z_k-x_\star\|^2_2.\qedhere
\end{aligned}
\end{equation*}
\end{proof}

The final worst-case guarantee is obtained by using the same chaining argument as in~\eqref{eq:chained_pot}, combined with an upper bound on $A_N$.

\begin{corollary}\label{cor:pot_fgm}
Let $f$ be an $L$-smooth convex function, and $x_\star\in\mathrm{argmin}_x\, f(x)$. For any $N\in\mathbb{N}$, the iterates of Algorithm~\ref{alg:FGM_1} satisfy
\[
f(x_N)-f_\star \leq \frac{2L\|x_0-x_\star\|^2_2}{N^2}.
\]
\end{corollary}
\begin{proof} Following the argument of~\eqref{eq:chained_pot}, we recursively use Theorem~\ref{thm:Nest_first} with $A_0=0$:
\[ A_N(f(x_N)-f_\star)\leq \phi_N\leq\hdots\leq\phi_0=\frac{L}{2}\|x_0-x_\star\|^2_2,\]
which yields
\[ f(x_N)-f_\star\leq \frac{L\|x_0-x_\star\|^2_2}{2A_N}\leq \frac{2L\|x_0-x_\star\|^2_2}{N^2},\]
where we used $A_N\geq N^2/4$ from~\eqref{eq:conv_Ak_fgm} to reach the last inequality.
\end{proof}

Before moving on, we emphasize that the rate of $O(N^{-2})$ matches that of lower bounds (see, e.g., Theorem~\ref{thm:smooth_LB}) up to absolute constants.

Finally, note that Nesterov's method is often written in a slightly different format, similar to that of Algorithm~\ref{alg:OGM_2}. The alternate formulation omits the third sequence $z_k$ and is provided in Algorithm~\ref{alg:FGM_2}. It is preferred in many references on the topic due to its simplicity. A third equivalent variant is provided in Algorithm~\ref{alg:FGM_3}; this variant turns out to be useful when generalizing the method beyond Euclidean spaces. The equivalence statements between Algorithm~\ref{alg:FGM_1}, Algorithm~\ref{alg:FGM_2}, and Algorithm~\ref{alg:FGM_3} are relatively simple and are provided in Appendix~\ref{s:eq_FGM}. Many references tend to favor one of these formulations, and we want to point out that they are equivalent in the base problem setup of unconstrained smooth convex minimization. Although the expression of the different formats in terms of the same external sequence $\{A_k\}_k$ does not always correspond to their simplest forms (i.e., alternate parameterizations might be simpler, particularly in the strongly convex case which follows), we proceed with this sequence to avoid introducing too many variations on the same theme.

\begin{algorithm}[!ht]
  \caption{Nesterov's method, form II}
  \label{alg:FGM_2}
  \begin{algorithmic}[1]
    \REQUIRE An $L$-smooth convex function $f$ and initial point $x_0$.
    \STATE \textbf{Initialize} $y_0=x_0$ and $A_0=0$.
    \FOR{$k=0,\ldots$}
      \STATE $a_{k}=\frac12(1+\sqrt{4 A_k+1})$
      \STATE $A_{k+1}=A_k+a_k$
      \STATE $x_{k+1}=y_k-\tfrac{1}{L}\nabla f(y_k)$
      \STATE $y_{k+1}=x_{k+1}+\frac{a_{k}-1}{a_{k+1}}(x_{k+1}-x_k)$
    \ENDFOR
    \ENSURE Approximate solution $x_{k+1}$.
  \end{algorithmic}
\end{algorithm}

\begin{algorithm}[!ht]
  \caption{Nesterov's method, form III}
  \label{alg:FGM_3}
  \begin{algorithmic}[1]
    \REQUIRE An $L$-smooth convex function $f$ and initial point $x_0$.
    \STATE \textbf{Initialize} $z_0=x_0$ and $A_0=0$.
    \FOR{$k=0,\ldots$}
      \STATE $a_{k}=\frac12(1+\sqrt{4 A_k+1})$
      \STATE $A_{k+1}=A_k+a_k$
      \STATE $y_{k}= x_k+(1-\frac{A_k}{A_{k+1}}) (z_k-x_k)$
      \STATE $z_{k+1}=z_k- \frac{A_{k+1}-A_k}{L}\nabla f(y_{k})$
      \STATE $x_{k+1}=\frac{A_k}{A_{k+1}}x_k+(1-\frac{A_k}{A_{k+1}})z_{k+1}$
    \ENDFOR
    \ENSURE Approximate solution $x_{k+1}$.
  \end{algorithmic}
\end{algorithm}
\subsection{Estimate Sequence Interpretation}\label{s:estimate}
We now relate the potential function approach to \emph{estimate sequences}. That is, we relate acceleration to first-order methods maintaining a model of the function throughout the iterative procedure. This approach was originally developed in~\citep[Section 2.2]{Nest03a}, and it has since been used in numerous works to obtain accelerated first-order methods in various settings (see discussions in Section~\ref{s:notes_ref_chapt_nest}). We present a slightly modified version, related to those of~\citep{baes2009estimate,wilson2016lyapunov}, which simplifies our comparisons with the previous material.

\subsubsection{Estimate Sequences}
As we see below, the basic idea underlying estimate sequences is closely related to that of potential functions, but it has explicit interpretations in terms of models of the objective function $f$. More precisely a sequence of pairs $\{(A_k,\varphi_k(x))\}_k$, with $A_k\geq 0$ and $\varphi_k:\mathbb{R}^d\rightarrow\mathbb{R}$, is called an estimate sequence of a function $f$ if\\ 
(i) for all $k\geq 0$ and $x\in\mathbb{R}^d$ we have
\begin{equation}\label{eq:est_seq}
 \varphi_k(x)-f(x) \leq A_k^{-1}(\varphi_0(x)-f(x)),
\end{equation}
(ii) $A_k\rightarrow \infty$ as $k\rightarrow\infty$.\\ If in addition, an estimate sequence satisfies\\ (iii)~for all $k\geq 0$, there exists some $x_k$ such that $f(x_k)\leq \varphi_k(x_\star)$, then we can guarantee that $f(x_k)-f_\star=O(A_k^{-1})$.

The purpose of estimate sequences is to start from an initial model $\varphi_0(x)$ satisfying $\varphi_0(x)\geq f_\star$ for all $x\in\mathbb{R}^d$ and then to design a sequence of convex models $\varphi_k$ that are increasingly good approximations of $f$, in the sense of~\eqref{eq:est_seq}. We provide further comments on conditions (i) and (iii), assuming for simplicity that $\{A_k\}$ is monotonically increasing (as is the case for all methods treated in this section).
\begin{itemize}
    \item Regarding (i), for all $x\in\mathbb{R}^d$, we have to design $\varphi_k$ to be either (a) a lower bound on the function (i.e., $\varphi_k(x)-f(x)\leq 0$ for that $x$) or (b) an increasingly good upper approximation of $f(x)$ when $0\leq \varphi_k(x)-f(x) \leq A_k^{-1}(\varphi_0(x)-f(x))$ for that~$x$. That is, we require that the error $|f(x)-\varphi_k(x)|$, incurred when approximating $f(x)$ by $\varphi_k(x)$, gets smaller for all $x$ for which $\varphi_k(x)$ is an upper bound on $f(x)$.
    
    To develop such models and the corresponding methods, three sequences of points are commonly used: (a) minimizers of our models $\varphi_k$ that correspond to iterates $z_k$ of the corresponding method; (b)~a sequence $y_k$ of points, whose first-order information is used to update the model of the function; and (c) the iterates $x_k$, corresponding to the best possible $f(x_k)$ that we can form. (The iteratives often do not correspond to the minimum of the model, $\varphi_k$, which is not necessarily an upper bound on the function.) 
    \item Regarding (iii), this condition ensures that the models $\varphi_k$ remain upper bounds on the optimal value $f_\star$. That is, it ensures that $f_\star\leq \varphi_k(x_\star)$ (since $f_\star\leq f(x_k)$) and hence that $\varphi_k(x_\star)-f_\star\geq 0$. From previous bullet point, this ensures that the modeling error of $f_\star$ goes to $0$ asymptotically as $k$ increases. More formally, conditions (ii) and (iii) allow us to construct proofs similar to potential functions and to obtain convergence rates. That is, under (iii), we get that
\begin{equation}\label{eq:est_sec_cons}
 f(x_k)-f_\star\leq \varphi_k(x_\star)-f(x_\star) \leq A_k^{-1} (\varphi_0(x_\star)-f(x_\star)),  
\end{equation}
and that therefore $f(x_k)-f_\star\leq O(A_k^{-1})$. The convergence rate is thereby dictated by the rate of $A_k^{-1}$, which goes to $0$ by~(ii). 
\end{itemize}

Now, the game consists of picking appropriate sequences $\{(A_k,\varphi_k)\}$ that correspond to simple algorithms. We thus translate our potential function results in terms of estimate sequences.

\subsubsection{Potential Functions as Estimate Sequences}
One can observe that potential functions and estimate sequences are closely related. First, in both cases, the convergence speed is dictated by that of a scalar sequence $A_{k}^{-1}$. In fact, there is one subtle but important difference between the two approaches: whereas $\varphi_k(x)$ should be an increasingly good approximation of $f$ for all $x$ in the context of estimate sequences,  potential functions require a model to be an increasingly good approximation of only $f_\star$, which is less restrictive. Hence, estimate sequences are more general but may not effectively handle situations in which the analysis actually requires having a weaker model that holds only on $f_\star$, and not of $f(x)$, for all $x$. We make this discussion more concrete via three examples, namely gradient descent, Nesterov's method, and the optimized gradient method.

\begin{itemize}
    \item Gradient descent: the potential inequality from Theorem~\ref{thm:pot_GM} actually holds for all $x$, and not only $x_\star$, as the proof does not exploit the optimality of $x_\star$. That is, it is proved that:
    \begin{equation*}
    \begin{aligned}
    (A_k+1) (f(x_{k+1})-f(x))&+\frac{L}{2}\|x_{k+1}-x\|^2_2\\&\leq A_k (f(x_{k})-f(x))+\frac{L}{2}\|x_{k}-x\|^2_2
    \end{aligned}
    \end{equation*}
    for all $x\in\mathbb{R}^d$. Therefore, the pair $\{(A_k,\varphi_k(x))\}_k$ with
    \[ \varphi_k(x)= f(x_k)+\frac{L}{2 A_k}\|x_k-x\|^2_2\]
    and $A_k=A_0+k$ (with $A_0>0$) is an estimate sequence for gradient descent. 
    \item Nesterov's first method: the potential inequality from Theorem~\ref{thm:Nest_first} also holds for all $x\in\mathbb{R}^d$, not only $x_\star$, as the proof does not exploit the optimality of $x_\star$. That is, it is proved that:
    \begin{equation*}
    \begin{aligned}
    A_{k+1} (f(x_{k+1})-f(x))&+\frac{L}{2}\|z_{k+1}-x\|^2_2\\&\leq A_k (f(x_{k})-f(x))+\frac{L}{2}\|z_{k}-x\|^2_2
    \end{aligned}
    \end{equation*}
    for all $x\in\mathbb{R}^d$. Hence, the pair $\{(A_k,\varphi_k(x))\}_k$ with
    \[ \varphi_k(x)= f(x_k)+\frac{L}{2 A_k}\|z_k-x\|^2_2\]
    and $A_k=A_{k-1}+\tfrac{1+\sqrt{4 A_{k-1}+1}}{2}$ (with $A_0>0$) is an estimate sequence for Nesterov's method. 
    \item Optimized gradient method: the potential inequality from Theorem~\ref{thm:pot_OGM} exploits the fact that $x_\star$ is an optimal point. Indeed, the proof relies on
    \[ f(x_\star)\geq f(y_k)+\langle \nabla f(y_k);x_\star-y_k\rangle +\frac1{2L}\|\nabla f(y_k)\|^2_2,\]
    which is an instance of Equation~\eqref{eq:interpolation_strcvx} exploiting $\nabla f(x_\star)=0$. This does not mean that there is no \emph{estimate sequence}-type model of the function as the algorithm proceeds, but the potential does not directly correspond to one. Alternatively, one can interpret 
    \[ \varphi_k(x)=f(y_k)-\frac{1}{2L}\|\nabla f(y_k)\|^2_2+\frac{L}{4\theta_{k,N}^2}\|z_{k+1}-x\|^2_2\]
    as an increasingly good model of $f_\star$ (i.e., it is an increasingly good approximation of $f(x)$ for all $x$ such that $\nabla f(x)=0$).
    
    A similar conclusion holds for the conjugate gradient method (CG), from Appendix~\ref{s:conj_methods}. We are not aware of any estimate sequence that can be used to prove that CG reaches the lower bound from Theorem~\ref{thm:smooth_LB}.
\end{itemize}
These discussions can be extended to the strongly convex setting, which we now address.

\section{Acceleration under Strong Convexity}
Before designing faster methods that exploit strong convexity, we briefly describe the benefits and limitations of this additional assumption. Roughly speaking, strong convexity guarantees that the gradient gets larger further away from the optimal solution. One way of looking at it is as follows: a function $f$ is $L$-smooth and $\mu$-strongly convex if and only if there exists some $(L-\mu)$-smooth convex function $\tilde{f}$ such that
\[ f(x) = \tilde{f}(x)+\frac{\mu}{2}\|x-x_\star\|^2_2,\]
where $x_\star$ is an optimal point for both $f$ and $\tilde{f}$. Therefore, one iteration of gradient descent can be described as follows:
\begin{equation*}
\begin{aligned}
x_{k+1}-x_\star&=x_k-x_\star-\gamma \nabla f(x_k)\\
&=x_k-x_\star-\gamma(\nabla\tilde{f}(x_k)+\mu (x_k-x_\star))\\
&=(1-\gamma \mu)(x_k-x_\star)-\gamma\nabla\tilde{f}(x_k).
\end{aligned}
\end{equation*}
We see that for sufficiently small step sizes $\gamma$, there is an additional \emph{contraction} effect due to the factor $(1-\gamma\mu)$, as compared to the effect that gradient descent has on smooth convex functions such as $\tilde{f}$. In what follows, we adapt our proofs to develop accelerated methods in the strongly convex case. Because the smooth strongly convex functions are sandwiched between two quadratic functions, these assumptions are of course much more restrictive than smoothness alone. 

\subsection{Gradient Descent and Strong Convexity}\label{s:gd_str_cvx}
As in the smooth convex case, the smooth strongly convex case can be studied through potential functions. There are many ways to prove convergence rates for this setting, but we only consider one that allows us to recover the $\mu=0$ case as its limit such that the results are well-defined even in degenerate cases. The next proof is essentially the same as that for the smooth convex case in Theorem~\ref{thm:pot_GM}, and the same inequalities are used, with strong convexity instead of convexity. The potential is only slightly modified, thereby allowing $A_k$ to have a geometric growth rate:
\[\phi_k \defeq A_k (f(x_k)-f_\star)+\frac{L+\mu A_k}{2}\|x_k-x_\star\|^2_2.\]
For notational convenience, we use $\cond=\tfrac{\mu}{L}$ to denote the inverse condition ratio. This quantity plays a key role in the geometric convergence of first-order methods in the presence of strong convexity.

\begin{theorem}\label{thm:pot_GM_strcvx}
Let $f$ be an $L$-smooth $\mu$-strongly (possibly with $\mu=0$) convex function, $x_\star\in\mathrm{argmin}_x\, f(x)$, and $k\in\mathbb{N}$. For any $A_k\geq 0$ and any $x_k$, it holds that
\begin{equation*}
\begin{aligned}
A_{k+1}(f(x_{k+1})-f_\star)&+\frac{L+\mu A_{k+1}}{2}\|x_{k+1}-x_\star\|^2_2\\
&\leq A_{k}(f(x_{k})-f_\star)+\frac{L+\mu A_k}{2}\|x_{k}-x_\star\|^2_2,
\end{aligned}
\end{equation*}
with $x_{k+1}=x_k-\tfrac1L \nabla f(x_k)$, $A_{k+1}=(1+A_k)/(1-\cond)$, and $\cond=\tfrac{\mu}{L}$.
\end{theorem}
\begin{proof} The proof consists of performing a weighted sum of the following inequalities.
\begin{itemize}
    \item Strong convexity of $f$ between $x_k$ and $x_\star$, with weight $\lambda_1=A_{k+1}-A_k$:
    \[ 0\geq f(x_k)-f_\star+\langle \nabla f(x_k); x_\star-x_k\rangle+\frac{\mu}{2}\|x_\star-x_k\|^2_2.\]
    \item Smoothness of $f$ between $x_k$ and $x_{k+1}$ with weight $\lambda_2=A_{k+1}$
    \[ 0\geq f(x_{k+1})-f(x_k)-\langle \nabla f(x_k);x_{k+1}-x_k\rangle-\frac{L}{2}\lVert x_k-x_{k+1}\rVert^2_2.\]
\end{itemize}
This weighted sum yields a valid inequality:
\begin{equation*}
\begin{aligned}
0\geq &\lambda_1[f(x_k)-f_\star+\langle \nabla f(x_k); x_\star-x_k\rangle+\frac{\mu}{2}\|x_\star-x_k\|^2_2]\\&+\lambda_2[f(x_{k+1})-f(x_k)-\langle \nabla f(x_k);x_{k+1}-x_k\rangle-\frac{L}{2}\lVert x_k-x_{k+1}\rVert^2_2].
\end{aligned}
\end{equation*}
Using $x_{k+1}=x_k-\tfrac1L \nabla f(x_k)$, this inequality can be rewritten exactly as
\begin{equation*}
\begin{aligned}
A_{k+1}&(f(x_{k+1})-f_\star)+\frac{L+\mu A_{k+1}}{2}\lVert x_{k+1}-x_\star\rVert^2_2\\
\leq& A_k(f(x_{k})-f_\star)+\frac{L+\mu A_k}{2}\lVert x_{k}-x_\star\rVert^2_2\\
&-\frac{(1-\cond)A_{k+1}-1}{2L}\lVert \nabla f(x_k)\rVert^2_2\\
&+((1-\cond) A_{k+1}- A_k-1)\langle \nabla f(x_k); x_k-x_\star\rangle.
\end{aligned}
\end{equation*}
The desired inequality follows from $A_{k+1}=(1+A_k)/(1-\cond)$ and the sign of $A_k$, making one of the last two terms nonpositive and the other equal to zero, thus reaching
\begin{equation*}
\begin{aligned}
A_{k+1}(f(x_{k+1})-f_\star)&+\frac{L+\mu A_{k+1}}{2}\|x_{k+1}-x_\star\|^2_2\\
&\leq A_{k}(f(x_{k})-f_\star)+\frac{L+\mu A_k}{2}\|x_{k}-x_\star\|^2_2.\qedhere
\end{aligned}
\end{equation*}
\end{proof}

From this theorem, we observe that adding strong convexity to the problem allows $A_k$ to follow a geometric rate given by $(1-\cond)^{-1}$ (where we again denote by $\cond=\tfrac{\mu}{L}$ the inverse condition number). The corresponding iteration complexity of gradient descent to find an approximate solution $f(x_k)-f_\star\leq \epsilon$ for  smooth strongly convex minimization is therefore $O(\tfrac{L}{\mu}\log\tfrac1\epsilon)$. This rate is essentially tight, as can be verified on quadratic functions (see, e.g., Section~\ref{c-Cheb}), and it follows from the following corollary whose result can be translated to iteration complexity using the same arguments as in Corollary~\ref{cor:iter_comp}.

\begin{corollary}\label{cor:pot_gm_strcvx}
Let $f$ be an $L$-smooth $\mu$-strongly convex function, and $x_\star\in\mathrm{argmin}_x\, f(x)$. For any $N\in\mathbb{N}$, the iterates of gradient descent with step size $\gamma_0=\gamma_1=\hdots=\gamma_N=\tfrac1L$ satisfy
\[f(x_N)-f_\star \leq \frac{\mu\|x_0-x_\star\|^2_2}{2((1-\cond)^{-N}-1)},\]
with the inverse condition number $\cond=\tfrac{\mu}{L}$.
\end{corollary}
\begin{proof}
Following the reasoning of~\eqref{eq:chained_pot}, we recursively use Theorem~\ref{thm:pot_GM} starting with $A_0=0$; that is,
\[ A_N (f(x_N)-f_\star)\leq \phi_N\leq \hdots\leq \phi_0=\frac{L}{2}\|x_0-x_\star\|^2_2,\]
and we notice that the recurrence equation $A_{k+1}=(A_k+1)/(1-\cond)$ has the solution $A_k=((1-\cond)^{-k}-1)/\cond$. The final bound is obtained by using $f(x_N)-f_\star\leq\tfrac{L\|x_0-x_\star\|^2_2}{2A_N}$ again.
\end{proof}
Note that as $\mu\rightarrow0$, the result of Corollary~\ref{cor:pot_gm_strcvx} tends to that of Corollary~\ref{cor:pot_gm}.
\begin{remark}[Lower bounds] 
As in the smooth convex case, one can derive lower complexity bounds for smooth strongly convex optimization. Using the lower bounds from smooth strongly convex quadratic minimization (for which Chebyshev's methods have optimal iteration complexity), one can conclude that no black-box first-order method can behave better than $f(x_k)-f_\star=O(\rho^k)$ with $\rho=\tfrac{(1-\sqrt{\cond})^2}{(1+\sqrt{\cond})^2}$ (see Section~\ref{c-Cheb}). In other words, lower complexity bounds from the quadratic optimization setting have the form $f(x_k)-f_\star=\Omega(\rho^k)$. We refer the reader to~\citet{Nest03a,nemirovsky1992information} for more details.

For smooth strongly convex problems beyond quadratics, this lower bound can be improved to $f(x_k)-f_\star=\Omega((1-\sqrt{\cond})^{2k})$ as provided in~\citep[Corollary~4]{drori2021exact}. In this context, we see that Nesterov's acceleration satisfies \[f(x_k)-f_\star=O((1-\sqrt{\cond})^k).\] That is, it has an $O(\sqrt{\tfrac{L}{\mu}}\log \tfrac1\epsilon)$ iteration complexity (using similar simplifications as those of Corollary~\ref{cor:iter_comp}), reaching the lower complexity bound up to a constant factor. As for the optimized gradient method provided in Section~\ref{s:OGM}, an optimal method for the smooth strongly convex case is detailed in Section~\ref{s:ogm_sc}, and it can be shown to match exactly the corresponding worst-case lower complexity bound.
\end{remark}

\subsection{Acceleration for Smooth Strongly Convex Objectives}
To adapt our proofs of convergence of accelerated methods to the strongly convex case, we need to make a small adjustment to the shape of the previous accelerated method
\begin{equation}\label{eq:FGM_1_strcvx}
\begin{aligned}
y_{k}&= x_k+\tau_k (z_k-x_k)\\
x_{k+1}&=y_{k}-\alpha_k \nabla f(y_{k})\\
z_{k+1}&=(1-\tfrac{\mu}{L} \delta_k) z_k+\tfrac{\mu}{L} \delta_k y_k - \gamma_k \nabla f(y_{k}).
\end{aligned}
\end{equation}
As discussed below, there is an optimized gradient method for smooth strongly convex minimization, similar to OGM for the smooth convex setting (see Section~\ref{s:OGM}), with this structure (details in Section~\ref{s:ogm_sc}). Following this scheme, Nesterov's method for strongly convex problems is presented in Algorithm~\ref{alg:FGM_1_strconvex}. As in the smooth convex case, we detail several of its convenient reformulations in Algorithm~\ref{alg:FGM_2_strconvex} and Algorithm~\ref{alg:FGM_3_strconvex}. The corresponding equivalences are established in Appendix~\ref{s:eq_FGM_strcvx}.

\begin{algorithm}[!ht]
  \caption{Nesterov's method, form I}
  \label{alg:FGM_1_strconvex}
  \begin{algorithmic}[1]
    \REQUIRE An $L$-smooth $\mu$-strongly (possibly with $\mu=0$) convex function~$f$ and initial $x_0$.
    \STATE \textbf{Initialize} $z_0=x_0$ and $A_0=0$; $\cond=\mu/L$ (inverse condition ratio).
    \FOR{$k=0,\ldots$}
      \STATE $A_{k+1}=\frac{2 A_k+1+\sqrt{4 A_k+4 \cond A_k^2+1}}{2 \left(1-\cond\right)}$ \COMMENT{$A_{k+1}$ solution to $(A_{k}-A_{k+1})^2-A_{k+1}-\cond A_{k+1}^2 =0$}
      \STATE Set $\tau_k=\frac{(A_{k+1}-A_k) (1+\cond A_k)}{A_{k+1}+2\cond A_k A_{k+1}-\cond A_k^2 }$ and $\delta_{k}=\frac{A_{k+1}-A_{k}}{1+\cond A_{k+1}}$
      \STATE $y_{k}=  x_k+\tau_k (z_k-x_k)$
      \STATE $x_{k+1}=y_{k}-\frac1L \nabla f(y_{k})$
      \STATE $z_{k+1}=(1-\cond\delta_k)z_k+ \cond\delta_k y_k-\frac{\delta_k}{L}\nabla f(y_{k})$
    \ENDFOR
    \ENSURE Approximate solution $x_{k+1}$.
  \end{algorithmic}
\end{algorithm}

Regarding the potential, we make the same adjustment as for gradient descent, arriving to the following theorem.

\begin{theorem}\label{thm:Nest_first_strcvx}
Let $f$ be an $L$-smooth $\mu$-strongly (possibly with $\mu=0$) convex function, $x_\star\in\mathrm{argmin}_x\, f(x)$, and $k\in\mathbb{N}$. For all $x_k,z_k\in\mathbb{R}^d$ and $A_k\geq 0$, the iterates of Algorithm~\ref{alg:FGM_1_strconvex} satisfy
\begin{equation*}
\begin{aligned}
A_{k+1}&(f(x_{k+1})-f_\star)+\frac{L+\mu A_{k+1}}{2}\|z_{k+1}-x_\star\|^2_2\\&\leq A_k (f(x_k)-f_\star)+\frac{L+\mu A_{k}}2 \|z_k-x_\star\|^2_2,
\end{aligned}
\end{equation*}
with $A_{k+1}=\frac{2 A_k+1+\sqrt{4 A_k+4\cond A_k^2+1}}{2 \left(1-\cond\right)}$ and $\cond=\tfrac{\mu}{L}$.
\end{theorem}
\begin{proof}
The proof consists of a weighted sum of the following inequalities.
\begin{itemize}
    \item Strong convexity between $x_\star$ and $y_{k}$ with weight $\lambda_1=A_{k+1}-A_k$:
    \[ f_\star \geq f(y_{k})+\langle \nabla f(y_{k});x_\star-y_k\rangle+\frac{\mu}{2}\|x_\star-y_k\|^2_2.\]
    \item Convexity between $x_k$ and $y_{k}$ with weight $\lambda_2=A_k$:
    \[ f(x_k) \geq f(y_{k})+\langle \nabla f(y_{k});x_k-y_{k}\rangle.\]
    \item Smoothness between $y_{k}$ and $x_{k+1}$ (\emph{descent lemma}) with weight $\lambda_3=A_{k+1}$
    \[ f(y_{k}) + \langle \nabla f(y_{k}); x_{k+1}-y_{k}\rangle +\frac{L}{2}\|x_{k+1}-y_{k}\|^2_2 \geq f(x_{k+1}).\]
\end{itemize}
We therefore arrive at the following valid inequality:
\begin{equation*}
\begin{aligned}
0\geq& \lambda_1 [ f(y_{k})-f_\star +\langle \nabla f(y_{k});x_\star-y_{k}\rangle+\frac{\mu}{2}\|x_\star-y_k\|^2_2]\\&+ \lambda_2[f(y_{k})-f(x_k)+\langle \nabla f(y_{k});x_k-y_{k}\rangle]\\&+\lambda_3[f(x_{k+1})-f(y_{k}) - \langle \nabla f(y_{k}); x_{k+1}-y_{k}\rangle -\frac{L}{2}\|x_{k+1}-y_{k}\|^2_2].
\end{aligned}
\end{equation*}
For the sake of simplicity, we do not substitute $A_{k+1}$ by its expression until the last stage of the reformulation. After substituting $x_{k+1}$, $z_{k+1}$ by their expressions in~\eqref{eq:FGM_1_strcvx} along with $\tau_k=\tfrac{(A_{k+1}-A_k) (1+\cond A_k)}{A_{k+1}+2\cond A_k A_{k+1}-\cond A_k^2}$, $\alpha_k=\frac1L$, $\delta_k=\frac{A_{k+1}-A_{k} }{1+\cond A_{k+1}}$, and $\gamma_k=\frac{\delta_k}{L}$, basic algebra shows that the previous inequality can be reorganized as 
\begin{equation*}
\begin{aligned}
A_{k+1}&(f(x_{k+1})-f_\star)+\frac{L+\mu A_{k+1}}{2}\|z_{k+1}-x_\star\|^2_2\\\leq& A_k (f(x_k)-f_\star)+\frac{L+\mu A_k}2 \|z_k-x_\star\|^2_2\\
&+\frac{ (A_{k}-A_{k+1})^2-A_{k+1}-\cond A_{k+1}^2 }{1+\cond A_{k+1}} \frac{1}{2L}\|\nabla f(y_k)\|^2_2\\
&-A_{k}^2  \frac{ (A_{k+1}-A_{k})  (1+ \cond A_{k})(1+ \cond A_{k+1})}{(A_{k+1}+2\cond A_k A_{k+1}-\cond  A_k^2)^2} \frac{ \mu}{2}\|x_k-z_k\|^2_2.
\end{aligned}
\end{equation*}
The desired statement follows from selecting $A_{k+1}\geq A_k\geq0$ such that \[(A_{k}-A_{k+1})^2-A_{k+1}-\cond A_{k+1}^2 =0,\]
thus yielding
\begin{equation*}
\begin{aligned}
A_{k+1}&(f(x_{k+1})-f_\star)+\frac{L+\mu A_{k+1}}{2}\|z_{k+1}-x_\star\|^2_2\\\leq& A_k (f(x_k)-f_\star)+\frac{L+\mu A_k}2 \|z_k-x_\star\|^2_2.\qedhere
\end{aligned}
\end{equation*}
\end{proof}

The final worst-case guarantee is obtained by using the same reasoning as before, together with a simple bound on $A_{k+1}$:
\begin{equation}\label{eq:conv_Ak_fgm_strcvx}
\begin{aligned}
A_{k+1}&=\frac{2 A_k+1+\sqrt{4 A_k+4 \cond A_k^2+1}}{2 \left(1-\cond\right)}\\&\geq \frac{2 A_{k}+\sqrt{\left(2 A_{k} \sqrt{\cond }\right)^2}}{2 \left(1-\cond \right)}=\frac{A_k}{1-\sqrt{\cond }},
\end{aligned}
\end{equation}
which means $f(x_k)-f_\star=O((1-\sqrt{\cond })^k)$ when $\mu>0$, or alternatively that $O\left(\sqrt{\tfrac{L}{\mu}}\log\tfrac1\epsilon\right)$ is the iteration complexity of obtaining an approximate solution $f(x_N)-f_\star\leq \epsilon$ (using similar simplifications as those of Corollary~\ref{cor:iter_comp}). The following corollary summarizes our result for Nesterov's method.
\begin{corollary}\label{cor:pot_fgm_strcvx}
Let $f$ be an $L$-smooth $\mu$-strongly (possibly with $\mu=0$) convex function and $x_\star\in\mathrm{argmin}_x\, f(x)$. For all $N\in\mathbb{N}$, $N\geq 1$, the iterates of Algorithm~\ref{alg:FGM_1_strconvex} satisfy
\[f(x_N)-f_\star \leq \min\left\{\frac{2}{N^2},\left(1-\sqrt{\cond}\right)^N\right\}L\|x_0-x_\star\|^2_2,\]
with $\cond=\tfrac{\mu}{L}$.
\end{corollary}
\begin{proof} Following the argument of~\eqref{eq:chained_pot}, we recursively use Theorem~\ref{thm:Nest_first_strcvx} with $A_0=0$, together with the bounds on $A_N$ for the smooth convex case~\eqref{eq:conv_Ak_fgm} and for the smooth strongly convex one~\eqref{eq:conv_Ak_fgm_strcvx}. (Note that $A_{k+1}$ is an increasing function of $\mu$, and hence the bound for the smooth case remains valid in the smooth strongly convex one.) We have $A_1=\tfrac{1}{1-\cond}=\tfrac1{(1-\sqrt{\cond})(1+\sqrt{\cond})}\geq \tfrac12 (1-\sqrt{\cond})^{-1}$, thus reaching $A_N\geq \tfrac12 (1-\sqrt{\cond})^{-N}$.
\end{proof}
\begin{remark}Before moving to the next section, we mention that another direct consequence of the potential inequality above (Theorem~\ref{thm:Nest_first_strcvx}) is that $z_k$ may also serve as an approximate solution to $x_\star$ when $\mu>0$. Indeed, by using the inequality
\[ \frac{L+\mu A_{N}}2 \|z_N-x_\star\|^2_2 \leq \phi_N \leq \hdots\leq\phi_0= \frac{L}{2}\|x_0-x_\star\|^2_2,\]
it follows that \[\|z_N-x_\star\|^2_2\leq \frac{1}{1+\cond A_N}\|x_0-x_\star\|^2_2\leq\frac{2 (1-\sqrt{\cond})^N}{2(1-\sqrt{\cond})^N+\cond}\|x_0-x_\star\|^2_2 ,\]
and hence that $\|z_N-x_\star\|^2_2=O((1-\sqrt{\cond})^N)$. Therefore it also follows that \[f(z_N)-f_\star\leq \frac{L}{2}\|z_N-x_\star\|^2_2=O((1-\sqrt{\cond})^N).\]
In addition, since $y_N$ is a convex combination of $x_N$ and $z_N$, the same conclusion holds for $\|y_N-x_\star\|^2_2$ and $f(y_N)-f_\star$. Similar observations also apply to other variants of accelerated methods when $\mu>0$.
\end{remark}
\subsection{A Simplified Stationary Method with Constant Momentum}\label{s:const_momentum}
Important simplifications are often made to the Nesterov's method in the strongly convex case where $\mu>0$. Several approaches produce the same method, known as the ``constant momentum'' version of Nesterov's accelerated gradient. We derive this version by observing that the asymptotic (or stationary) behavior of Algorithm~\ref{alg:FGM_1_strconvex} can be characterized explicitly. In particular, when $k\rightarrow\infty$, it is clear that $A_k\rightarrow\infty$ as well. We can thus take the limits of all parameters as $A_k\rightarrow\infty$, to obtain a corresponding ``limit/stationary method.'' This is similar in spirit to the result showing that Polyak's heavy-ball method is the asymptotic version of Chebyshev's method, discussed in Section~\ref{s:HeavyBall}. First, the convergence rate is obtained as
\[ \lim_{A_k\rightarrow\infty} \frac{A_{k+1}}{A_k}=\left(1-\sqrt{\cond}\right)^{-1}.\]
By taking the limits of all the algorithmic parameters, that is,
\[ \lim_{A_k\rightarrow\infty} \tau_k=\frac{\sqrt{\cond}}{1+\sqrt{\cond}} ,\quad \lim_{A_k\rightarrow\infty} \delta_k=\frac1{\sqrt{\cond}}, \]
we obtain Algorithm~\ref{alg:FGM_1_strcvx_constmomentum} and its equivalent, probably most well-known, second form, provided as Algorithm~\ref{alg:FGM_2_strcvx_constmomentum}.
\begin{algorithm}[!ht]
  \caption{Nesterov's method, form I, constant momentum}
  \label{alg:FGM_1_strcvx_constmomentum}
  \begin{algorithmic}[1]
    \REQUIRE $L$-smooth $\mu$-strongly convex function $f$ and initial point $x_0$.
    \STATE \textbf{Initialize} $z_0=x_0$ and $A_0>0$; $\cond=\mu/L$ (inverse condition ratio).
    \FOR{$k=0,\ldots$}
      \STATE $A_{k+1}=\frac{A_k}{1-\sqrt{\cond}}$\COMMENT{Only for the proof/relation to previous methods.}
      \STATE $y_{k}=  x_k+\frac{\sqrt{\cond}}{1+\sqrt{\cond}}(z_k-x_k)$
      \STATE $x_{k+1}=y_{k}-\frac1L \nabla f(y_{k})$
      \STATE $z_{k+1}=\left(1-\sqrt{\cond}\right)z_k+\sqrt{\cond}\left(y_k- \frac1{\mu}\nabla f(y_{k})\right)$
    \ENDFOR
    \ENSURE Approximate solutions $(y_k,x_{k+1},z_{k+1})$.
  \end{algorithmic}
\end{algorithm}
\begin{algorithm}[!ht]
  \caption{Nesterov's method, form II, constant momentum}
  \label{alg:FGM_2_strcvx_constmomentum}
  \begin{algorithmic}[1]
    \REQUIRE $L$-smooth $\mu$-strongly convex function $f$ and initial point $x_0$.
    \STATE \textbf{Initialize} $y_0=x_0$ and $A_0>0$; $\cond=\mu/L$ (inverse condition ratio).
    \FOR{$k=0,\ldots$}
      \STATE $A_{k+1}=\frac{A_k}{1-\sqrt{\cond}}$\COMMENT{Only for the proof/relation to previous methods.}
      \STATE $x_{k+1}=y_k-\tfrac{1}{L}\nabla f(y_k)$
      \STATE $y_{k+1}=x_{k+1}+\frac{1-\sqrt{\cond}}{1+\sqrt{\cond}}(x_{k+1}-x_k)$
    \ENDFOR
    \ENSURE Approximate solutions $(y_k,x_{k+1})$.
  \end{algorithmic}
\end{algorithm}

\enlargethispage{\baselineskip}
From a worst-case analysis perspective, these simplifications correspond to using a Lyapunov function obtained by dividing the potential function of Theorem~\ref{thm:Nest_first_strcvx} by $A_k$ and then taking the limit of the inequality:
\begin{equation*}
\begin{aligned}
\rho^{-1}\left(f(x_{k+1})-f_\star+\frac{\mu}{2}\|z_{k+1}-x_\star\|^2_2\right)\leq f(x_{k})-f_\star+\frac{\mu}{2}\|z_{k}-x_\star\|^2_2,
\end{aligned}
\end{equation*}
with $\rho=(1-\sqrt{\cond})$.

\begin{theorem}\label{thm:nest_constmom} Let $f$ be an $L$-smooth $\mu$-strongly convex function, $x_\star=\mathrm{argmin}_x\, f(x)$ and $k\in\mathbb{N}$. For any $x_k,z_k\in\mathbb{R}^d$, the iterates of Algorithm~\ref{alg:FGM_1_strcvx_constmomentum} (or equivalently those of Algorithm~\ref{alg:FGM_2_strcvx_constmomentum}) satisfy
\begin{equation*}
\begin{aligned}
\rho^{-1}\left(f(x_{k+1})-f_\star+\frac{\mu}{2}\|z_{k+1}-x_\star\|^2_2\right)\leq f(x_{k})-f_\star+\frac{\mu}{2}\|z_{k}-x_\star\|^2_2.
\end{aligned}
\end{equation*}
\end{theorem}
\begin{proof}
Let $A_k>0$ and $A_{k+1}=A_k/(1-\sqrt{\cond})$. The proof is essentially the same as that of as Theorem~\ref{thm:Nest_first_strcvx}. That is, the weights used in this proof are those used in Theorem~\ref{thm:Nest_first_strcvx} divided by $A_k$, leading to a slight variation in the reformulation of the weighted sum, and the following valid inequality:
\begin{equation*}
\begin{aligned}
0\geq& \lambda_1 [ f(y_{k})-f_\star +\langle \nabla f(y_{k});x_\star-y_{k}\rangle+\frac{\mu}{2}\|x_\star-y_k\|^2_2]\\&+ \lambda_2[f(y_{k})-f(x_k)+\langle \nabla f(y_{k});x_k-y_{k}\rangle]\\&+\lambda_3[f(x_{k+1})-f(y_{k}) - \langle \nabla f(y_{k}); x_{k+1}-y_{k}\rangle -\frac{L}{2}\|x_{k+1}-y_{k}\|^2_2],
\end{aligned}
\end{equation*}
with weights
\begin{align*}
\lambda_1=&\frac{A_{k+1}-A_k}{A_k}=\frac{\sqrt{\cond}}{1-\sqrt{\cond}},\,\,\, \\
\lambda_2=&\frac{A_k}{A_k}=1,\, \text{ and }\\
\lambda_3=&\frac{A_{k+1}}{A_k}=(1-\sqrt{\cond})^{-1}.
\end{align*}
By substituting 
\begin{equation*}
\begin{aligned}
y_k&=x_k+\frac{\sqrt{\cond}}{1+\sqrt{\cond}}(z_k-x_k)\\
x_{k+1}&=y_k-\frac1L \nabla f(y_k)\\
z_{k+1}&=\left(1-\sqrt{\cond}\right)z_k+\sqrt{\cond}\left(y_k-\frac1\mu \nabla f(y_k)\right),
\end{aligned}
\end{equation*}
we arrive at the following valid inequality:
\begin{equation*}
\begin{aligned}
\rho^{-1}&\left(f(x_{k+1})-f_\star+\frac{\mu}{2}\|z_{k+1}-x_\star\|^2_2\right)\\&\leq f(x_{k})-f_\star+\frac{\mu}{2}\|z_{k}-x_\star\|^2_2-\frac{\sqrt{\cond}}{(1+\sqrt{\cond})^2}\frac{\mu}{2}\|x_k-z_k\|^2_2.
\end{aligned}
\end{equation*}
We reach the desired statement from the last term being nonpositive:
\[ \rho^{-1}\left(f(x_{k+1})-f_\star+\frac{\mu}{2}\|z_{k+1}-x_\star\|^2_2\right)\leq f(x_{k})-f_\star+\frac{\mu}{2}\|z_{k}-x_\star\|^2_2.\qedhere\]
\end{proof}
\begin{corollary}\label{cor:pot_fgm_strcvx_cst_mom}
Let $f$ be an $L$-smooth $\mu$-strongly convex function and $x_\star\in\mathrm{argmin}_x\, f(x)$. For all $N\in\mathbb{N}$, $N\geq 1$, the iterates of Algorithm~\ref{alg:FGM_1_strcvx_constmomentum} (or equivalently of Algorithm~\ref{alg:FGM_2_strcvx_constmomentum}) satisfy
\[f(x_N)-f_\star \leq (1-\sqrt{\cond})^N\left(f(x_0)-f_\star+\frac{\mu}{2}\|x_0-x_\star\|^2_2\right),\]
with $\cond=\tfrac{\mu}{L}$.
\end{corollary}
\begin{proof} The desired result directly follows from Theorem~\ref{thm:nest_constmom} with 
\[ \phi_k \defeq \rho^{-k}\left(f(x_k)-f_\star+\frac{\mu}{2}\|z_k-x_\star\|^2_2\right),\]
$\rho=1-\sqrt{\cond}$, $z_0=x_0$, and $\rho^{-N} (f(x_N)-f_\star)\leq\phi_N$.
\end{proof}
\begin{remark} In view of Section~\ref{s:estimate}, one can also find estimate sequence interpretations of Algorithms~\ref{alg:FGM_1_strconvex} and~\ref{alg:FGM_1_strcvx_constmomentum} from their respective potential functions.
\end{remark}
\begin{remark}
A few works on accelerated methods focus on understanding this particular instance of Nesterov's method. Our analysis here is largely inspired by that of~\citet{Nest03a}, but such potentials can be obtained in different ways, see, for example~\citep{wilson2016lyapunov,shi2018understanding,bansal2019potential}.
\end{remark}

\section{Recent Variants of Accelerated Methods}
In this section, we first push the reasoning in terms of \emph{potential functions} to its limit. We present the \emph{information-theoretic exact method}~\citep{drori2021optimal}, which generalizes the optimized gradient descent in the strongly convex case. Similar to Nesterov's method with constant momentum, the \emph{information-theoretic exact method} has a limit case that is known as the \emph{triple momentum method}~\citep{van2017fastest}. We then discuss a more geometric variant, known as \emph{geometric descent}~\citep{bubeck2015geometric} or \emph{quadratic averaging}~\citep{drusvyatskiy2018optimal}.

\subsection{An Optimal Method for Strongly Convex Minimization}\label{s:ogm_sc}

It turns out that there also exist optimal gradient methods for smooth strongly convex minimization that are similar to the optimized gradient method for smooth convex minimization. Such methods can be obtained by solving a minimax problem similar to~\eqref{eq:pep_opt_h} with different objectives. 

The following scheme is optimal for the criterion $\tfrac{\|z_N-x_\star\|^2_2}{\|x_0-x_\star\|^2_2}$, reaching the exact worst-case lower complexity bound for this criterion, as discussed below. In addition, this method reduces to the OGM (see Section~\ref{s:OGM}) when $\mu=0$ by using the correspondence $A_{k+1}=4\theta_{k,N}^2$ (for $k<N$). Therefore, this method is \emph{doubly optimal}, i.e., optimal according to two criteria, in the sense that it also achieves the lower complexity bound for $\tfrac{f(y_N)-f_\star}{\|x_0-x_\star\|^2_2}$ when $\mu=0$, using the last iteration adjustment from Lemma~\ref{thm:pot_OGM_final}.

The following analysis is reminiscent of Nesterov's method in Algorithm~\ref{alg:FGM_1_strconvex} but also of the optimized gradient method and its proof (see Theorem~\ref{thm:pot_OGM}). That is, the known potential function for the information-theoretic exact method (ITEM) relies on inequality~\eqref{eq:interpolation_strcvx}, not only for its proof but also simply to simply show that it is nonnegative, which follows from instantiating~\eqref{eq:interpolation_strcvx} at $y=x_\star$. The following analyses can be found almost verbatim in~\citep{drori2021optimal}. The main proof of this section is particularly algebraic, but it can be reasonably skipped as it follows from similar ideas found in previous developments.

\begin{algorithm}[!ht]
  \caption{Information-theoretic exact method (ITEM)}
  \label{alg:ITEMS}
  \begin{algorithmic}[1]
    \REQUIRE $L$-smooth $\mu$-strongly (possibly with $\mu=0$) convex function $f$ and initial $x_0$.
    \STATE \textbf{Initialize} $z_0=x_0$ and $A_0=0$; $\cond=\mu/L$ (inverse condition ratio).
    \FOR{$k=0,\ldots$}
      \STATE $A_{k+1}= \frac{(1+\cond) A_k+2 \left(1+\sqrt{(1+A_k) (1+\cond  A_k)}\right)}{(1-\cond )^2}$
      \STATE Set $\tau_k=1-\frac{A_k}{(1-\cond)A_{k+1}}\text{, and }  \delta_k=\frac1{2}\frac{(1-\cond )^2 A_{k+1}-(1+\cond ) A_k}{1+\cond+\cond  A_k}$
      \STATE ${y_{k}}= x_{k} + \tau_k (z_k-x_k)$
      \STATE $x_{k+1}={y_{k}-\tfrac1L \nabla f(y_{k})}$
    \STATE $z_{k+1}=(1-\cond\delta_k)z_k+\cond \delta_k y_{k}-\tfrac{\delta_k}L  \nabla f(y_{k})$
    \ENDFOR
    \ENSURE Approximate solutions $(y_k,x_{k+1},z_{k+1})$.
  \end{algorithmic}
\end{algorithm}

\begin{theorem}\label{thm:ITEMS_pot}
Let $f$ be an $L$-smooth $\mu$-strongly (possibly with $\mu=0$) convex function, $x_\star\in\mathrm{argmin}_x\, f(x)$ and $k\in\mathbb{N}$. For all $y_{k-1},z_k\in\mathbb{R}^d$ and $A_k\geq 0$, the iterates of Algorithm~\ref{alg:ITEMS} satisfy
\[ \phi_{k+1}\leq \phi_k,\]
with
\begin{equation*}
\begin{aligned}
\phi_k \defeq &A_{k}\bigg[f(y_{k-1})-f_\star-\frac{1}{2L}\|\nabla f(y_{k-1})\|^2_2\\
& \quad-\frac{\mu}{2(1-\cond)}\|y_{k-1}-\tfrac1L \nabla f(y_{k-1})-x_\star\|^2_2\bigg]\\
&+\frac{L+\mu A_{k}}{1-\cond}\|z_{k}-x_\star\|^2_2
\end{aligned}
\end{equation*}
and $A_{k+1}=\frac{(1+\cond) A_k+2 \left(1+\sqrt{(1+A_k) (1+\cond  A_k)}\right)}{(1-\cond )^2}$.
\end{theorem}

\begin{proof}
We first perform a weighted sum of two inequalities from Theorem~\ref{eq:interpolation_strcvx}.
\begin{itemize}
    \item Smoothness and strong convexity between ${y_{k-1}}$ and ${y_{k}}$ with weight $\lambda_1=A_k$:
    \begin{equation*}
    \begin{aligned}
    {f(y_{k-1})} \geq& f(y_{k})+\langle \nabla f(y_{k});y_{k-1}-y_{k}\rangle\\&+\frac1{2L}\|\nabla f(y_{k})-\nabla f(y_{k-1})\|^2_2\\&+\frac{\mu}{2(1-\cond)}\|y_{k}-y_{k-1}-\tfrac1L (\nabla f(y_{k})-\nabla f(y_{k-1}))\|^2_2.
    \end{aligned}
    \end{equation*}
    \item Smoothness and strong convexity of $f$ between $x_\star$ and {$y_{k}$} with weight $\lambda_2=A_{k+1}-A_k$:
    \begin{equation*}
    \begin{aligned}
    f_\star \geq& f(y_{k})+\langle \nabla f(y_{k});x_\star-y_{k}\rangle+\frac1{2L}\|\nabla f(y_{k})\|^2_2\\&+\frac{\mu}{2(1-\cond)}\|y_{k}-x_\star-\tfrac1L \nabla f(y_{k})\|^2_2.
    \end{aligned}
    \end{equation*}
\end{itemize}
By summing up and reorganizing these two inequalities (without substituting $A_{k+1}$ by its expression, for simplicity), we arrive at the following valid inequality:
\begin{equation*}
\begin{aligned}
0\geq \lambda_1 &\bigg[ f(y_{k})-f(y_{k-1})+\langle \nabla f(y_{k});y_{k-1}-y_{k}\rangle\\
& \quad +\frac1{2L}\|\nabla f(y_{k})-\nabla f(y_{k-1})\|^2_2\\
&\quad+ {\frac{\mu}{2(1-\cond)}\|y_{k}-y_{k-1}-\tfrac1L (\nabla f(y_{k})-\nabla f(y_{k-1}))\|^2_2}\bigg]\\
+\lambda_2 &\bigg[f(y_{k})-f_\star+\langle \nabla f(y_{k});x_\star-y_{k}\rangle+\frac1{2L}\|\nabla f(y_{k})\|^2_2\\&\quad+\frac{\mu}{2(1-\cond)}\|y_{k}-x_\star-\tfrac1L \nabla f(y_{k})\|^2_2\bigg].
\end{aligned}
\end{equation*}
By substituting the expressions of $y_k$ and $z_{k+1}$ with
\begin{equation*}
\begin{aligned}
y_k&=(1-\tau_k)\left(y_{k-1}-\tfrac1L\nabla f(y_{k-1})\right)+\tau_k z_k\\
z_{k+1}&=(1-q\delta_k)z_k+q\delta_k y_k -\tfrac{\delta_k}{L}\nabla f(y_k)\
\end{aligned}
\end{equation*}
(noting that this substitution is valid even for $k=0$ since $A_0=0$ in that case and hence, $\tau_0=1$ and $y_0=z_0$), the previous inequality can be reformulated exactly as
\begin{equation*}
\begin{aligned}
\phi_{k+1}
\leq& \phi_k-\frac{LK_1 }{1-\cond} P(A_{k+1}) \|z_k-x_\star\|^2_2\\
&+\frac{K_2}{4L(1-\cond)} P(A_{k+1})\\&\quad\times\| (1-\cond) A_{k+1}{\nabla f(y_{k})}-\mu A_k ({x_{k}}-x_\star)+K_3\mu (z_k-x_\star)\|^2_2,
\end{aligned}
\end{equation*}
with the three parameters (which are well-defined given that $0\leq\mu<L$ and $A_k,A_{k+1}\geq 0$)
\begin{equation*}
\begin{aligned}
K_1&=\frac{ \cond ^2 }{(1+\cond)^2+(1-\cond)^2\cond A_{k+1}}\\
K_2&=\frac{(1+\cond)^2+(1-\cond )^2 \cond  A_{k+1}}{(1-\cond)^2\left(1+\cond+\cond  A_k\right)^2 A_{k+1}^2}\\
K_3&=(1+\cond)\frac{( 1+\cond) A_k- (1-\cond)(2+\cond A_k)A_{k+1}}{(1+\cond)^2+(1-\cond )^2 \cond  A_{k+1}},
\end{aligned}
\end{equation*}
as well as \[P(A_{k+1})=(A_k-(1-\cond ) A_{k+1})^2-4 A_{k+1} (1+\cond  A_k).\] To obtain the desired inequality, we select $A_{k+1}$ such that $A_{k+1}\geq A_k$ and $P(A_{k+1})=0$, thereby demonstrating the claim $\phi_{k+1}\leq \phi_k$ and the choice of $A_{k+1}$.
\end{proof}

The final bound for this method is obtained after the usual growth analysis of the sequence $A_k$, as follows. When $\mu=0$, we have
\[A_{k+1}=2+A_k+2\sqrt{1+A_k}\geq 2+A_k+2\sqrt{A_k}\geq (1+\sqrt{A_k})^2,\]
reaching $\sqrt{A_{k+1}}\geq 1+\sqrt{A_k}$ and hence $\sqrt{A_k}\geq k$ and $A_k\geq k^2$. When $\mu>0$, we can use an alternate bound:
\begin{align*}
A_{k+1}=&\frac{(1+\cond) A_k+2 \left(1+\sqrt{(1+A_k) (1+\cond  A_k)}\right)}{(1-\cond )^2}\\
\geq& \frac{(1+\cond)A_k+2 \sqrt{\cond  A_k^2}}{(1-\cond)^2}=\frac{A_k}{(1-\sqrt{\cond})^2},
\end{align*} 
therefore achieving similar bounds as before. In this case, we only emphasize the convergence result for $\|z_N-x_\star\|$ since it corresponds to the lower complexity bound for smooth strongly convex minimization (provided below).
\begin{corollary}\label{corr:finalbound}
Let $f\in\mathcal{F}_{\mu,L}$, and denote $\cond=\mu/L$. For any $x_0=z_0\in\mathbb{R}^d$ and $N\in\mathbb{N}$ with $N\geq 1$, the iterates of Algorithm~\ref{alg:ITEMS} satisfy
\begin{equation*}
\begin{aligned}
\|z_N-x_\star\|^2_2 &\leq \frac{1}{1+\cond A_N}\|z_0-x_\star\|^2_2\leq  \frac{(1-\sqrt{\cond})^{2N}}{(1-\sqrt{\cond})^{2N}+\cond}\|z_0-x_\star\|^2_2.
\end{aligned}
\end{equation*}
\end{corollary}
\begin{proof} From Theorem~\ref{thm:ITEMS_pot}, we get
\[ \phi_N\leq\phi_{N-1}\leq \hdots\leq \phi_0=\frac{L}{1-q}\|z_0-x_\star\|^2_2.\]
From~\eqref{eq:interpolation_strcvx}, we have that $\phi_N\geq \frac{L+A_N\mu}{1-q}\|z_N-x_\star\|^2_2$. It remains to use the bounds on $A_N$. That is, by using $A_1=\tfrac{4}{(1-\cond)^2}=\tfrac{4}{(1+\sqrt{\cond})^2(1-\sqrt{\cond})^2}\geq (1-\sqrt{\cond})^{-2}$, we have $A_N\geq (1-\sqrt{\cond})^{-2N},$ which concludes the proof.
\end{proof}

Before concluding, we mention that the algorithm is non-improvable when minimizing large-scale smooth strongly convex functions in the following sense.

\begin{theorem}\citep[Corollary 4]{drori2021exact}\label{thm:smooth_sc_LB} Let $0\leq \mu<L<\infty$, $d,N\in\mathbb{N}$ with $d\geq 2N+1$. For any black-box first-order method that performs at most $N$ calls to the first-order oracle $(f(\cdot),\nabla f(\cdot))$, there exists a function $f\in\mathcal{F}_{\mu,L}(\mathbb{R}^d)$ and $x_0\in\mathbb{R}^d$ such that
\[ \|x_N-x_\star\|^2_2 \geq \frac{1}{1+\cond A_N}\lVert x_0-x_\star\rVert^2_2,\]
where $x_\star\in\mathrm{argmin}_x\, f(x)$, $x_N$ is the output of the method under consideration, $x_0$ is its input, and $A_N$ is defined as in Algorithm~\ref{alg:ITEMS}.
\end{theorem}

\begin{remark} Just as for the optimized gradient method from Section~\ref{s:OGM}, the ITEM might serve as a template for the design of other accelerated schemes. However, it has the same caveats as the optimized gradient method, which are also similar to those of the triple momentum method, presented in the next section. As emphasized in Section~\ref{sec:LO_OGM}, it is unclear how to generalize the ITEM to broader classes of problems, e.g., problems involving constraints.
\end{remark}

\subsection{The Triple Momentum Method}\label{s:tmm}
The \emph{triple momentum method} (TMM) is due to~\citet{van2017fastest} and is reminiscent of Nesterov's method with constant momentum, provided as Algorithm~\ref{alg:FGM_1_strcvx_constmomentum}. It corresponds to the asymptotic behavior of the information-theoretic exact method, just as Nesterov's accelerated method with constant momentum is the limit case of Nesterov's method (see Section~\ref{s:const_momentum}) and as Polyak's heavy-ball is the limit case of Chebyshev's method (see Section~\ref{s:HeavyBall}). Indeed, considering Algorithm~\ref{alg:ITEMS}, one can explicitly compute
\[ \lim_{A_k\rightarrow\infty} \frac{A_{k+1}}{A_k}=\left(1-\sqrt{\cond}\right)^{-2}\] as well as
\[ \lim_{A_k\rightarrow\infty} \tau_k=1-\frac{1-\sqrt{\cond}}{1+\sqrt{\cond}} ,\quad \lim_{A_k\rightarrow\infty} \delta_k=\frac1{\sqrt{\cond}}. \]

\begin{algorithm}[!ht]
  \caption{Triple momentum method (TMM)}
  \label{alg:TMM}
  \begin{algorithmic}[1]
    \REQUIRE $L$-smooth $\mu$-strongly convex function $f$ and initial point $x_0$.
    \STATE \textbf{Initialize} $y_{-1}=z_0=x_0$; $\cond=\mu/L$ (inverse condition ratio).
    \FOR{$k=0,\ldots$}
      \STATE $A_{k+1}=\frac{A_k}{(1-\sqrt{\cond})^2}$\COMMENT{Only for the proof/relation to previous methods.}
      \STATE $y_{k}=\frac{1-\sqrt{\cond}}{1+\sqrt{\cond}} \left(y_{k-1}-\tfrac1L \nabla f(y_{k-1})\right) + \left(1-\frac{1-\sqrt{\cond}}{1+\sqrt{\cond}}\right) z_k $
      \STATE $z_{k+1}=\sqrt{\cond} \left(y_{k}-\tfrac1{\mu} \nabla f(y_{k})\right)+\left(1-\sqrt{\cond}\right)z_k$
    \ENDFOR
    \ENSURE Approximate solutions $(y_{k},z_{k+1})$.
  \end{algorithmic}
\end{algorithm}

As for the information-theoretic exact method, the known potential function for the triple momentum method relies on inequality~\eqref{eq:interpolation_strcvx}, not only for its proof but to show that it is nonnegative, which follows from instantiating~\eqref{eq:interpolation_strcvx} at $y=x_\star$.

\begin{theorem}\label{thm:tmm} Let $f$ be an $L$-smooth $\mu$-strongly convex function, $x_\star=\mathrm{argmin}_x\, f(x)$, and $k\in\mathbb{N}$. For any $x_k,z_k\in\mathbb{R}^d$, the iterates of Algorithm~\ref{alg:TMM} satisfy
\begin{equation*}
\begin{aligned}
f&(y_{k})-f_\star-\frac1{2L}\|\nabla f(y_{k})\|^2_2-\frac{\mu }{2(1-\cond)}\|y_{k}-x_\star-\tfrac1L \nabla f(y_{k})\|^2_2\\
&\quad+\frac{\mu}{1-\cond}\|z_{k+1}-x_\star\|^2_2\\
\leq& \rho^2 \bigg(f(y_{k-1})-f_\star-\frac1{2L}\|\nabla f(y_{k-1})\|^2_2\\&\quad\quad-\frac{\mu}{2(1-\cond)}\|y_{k-1}-x_\star-\tfrac1L \nabla f(y_{k-1})\|^2_2+\frac{\mu}{1-\cond}\|z_{k}-x_\star\|^2_2\bigg),
\end{aligned}
\end{equation*}
with $\rho=1-\sqrt{\cond}$.
\end{theorem}
\enlargethispage{\baselineskip}
\begin{proof} For simplicity, we consider Algorithm~\ref{alg:TMM} in the following form, parameterized by~$\rho$
\begin{equation*}
\begin{aligned}
y_{k}&=\frac{\rho}{2-\rho } \left(y_{k-1}-\tfrac1{L}\nabla f(y_{k-1})\right)+\left(1-\frac{\rho}{2-\rho }\right)z_k\\
z_{k+1}&=(1-\rho ) \left(y_{k}-\tfrac1\mu \nabla f(y_{k})\right)+\rho  z_k.
\end{aligned}
\end{equation*}
We combine the following two inequalities.
\begin{itemize}
    \item Smoothness and strong convexity between $y_{k-1}$ and $y_{k}$ with weight~$\lambda_1=\rho^2$:
    \begin{equation*}
    \begin{aligned}
     f(y_{k-1})\geq& f(y_{k})+\langle \nabla f(y_{k});y_{k-1}-y_{k}\rangle\\
     & +\frac{1}{2L}\|\nabla f(y_{k})-\nabla f(y_{k-1})\|^2_2\\
     & +\frac{\mu }{2(1-\cond)}\|y_{k}-y_{k-1}-\tfrac1L( \nabla f(y_{k})-\nabla f(y_{k-1}))\|^2_2.
    \end{aligned}
    \end{equation*}
    \item Smoothness and strong convexity between $x_\star$ and $y_{k}$ with weight~$\lambda_2$ $=1-\rho^2$:
    \begin{equation*}
    \begin{aligned}
    f_\star \geq& f(y_{k})+\langle \nabla f(y_{k});x_\star-y_{k}\rangle +\frac{1}{2L}\|\nabla f(y_{k})\|^2_2\\&+\frac{\mu }{2(1-\cond)}\|y_{k}-x_\star-\tfrac1L \nabla f(y_{k})\|^2_2,
    \end{aligned}
    \end{equation*}
\end{itemize}
After some algebra, the weighted sum can be reformulated exactly as follows (it is simpler not to use the expression of $\rho$ to verify this):
\begin{equation*}
\begin{aligned}
f&(y_{k})-f_\star-\frac1{2L}\|\nabla f(y_{k})\|^2_2-\frac{\mu }{2(1-\cond)}\|y_{k}-x_\star-\tfrac1L \nabla f(y_{k})\|^2_2\\
&+\frac{\mu}{1-\cond}\|z_{k+1}-x_\star\|^2_2\\
\leq& \rho^2 \bigg(f(y_{k-1})-f_\star-\frac1{2L}\|\nabla f(y_{k-1})\|^2_2\\
& \quad -\frac{\mu }{2(1-\cond)}\|y_{k-1}-x_\star-\tfrac1L \nabla f(y_{k-1})\|^2_2+\frac{\mu}{1-\cond}\|z_{k}-x_\star\|^2_2\bigg)\\
&-\frac{\cond-(\rho -1)^2}{(\rho -2) \left(1-\cond\right)}\\&\quad\quad\times \langle \nabla f(y_{k});\rho  (y_{k-1}-x_\star)-2 (\rho -1) (z_k-x_\star) -\tfrac{ \rho }{L}\nabla f(y_{k-1}) \rangle\\
&-\frac{\cond-(\rho -1)^2}{ \left(1-\cond\right)\mu}  \|\nabla f(y_{k})\|^2_2.
\end{aligned}
\end{equation*}
Using the expression $\rho=1-\sqrt{\cond}$, the last two terms cancel, and we arrive at the desired result:
\begin{equation*}
\begin{aligned}
f&(y_{k})-f_\star-\frac1{2L}\|\nabla f(y_{k})\|^2_2-\frac{\mu }{2(1-\cond)}\|y_{k}-x_\star-\tfrac1L \nabla f(y_{k})\|^2_2 \\
& \quad +\frac{\mu}{1-\cond}\|z_{k+1}-x_\star\|^2_2\\
\leq& \rho^2 \bigg(f(y_{k-1})-f_\star-\frac1{2L}\|\nabla f(y_{k-1})\|^2_2\\&\quad\quad-\frac{\mu }{2(1-\cond)}\|y_{k-1}-x_\star -\tfrac1L \nabla f(y_{k-1})\|^2_2+\frac{\mu}{1-\cond}\|z_{k}-x_\star\|^2_2\bigg).\qedhere
\end{aligned}
\end{equation*}
\end{proof}

\begin{remark} The triple momentum method was proposed in~\citep{van2017fastest}, heavily relying on the control-theoretic framework developed by~\citep{lessard2016analysis} (which is discussed in Appendix~\ref{a-WC_FO}). It was further studied in~\citet{cyrus2018robust} from a robust control perspective and by~\citet{zhou2020boosting} as an accelerated method for a different objective. The triple momentum method can also be obtained as a time-independent optimized gradient method~\citep{lessard2020direct}. Of course, all of the drawbacks of the OGM and ITEM also apply to the triple momentum method, so the same questions related to generalizations of this scheme remain open. Furthermore, this method is defined only for $\mu>0$, like Nesterov's method with constant momentum.
\end{remark}

\subsection{Quadratic Averaging and Geometric Descent}
Accelerated methods tailored for the strongly convex setting, such as Algorithms~\ref{alg:FGM_1_strcvx_constmomentum} and~\ref{alg:TMM}, make use of two kinds of step sizes to update the different sequences. First, they perform \emph{small gradient steps} with the step size $1/L$. Such steps correspond to minimizing quadratic upper bounds~\eqref{eq:smooth}. Second, they use so-called \emph{large steps} $1/\mu$, which correspond to minimizing quadratic lower bounds~\eqref{eq:strcvx}. This algorithmic structure is provided with an interpretation which was further exploited by \emph{geometric descent}~\citep{bubeck2015geometric} and \emph{quadratic averaging}~\citep{drusvyatskiy2018optimal}. These two methods produce the same sequences of iterates (see~\citep[Theorem 4.5]{drusvyatskiy2018optimal}), and we therefore take the stand of presenting them through the quadratic averaging viewpoint, which relates more directly to previous sections.

\citet{drusvyatskiy2018optimal} propose a method based on \emph{quadratic averaging}. It is similar in shape to Algorithm~\ref{alg:FGM_1_strcvx_constmomentum} except that the last sequence $z_k$ is explicitly computed as the minimum of a quadratic lower bound on the function. (The coefficients arising in the computation of $z_k$ are dynamically selected to maximize this lower bound). To construct the new lower bound at iteration $k$, the algorithm combines the lower bound constructed at iteration $k-1$, whose minimum is achieved at $z_k$, with a new lower bound constructed using the strong convexity assumption along with the first-order information of the current iterate (more precisely this second lower bound on $f(x)$ is $f(y_k)+\langle \nabla f(y_k),x-y_k\rangle+\tfrac{\mu}{2}\|x-y_k\|^2_2$), whose minimum is $y_k-\tfrac1\mu \nabla f(y_k)$. Because of the specific shape of these two lower bounds, their convex combinations has a minimum that is the convex combination of their respective minima, with the same weights, thereby motivating the update rule $z_{k}=\lambda (y_k-\tfrac1\mu \nabla f(y_k)) +(1-\lambda) z_{k-1}$, with dynamically chosen $\lambda$ (to maximize the minimum value of the new under-approximation).

Alternatively, the sequence $z_k$ can be interpreted in terms of a localization method for tracking $x_\star$, using intersections of balls. In this case, the sequence $z_k$ corresponds to the centers of the balls containing $x_\star$. This alternate viewpoint is referred to as \emph{geometric descent}~\citep{bubeck2015geometric}, and the $\lambda$ are chosen to minimize the radius of the ball, centered at $z_{k}=\lambda (y_k-\tfrac1\mu \nabla f(y_k)) +(1-\lambda) z_{k-1}$ while ensuring that the new ball contains~$x_\star$. 

Geometric descent is detailed in~\citep{bubeck2015geometric} and~\citep[Section 3.6.3]{bubeck2015convex}. It has been extended to handle constraints~\citep{chen2017geometric} and has been studied using the same Lyapunov function as that used in Theorem~\ref{thm:nest_constmom}~\citep{karimi2017single}.

\section{Practical Extensions}\label{sec:nest_practice}\label{s:backtracking} 
The goal of this section is to detail a few of the many extensions of Nesterov's accelerated methods. We see below that additional elements can be incorporated into the accelerated methods while maintaining the same proof structures. That is, we perform weighted sums involving the same inequalities that we used for the smooth (possibly strongly) convex functions and only introduce a few additional inequalities to account for the new elements.

We also seek to provide intuition, along with bibliographical pointers for going further. The following scenarios are particularly important in practice.

\paragraph{Constraints.} In the presence of constraints or nonsmooth functions, a common approach is to resort to (proximal) splitting techniques. This idea is not recent; see, e.g.,~\citep{douglas1956numerical,glowinski1975approximation,lions1979}. However, it remains highly relevant in signal processing, computer vision, and statistical learning~\citep{peyre2011numerical,parikh2014proximal,chambolle2016introduction,fessler2020optimization}.
\paragraph{Adaptation.} Problem constants, such as smoothness and strong convexity parameters, are generally unknown. Furthermore, their \emph{local} values tend to be much more favorable than their, typically conservative, global values. In general, smoothness constants are estimated on the fly using backtracking line-search strategies; see, e.g.~\citep{goldstein1962cauchy,armijo1966minimization,Nest83}. Strong convexity constants, or the more general H\"olderian error bounds (see Section~\ref{c-restart}), on the other hand, are more difficult to estimate, and \emph{restart} schemes are often used to adapt to these additional regularity properties; see, e.g.~\citep{Nest13,Beck12,ODon15,roulet2017sharpness}. Such schemes are the workhorse of Section~\ref{c-restart}.
\paragraph{Non-Euclidean settings.} Although we only briefly mention this topic, accounting for the geometry of the problem at hand is generally key to obtaining good empirical performance. In particular, optimization problems can often be formulated more naturally in a non-Euclidean space, with non-Euclidean norms producing better implicit models for the function. A popular method in this setting is commonly known as mirror descent~\citep{Nemi83}---see also~\citep{Bentc01,juditsky2011first}---which we do not to detail at length here. Good surveys are provided by~\citet{beck2003mirror,bubeck2015convex}. However, we do describe an accelerated method in this setting in Section~\ref{s:acc_MD}.

\paragraph{Numerical stability and monotone accelerated methods.} Gradient descent guarantees the iterates to be monotonically good approximations of an optimal solution (i.e., $f(x_{k+1})\leq f(x_k)$ for all $k$). This desirable property is generally not true for accelerated methods. Although the worst-case guarantees on $f(x_k)-f_\star$ are indeed monotonically decreasing functions of the iteration counter (such methods are sometimes referred to as \emph{quasi-monotone methods}~\citep{nesterov2015quasi}), accelerated methods are in general not descent schemes. Monotonicity is a desirable feature for improving numerical stability of algorithms, and we show in Section~\ref{s:monotone_variants} that simple modifications allows enforcing monotonicity of accelerated methods at low technical and computational cost. Albeit with a different presentation, such developments can be found in, e.g.,~\citep{tseng2008accelerated,beck2009fastmonotone}. The technique is particularly simple to incorporate within the potential function-based analyses of this section.

\subsection{Handling Nonsmooth Terms/Constraints} \label{ssec:composite_operator}
In this section, we consider the problem of minimizing a sum of two convex functions:
\begin{equation}\label{eq:opt_composite}
F_\star=\min_{x\in\mathbb{R}^d} \{ F(x)\defeq f(x)+h(x)\},
\end{equation}
where $f$ is $L$-smooth and (possibly $\mu$-strongly) convex and where $h$ is convex, closed, and proper (CCP), which we denote by $f\in\mathcal{F}_{\mu,L}$ and $h\in\mathcal{F}_{0,\infty}$ (These technical conditions ensure that the proximal operator, defined hereafter, is well defined everywhere on $\mathbb{R}^d$. We refer to the clear introduction by~\citet{ryu2016primer} and the references therein for further details.) In addition, we assume a proximal operator of $h$ to be readily available, so
\begin{equation}\label{eq:prox_chapt_Nest}
    \mathrm{prox}_{\gamma h}(x)\defeq \underset{y}{\mathrm{argmin}}\{\gamma h(y)+\tfrac12\|x-y\|^2_2\},
\end{equation}
can be evaluated efficiently. (Section~\ref{c-prox} deals with some cases where this operator is approximated; see also the discussions in Section~\ref{s:notes_ref_chapt_nest}.) The proximal operator can be seen as an \emph{implicit (sub)gradient} step on $h$, as dictated by the optimality conditions of the proximal operation 
\[ x_+=\mathrm{prox}_{\gamma h}(x) \Leftrightarrow x_+=x-\gamma g_h(x_+) \text{ with }g_h(x_+)\in\partial h(x_+).\]
In particular, when $h(x)$ is the indicator function of a closed convex set $Q$, the proximal operation corresponds to the orthogonal projection onto $Q$. There are a few commonly used functions for which the proximal operation has an analytical solution such as $h(x)=\|x\|_1$; see, for instance, the list provided by~\citep{chierchia2020proximity}. In the proofs below, we incorporate $h$ using inequalities that characterize convexity, that is,
\[ h(x)\geq h(y)+\langle g_h(y);x-y\rangle,\]
where $g_h(y)\in\partial h(y)$ is some subgradient of $h$ at $y$. The proximal step (sometimes referred to as backward, or implicit, step) is a base algorithmic tool in the first-order optimization toolbox.

In this setting, classical methods for solving~\eqref{eq:opt_composite} involve a \emph{forward-backward splitting} strategy (in other words, forward steps (a.k.a. gradient steps) on $f$ and backward steps (a.k.a. proximal steps) on $h$), introduced by~\citet{passty1979ergodic}. This topic is addressed in many references, and we refer to~\citep{parikh2014proximal,ryu2016primer} and the references therein for further details. In the context of accelerated methods, forward-backward splitting was introduced by Nesterov (\citeyear{Nest03a,Nest13}) through the concept of \emph{gradient mapping}; see also~\citet{tseng2008accelerated,Beck09}. Problem~\eqref{eq:opt_composite} is also sometimes referred to as the \emph{composite convex optimization setting}~\citep{Nest13}.
Depending on the assumptions made on $f$ and $h$, there are alternate ways of solving this problem---for example, when the proximal operator is available for both, one can use the Douglas-Rachford splitting~\citep{douglas1956numerical,lions1979}. However, this is beyond the scope of this section and we refer to~\citep{ryu2016primer,condat2019proximal} and the references therein for further discussions on this topic.

\subsection{Adaptation to Unknown Regularity Parameters}
In previous sections, we assumed $f$ to be $L$-smooth and possibly $\mu$-strongly convex. Moreover, in the previous algorithms, we explicitly used the values of both $L$ and $\mu$ to design the methods. However, this is not a desirable feature. First, it means that we need to be able to estimate valid values for $L$ and $\mu$. Second, it means that the methods are not adaptive to potentially better (local) parameter values. That is, we do not benefit from the problems being simpler than specified, i.e., where the smallest valid $L$ is much smaller than our estimate and/or the largest valid $\mu$ is much larger than our estimation. Furthermore, we want to benefit from the typically better local properties of the function at hand, along the path taken by the method, rather than relying on the global properties. The difference between local and global regularity properties is often significant, and adaptive methods often converge much faster in practice. 

We discuss below how adaptation is implemented for the smoothness constant, using line-search techniques. However, it remains an open question whether strong convexity parameters can be efficiently estimated while maintaining reasonable worst-case guarantees and without resorting to restart schemes (i.e., outer iterations) (see Section~\ref{c-restart}).

To handle unknown parameters, the key is to examine the inequalities used in the proofs of the desired method. It turns out that smoothness is usually only used in inequalities between pairs of iterates, which means that these inequalities can be tested at runtime, at each iteration. Therefore, for our guarantees to hold, we do not need the function to be $L$-smooth everywhere, but rather we only need a given inequality to hold for the value of $L$ that we are using (where the smoothness of the function ensures that such an $L$ exists). Conversely, strong convexity is typically only used in inequalities involving the optimal point (see, for example, the proof of Theorem~\ref{thm:Nest_first_strcvx}), which we do not know a priori. As a result, these inequalities cannot be tested as the algorithm proceeds, which complicates the estimation of strong convexity while running the algorithm. Adaptation to strong convexity is therefore typically accomplished via the use of \emph{restarts}.

These approaches are common, and are not new~\citep{goldstein1962cauchy,armijo1966minimization}; they were already used by~\citet{Nest83}. They were later adapted to the forward-backward setting in~\citet{Nest13,Beck09} and have been further exploited to improve performance in various settings; see, e.g.,~\citep{scheinberg2014fast}. The topic is further discussed in the next section as well as in the notes and references provided in Section~\ref{s:notes_ref_chapt_nest}.

\subsubsection{An Accelerated Forward-backward Methods with Backtracking}
As discussed above, the smoothness constant $L$ is used very sparsely in the proofs of both the gradient descent (Theorem~\ref{thm:pot_GM} and Theorem~\ref{thm:pot_GM_strcvx}) and the accelerated variants (Theorem~\ref{thm:Nest_first}, Theorem~\ref{thm:Nest_first_strcvx}, and Theorem~\ref{thm:nest_constmom}). Essentially, it is only used in three places: (i) to compute $A_{k+1}$ (only when $\mu>0$);  (ii) to compute $x_{k+1}=y_k-\tfrac1L \nabla f(y_k)$; and (iii) in the inequality
\begin{equation}\label{eq:backtracking_cond}
f(x_{k+1})\leq f(y_{k}) + \langle \nabla f(y_{k}); x_{k+1}-y_{k}\rangle +\frac{L}{2}\|x_{k+1}-y_{k}\|^2_2.
\end{equation}
(Recall that this is known as the \emph{descent lemma} since by substituting the gradient step $x_{k+1}=y_k-\tfrac1L \nabla f(y_k)$ it can be written as $f(x_{k+1})\leq f(y_k)-\tfrac1{2L}\|\nabla f(y_k)\|^2_2$.) Other than $L$,~\eqref{eq:backtracking_cond} only contains information that \emph{we observe}. Hence, we can simply \emph{check} whether this inequality holds for a given estimate of $L$. When it does not hold, we simply increase the current approximation of $L$ and then with this new estimate, recompute (i) $A_{k+1}$ (necessary only if $\mu>0$) and the corresponding $y_k$ and (ii) $x_{k+1}$, using the new step size. We then check again whether~\eqref{eq:backtracking_cond} is satisfied. If~\eqref{eq:backtracking_cond} is satisfied, then we can proceed (because the potential inequality of the desired method is then verified---see, e.g., Theorem~\ref{thm:pot_GM} or Theorem~\ref{thm:pot_GM_strcvx} for gradient descent, or Theorem~\ref{thm:Nest_first}, Theorem~\ref{thm:Nest_first_strcvx}, or Theorem~\ref{thm:nest_constmom} for Nesterov's methods), and otherwise we continue increasing our approximation of $L$ until the descent condition~\eqref{eq:backtracking_cond} is satisfied. Finally, to guarantee that we only perform a finite number of ``wasted'' gradient steps to estimate $L$, we need an appropriate rule for how to increase our approximation. It is common to simply multiply the current approximation by some constant $\alpha>1$, thereby guaranteeing that at most $\lceil\log_\alpha \tfrac{L}{L_0}\rceil$ gradient steps, where $L$ is the true smoothness constant and $L_0$ is our starting estimate, are wasted in the process. As we see below, both backtracking and nonsmooth terms require proofs very similar to those presented above, and potential-based analyses are suitable.

We present two extensions of Nesterov's first method that can handle nonsmooth terms, and that have a backtracking procedure on the smoothness parameter. The first, the fast iterative shrinkage-thresholding algorithm (FISTA), is particularly popular, while the second resolves one potential issue that can arise in the original FISTA.

\paragraph{FISTA.}
The fast iterative shrinkage-thresholding algorithm, due to~\citet{Beck09}, is a natural extension of~\citet{Nest83} in its first form (see Algorithm~\ref{alg:FGM_1_strconvex}), handling an additional nonsmooth term.  In this section, we present a strongly convex variant of FISTA, provided as Algorithm~\ref{alg:fista_1}. The proof contains the same ingredients as in the original work, and it can easily be compared to previous material.

\begin{algorithm}[!ht]
  \caption{Strongly convex FISTA, form I}
  \label{alg:fista_1}
  \begin{algorithmic}[1]
    \REQUIRE $L$-smooth $\mu$-strongly (possibly with $\mu=0$) convex function~$f$, a convex function~$h$ with proximal operator available, an initial point $x_0$, and an initial estimate $L_0>\mu$.
    \STATE \textbf{Initialize} $z_0=x_0$, $A_0=0$, and some $\alpha>1$.
    \FOR{$k=0,\ldots$}
        \STATE $L_{k+1}=L_k$\label{alg:init_bt}
        \LOOP 
            \STATE $\cond_{k+1}=\mu/L_{k+1}$
            \STATE $A_{k+1}=\frac{2 A_k+1+\sqrt{4 A_k+4\cond_{k+1} A_k^2 +1}}{2 \left(1-\cond_{k+1}\right)}$
            \STATE set $\tau_k=\frac{(A_{k+1}-A_k) (1+\cond_{k+1} A_k )}{A_{k+1}+2\cond_{k+1} A_k A_{k+1}-\cond_{k+1} A_k^2}$ and $\delta_k=\frac{A_{k+1}-A_{k}}{1+\cond_{k+1} A_{k+1}}$
            \STATE $y_{k}= x_k+\tau_k(z_k-x_k)$
            \STATE $x_{k+1}=\mathrm{prox}_{h/L_{k+1}}\left(y_{k}-\tfrac{1}{L_{k+1}}\nabla f(y_{k})\right)$
            \STATE $z_{k+1}=(1-\cond_{k+1}\delta_k)z_k+\cond_{k+1}\delta_ky_k+\delta_k \left(x_{k+1}-y_k\right)$
            \IF{\eqref{eq:backtracking_cond} holds}
                \STATE \textbf{break} \COMMENT{Iterates accepted; $k$ will be incremented.}
            \ELSE
                \STATE $L_{k+1}=\alpha L_{k+1}$ \COMMENT{Iterates  not accepted; recompute new $L_{k+1}$.}
            \ENDIF 
        \ENDLOOP
    \ENDFOR
    \ENSURE An approximate solution $x_{k+1}$.
  \end{algorithmic}
\end{algorithm}

\enlargethispage{\baselineskip}
The proof follows exactly the same steps as the proof of Theorem~\ref{thm:Nest_first_strcvx} (Nesterov's method for strongly convex functions), but it also accounts for the nonsmooth function $h$. (Observe that the potential is stated in terms of $F$ and not $f$.) Two additional inequalities, involving the convexity of $h$ between two different pairs of points, allow this nonsmooth term to be taken into account. In this case, $f$ is assumed to be smooth and convex over $\mathbb{R}^d$ (i.e., it has full domain, $\dom f=\mathbb{R}^d$), and we are therefore allowed to evaluate gradients of $f$ outside of domain of $h$.

\begin{theorem}\label{thm:FISTA_strcvx} Let $f\in\mathcal{F}_{\mu,L}$ (with full domain: $\dom f=\mathbb{R}^d$); $h$ be a closed convex proper function; $x_\star \in\mathrm{argmin}_x\, \{f(x)+h(x)\}$; and $k\in\mathbb{N}$. For any $x_k,z_k\in\mathbb{R}^d$ and $A_k\geq 0$, the iterates of Algorithm~\ref{alg:fista_1} that satisfy~\eqref{eq:backtracking_cond} also satisfy
\begin{equation*}
\begin{aligned}
A_{k+1}&(F(x_{k+1})-F_\star)+\frac{L_{k+1}+\mu A_{k+1}}{2}\|z_{k+1}-x_\star\|^2_2\\&\leq A_k (F(x_k)-F_\star)+\frac{L_{k+1}+\mu A_{k}}2 \|z_k-x_\star\|^2_2,
\end{aligned}
\end{equation*}
with $A_{k+1}=\frac{2 A_k+1+\sqrt{4 A_k+4 A_k^2 \cond_{k+1}+1}}{2 \left(1-\cond_{k+1}\right)}$ and $q_{k+1}=\mu/L_{k+1}$.
\end{theorem}
\bgroup
\addtolength{\jot}{-0.1em}
\begin{proof}
The proof consists of a weighted sum of the following inequalities.
\begin{itemize}
    \item Strong convexity of $f$ between $x_\star$ and $y_{k}$ with weight $\lambda_1=A_{k+1}-A_k$:
    \[ f(x_\star) \geq f(y_{k})+\langle \nabla f(y_{k});x_\star-y_k\rangle+\frac{\mu}{2}\|x_\star-y_k\|^2_2.\]
    \item Strong convexity of $f$ between $x_k$ and $y_{k}$ with weight $\lambda_2=A_k$:
    \[ f(x_k) \geq f(y_{k})+\langle \nabla f(y_{k});x_k-y_{k}\rangle.\]
    \item Smoothness of $f$ between $y_{k}$ and $x_{k+1}$ (\emph{descent lemma}) with weight $\lambda_3=A_{k+1}$:
    \[ f(y_{k}) + \langle \nabla f(y_{k}); x_{k+1}-y_{k}\rangle +\frac{L_{k+1}}{2}\|x_{k+1}-y_{k}\|^2_2 \geq f(x_{k+1}).\]
    \item Convexity of $h$ between $x_\star$ and $x_{k+1}$ with weight $\lambda_4=A_{k+1}-A_k$:
    \[ h(x_\star)\geq h(x_{k+1})+\langle g_h(x_{k+1});x_\star-x_{k+1}\rangle,\]
    with $g_h(x_{k+1})\in\partial h(x_{k+1})$ and $x_{k+1}=y_k-\tfrac1{L_{k+1}} (\nabla f(y_k)+g_h(x_{k+1}))$.
    \item Convexity of $h$ between $x_k$ and $x_{k+1}$ with weight $\lambda_5=A_k$:
    \[ h(x_k)\geq h(x_{k+1})+\langle g_h(x_{k+1});x_k-x_{k+1}\rangle.\]
\end{itemize}
We get the following inequality
\begin{equation*}
\begin{aligned}
0\geq &\lambda_1 [f(y_{k})-f(x_\star)+\langle \nabla f(y_{k});x_\star-y_k\rangle+\frac{\mu}{2}\|x_\star-y_k\|^2_2]\\&+\lambda_2[f(y_{k})-f(x_k)+\langle \nabla f(y_{k});x_k-y_{k}\rangle]\\
&+\lambda_3[f(x_{k+1})-(f(y_{k}) + \langle \nabla f(y_{k}); x_{k+1}-y_{k}\rangle \\
& \quad +\frac{L_{k+1}}{2}\|x_{k+1}-y_{k}\|^2_2)]\\
&+\lambda_4[ h(x_{k+1})-h(x_\star)+\langle g_h(x_{k+1});x_\star-x_{k+1}\rangle]\\&+\lambda_5[ h(x_{k+1})-h(x_k)+\langle g_h(x_{k+1});x_k-x_{k+1}\rangle].
\end{aligned}
\end{equation*}
By substituting $y_k$, $x_{k+1}$, and $z_{k+1}$ with
\begin{equation*}
\begin{aligned}
y_k&=x_k+\tau_k  (z_k-x_k)\\
x_{k+1}&=y_k-\tfrac{1}{L_{k+1}}(\nabla f(y_k)+g_h(x_{k+1})) \\
z_{k+1}&=(1-\cond_{k+1}\delta_k)z_k+\cond_{k+1}\delta_k y_k+\delta_k \left(x_{k+1}-y_k\right),
\end{aligned}
\end{equation*}
the previous weighted sum can be reformulated exactly as
\begin{equation*}
\begin{aligned}
A_{k+1}&(f(x_{k+1})+h(x_{k+1})-f(x_\star)-h(x_\star))\\
& \quad +\frac{L_{k+1}+\mu A_{k+1}}{2}\|z_{k+1}-x_\star\|^2_2\\
\leq& A_k (f(x_k)+h(x_k)-f(x_\star)-h(x_\star))+\frac{L_{k+1}+A_k\mu}{2}\|z_k-x_\star\|^2_2\\
&+\frac{(A_{k}-A_{k+1})^2-A_{k+1}-\cond_{k+1}  A_{k+1}^2}{1+\cond_{k+1}  A_{k+1}} \frac1{2L_{k+1}}\|\nabla f(y_k)+g_h(x_{k+1})\|^2_2\\
&-\frac{A_{k}^2 (A_{k+1}-A_{k}) (1+\cond_{k+1}  A_{k}) (1+\cond_{k+1}  A_{k+1})}{\left(A_{k+1} +2 \cond_{k+1}  A_{k}A_{k+1}-\cond_{k+1}  A_{k}^2\right)^2}\frac{\mu}{2}\|x_k-z_k\|^2_2.
\end{aligned}
\end{equation*}
Using $0\leq \cond_{k+1}\leq 1$ and selecting $A_{k+1}$ such that $A_{k+1}\geq A_k$ and
\[ (A_{k}-A_{k+1})^2-A_{k+1}-\cond_{k+1}  A_{k+1}^2 = 0,\]
yields the desired result:
\begin{equation*}
\begin{aligned}
A_{k+1}&(f(x_{k+1})+h(x_{k+1})-f(x_\star)-h(x_\star))\\
&\quad+\frac{L_{k+1}+\mu A_{k+1}}{2}\|z_{k+1}-x_\star\|^2_2\\
\leq& A_k (f(x_k)+h(x_k)-f(x_\star)-h(x_\star))+\frac{L_{k+1}+\mu A_k}{2}\|z_k-x_\star\|^2_2.\qedhere
\end{aligned}
\end{equation*}
\end{proof}
\egroup

Finally, we obtain a complexity guarantee by adapting the potential argument~\eqref{eq:chained_pot} and by noting that $A_{k+1}$ is a decreasing function of $L_{k+1}$ (whose maximal value is $\alpha L$ assuming $L_0<L$ and is otherwise $L_0$). The growth rate of $A_k$ in the smooth convex setting remains unchanged; see~\eqref{eq:conv_Ak_fgm}. However, when $L_0<L$, its geometric grow rate might actually be slightly degraded to
\[ A_{k+1} \geq \left(1-\sqrt{\frac{\mu}{\alpha L}}\right)^{-1} A_k,\]
which remains better than the worst-case $(1-\tfrac{\mu}{\alpha L})$ rate of gradient descent with backtracking, assuming in both cases that $L_0<L$. When $L_0>L$ the rates might respectively be degraded to $(1-\sqrt{\tfrac{\mu}{L_0}})$ and $(1-\tfrac{\mu}{L_0})$ instead.

\begin{corollary}\label{cor:FISTA} Let $f\in\mathcal{F}_{\mu,L}$ (with full domain: $\dom f=\mathbb{R}^d$); $h$ be a closed convex proper function; and $x_\star\in\mathrm{argmin}_x\, \{F(x)\defeq f(x)+h(x)\}$. For any $N\in\mathbb{N}$, $N\geq 1$,  and $x_0\in\mathbb{R}^d$, the output of Algorithm~\ref{alg:fista_1} satisfies
\[F(x_N)-F_\star\leq \min\left\{\frac{2}{N^2},\left(1-\sqrt{\frac{\mu}{\ell}}\right)^{N}\right\}\ell\|x_0-x_\star\|^2_2,\]
with $\ell=\max\{\alpha L,L_0\}$.
\end{corollary}
\begin{proof} We assume that $L>L_0$ since otherwise $f\in\mathcal{F}_{\mu,L_0}$ and the proof would directly follow the case without backtracking.
Define
\[ \phi_k\defeq A_k (F(x_k)-F_\star)+\frac{L_k+\mu A_{k}}2 \|z_k-x_\star\|^2_2.\]
Since $L_{k+1}/L_k\geq 1$, we have
\begin{equation*}
\begin{aligned}
\phi_{k+1}&\leq A_k (F(x_k)-F_\star)+\frac{L_{k+1}+\mu A_{k}}2 \|z_k-x_\star\|^2_2\leq \frac{L_{k+1}}{L_k}\phi_k.
\end{aligned}
\end{equation*}
The chained potential argument~\eqref{eq:chained_pot} can then be adapted to obtain
\[A_N (F(x_N)-F_\star)\leq \phi_N\leq \frac{L_{N}}{L_{N-1}} \phi_{N-1}\leq\frac{L_{N}}{L_{N-2}} \phi_{N-2}\leq \hdots \leq \frac{L_N}{L_0}\phi_0, \]
where we used Theorem~\ref{thm:FISTA_strcvx} and the fact that the output of the algorithm satisfies~\eqref{eq:backtracking_cond}.
Using $A_0=0$, we reach
\[F(x_N)-F_\star\leq \frac{L_N\|x_0-x_\star\|^2_2}{2 A_N}.\]
Using our previous bounds on $A_N$ (noting that $A_{k+1}$ is a decreasing function of $L_{k+1}$) in, e.g., Corollary~\ref{cor:pot_fgm_strcvx}, along with the fact that in the worst-case the estimated smoothness cannot be larger than the growth rate $\alpha$ times the true constant $L_N< \alpha L$ except if $L_0$ were already larger than the true $L$, in which case $L_N=L_0$. Therefore, we get $L_N\leq\ell=\max\{\alpha L,L_0\}$, yielding the desired result.
\end{proof}

\begin{remark}\label{rem:bt_strat} There are two common variations on the backtracking strategy presented in this section. One can, for example, reset $L_{k+1}\leftarrow L_0$ (in line~\ref{alg:init_bt} of Algorithm~\ref{alg:fista_1}) at each iteration, potentially using a total of $N\lceil\log_\alpha \tfrac{L}{L_0}\rceil$ additional gradient evaluations over all iterations. Another possibility is to pick some additional constant $0<\beta<1$ and to initiate  $L_{k+1}\leftarrow \beta L_k$ (in line~\ref{alg:init_bt} of Algorithm~\ref{alg:fista_1}). In the case $\beta=1/\alpha$, this strategy potentially costs $1$ additional gradient evaluation per iteration due to the backtracking strategy, thus potentially using a total of $N+\lceil\log_\alpha \tfrac{L}{L_0}\rceil$ additional gradient evaluations over all iterations. 

Such \emph{non-monotonic} estimations of $L$ can be incorporated at a low additional technical cost. The corresponding methods and their analyses are essentially the same as those of this section; they are provided in Appendix~\ref{s:Fista_nonmonotone_BT} and ~\ref{s:Acc_nonmonotone_BT} (see Algorithm~\ref{alg:fista_1_general_BT} and Algorithm~\ref{alg:MST_strconvex_general_BT}).
\end{remark}

\begin{remark} Variations on strongly convex extensions of FISTA, involving backtracking line-searches, can be found in, e.g.,~\citep{chambolle2016introduction,calatroni2019backtracking,florea2018accelerated,florea2020generalized}, together with practical improvements. The method presented in this section was designed for easy comparison with the previous material.
\end{remark}

\subsubsection{Another Accelerated Proximal Gradient Method}
FISTA potentially evaluates gradients outside of the domain of $h$, and it therefore implicitly assumes that $f$ is defined even outside this region. In many situations, this is not an issue, such as when $f$ is quadratic. In this section, we instead assume that $f$ is continuously differentiable and satisfies smoothness condition~\eqref{eq:smooth_res} only for all $x,y\in\dom h$. 

\begin{definition}\label{def:smoothstrconvex_constrained} Let $0\leq\mu< L\leq+\infty$ and $C\subseteq\mathbb{R}^d$. A closed convex proper function $f:\mathbb{R}^d\rightarrow\mathbb{R}\cup\{+\infty\}$ is $L$-smooth and $\mu$-strongly convex on $C$ (written $f\in\mathcal{F}_{\mu,L}(C)$) if and only if 
\begin{itemize}
    \item ($L$-smoothness) there exists an open set $C'$ such that $C\subseteq C'$ and $f$ is continuously differentiable on $C'$, and for all $x,y\in C$, it holds that
    \begin{equation}\label{eq:smooth_res}
    f(x)\leq f(y)+\langle \nabla f(y);x-y\rangle+\frac{L}{2}\|x-y\|^2_2;
    \end{equation}
    \item ($\mu$-strong convexity)  for all $x,y\in C$, it holds that
    \begin{equation}\label{eq:strcvx_res}
    f(x)\geq f(y)+\langle g_f(y);x-y\rangle+\frac{\mu}{2}\|x-y\|^2_2,
    \end{equation}
    where $g_f(y)\in\partial f(y)$ is a subgradient of $f$ at $y$. (Note that $g_f(y)=\nabla f(y)$ when $f$ is differentiable.)
\end{itemize}
By extension, $\mathcal{F}_{\mu,\infty}(C)$ denotes the set of closed $\mu$-strongly convex proper functions whose domain contains $C$, and $\mathcal{F}_{0,\infty}$ denotes the set of closed convex proper functions.
\end{definition}

There exist different ways of handling this situation. The method presented in this section relies on using the proximal operator on the sequence $z_k$ and on formulating Nesterov's method in form III (see Algorithm~\ref{alg:FGM_3}). In this situation, assuming the initial point is feasible ($F(x_0)<\infty$) implies both the $x_k$ and $y_k$ are obtained from convex combinations of feasible points and hence are feasible.

A wide variety of accelerated methods exists; most variants solve the issue of FISTA using two proximal operations per iteration (on both of the sequences $x_k$ and $z_k$). The variant in this section performs only one projection per iteration, while also fixing the infeasibility issue of $y_k$ in FISTA. Variations on this theme can found in a number of references; see, for example,~\citep[``Improved interior gradient algorithm'']{auslender2006interior},~\citep[Algorithm 1]{tseng2008accelerated}, or more recently~\citep[``Method of similar triangles'']{gasnikov2018universal}.

\begin{algorithm}[!h]
  \caption{A proximal accelerated gradient method}
  \label{alg:MST_strconvex}
  \begin{algorithmic}[1]
    \REQUIRE $h\in\mathcal{F}_{0,\infty}$ with proximal operator available, $f\in\mathcal{F}_{\mu,L}(\dom h)$, an initial point $x_0\in\dom h$, and an initial estimate $L_0>\mu$.
    \STATE \textbf{Initialize} $z_0=x_0$, $A_0=0$, and some $\alpha>1$.
    \FOR{$k=0,\ldots$}
        \STATE $L_{k+1}=L_k$
        \LOOP
            \STATE $\cond_{k+1}=\mu/L_{k+1}$
            \STATE $A_{k+1}=\frac{2 A_k+1+\sqrt{4 A_k+4\cond_{k+1} A_k^2 +1}}{2 \left(1-\cond_{k+1}\right)}$
            \STATE set $\tau_k=\frac{(A_{k+1}-A_k) (1+\cond_{k+1} A_k)}{A_{k+1}+2\cond_{k+1} A_k A_{k+1}-\cond_{k+1} A_k^2 }$ and $\delta_{k}=\frac{A_{k+1}-A_{k}}{1+\cond_{k+1} A_{k+1}}$
            \STATE $y_{k}=  x_k+\tau_k (z_k-x_k)$
            \STATE $z_{k+1}=\mathrm{prox}_{\delta_k h/L_{k+1}}\left((1-\cond_{k+1}\delta_k)z_k+ \cond_{k+1}\delta_k y_k- \frac{\delta_k}{L_{k+1}}\nabla f(y_{k})\right)$
            \STATE $x_{k+1}=\frac{A_k}{A_{k+1}}x_k+(1-\frac{A_k}{A_{k+1}})z_{k+1}$
            \IF{\eqref{eq:backtracking_cond} holds}
                \STATE \textbf{break} \COMMENT{Iterates accepted; $k$ will be incremented.}
            \ELSE
                \STATE $L_{k+1}=\alpha L_{k+1}$ \COMMENT{Iterates  not accepted; recompute new $L_{k+1}$}
            \ENDIF 
        \ENDLOOP
    \ENDFOR
    \ENSURE Approximate solution $x_{k+1}$.
  \end{algorithmic}
\end{algorithm}

\begin{theorem}\label{thm:MST_strcvx} Let $h\in\mathcal{F}_{0,\infty}$, $f\in\mathcal{F}_{\mu,L}(\dom h)$; $x_\star\in\mathrm{argmin}_x\, \{F(x)\defeq f(x)+h(x)\}$; and $k\in\mathbb{N}$. For any $x_{k},z_k\in\mathbb{R}^d$ and $A_k\geq0$, the iterates of Algorithm~\ref{alg:MST_strconvex}  that satisfy~\eqref{eq:backtracking_cond} also satisfy
\begin{equation*}
\begin{aligned}
A_{k+1}&(F(x_{k+1})-F_\star)+\frac{L_{k+1}+\mu A_{k+1}}{2}\|z_{k+1}-x_\star\|^2_2\\&\leq A_k (F(x_k)-F_\star)+\frac{L_{k+1}+\mu A_{k}}2 \|z_k-x_\star\|^2_2,
\end{aligned}
\end{equation*}
with $A_{k+1}=\frac{2 A_k+1+\sqrt{4 A_k+4 A_k^2 \cond_{k+1}+1}}{2 \left(1-\cond_{k+1}\right)}$ and $q_{k+1}=\mu/L_{k+1}$.
\end{theorem}
\begin{proof} First, $z_k$ is in $\dom h$ by construction---it is the output of the proximal/projection step. Furthermore, we have $0\leq \tfrac{A_{k}}{A_{k+1}}\leq 1$ given that $A_{k+1}\geq A_k\geq 0$. A direct consequence is that since $z_0=x_0\in\dom h$, all subsequent $\{y_k\}$ and $\{x_k\}$ are also in $\dom h$ (as they are obtained from convex combinations of feasible points).

The rest of the proof consists of a weighted sum of the following inequalities (which are valid due to the feasibility of the iterates).
\begin{itemize}
    \item Strong convexity of $f$ between $x_\star$ and $y_{k}$ with weight $\lambda_1=A_{k+1}-A_k$:
    \[ f(x_\star) \geq f(y_{k})+\langle \nabla f(y_{k});x_\star-y_k\rangle+\frac{\mu}{2}\|x_\star-y_k\|^2_2.\]
    \item Convexity of $f$ between $x_k$ and $y_{k}$ with weight $\lambda_2=A_k$:
    \[ f(x_k) \geq f(y_{k})+\langle \nabla f(y_{k});x_k-y_{k}\rangle.\]
    \item Smoothness of $f$ between $y_{k}$ and $x_{k+1}$ (\emph{descent lemma}) with weight $\lambda_3=A_{k+1}$:
    \[ f(y_{k}) + \langle \nabla f(y_{k}); x_{k+1}-y_{k}\rangle +\frac{L_{k+1}}{2}\|x_{k+1}-y_{k}\|^2_2 \geq f(x_{k+1}).\]
    \item Convexity of $h$ between $x_\star$ and $z_{k+1}$ with weight $\lambda_4=A_{k+1}-A_k$:
    \[ h(x_\star)\geq h(z_{k+1})+\langle g_h(z_{k+1});x_\star-z_{k+1}\rangle,\]
    with $g_h(z_{k+1})\in\partial h(z_{k+1})$ and $z_{k+1}=(1-\cond\delta_k)z_k+ \cond\delta_k y_k- \frac{\delta_k}{L_{k+1}}(\nabla f(y_{k})+g_h(z_{k+1}))$.
    \item Convexity of $h$ between $x_k$ and $x_{k+1}$ with weight $\lambda_5=A_k$:
    \[ h(x_k)\geq h(x_{k+1})+\langle g_h(x_{k+1});x_k-x_{k+1}\rangle,\]
    with $g_h(x_{k+1})\in\partial h(x_{k+1})$.
    \item Convexity of $h$ between $z_{k+1}$ and $x_{k+1}$ with weight $\lambda_6=A_{k+1}-A_k$:
    \[ h(z_{k+1})\geq h(x_{k+1})+\langle g_h(x_{k+1});z_{k+1}-x_{k+1}\rangle.\]
\end{itemize}
We this obtain the following inequality:
\begin{equation*}
\begin{aligned}
0\geq &\lambda_1 [f(y_{k})-f(x_\star)+\langle \nabla f(y_{k});x_\star-y_k\rangle+\frac{\mu}{2}\|x_\star-y_k\|^2_2]\\
&+\lambda_2[f(y_{k})-f(x_k)+\langle \nabla f(y_{k});x_k-y_{k}\rangle]\\
&+\lambda_3[f(x_{k+1})-(f(y_{k}) + \langle \nabla f(y_{k}); x_{k+1}-y_{k}\rangle\\
&+\frac{L_{k+1}}{2}\|x_{k+1}-y_{k}\|^2_2)]\\
&+\lambda_4[ h(z_{k+1})-h(x_\star)+\langle g_h(z_{k+1});x_\star-z_{k+1}\rangle]\\
&+\lambda_5[ h(x_{k+1})-h(x_k)+\langle g_h(x_{k+1});x_k-x_{k+1}\rangle]\\
&+\lambda_6[ h(x_{k+1})-h(z_{k+1})+\langle g_h(x_{k+1});z_{k+1}-x_{k+1}\rangle].
\end{aligned}
\end{equation*}
By substituting $y_k$, $z_{k+1}$, and $x_{k+1}$ with
\begin{equation*}
\begin{aligned}
y_k&=x_k+\tau_k  (z_k-x_k)\\
z_{k+1}&= (1- \cond_{k+1} \delta_k)z_k +\cond_{k+1}\delta_k y_k  - \frac{ \delta_k}{L_{k+1}} (\nabla f(y_k) + g_h(z_{k+1}))\\
x_{k+1}&= \frac{A_k}{A_{k+1}} x_k + \left(1 - \frac{A_k}{A_{k+1}}\right) z_{k+1},
\end{aligned}
\end{equation*}
we reformulate the previous inequality as
\begin{equation*}
\begin{aligned}
A_{k+1}&(f(x_{k+1})+h(x_{k+1})-f(x_\star)-h(x_\star))\\
& \quad +\frac{L_{k+1}+A_{k+1}\mu}{2}\|z_{k+1}-x_\star\|^2_2\\
\leq& A_k (f(x_k)+h(x_k)-f(x_\star)-h(x_\star))+\frac{L_{k+1}+A_k\mu}{2}\|z_k-x_\star\|^2_2\\
& \quad +\frac{(A_k-A_{k+1})^2 \left((A_{k}-A_{k+1})^2-A_{k+1}-\cond_{k+1}  A_{k+1}^2\right)}{ A_{k+1} (1+\cond_{k+1}  A_{k+1})^2} \\
& \quad\quad\quad  \times \frac1{2L_{k+1}}\|\nabla f(y_k)+g_h(z_{k+1})\|^2_2\\
& \quad -\frac{A_k^2 (A_{k+1}-A_k) (1+ \cond_{k+1} A_k) (1+\cond_{k+1} A_{k+1})}{\left(A_{k+1}+2 \cond_{k+1}  A_k A_{k+1}-\cond_{k+1}  A_k^2\right)^2} \frac{\mu}{2}\|x_k-z_k\|^2_2.
\end{aligned}
\end{equation*}
Then selecting $A_{k+1}\geq A_k$ such that
\[ (A_{k}-A_{k+1})^2-A_{k+1}-\cond_{k+1}  A_{k+1}^2 = 0,\]
yields the desired result:
\begin{equation*}
\begin{aligned}
A_{k+1}&(f(x_{k+1})+h(x_{k+1})-f(x_\star)-h(x_\star))\\
& \quad +\frac{L_{k+1}+\mu A_{k+1}}{2}\|z_{k+1}-x_\star\|^2_2\\
\leq& A_k (f(x_k)+h(x_k)-f(x_\star)-h(x_\star))+\frac{L_{k+1}+\mu A_k}{2}\|z_k-x_\star\|^2_2.\qedhere
\end{aligned}
\end{equation*}
\end{proof}

We have the following corollary.

\begin{corollary} Let $h\in\mathcal{F}_{0,\infty}$, $f\in\mathcal{F}_{\mu,L}(\dom h)$, and $x_\star\in\mathrm{argmin}_x$ $\{F(x)\defeq f(x)+h(x)\}$. For any $N\in\mathbb{N}$, $N\geq 1$, and $x_0\in\mathbb{R}^d$, the output of Algorithm~\ref{alg:MST_strconvex} satisfies
\[F(x_N)-F_\star\leq \min\left\{\frac{2}{N^2},\left(1-\sqrt{\frac{\mu}{\ell}}\right)^{-N}\right\}\ell\|x_0-x_\star\|^2_2,\]
with $\ell=\max\{\alpha L,L_0\}$.
\end{corollary}
\begin{proof} The proof follows the same arguments as those for Corollary~\ref{cor:FISTA}, using the potential from Theorem~\ref{thm:MST_strcvx} with the fact the output of the algorithm satisfies~\eqref{eq:backtracking_cond}.
\end{proof}
\begin{remark}
In this section, we introduced backtracking techniques by examining how the inequalities are used in previous proofs. In particular, because smoothness is used only through the \emph{descent lemma} in which the only \emph{unknown} value is $L$, one can simply check this inequality at runtime. Another way to exploit the observation of which inequalities are needed in a proof is to identify \emph{minimal} assumptions on the class of functions under which it is possible to prove accelerated rates; this topic is explored by, e.g.,~\citet{necoara2019linear,hinder20a}. More generally the same question holds for the ability to prove the convergence rates of simpler methods, such as gradient descent~\citep{bolte2017error,necoara2019linear}.
\end{remark}

\subsection{Monotone Accelerated Methods}\label{s:monotone_variants}
As emphasized in previous sections, accelerated methods are \emph{quasi-monotone}, meaning that their worst-case guarantees are decreasing functions of the number of iteration. However, they are generally not monotone, as the function values are not guaranteed to be improved from iteration to iteration.

In this section, we introduce a simple trick for making common accelerated methods monotone. The technique stems from a simple observation that the iterates $\{x_k\}_k$ can be slightly changed while maintaining the worst-case guarantees (see, e.g.,~\citep{tseng2008accelerated}). That is, we introduce an additional sequence, denoted by $\{\tilde{x}_k\}_k$ for which $\{F(\tilde{x}_k)\}_k$ is monotonically decreasing. For instance, for all problems on which $\{F(x_k)\}_k$ is already monotonically decreasing, we have $x_k=\tilde{x}_k$ for all $k$.

As it is, the trick does not apply to the optimized gradient method (Algorithm~\ref{alg:OGM_1}), the information theoretic exact method (Algorithm~\ref{alg:ITEMS}) and the triple momentum method (Algorithm~\ref{alg:TMM}) due to the slightly different structure of their potential functions. However, this trick is valid for all other methods presented in this section, as well as those presented in Appendix~\ref{a-backandforth} and the proximal accelerated methods from Section~\ref{c-prox}.

\begin{algorithm}[!ht]
  \caption{Wrapper: monotone accelerated methods}
  \label{alg:accel_monotone}
  \begin{algorithmic}[1]
    \REQUIRE Pick an algorithm $\mathcal{A}$ among Algorithms $\{$\ref{alg:FGM_1},
\ref{alg:FGM_2}, \ref{alg:FGM_3}, \ref{alg:FGM_1_strconvex}, \ref{alg:FGM_1_strcvx_constmomentum}, \ref{alg:FGM_2_strcvx_constmomentum}, \ref{alg:fista_1}, \ref{alg:MST_strconvex}, \ref{alg:FGM_2_strconvex}, \ref{alg:FGM_3_strconvex}, \ref{alg:fista_1_general_BT}, \ref{alg:MST_strconvex_general_BT}$\}$ and use the same inputs the same associated problem inputs, including an initial guess $x_0\in\mathbb{R}^d$.
    \STATE {\bf Initialize} Execute initialization step from the chosen algorithm, as well as $\tilde{x}_0=x_0$.
    \FOR{$k=0,\ldots$}
      \STATE $x_k=\tilde{x}_k$
      \STATE Execute one iteration of the chosen algorithm (one update for each sequence).
      \STATE $\tilde{x}_{k+1}=\argmin_x\{F(x):\,x\in\{x_{k+1},\tilde{x}_k\}\}$
    \ENDFOR
    \ENSURE Approximate solution $x_{k+1}$.
  \end{algorithmic}
\end{algorithm}

The following fact summarizes the result for this method. In a nutshell, the desired result (i.e., same worst-case guarantees as previous algorithms while having a monotone sequence $\{F(\tilde{x}_k)\}_k$) is achieved due to two fact. First, the sequence $\{\tilde{x}_k\}_k$ is constructed for satisfying $F(\tilde{x}_{k+1})\leq F(x_{k+1})$, and hence the potential function analyses of all the algorithmic schemes is preserved. Secondly, the sequence is built for satisfying $F(\tilde{x}_{k+1})\leq F(\tilde{x_k})$ and hence it is monotonically decreasing.

\begin{theorem}\label{thm:monotone_convergence} \looseness=-1 Let $\mathcal{A}$ be an algorithm in \{\ref{alg:FGM_1},
\ref{alg:FGM_2}, \ref{alg:FGM_3}, \ref{alg:FGM_1_strconvex}, \ref{alg:FGM_1_strcvx_constmomentum}, \ref{alg:FGM_2_strcvx_constmomentum}, \ref{alg:fista_1}, \ref{alg:MST_strconvex}, \ref{alg:FGM_2_strconvex}, \ref{alg:FGM_3_strconvex}, \ref{alg:fista_1_general_BT}, \ref{alg:MST_strconvex_general_BT}\}, let $F\in\mathcal{F}_{0,\infty}$ be the convex function on which this algorithm is applied (possibly $F\defeq f+h$ with $f$ and $h$ satisfying the input requirements under which $\mathcal{A}$ operates), and $x_0\in\mathbb{R}^d$. Further let $\{\tilde{x}_k\}_k$ be the sequence generated by Algorithm~\ref{alg:accel_monotone} with $\mathcal{A}$ on $F$ and $x_0$. It follows that $F(\tilde{x}_N)$ satisfies the same worst-case guarantee as that of $\mathcal{A}$:
\[ F(\tilde{x}_N)-F_\star\leq \frac{L\|x_0-x_\star\|^2_2}{2 A_N},\]
with $A_N$ being defined in $\mathcal{A}$. In addition, the sequence $\{F(\tilde{x}_k)\}_k$ is monotonically decreasing.
\end{theorem}
\begin{proof} The proof follows from the fact that worst-case guarantees for all algorithms under consideration rely on the same potential function:
\[ \phi_k\defeq A_k(F(x_k)-F_\star)+\frac{L+\mu A_k}{2}\|z_k-x_\star\|^2_2.\]
Defining the alternate potential $\tilde{\phi}_k$ as \[ \tilde{\phi}_k\defeq A_k(F(\tilde{x}_k)-F_\star)+\frac{L+\mu A_k}{2}\|z_k-x_\star\|^2_2,\]
it follows from $F(\tilde{x}_k)\leq F(x_k)$ and $A_k\geq 0$ that
\[ \tilde{\phi}_k\leq \phi_k \]
for all $k\geq 0$.
Furthermore, it follows from the fact that 
\[ \tilde{\phi}_{k+1}\leq \phi_{k+1}\leq \phi_k\]
for all $x_k\in\mathbb{R}^d$ that we can choose $x_k=\tilde{x}_k$, thereby reaching $\tilde{\phi}_{k+1}\leq \tilde{\phi}_k$. Therefore, it holds that
\[ A_N (F(\tilde{x}_N)-F_\star)\leq \tilde{\phi}_{N}\leq \tilde{\phi}_0=\frac{L\|x_0-x_\star\|^2_2}{2}.\]
Hence, the same worst-case guarantees are achieved.

Finally, monotonicity is obtained by construction. That is, the sequence satisfies $F(\tilde{x}_{k+1})\leq F(\tilde{x}_k)$  (for all $k\geq 0$) and hence is monotonically decreasing.
\end{proof}
\subsection{Beyond Euclidean Geometries using Mirror Maps}\label{s:acc_MD}
In this section, we put ourselves in a slightly different scenario, often referred to as the \emph{non-Euclidean setting} or the \emph{mirror descent setup}. We consider the convex minimization problem:
\begin{equation}\label{eq:Comp_MD}
    F_\star=\min_{x\in\mathbb{R}^d} \{F(x)\defeq f(x)+h(x)\},
\end{equation}
with $h,f:\mathbb{R}^d\rightarrow\mathbb{R}\cup\{+\infty\}$ closed, convex, and proper. Furthermore, we assume $f$ to be differentiable and to have Lipschitz gradients with respect to some (possibly non-Euclidean) norm $\|.\|$. That is, with the corresponding dual norm denoted by $\|s\|_*=\sup_x\{\langle s;x\rangle\,:\, \|x\|\leq 1\}$, we require
\[ \| \nabla f(x)-\nabla f(y)\|_*\leq L \|x-y\|\]
for all $x,y\in\dom h$. In this setting, inequality~\eqref{eq:smooth_res} also holds (see Appendix~\ref{app:inequalities_nonEuclid}), and we, perhaps abusively, also denote $f\in\mathcal{F}_{0,L}(\dom h)$.

To solve~\eqref{eq:Comp_MD}, we define a few additional ingredients. First, we pick a $1$-strongly convex, closed proper function $w:\mathbb{R}^d\rightarrow\mathbb{R}\cup\{+\infty\}$ such that $\dom h\subseteq \dom w$. (Recall that by assumption on $h$,  $\dom h\neq \emptyset$, and therefore $\dom w\neq\emptyset$.) These assumptions ensure that the proximal operations below are well defined; $w$ is commonly referred to as the \emph{distance generating function}. Under additional technical assumptions, $w$ could be chosen as a strictly convex function instead, and similar algorithms can be used, but we focus on the strongly convex case here.

Finally, pick $g_w(y)\in\partial w(y)$, and define the Bregman divergence generated by $w$ as
\begin{equation}\label{eq:Breg_dist}
D_w(x;y)=w(x)-w(y)-\langle g_w(y);x-y\rangle,
\end{equation}
which we use below as a notion of distance to generalize the previous proximal operator~\eqref{eq:prox_chapt_Nest}. Note that the Bregman divergence $D_w(.;.)$ generated by any subgradient of $w$ at $z_k$ is considered valid here.

The base ingredient we use to solve~\eqref{eq:Comp_MD} is the Bregman proximal gradient step, with step size $\tfrac{a_k}{L}$:
\begin{equation}\label{eq:breg_ForwardBackward}
z_{k+1}= \argmin_y\left\{ \frac{a_k}{L} \left(h(y)+\langle \nabla f(y_k);y-y_k\rangle\right) +D_w(y;z_k)\right\},
\end{equation}
which corresponds to the usual Euclidean proximal gradient step when $w(x)=\tfrac12\|x\|^2_2$. Under previous assumptions,~\eqref{eq:breg_ForwardBackward} is well defined and we can explicitly write
\[ g_w(z_{k+1})\ni g_w(z_k)-\frac{a_k}{L} \left(\nabla f(y_k)+g_h(z_{k+1})\right),\]
with some $g_w(z_{k+1})\in\partial w(z_{k+1})$, $g_w(z_{k})\in\partial w(z_{k})$, and $g_h(z_{k+1})\in\partial h(z_{k+1})$. 

Under this construction, one can rely on~\eqref{eq:breg_ForwardBackward} to solve~\eqref{eq:Comp_MD} using Algorithm~\ref{alg:IGA}. Note that when $w$ is differentiable (which is usually the case but which necessitates further discussions when requiring $w$ to be closed, convex, and proper), we often refer to $\nabla w$ as a \emph{mirror map}. This refers to a bijective mapping due to the strong convexity and differentiability of $w$. In this case, the iterations can be described as \[ \nabla w(z_{k+1})=\nabla w(z_k)-\frac{a_k}{L} \left(\nabla f(y_k)+g_h(z_{k+1})\right).\]

\begin{algorithm}[!ht]
  \caption{Proximal accelerated Bregman gradient method}
  \label{alg:IGA}
  \begin{algorithmic}[1]
    \REQUIRE $h\in\mathcal{F}_{0,\infty}$, $f\in\mathcal{F}_{0,L}(\dom h)$, $w\in\mathcal{F}_{1,\infty}$ with $\dom h\subseteq\dom w$, and $x_0\in\dom h$ (such that $\partial w(z_0)\neq\emptyset$).
    \STATE \textbf{Initialize} $z_0=x_0$ and $A_0=0$.
    \FOR{$k=0,\ldots$}
      \STATE $a_k=\frac{1+\sqrt{4A_k+1}}{2}$
      \STATE $A_{k+1}=A_k+a_k$
      \STATE $y_{k}=  \frac{A_k}{A_{k+1}}x_k+\left(1-\frac{A_k}{A_{k+1}}\right) z_k$
      \STATE $z_{k+1}= \argmin_y\left\{ \frac{a_k}{L} \left(h(y)+\langle \nabla f(y_k);y-y_k\rangle\right) +D_w(y;z_k)\right\}$
      \STATE $x_{k+1}=\frac{A_k}{A_{k+1}}x_k+\left(1-\frac{A_k}{A_{k+1}}\right)z_{k+1}$
    \ENDFOR
    \ENSURE Approximate solution $x_{k+1}$.
  \end{algorithmic}
\end{algorithm}

Theorem~\ref{thm:IGA} provides a convergence guarantee for Algorithm~\ref{alg:IGA} by a potential argument.
\begin{theorem}\label{thm:IGA} Let $h\in\mathcal{F}_{0,\infty}$, $f\in\mathcal{F}_{0,L}(\dom h)$, $w\in\mathcal{F}_{1,\infty}$ with $\dom h\subseteq\dom w$,  $x_\star\in\argmin_{x} \{ F(x)\defeq f(x)+h(x)\}$, and $k\in\mathbb{N}$. For any $x_k,z_k\in\dom h$ such that $\partial w(z_k)\neq\emptyset$ and $A_k\geq0$; the iterates of Algorithm~\ref{alg:IGA} satisfy
\begin{equation*}
\begin{aligned}
A_{k+1}(F(x_{k+1})-F_\star)&+L D_w(x_\star;z_{k+1})\\ &\leq A_k (F(x_k)-F_\star)+L D_w(x_\star;z_k),
\end{aligned}
\end{equation*}
where $A_{k+1}=A_k+\frac{1+\sqrt{4A_k+1}}{2}$ and $D_w(.;.)$ is a Bregman divergence~\eqref{eq:Breg_dist} with respect to $w$. Furthermore, $\partial w(z_{k+1})\neq\emptyset$.
\end{theorem}
\begin{proof} First, $z_k$ is feasible, i.e., $z_k\in\dom h$, by construction. Indeed, $z_0$ is feasible by assumption, and the following iterates $z_k$ ($k>0$) are obtained after proximal steps; hence, $\partial h(z_k)\neq\emptyset$, and therefore, $z_k\in\dom h$.

Second, it can be directly verified that $0\leq \tfrac{A_{k}}{A_{k+1}}\leq 1$ given that $A_{k+1}\geq A_k\geq 0$. It follows from $z_0=x_0\in\dom h$ that the elements of $\{y_k\}$ and $\{x_k\}$ are obtained as convex combinations of elements of $\{z_k\}$. Hence, the sequences $\{y_k\}$ and $\{x_k\}$ are also in $\dom h$, which is convex. The rest of the proof consists of a weighted sum of the following inequalities.
\begin{itemize}
    \item Convexity of $f$ between $x_\star$ and $y_{k}$ with weight $\lambda_1=A_{k+1}-A_k$:
    \[ f(x_\star) \geq f(y_{k})+\langle \nabla f(y_{k});x_\star-y_k\rangle.\]
    \item Convexity of $f$ between $x_k$ and $y_{k}$ with weight $\lambda_2=A_k$:
    \[ f(x_k) \geq f(y_{k})+\langle \nabla f(y_{k});x_k-y_{k}\rangle.\]
    \item Smoothness of $f$ between $y_{k}$ and $x_{k+1}$ (\emph{descent lemma}) with weight $\lambda_3=A_{k+1}$:
    \[ f(y_{k}) + \langle \nabla f(y_{k}); x_{k+1}-y_{k}\rangle +\frac{L}{2}\|x_{k+1}-y_{k}\|^2 \geq f(x_{k+1}).\]
    \item Convexity of $h$ between $x_\star$ and $z_{k+1}$ with weight $\lambda_4=A_{k+1}-A_k$:
    \[ h(x_\star)\geq h(z_{k+1})+\langle g_h(z_{k+1});x_\star-z_{k+1}\rangle,\]
    with $g_h(z_{k+1})\in\partial h(z_{k+1})$.
    \item Convexity of $h$ between $x_k$ and $x_{k+1}$ with weight $\lambda_5=A_k$:
    \[ h(x_k)\geq h(x_{k+1})+\langle g_h(x_{k+1});x_k-x_{k+1}\rangle,\]
    with $g_h(x_{k+1})\in\partial h(x_{k+1})$.
    \item Convexity of $h$ between $z_{k+1}$ and $x_{k+1}$ with weight $\lambda_6=A_{k+1}-A_k$
    \[ h(z_{k+1})\geq h(x_{k+1})+\langle g_h(x_{k+1});z_{k+1}-x_{k+1}\rangle.\]
    \item Strong convexity of $w$ between $z_{k+1}$ and $z_k$ with weight $\lambda_7=L$
    \[ w(z_{k+1})\geq w(z_k)+\langle g_w(z_{k});z_{k+1}-z_k\rangle+\frac12\|z_k-z_{k+1}\|^2,\]
    with $g_w(z_k)\in\partial w(z_k)$.
\end{itemize}
We thus obtain the following inequality:
\begin{equation*}
\begin{aligned}
0\geq &\lambda_1 [f(y_{k})-f_\star+\langle \nabla f(y_{k});x_\star-y_k\rangle]\\
&+\lambda_2[f(y_{k})-f(x_k)+\langle \nabla f(y_{k});x_k-y_{k}\rangle]\\
&+\lambda_3[f(x_{k+1})-(f(y_{k}) + \langle \nabla f(y_{k}); x_{k+1}-y_{k}\rangle +\frac{L}{2}\|x_{k+1}-y_{k}\|^2)]\\
&+\lambda_4[ h(z_{k+1})-h(x_\star)+\langle g_h(z_{k+1});x_\star-z_{k+1}\rangle]\\
&+\lambda_5[ h(x_{k+1})-h(x_k)+\langle g_h(x_{k+1});x_k-x_{k+1}\rangle]\\
&+\lambda_6[ h(x_{k+1})-h(z_{k+1})+\langle g_h(x_{k+1});z_{k+1}-x_{k+1}\rangle]\\
&+\lambda_7[ w(z_k)-w(z_{k+1})+\langle g_w(z_{k});z_{k+1}-z_k\rangle+\frac12\|z_k-z_{k+1}\|^2].
\end{aligned}
\end{equation*}
Now, by substituting (using $g_w(z_{k+1})\in\partial w(z_{k+1})$)
\begin{equation*}
\begin{aligned}
y_k&=\frac{A_k}{A_{k+1}}x_k+\left(1-\frac{A_k}{A_{k+1}}\right)z_k\\
g_w(z_{k+1})&=g_w(z_{k})-\frac{A_{k+1}-A_k}{L} \left(\nabla f(y_k)+g_h(z_{k+1})\right)\\
x_{k+1}&=\frac{A_{k}}{A_{k+1}}x_k+\left(1-\frac{A_k}{A_{k+1}}\right)z_{k+1},
\end{aligned}
\end{equation*}
we obtain exactly $\|x_{k+1}-y_{k}\|^2=\tfrac{(A_{k+1}-A_k)^2}{A_{k+1}^2}\|z_k-z_{k+1}\|^2$. The weighted sum can then be reformulated exactly as
\begin{equation*}
\begin{aligned}
A_{k+1}&(F(x_{k+1})-F_\star)+L D_w(x_\star,z_{k+1})\\ &\leq A_k (F(x_k)-F_\star)+L D_w(x_\star,z_k)\\
& \quad +\frac{ (A_k-A_{k+1})^2-A_{k+1} }{A_{k+1}} \frac{L}{2}\|z_k-z_{k+1}\|^2,
\end{aligned}
\end{equation*}
and we obtain the desired inequality from selecting $A_{k+1}$ that satisfies $A_{k+1}\geq A_k$ and
\[ (A_k-A_{k+1})^2-A_{k+1} =0.\qedhere\]
\end{proof}

We conclude this section by providing a final corollary to describe the worst-case performance of the method.
\begin{corollary}  Let $h\in\mathcal{F}_{0,\infty}$, $f\in\mathcal{F}_{0,L}(\dom h)$, $w\in\mathcal{F}_{1,\infty}$ with $\dom h\subseteq\dom w$, and $x_\star\in\argmin_{x} \{ F(x)\defeq f(x)+h(x)\}$. For any $x_0\in\dom h$ such that $\partial w(x_0)\neq\emptyset$ and any $N\in\mathbb{N}$; the iterates of Algorithm~\ref{alg:IGA} satisfy
\[F(x_N)-F_\star\leq \frac{L D_w(x_\star;z_0)}{A_N}\leq\frac{4L D_w(x_\star;z_0)}{N^2} .\]
\end{corollary}
\begin{proof} The claim directly follows from previous arguments using the potential $\phi_k\defeq A_k (F(x_k)-F_\star)+L D_w(x_\star;z_k)$, along with $A_N\geq \frac{N^2}{4}$ from~\eqref{eq:conv_Ak_fgm}.
\end{proof}
\begin{remark} \emph{Bregman} first-order methods are often split into two different families: \emph{mirror descent} and \emph{dual averaging} (which we did not explicitly mention here), for which we refer the reader to the discussions in~\citep[Chapter 4]{bubeck2015convex}.
The method presented in this section is essentially a case of~\citep[Algorithm 1]{tseng2008accelerated}, and it corresponds to~\citep[``Improved interior gradient algorithm'']{auslender2006interior} when the norm is Euclidean. It is also similar to~\citep[Algorithm 3]{tseng2008accelerated} and to~\citep[``Method of similar triangles'']{gasnikov2018universal}, which are ``dual-averaging'' versions of the same algorithm: they are essentially equivalent in the Euclidean setup without constraints. The method presented here enjoys a number of variants (see, e.g., discussions in \citep{tseng2008accelerated}), some of which may involve two projections per iteration, as in~\citep[Section 3]{Nest03},~\citep[Section 3]{lan2011primal}. The method can also naturally be embedded with a backtracking procedure, exactly as in previous sections. Finally, such methods can be adapted to strong convexity, either on $f$, as in~\citep{gasnikov2018universal} or on $h$, as in~\citep{diakonikolas2021complementary}.
\end{remark}
\begin{remark} Note that the technique for rendering accelerated methods monotone (see Section~\ref{s:monotone_variants}) also directly applies to Algorithm~\ref{alg:IGA}.
\end{remark}
\begin{remark} Beyond the Euclidean setting where $w(x)=\tfrac12\|.\|^2_2$, a classical example of the impact of the non-Euclidean setup optimizes a simple function over the simplex. In this case, writing $x^{(i)}$ for the $i$th component of some $x\in\mathbb{R}^d$, we consider the situation where $h$ is the indicator function of the simplex
\[ h(x)=\left\{\begin{array}{ll}
    0 & \text{if } \sum_{i=1}^d x^{(i)}=1,\, x^{(i)}\geq 0 \, (i=1,\hdots,d)  \\
    +\infty & \text{otherwise,} 
\end{array}\right.\]
and $w$ is the \emph{entropy} (which is closed, convex, and proper, and $1$-strongly convex over the simplex for $\|.\|=\|.\|_1$; this is known as Pinsker's inequality). That is, define some $w_i:\mathbb{R}\rightarrow\mathbb{R}$:
\[ w_i(x)=\left\{\begin{array}{ll}
    x \log x & \text{ if $x>0$,}  \\
    0 & \text{if $x=0$,}\\
    +\infty & \text{ otherwise,}
\end{array}\right.\]
and set $w(x)=\sum_{i=1}^d w_i(x^{(i)})$. In this case, the expression for the Bregman proximal gradient step in Algorithm~\ref{alg:IGA} can be computed exactly, assuming $y_k^{(i)}\neq 0$ (for all $i=1,...,d$):
\begin{equation*}
\begin{aligned}
z_{k+1}^{(i)}&=\frac{y_k^{(i)}\exp\left[-\frac{a_k}{L} [\nabla f(y_k)]^{(i)}\right]}{\sum_{i=1}^d y_k^{(i)}\exp\left[-\frac{a_k}{L} [\nabla f(y_k)]^{(i)}\right]}.
\end{aligned}
\end{equation*}
Hence, we also have that $y_{k}^{(i)}\neq 0$ as long as $y_{0}^{(i)}\neq 0$; a common technique is to instantiate $x_0^{(i)}=\tfrac1d$. 

In this setup, a non-Euclidean geometry often provides a significant practical advantage when optimizing large-scale functions by improving the dependence from $d$ to $\ln d$ in the final complexity bound. In fact, here we have $D_w(x_\star,x_0)\leq \ln d$ compared to $D_{\tfrac12\|.\|_2^2}(x_\star,x_0)\leq \tfrac12$ in the Euclidean case, so the dependence in $d$ is seemingly better in the Euclidean case. However, the choice of norms has a very significant impact. When the gradient has a small Lipschitz constant measured with respect to $\|.\|=\|.\|_1$, that is,
\[\|\nabla f(x)-\nabla f(y)\|_{\infty}\leq L_1 \|x-y\|_1,\]
the Lipschitz constant might be up to $d$ times smaller than the constant computed using the Euclidean norm $\|.\|=\|.\|_2$ (using norm equivalences), i.e.,
\[\|\nabla f(x)-\nabla f(y)\|_{2}\leq L_2 \|x-y\|_2,\]
with $L_1\sim L_2/d$. The final complexity bound using a Euclidean geometry then reads 
\[
    F(x_N)-F_\star\leq \tfrac{2L_1 d}{N^2},
\] 
as compared to 
\[
    F(x_N)-F_\star\leq \tfrac{4L_1 \ln d}{N^2}
\]
using the geometry induced by the entropy. The impact of the choice of norm is discussed extensively in, e.g.~\citep[Example 2.1]{Judic07} and \citep{d2018optimal}.

Another related example is optimizing over a spectrahedron; see, for example, the nice introduction by~\citet[Chapter 4]{bubeck2015convex}. This setup is largely motivated in~\citep{Nest03}. We refer to~\citep{allen2017linear,diakonikolas2021complementary} and the references therein for further discussions on this topic.
\end{remark}

\section{Continuous-time Interpretations}

Before concluding this section, we survey another last popular approach to Nesterov's acceleration. The idea is to study continuous versions of first-order schemes, which enables simpler (less technical) proofs to emerge. One notorious caveat of such approaches is that they usually defer implementation details (i.e., discretization) to integration solvers (such as explicit Euler's scheme), but the simplicity of the proofs arguably renders this approach worth investigating. In particular, we see later that one usually does not use \emph{smoothness} when computing worst-case convergence speeds. That is, smoothness is intrinsically linked to the discretization procedure, and not at all to the convergence speed of the continuous-time processes.

The section is divided in three parts. We start by reviewing results related to the \emph{gradient flow}, a natural ordinary differential equation (ODE) for modelling first-order methods. Then, we continue with Nesterov's ODE, the continuous-time limit of Nesterov's accelerated gradient method as the step size vanishes~\citep{su2014differential}. We conclude with discussions and pointers to different results and research directions relying on ODE interpretations of Nesterov's method. 

For simplicity purposes, we make the choice of only presenting \emph{informal} arguments for obtaining the continuous-time versions of gradient descent and of Nesterov's acceleration.

\subsection{Gradient Flow}

A starting point for linking first-order methods to continuous-time processes is to minimize a convex function $f$ by following its \emph{gradient flow}
\begin{equation}\label{eq:gradient_flow}
\left\{ \begin{array}{l}\dot{x}(t)=-\nabla f(x(t)),\\
x(0)=x_0\end{array}\right.
\end{equation}
where $\dot{x}(t)=\tfrac{dx(t)}{dt}$ denotes the usual time derivative of $x$, and where the time $t\in\mathbb{R}$ takes the role of the usual iteration counter $k\in\mathbb{N}$. For technical reasons, we also assume throughout that $f$ is $L$-smooth, as it ensures the existence of a solution to~\eqref{eq:gradient_flow}. 

From~\eqref{eq:gradient_flow}, one can recover usual first-order methods by playing with different integration schemes. Hence, first-order methods can typically be interpreted in the light of integration schemes, through specific notions of numerical stability, see, e.g.,~\citep{scieur2017integration}. In short, whereas integration schemes typically aim at tracking the full trajectory of an ODE, optimization methods only target tracking the stationary point of~\eqref{eq:gradient_flow} thereby requiring less stringent notions of numerical stability.

As a particular case, the classical explicit Euler integration scheme applied to~\eqref{eq:gradient_flow} boils down to gradient descent for minimizing $f(x)$. On the other hand,~\eqref{eq:gradient_flow} can be obtained as a natural continuous-time counterpart to gradient descent with vanishing step size. That is, gradient descent with step size $\Delta$ can be written as
\[ \frac{x_{k+1}-x_k}{\Delta}=-\nabla f(x_k).\]
Informally, assume that the sequence $\{x_k\}_{k}$ is obtained as an approximation to a solution $x(t)$ to an ODE, with $x(0)=x_0$ and $x(t)\approx x(k\Delta)$. It is clear that 
\[x_{k+1}\approx x(t+\Delta)=x(t)+\Delta\dot{x}(t)+o(\Delta),\] 
and hence
\[ -\nabla f(x(t))=\frac{x(t+\Delta)-x(\Delta)}{\Delta}=\dot{x}(t)+\tfrac{o(\Delta)}{\Delta},\]
leading to~\eqref{eq:gradient_flow} by taking the limit $\Delta\rightarrow0$ on both sides.

Before moving to Nesterov's ODE, we glance at the convergence speed of the gradient flow towards an optimum when $f$ is a convex function. For doing that, we use a similar potential function as that used in the discrete setup.

\paragraph{Convergence speed.} We proceed essentially as in the previous section, but we see below that the proofs are much less technical. We introduce the (continuous) potential function:
\[\phi(t)\defeq a(t)\, (f(x(t))-f(x_\star))+\|x(t)-x_\star\|^2_2.\]
The analysis simply consists in showing that $\dot{\phi}(t)\leq 0$. Indeed, just as in the discrete setup, we can then write:
\[ \phi(t)\leq \phi(0),\]
with $\phi(t)\geq a(t) \, (f(x(t))-f(x_\star))$ and $\phi(0)=\|x(t)-x_\star\|^2_2$. Thereby, we reach
\[ f(x(t))-f(x_\star) \leq \frac{\phi(0)}{a(t)}=\frac{\|x(t)-x_\star\|^2_2}{a(t)},\]
and the worst-case speed of convergence of $f(x(t))$ towards $f(x_\star)$ is dictated by the growth rate of $a(t)$.
\begin{theorem} \label{thm:cont_gd} Let $f$ be a closed proper convex function, $x_\star\in\argmin_x$ $f(x)$, and let $a(t)=t$. For any $x(\cdot)$ solution to~\eqref{eq:gradient_flow} and any $t\geq 0$, it holds that
\[\dot{\phi}(t)\leq 0,\]
with
\[\phi(t)\defeq a(t)\,(f(x(t))-f(x_\star))+\|x(t)-x_\star\|^2_2.\]
\end{theorem}\begin{proof}
The proof simply uses a convexity inequality between $x(t)$ and $x_\star$. That is, explicit computations allow obtaining:
\[\dot{\phi}(t)=\dot{a}(t)\,(f(x(t))-f(x_\star))-a(t)\|\nabla f(x(t))\|^2_2-\langle \nabla f(x(t));x(t)-x_\star\rangle.\]
Then, a convexity inequality
 \[  f(x(t))-f(x_\star)\leq \langle \nabla f(x(t)); x(t)-x_\star\rangle\]
allows writing
\[ \dot{\phi}(t)\leq -a(t)\|\nabla f(x(t))\|^2_2+(\dot{a}(t)-1)\langle \nabla f(x(t));x(t)-x_\star\rangle,\]
and we reach the desired inequality using $\dot{a}(t)=1$: \[\dot{\phi}(t)\leq -a(t) \|\nabla f(x(t))\|^2_2\leq 0.\qedhere\]
\end{proof}
\begin{corollary}\label{cor:cont_gd} Let $f$ be a closed proper convex function and $x_\star\in\argmin_x f(x)$. For any $x(\cdot)$ solution to~\eqref{eq:gradient_flow} and any $t\geq 0$, it holds that
\[ f(x(t))-f(x_\star) \leq \frac{\|x(0)-x_\star\|^2_2}{t}.\]
\end{corollary}
\begin{proof}
The proof follows from $\dot{\phi}(t)\leq 0$ (Theorem~\ref{thm:cont_gd}) and \[f(x(t))-f(x_\star) \leq \frac{\phi(0)}{a(t)}=\frac{\|x(0)-x_\star\|^2_2}{a(t)},\]
with $a(t)=t$.
\end{proof}

\subsection{An Ordinary Differential Equation for Nesterov's Method}
In this section, we briefly study a different ODE, namely
\begin{equation}\label{eq:ODE_NEST}
    \left\{\begin{array}{l}
    \ddot{y}(t)+\frac{3}{t+T}\dot{y}(t)+\nabla f(y(t))=0,\\
    y(0)=x_0,\\
    \dot{y}(0)=0,
    \end{array}\right.
\end{equation}
where $T\geq 0$ is some constant. This ODE with $T=0$ was proposed in~\citep{su2014differential} as the continuous-time interpretation of Nesterov's method. The case $T>0$ might be considered instead for simplicity of the exposition, rendering the ODE directly well-defined even for $t=0$. For simplicity, we also consider Algorithm~\ref{alg:FGM_2} with $a_k=\tfrac{k}{2}$ for $k=0,1,\hdots$ (leading to $A_k=\frac{k^2}{4}$), that is,
\begin{equation}\label{eq:nest_simp}
\begin{aligned}
x_{k+1}&=y_k-\Delta\nabla f(y_k)\\
y_{k+1}&=x_{k+1}+\frac{k}{k+3}(x_{k+1}-x_k),
\end{aligned}
\end{equation}
for $k=0,1,\hdots$. In this setup, the method can also be described in terms of a single sequence:
\begin{equation}\label{eq:nest_simp1}
\begin{aligned}
y_{k+1}=y_k&-\Delta\nabla f(y_k)+\frac{k}{k+3}(y_k-y_{k-1})\\&-\frac{k}{k+3}\Delta(\nabla f(y_k)-\nabla f(y_{k-1})).
\end{aligned}
\end{equation}
As shown in~\citep{su2014differential}, it turns out that~\eqref{eq:ODE_NEST} can be seen as the continuous-time limit of~\eqref{eq:nest_simp} as $\Delta\rightarrow 0$. Using a similar informal development to that of~\citep{su2014differential}, a few approximations allow to arrive to the desired ODE. For doing that, let us denote by $y(t)$ a trajectory of the limiting ODE, and let us assume that it is approximated by Nesterov's method (acting as a numerical integrator) as $y(t)\approx y_{t/\sqrt{\Delta}}$ (i.e., we make the correspondence $y\left(k\sqrt{\Delta}\right)\approx y_k$). Note the time scaling in $\sqrt{\Delta}$ instead of $\Delta$ due to the second-order dynamics of the system. Next, we use the following Taylor expansions
\begin{equation*}
    \begin{aligned}
    \frac{y\left(t+\sqrt{\Delta}\right)-y(t)}{\sqrt{\Delta}}&=\dot{y}(t)+\frac{\sqrt{\Delta}}{2}\ddot{y}(t)+o\left(\sqrt{\Delta}\right)\\
    \frac{y\left(t-\sqrt{\Delta}\right)-y(t)}{\sqrt{\Delta}}&=-\dot{y}(t)+\frac{\sqrt{\Delta}}{2}\ddot{y}(t)+o\left(\sqrt{\Delta}\right)\\
    \sqrt{\Delta}\nabla f\left(y\left(t+\sqrt{\Delta}\right)\right)&=\sqrt{\Delta}\nabla f(y(t))+o\left(\sqrt{\Delta}\right)\\
    \frac{t}{t+3\sqrt{\Delta}}&=1-3\frac{\sqrt{\Delta}}{t}+o\left(\sqrt{\Delta}\right)\\
    \end{aligned}
\end{equation*}
for writing
\begin{align*}
\dot{y}(t)+&\frac{\sqrt{\Delta}}{2}\ddot{y}(t)+o\left(\sqrt{\Delta}\right)\\=&\left(1-3\frac{\sqrt{\Delta}}{t}\right)\left(\dot{y}(t)-\frac{\sqrt{\Delta}}{2}\ddot{y}(t)+o\left(\sqrt{\Delta}\right)\right)-\sqrt{\Delta}\nabla f(y(t)).\end{align*}
Simplifying these expressions, dividing all terms by $\sqrt{\Delta}$ and taking the limit $\Delta\rightarrow 0$ leads to the desired
\[ \ddot{y}(t)+\frac{3}{t}\dot{y}(t)+\nabla f(y(t))=0.\]

\paragraph{Convergence speed.} As for the gradient flow, the analysis only requires an appropriate (continuous) potential function $\phi(t)$. The desired conclusion follows from showing that $\dot{\phi}(t)<0$, as provided by the following theorem.
\begin{theorem} \label{thm:cont_nest} Let $f$ be a closed proper convex function and $x_\star\in\argmin_x f(x)$. For any $x(\cdot)$ solution to~\eqref{eq:ODE_NEST} and any $t\geq 0$, it holds that
\[\dot{\phi}(t)\leq 0,\]
with
\[\phi(t)\defeq (t+T)^2 (f(x(t))-f_\star)+2\| x(t)+\tfrac{t+T}{2}\dot{x}(t)-x_\star\|^2_2.\]
\end{theorem}
\begin{proof} From the expression of $\phi(t)$, it is relatively straightforward to obtain that
\begin{equation*}
\begin{aligned}
\dot{\phi}(t)=& 2(t+T)(f(x(t))-f(x_\star))+(t+T)^2\langle \nabla f(x(t));\dot{x}(t)\rangle\\&+4\langle (x(t)+\frac{t+T}{2}\dot{x}(t)-x_\star);\frac{3}{2}\dot{x}(t)+\frac{t+T}{2}\ddot{x}(t)\rangle.
\end{aligned}
\end{equation*}
Using~\eqref{eq:ODE_NEST}, one can substitute the expression of $\ddot{x}(t)$, leading to
\begin{equation*}
\begin{aligned}
\dot{\phi}(t)&= 2(t+T)(f(x(t))-f(x_\star))+2(t+T)\langle \nabla f(x(t));x_\star-x(t)\rangle
\end{aligned}
\end{equation*}
and it follows from convexity that $\dot{\phi}(t)\leq 0$, as desired.
\end{proof}
\begin{corollary}\label{cor:cont_agd} Let $f$ be a closed proper convex function and $x_\star\in\argmin_x f(x)$. For any $x(\cdot)$ solution to~\eqref{eq:ODE_NEST} and any $t\geq 0$, it holds that
\[ f(x(t))-f(x_\star) \leq \frac{T^2(f(x(0))-f(x_\star))+\|x(0)-x_\star\|^2_2}{(t+T)^2}.\]
\end{corollary}
\begin{proof}
The proof follows from $\dot{\phi}(t)\leq 0$ (Theorem~\ref{thm:cont_nest}), implying $\phi(t)\leq\phi(0)$ and thereby \[f(x(t))-f(x_\star) \leq \frac{\phi(0)}{(t+T)^2}.\qedhere\]
\end{proof}

\subsection{Continuous-time Approaches to Acceleration: Summary}

In this section, we saw that some ordinary differential equations can be used for modelling gradient and accelerated gradient-type methods. The corresponding convergence proofs are much simpler, as they only involve using a single inequality, namely convexity between two points: $x(t)$ and $x_\star$. However, convergence speeds of continuous-time versions of algorithms might not be representative of their behaviors, as the corresponding ODE might be complicated to integrate, and as using numerical integration solvers might break the potentially nice convergence properties of the continuous-time dynamics.

The content of this section is explored at length in many references, see e.g.,~\citep{su2014differential,krichene2015accelerated,wibisono2016variational,attouch2018fast}. We discuss a few extensions and limitations below, before concluding the section.

\paragraph{Integration methods and optimization.} Continuous-time formulations of gradient-based methods cannot be implemented as is on digital computers. Therefore, continuous-time analyses cannot provide a complete picture on the topic without incorporating numerical integration methods into the analyses. Another symptom of this incompleteness is that of different optimization methods giving rise to the same limiting ODEs. For instance, the same limiting ODE is obtained from Polyak's heavy-ball, from Nesterov's accelerated methods, and from the triple momentum methods; see~\citep{shi2018understanding} and~\citep{sun2020high}.

\enlargethispage{\baselineskip}
This observation motivates different lines of works. In~\citep{scieur2017integration}, the authors focus on the gradient flow~\eqref{eq:gradient_flow} and propose different integration methods for recovering classical first-order methods in the quadratic minimization setup. In~\citep{su2014differential}, it is shown that forward Euler integration of Nesterov's ODE~\eqref{eq:ODE_NEST} leads to a heavy-ball type method, close to Nesterov's acceleration. In~\citep{shi2019acceleration}, the authors obtain accelerated methods by integrating ``high resolution'' variants of the \emph{accelerated ODEs} via symplectic methods (partially implicit, and partially explicit integration rules). Various discussions, developments, and connections between continuous-time systems and their discrete counterparts are further presented in~\citet{diakonikolas2019approximate,siegel2019accelerated,sanz2021connections}.

\paragraph{Strongly convex ODEs.} In the strongly convex case, limiting ODEs for ``stationnary'' accelerated methods such as Nesterov's method with constant momentum (Algorithm~\ref{alg:FGM_1_strcvx_constmomentum}) or triple momentum method (Algorithm~\ref{alg:TMM}) are also presented in different works; see, e.g.,~\citep{shi2018understanding,sun2020high}.

\paragraph{Continuized methods.} As previously discussed, it might not be simple to use classical continuous-time methods on digital computers (nontrivial integration schemes must be deployed). A family of so-called ``\emph{continuized methods}'' are directly implementable while keeping the benefits of the continuous-time approaches. Those methods rely on \emph{randomized} discretizations of the continuous-time process~\citep{even2021continuized}.

\section{Notes and References}\label{s:notes_ref_chapt_nest}

\paragraph{Estimate sequences, potential functions, and differential equations.}
Potential functions were already used in the original paper by~\citet{Nest83} to develop accelerated methods. \citet{Nest03a} developed estimate sequences as an alternate, more constructive, approach to obtaining optimal first-order methods. Since then, both approaches have been used in many references on this topic, in a variety of settings. \citet{tseng2008accelerated} provides a helpful unified view of accelerated methods. Estimate sequences have been extensively studied by, e.g.,~\citet{Nest13,baes2009estimate,devolder2011stochastic,kulunchakov2019estimate}. Another related approach is that of the \emph{approximate duality gap}~\citet{diakonikolas2019approximate} which is a constructive approach to estimate sequences/potential functions with a continuous-time counterpart.

\paragraph{Beyond Euclidean geometries.} 
\emph{Mirror descent} dates back to the work of~\citet{Nemi83}. It has been further developed and used in many subsequent works~\citep{Bentc01,Nest03,Nest09,xiao2010dual,juditsky2014deterministic,diakonikolas2021complementary}. Sound pedagogical surveys can be found in~\citep{beck2003mirror,juditsky2011first,juditsky2011first2,bubeck2015convex}.

Beyond the setting described in this section, mirror descent has also been studied in the \emph{relative smoothness} setting, introduced by~\citet{Baus16}---see also~\citep{teboulle2018simplified}---, and extended to the notion of \emph{relative strong convexity} by~\citet{Lu18a}. However, acceleration remains an open issue in the context of relative smoothness and relative strong convexity, and it is generally unclear which additional assumptions allow accelerated rates. It is, however, clear that additional assumptions are required, as emphasized by the lower bound provided by~\citet{dragomir2019optimal}. In particular, accelerated schemes are known under an additional \emph{triangle scaling inequality}~\citep{hanzely2021accelerated,gutman2018unified}.

\paragraph{Lower complexity bounds.} Lower complexity bounds have been studied in a variety of settings to establish limits on the worst-case performance of black-box methods. The classical reference on this topic is the book by~\citet{Book:NemirovskyYudin}.

Of particular interest to us,~\citet{nemirovsky1991optimality,nemirovsky1992information} establish the optimality of the Chebyshev and of the conjugate gradient methods for convex quadratic minimization; see, also, a complete picture provided in the course notes by~\citet{nemirovskinotes1995}. Lower bounds for black-box first-order methods in the context of smooth convex and smooth strongly convex optimization can be found in~\citep{Nest03a}. The final lower bound for black-box smooth convex minimization was obtained by~\citet{drori2017exact}; it demonstrates the optimality of the optimized gradient method, as well as that of conjugate gradients, as discussed earlier in this section. Lower bounds for $\ell_p$ norms in the mirror descent setup are constructed in~\citet{guzman2015lower}, whereas a lower bound for mirror descent in the relative smoothness setup is provided by~\citet{dragomir2019optimal}.

\paragraph{Changing the performance measure.} 
Obtaining (practical) accelerated method for other types of convergence criteria, such as gradient norms, is still not a fully settled issue. These criterion are important in other contexts, including dual methods, and can be used to draw links between methods intrinsically designed to solve convex problems and those used in nonconvex settings, where the goal is to find stationary points. There are a few \emph{tricks} that make it possible to pass from a guarantee in one context to another one. For example, a regularization trick was proposed in~\citep{nesterov2012make} that yields approximate solutions with small gradient norm. Beyond that, in the context of smooth convex minimization, recent progresses have been made by~\citet{kim2018optimizing}, who designed an optimized method for minimizing the gradient norm after a given number of iterations. This method was analyzed through potential functions in~\citep{diakonikolas2021potential}, and its geometric structure was further explored and exploited in~\citep{lee2021geometric}. Corresponding lower bounds, based on quadratic minimization, for a variety of performance measures can be found in~\citep{nemirovsky1992information}.

\paragraph{Adaptation and backtracking line-searches.} 
The idea of using backtracking line-searches is classical and is attributed to~\citet{goldstein1962cauchy} and~\citet{armijo1966minimization}; see also discussions in~\citep{nocedal2006numerical,bonnans2006numerical}. It was already incorporated in the original work of~\citet{Nest83} to estimate the smoothness constant within an accelerated gradient method. Since then, many works on the topic have relied heavily on this technique, which is often adapted to obtain better practical performance; see, for example,~\citep{scheinberg2014fast,chambolle2016introduction,florea2018accelerated,calatroni2019backtracking}. A more recent adaptive step size strategy (without line-search) can be found in~\citep{malitsky20a}.

\paragraph{Numerical stability, inexactness, stochasticity, and randomness.} 

The ability to use approximate first-order information, be it stochastic or deterministic, is key for tackling certain problems for which computing the exact gradient is expensive. Deterministic (or adversarial) error models are studied in, e.g.,~\citep{dAsp05,schmidt2011convergence,devolder2014first,devolder2013exactness,devolder2013intermediate,aybat2020robust} through different noise models. Such approaches can also be deployed when the projection/proximal operation is computed approximately~\citep{schmidt2011convergence,villa2013accelerated} (see also Section~\ref{c-prox} and the references therein).

Similarly, stochastic approximations and incremental gradient methods are key in many statistical learning problems, where samples are accessed one at a time and for which it is not desirable to optimize beyond the data accuracy~\citep{bottou2007tradeoffs}. For this reason, the old idea of stochastic approximations~\citep{robbins1951stochastic} is still widely used and remains an active area of research. The ``optimal'' variants of stochastic approximations were developed much later~\citep{lanstoch2008} with the rise of machine learning applications. In this context, it is not possible to asymptotically accelerate convergence rates but only to accelerate the transient phase toward a purely stochastic regime; see also~\citep{hu2009accelerated,xiao2010dual,devolder2011stochastic,lan2012optimal,dvurechensky2016stochastic,aybat2019universally,gorb20}---in particular, we note that ``stochastic'' estimate sequences were developed in~\citep{devolder2011stochastic,kulunchakov2019estimate}. The case of stochastic noise arising from sampling an objective function that is a \emph{finite sum} of smooth components attracted substantial attention in the 2010s, starting with ~\citep{schmidt2017minimizing,johnson2013accelerating,shalev2013stochastic,defazio2014saga,defazio2014finito,mairal2015incremental} and was then extended to feature acceleration techniques~\citep{shalev2014accelerated,allen2017katyusha,zhou2018simple,zhou2019direct}. Acceleration techniques also apply in the context of randomized block coordinate descent; see, for example,~\citet{nesterov2012efficiency,lee2013efficient,fercoq2015accelerated,nesterov2017efficiency}.

\paragraph{Higher-order methods.} Acceleration mechanisms have also been proposed in the context of higher-order methods. This line of work started with the cubic regularized Newton method introduced in~\citep{nesterov2006cubic} and its acceleration using estimate sequence mechanisms~\citep{nesterov2008accelerating}; see also~\citep{baes2009estimate,wilson2016lyapunov} and~\citep{monteiro2013accelerated} (which we also discuss in the next section). Optimal higher-order methods were presented by~\citep{gasnikov2019optimal}. It was not clear before the work of~\citet{nesterov2019implementable} that intermediate subproblems arising in the context of higher-order methods were tractable. The fact that tractability is not an issue has attracted significant attention to these methods.

\paragraph{Optimized methods.} 
Optimized gradient methods were discovered by~\citet{kim2016optimized}, based on the work by~\citet{Dror14}. Since then, optimized methods have been studied in various settings: incorporating constraints/proximal terms~\citep{kim2018another,taylor2017exact}; optimizing gradient norms~\citep{kim2018generalizing,kim2018optimizing,diakonikolas2021potential} (as an alternative to~\citep{nesterov2012make}); adapting to unknown problem parameters using exact line-searches~\citep{drori2019efficient} or restarts~\citep{kim2018adaptive}; and in the strongly convex case~\citep{van2017fastest,cyrus2018robust,park2021factor,drori2021optimal}. Such methods have also appeared in the context of fixed-point iterations~\citep{lieder2020convergence} and proximal methods~\citep{kim2019accelerated,barre2020principled}.

\paragraph{On obtaining proofs from this section.} The
worst-case performance of first-order methods can often be computed numerically. This has been shown in~\citep{Dror14,drori2014contributions,drori2016optimal,taylor2017smooth} through the introduction of performance estimation problems. Such techniques might be framed in different ways, e.g., from a purely optimization-based point of view~\citep{Dror14,taylor2017smooth} or from a control-theoretical perspective~\citep{lessard2016analysis,fazlyab2018analysis}. We provide a brief summary in the following lines with more details in Appendix~\ref{a-WC_FO}.

The performance estimation approach was shown to provide \emph{tight} certificates, from which one can recover both worst-case certificates and matching worst-case problem instances in~\citep{taylor2017smooth,taylor2017exact}. One consequence is that worst-case guarantees for first-order methods such as those detailed in this section can \emph{always} be obtained as a weighted sum of the appropriate inequalities characterizing the problem at hand; see, for instance,~\citep{de2017worst,dragomir2019optimal}. A similar approach framed in control theoretic terms, and originally tailored to obtain geometric convergence rates, was developed by~\citet{lessard2016analysis} and can also be used to form potential functions~\citep{hu2017dissipativity} as well as optimized methods such as the triple momentum method~\citep{van2017fastest,cyrus2018robust,lessard2020direct}.

The proofs in this section were obtained by using the performance estimation approach tailored for potential functions~\citep{taylor19bach} together with the performance estimation toolbox~\citep{pesto2017}. In particular, the potential function for the optimized gradient method can be found in~\citep[Theorem 11]{taylor19bach} (see also~\citep{drori2021optimal,park2021factor}). These techniques can be used for to either validate or rediscover the proofs in this section numerically, through semidefinite programming. More details are provided in Appendix~\ref{a-WC_FO}.

For the purpose of reproducibility, we provide the corresponding code, as well as notebooks for numerically and symbolically verifying the algebraic reformulations in this section at \togglecodeurl.

\chapter{Proximal Acceleration and Catalysts}\label{c-prox}

In this section, we present simple methods based on approximate proximal operations that produce accelerated gradient-based methods. This idea is exploited for example in the Catalyst~\citep{lin2015universal,lin2017catalyst} and Accelerated Hybrid Proximal Extragradient (A-HPE)~\citep{monteiro2013accelerated} frameworks. In essence, the idea is to develop (conceptual) accelerated proximal point algorithms and to use classical iterative methods to approximate the proximal point. In particular, these frameworks produce accelerated gradient methods (in the same sense as Nesterov's acceleration) when the approximate proximal points are computed using linearly converging gradient-based optimization methods.

\section{Introduction}
We review acceleration from the perspective of proximal point algorithms (PPA). The key concept here, called \emph{proximal operation}, dates back to the 1960s, with the works of Moreau~(\citeyear{moreau1962proximite,moreau1965proximite}).  Its introduction to optimization is attributed to Martinet~(\citeyear{martinet1970breve,martinet1972det}) and was primarily motivated by its link with augmented Lagrangian techniques. In contrast with previous sections, where information about the functions to be minimized was obtained through their gradients, the following sections deal with the case in which information is gathered through a \emph{proximal operator} or an approximation of that operator.

The proximal point algorithm and its use in the development of optimization schemes are surveyed in~\citep{parikh2014proximal}. We aim to go in a slightly different direction here and describe the use of the PPA in an outer loop to obtain improved convergence guarantees in the spirit of the Accelerated Hybrid Proximal Extragradient (A-HPE) method~\citep{monteiro2013accelerated} and of Catalyst~\citep{lin2015universal,lin2017catalyst}.

In this section, we focus on the problem of solving 
\begin{equation}\label{eq:opt_prob_ppa}
    f_\star=\min_{x\in\mathbb{R}^d} f(x),
\end{equation}
where $f$ is closed, convex, and proper (it has a closed convex non-empty epigraph), which we denote by $f\in\mathcal{F}_{0,\infty}$ in line with Definition~\ref{def:smoothstrconvex} from Section~\ref{c-Nest}. We denote by $\partial f(x)$ the subdifferential of $f$ at $x\in\mathbb{R}^d$ and by $g_f(x)\in\partial f(x)$ some element of the subdifferential at $x$, irrespective of whether $f$ is continuously differentiable. We aim to find an $\epsilon$-approximate solution $x$ such that $f(x)-f_\star\leq \epsilon$.

It is possible to develop optimized proximal methods in the spirit of the optimized gradient methods presented in Section~\ref{c-Nest}. That is, given a computational budget---in the proximal setting, this consists of a number of iterations and a sequence of step sizes---one can choose the algorithmic parameters to optimize the worst-case performance. The proximal equivalent to the optimized gradient method is G\"uler's second method~\citep[Section 6]{guler1992new} (see the discussions in Section~\ref{s:notes_ref_prox}). We do not spend time on this method here and directly aim for methods designed from simple potential functions, in the same spirit our approach to Nesterov's accelerated gradient methods in Section~\ref{c-Nest}.

\section{Proximal Point Algorithm and Acceleration}
Whereas the base method for minimizing a function using its gradient is gradient descent:
\[ x_{k+1} = x_k - \lambda  g_f(x_k),\]
the base method for minimizing a function using its proximal oracle is the proximal point algorithm:
\BEQ x_{k+1}= \mathrm{prox}_{\lambda f}(x_k), \label{eq:prox_point_algo} \EEQ
where the proximal operator is given by
\[ \mathrm{prox}_{ \lambda f }(x)\defeq \argmin_y\{\Phi(y;x)\defeq \lambda f(y)+\frac{1}{2}\rVert y-x\lVert^2_2\}.\]
The proximal point algorithm has a number of intuitive interpretations, with two of them being particularly convenient for our purposes.
\begin{itemize}
    \item Optimality conditions of the proximal subproblem reveal that a proximal step corresponds to an implicit (sub)gradient method:
    \[ x_{k+1} = x_k - \lambda   g_f(x_{k+1}),\]
    where $ g_f(x_{k+1})\in \partial f(x_{k+1})$.
    \item Using the proximal point algorithm is equivalent to applying gradient descent to the Moreau envelope of $f$, where the Moreau envelope, denoted $F_\lambda$, is provided by 
    \[
        F_\lambda(x)\defeq \min_{y}\{ f(y)+\frac{1}{2\lambda}\rVert y-x\lVert^2_2\}.
    \] 
    The Moreau envelope has the same set of optimal solutions as $f$, while enjoying attractive additional regularity properties (it is $1/\lambda$-smooth and convex; see Definition~\ref{def:smoothstrconvex}). More precisely, its gradient is given by
    \[\nabla F_{\lambda}(x)=\left(x-\mathrm{prox}_{ \lambda f }(x)\right)/\lambda,\]
    (see~\citep{lemarechal1997practical} for more details). This allows us to write \[x_{k+1}=\mathrm{prox}_{ \lambda f }(x_k)=x_k-\lambda \nabla F_\lambda(x_k),\]  and hence to write proximal minimization methods (as well as their inexact and accelerated variants) applied to $f$ as classical gradient methods (and their inexact and accelerated variants) applied to $F_{\lambda}$.
\end{itemize}
In general, proximal operations are expensive, sometimes nearly as expensive as minimizing the function itself. However, there are many cases, especially in the context of composite optimization problems, where one can isolate parts of the objective for which proximal operators actually have analytical solutions; see, e.g.~\citep{chierchia2020proximity} for a list of such examples. 

In the following sections, we start by analyzing such proximal point methods, and then at the end of the section we show how proximal methods can be used in \emph{outer loops}, where proximal subproblems are solved approximately using a classical iterative method (in inner loops). In particular, we describe how this combination produces accelerated numerical schemes. 

\subsection{Convergence Analysis}
Given the links between proximal operations and gradient methods, it is probably not surprising that proximal point methods for convex optimization can be analyzed using potential functions similar to those used for gradient methods.  

However, there is a huge difference between gradient and proximal steps, as the latter can be made arbitrarily ``powerful'' by taking large step sizes. In other words, a single proximal operation can produce an arbitrarily good approximate solution by picking an arbitrarily large step size. This contrasts with gradient descent, where large step sizes make the method diverge. This fact is clarified later by Corollary~\ref{cor:rate_prox}. However, this nice property of proximal operators comes at a cost: we may not be able to efficiently compute the proximal step.

\enlargethispage{\baselineskip}
As emphasized by the next theorem, proximal point methods for solving~\eqref{eq:opt_prob_ppa} can be analyzed by using similar potentials as those of gradient-based methods. We use
\begin{equation}\label{eq:pot_PPA}
    \phi_k\defeq A_k(f(x_{k})-f(x_\star))+\frac{1}{2}\|x_{k}-x_\star\|^2_2
\end{equation}
and show that $\phi_{k+1}\leq \phi_k$. As before, this type of reasoning can be used recursively:
\begin{equation}\label{eq:pot_chained_PPA}
\begin{aligned}
    A_N(f(x_{N})-f(x_\star))\leq \phi_N\leq \phi_{N-1}\leq\hdots\leq \phi_0=&A_0(f(x_{0})-f(x_\star))\\
    &+\frac{1}{2}\|x_{0}-x_\star\|^2_2,
\end{aligned}
\end{equation}
thereby reaching bounds of the type $f(x_N)-f_\star\leq \tfrac{1}{2A_N}\|x_0-x_\star\|^2_2=O(A_N^{-1})$, assuming $A_0=0$. Since the convergence rates are dictated by the growth rate of the scalar sequence $\{A_k\}_k$, the proofs are designed to increase $A_k$ as fast as possible.

\begin{theorem}\label{thm:ppa_base} Let $f \in \mathcal{F}_{0, \infty}$. For any $k\in\mathbb{N}$, $A_k,\lambda_k\geq0$ and any $x_k$, it holds that
\begin{equation*}
\begin{aligned}
A_{k+1}(f(x_{k+1})-f(x_\star))&+\frac{1}{2}\|x_{k+1}-x_\star\|^2_2\\&\leq A_k(f(x_{k})-f(x_\star))+\frac{1}{2}\|x_{k}-x_\star\|^2_2,
\end{aligned}
\end{equation*}
with $x_{k+1}=\mathrm{prox}_{\lambda_k f}(x_k)$ and $A_{k+1}=A_k+\lambda_k$.
\end{theorem}
\begin{proof}
We perform a weighted sum of the following valid inequalities originating from our assumptions.
\begin{itemize}
    \item Convexity between $x_{k+1}$ and $x_\star$ with weight $\lambda_k$:
    \[ f(x_\star)\geq f(x_{k+1})+\langle  g_f(x_{k+1});x_\star-x_{k+1}\rangle,\]
    with some $g_f(x_{k+1})\in\partial f(x_{k+1})$.
    \item Convexity between $x_{k+1}$ and $x_k$ with weight $A_k$:
    \[f(x_k)\geq f(x_{k+1})+\langle  g_f(x_{k+1});x_k-x_{k+1}\rangle,\]
    with the same $g_f(x_{k+1})\in\partial f(x_{k+1})$ as before.
\end{itemize}
By performing a weighted sum of these two inequalities with their respective weights, we obtain the following valid inequality:
\begin{equation*}
    \begin{aligned}
    0\geq &\lambda_k \left[f(x_{k+1})-f(x_\star)+\langle  g_f(x_{k+1});x_\star-x_{k+1}\rangle\right]\\&+A_k \left[f(x_{k+1})-f(x_k)+\langle  g_f(x_{k+1});x_k-x_{k+1}\rangle \right].
    \end{aligned}
\end{equation*}
By matching the expressions term by term and by substituting $x_{k+1}=x_k-\lambda_k  g_f(x_{k+1})$, one can easily check that the previous inequality can be rewritten exactly as
\begin{equation*}
\begin{aligned}
(A_k+\lambda_k)&(f(x_{k+1})-f(x_\star))+\frac12 \|x_{k+1}-x_\star\|^2_2\\
\leq & A_k(f(x_k)-f(x_\star))+\frac12\|x_k-x_\star\|^2_2 - \lambda_k\frac{2A_k+\lambda_k}{2}\| g_f(x_{k+1})\|^2_2.
\end{aligned}
\end{equation*}
By omitting the last term on the right hand-side (which is nonpositive), we reach the desired statement.
\end{proof}

The first proof of the following worst-case guarantee is due to~\citet{guler1991convergence} and directly follows from the previous potential. 

\begin{corollary} \label{cor:rate_prox}
Let $f \in \mathcal{F}_{0, \infty}$, $\{\lambda_i\}_{i\geq 0}$ be a sequence of nonnegative step sizes, and $\{x_i\}_{i\geq 0}$ be the sequence of iterates from the corresponding proximal point algorithm \eqref{eq:prox_point_algo}. For all $k\in\mathbb{N}$, $k\geq 1$, it holds that
\[ f(x_k)-f_\star \leq \frac{\|x_0-x_\star\|^2_2}{2\sum_{i=0}^{k-1}\lambda_{i}}.\]
\end{corollary}
\begin{proof} It follows directly from the potential with the choice $A_0=0$. That is, we use the potential $\phi_k$ defined in~\eqref{eq:pot_PPA} with Theorem~\ref{thm:ppa_base}, which allows using the chaining argument from~\eqref{eq:pot_chained_PPA}. We obtain:
\[f(x_k)-f(x_\star)\leq \frac1{2A_k} \lVert x_0-x_\star\rVert^2_2,\]
where $A_k=\sum_{i=0}^{k-1}\lambda_i$ and the claim directly follows.
\end{proof}

Note again that we can make this bound arbitrarily good by simply increasing the value of the $\lambda_k$. There is no contradiction here because the proximal oracle is massively stronger than the usual gradient step, as previously discussed. However, solving even a single proximal step is usually (nearly) as hard as solving the original optimization problem, so the proximal method, as detailed here, is a purely conceptual algorithm.

Note that the choice of a constant step size $\lambda_k=\lambda$ results in $f(x_N)-f_\star=O(N^{-1})$ convergence, reminiscent of gradient descent. It turns out that as for gradient-based optimization of smooth convex functions, it is possible to improve this result to $O(N^{-2})$ by using information from previous iterations. This idea was proposed by~\citet{guler1992new}. In the case of a constant step size $\lambda_k=\lambda$, one possible way to obtain this improvement is to apply Nesterov's method (or any other accelerated variant) to the Moreau envelope of $f$.  For varying step sizes, the corresponding bound has the form
\[ f(x_k)-f_\star \leq \frac{2\|x_0-x_\star\|^2_2}{\left(\sum_{i=0}^{k-1}\sqrt{\lambda_{i}}\right)^2}.\]
In addition, G\"uler's acceleration is actually robust to computation errors (as described in the next sections), while allowing for varying step size strategies. 

These two key properties allow using G\"uler's acceleration to design improved numerical optimization schemes using
\begin{enumerate}
    \item Approximate proximal steps, for example by approximately solving the proximal subproblems via iterative methods; and
    \item The step size $\lambda_i$'s can be increased from one iteration to the next, allowing for arbitrarily fast convergence rates to be achieved (assuming that the proximal subproblems can be solved efficiently).
\end{enumerate}
It is important to note that the classical lower bounds for gradient-type methods do not apply here, as we use the much stronger, and more expensive, proximal oracle. It is therefore not a surprise that such techniques (i.e., increasing the step sizes $\lambda_k$ from iteration to iteration) might beat the $O(N^{-2})$ bound obtained through Nesterov's acceleration. Such increasing step size rules can be used, for example, when solving the proximal subproblem via Newton's method, as proposed by~\citet{monteiro2013accelerated}.

\section{G\"uler and Monteiro-Svaiter Acceleration}
In this section, we describe an accelerated version of the proximal point algorithm which may involve inexact proximal evaluations. The method detailed below is a simplified version, sufficient for our purposes, of that of~\citet{monteiro2013accelerated}, and we provide a simple convergence proof for it. The method essentially boils down to that of~\citet{guler1992new} when exact proximal evaluations are used (note that the inexact analysis of~\citet{guler1992new} has a few gaps).

Before proceeding, we mention that there exist quite a few natural notions of inexactness for proximal operations. In this section, we focus on approximately satisfying the first-order optimality conditions of the proximal problem
\[ x_{k+1}=\argmin_x \{ \Phi(x;y_k)\defeq f(x)+\frac{1}{2\lambda_k}\|x-y_k\|^2_2\}.\]
In other words, optimality conditions of the proximal subproblem are
\[0= \lambda_k  g_f(x_{k+1})+ x_{k+1}-y_k,\] for some $ g_f(x_{k+1})\in \partial f(x_{k+1})$. In the following lines, we instead tolerate an error $e_k$:
\[ e_k = \lambda_k  g_f(x_{k+1})+x_{k+1}-y_k,\] 
and require $\|e_k\|_2$ to be small enough to guarantee convergence---even starting at an optimal point does not imply staying at it, without proper assumptions on $e_k$. One possibility is to require $\|e_k\|_2$ to be small with respect to the distance between the starting point $y_k$ and the approximate solution to the proximal subproblem $x_{k+1}$.

Formally, we use the following definition for an approximate solution with relative inaccuracy $0\leq \delta \leq 1$:
\BEQ\label{eq:error_crit}
x_{k+1}\approx_{\delta} \mathrm{prox}_{\lambda_k f}(y_k) \Longleftrightarrow \left[\begin{array}{l} \| e_k\|_2\leq \delta\| x_{k+1}-y_k\|_2\\
\text{with } e_k\defeq x_{k+1}-y_k+\lambda_k  g_f(x_{k+1})\\
\text{for some } g_f(x_{k+1})\in  \partial f(x_{k+1})
      \end{array}\right]
\EEQ
Intuitively, this notion allows for tolerance of relatively large errors when the solution of the proximal subproblem is far from $y_k$ (meaning that $y_k$ is also far away from a minimum of $f$), while demanding relatively small errors when  approaching a solution. On the other side, if $y_k$ is an optimal point for $f$, then so is $x_{k+1}$, as shown by the following proposition.

\begin{proposition}
    Let $y_k\in\argmin_x f(x)$. For any $\delta\in[0,1]$ and any $x_{k+1}\approx_{\delta} \mathrm{prox}_{\lambda_k f}(y_k)$, 
    it holds that $x_{k+1}\in\argmin_x f(x)$.
\end{proposition}
\begin{proof} We only consider the case $\delta=1$ since without loss of generality  $x_{k+1}\approx_{1} \mathrm{prox}_{\lambda_k f}(y_k)\Rightarrow x_{k+1}\approx_{\delta} \mathrm{prox}_{\lambda_k f}(y_k)$ for any $\delta\in[0,1]$. By definition of $x_{k+1}$, we have
    \BEA
    \lVert x_{k+1}-y_k+\lambda_k  g_f(x_{k+1})\rVert^2_2 &\leq & \lVert x_{k+1}-y_k\rVert^2_2 \nonumber\\
    \Leftrightarrow \quad 2\lambda_k\langle  g_f(x_{k+1});\,x_{k+1}-y_k\rangle &\leq& -\lambda_k^2\lVert  g_f(x_{k+1})\rVert^2_2, \label{eq:scalar_product_approx_prox_operator}
    \EEA
    for some $g_f(x_{k+1})\in\partial f(x_{k+1})$. (The second inequality follows from base algebraic manipulations of the first one.) In addition, optimality of $y_k$ implies that
    \[
        \langle  g_f(x_{k+1});\,x_{k+1}-y_k\rangle=\langle  g_f(x_{k+1})- g_f(y_k);\,x_{k+1}-y_k\rangle\geq 0,
    \]
    with $g_f(y_k)=0\in\partial f(y_k)$, where the second inequality follows from the convexity of $f$ (see, e.g., Section~\ref{a-inequalities}). Therefore, condition \eqref{eq:scalar_product_approx_prox_operator} can be satisfied only when $ g_f(x_{k+1})=0$, meaning that $x_{k+1}$ is a minimizer of~$f$.
\end{proof}

Now, assuming that it is possible to find an approximate solution to the proximal operator, one can use Algorithm~\ref{alg:inexact-AccProx}, originally from~\citep{monteiro2013accelerated}, to minimize the convex function~$f$. For simplicity, the parameters $A_k, a_k$ in the algorithm are optimized for $\delta=1$; they can be slightly improved by exploiting the case $0\leq\delta<1$.
\begin{algorithm}[!ht]
  \caption{An inexact accelerated proximal point method~\citep{monteiro2013accelerated}}
  \label{alg:inexact-AccProx}
  \begin{algorithmic}[1]
    \REQUIRE
     A convex function $f$ and an initial point $x_0$.
    \STATE \textbf{Initialize} $z_0=x_0$ and $A_0=0$.
    \FOR{$k=0,\ldots$}
      \STATE Pick $a_{k}= \tfrac{\lambda_k+\sqrt{\lambda_k^2+4A_k\lambda_k}}{2}$
      \STATE $A_{k+1} = A_k + a_k.$
      \STATE $y_{k}= \tfrac{A_k}{A_k+a_{k}} x_k + \tfrac{a_{k}}{A_k+a_k} z_k$
      \STATE $x_{k+1}\approx_{\delta} \mathrm{prox}_{\lambda_k f}(y_k)$ \hspace{1cm}(see Eq.~\eqref{eq:error_crit}, for some $\delta\in\,[0,\,1]$)
      \STATE $z_{k+1} = z_k - a_{k}  g_f(x_{k+1})$
    \ENDFOR
    \ENSURE Approximate solution $x_{k+1}$.
  \end{algorithmic}
\end{algorithm}

Perhaps surprisingly, this method can be analyzed with the same potential as the proximal point algorithm, despite the presence of computation errors.

\begin{theorem}\label{thm:pot_monteiro} Let $f \in \mathcal{F}_{0, \infty}$. For any $k\in\mathbb{N}$, $A_k,\lambda_k\geq0$ and any $x_k,z_k\in\mathbb{R}^d$, it holds that
\begin{equation*}
\begin{aligned}
A_{k+1}(f(x_{k+1})-f(x_\star))&+\frac{1}{2}\|z_{k+1}-x_\star\|^2_2\\&\leq A_k(f(x_{k})-f(x_\star))+\frac{1}{2}\|z_{k}-x_\star\|^2_2,
\end{aligned}
\end{equation*}
where $x_{k+1}$ and $z_{k+1}$ are generated by one iteration of Algorithm~\ref{alg:inexact-AccProx}, and $A_{k+1}=A_k+a_k$.
\end{theorem}
\begin{proof}
We perform a weighted sum of the following valid inequalities, which stem from our assumptions.
\begin{itemize}
    \item Convexity between $x_{k+1}$ and $x_\star$ with weight $A_{k+1}-A_k$:
    \[ f(x_\star)\geq f(x_{k+1})+\langle  g_f(x_{k+1});x_\star-x_{k+1}\rangle,\]
    for some $g_f(x_{k+1})\in\partial f(x_{k+1})$, which we also use below.
    \item Convexity between $x_{k+1}$ and $x_k$ with weight $A_k$:
    \[f(x_k)\geq f(x_{k+1})+\langle  g_f(x_{k+1});x_k-x_{k+1}\rangle.\]
    \item Error magnitude with weight $(A_{k+1})/(2\lambda_k)$:
    \[ \|e_k\|^2_2 \leq \|x_{k+1}-y_k\|^2_2,\]
    which is valid for all $\delta\in[0,1]$ in~\eqref{eq:error_crit}.
\end{itemize}
By performing the weighted sum of these three inequalities, we obtain the following valid inequality:
\begin{equation*}
    \begin{aligned}
    0\geq & (A_{k+1}-A_k) \left[f(x_{k+1})-f(x_\star)+\langle  g_f(x_{k+1});x_\star-x_{k+1}\rangle\right]\\&+A_k \left[f(x_{k+1})-f(x_k)+\langle  g_f(x_{k+1});x_k-x_{k+1}\rangle \right]\\&+\frac{A_{k+1}}{2\lambda_k}\left[ \|e_k\|^2_2 -\|x_{k+1}-y_k\|^2_2\right].
    \end{aligned}
\end{equation*}
After substituting $A_{k+1}=A_k+a_k$, $x_{k+1}=A_k/(A_k+a_k) x_k+a_k/(A_k+a_k) z_k-\lambda_k  g_f(x_{k+1})+e_k$ and $z_{k+1}=z_k-a_k g_f(x_{k+1})$, one can easily check that the previous inequality can be rewritten as
\begin{equation*}
\begin{aligned}
(A_k&+a_k)(f(x_{k+1})-f(x_\star))+\frac12 \|z_{k+1}-x_\star\|^2_2\\
\leq & A_k(f(x_k)-f(x_\star))+\frac12\|z_k-x_\star\|^2_2 \\
& - \frac{\lambda_k(a_k+A_k)-a_k^2}{2}     \| g_f(x_{k+1})\|^2_2,
\end{aligned}
\end{equation*}
either by comparing the expressions on a term-by-term basis or by using an appropriate ``complete the squares'' strategy. We obtain the desired statement by enforcing $\lambda_k(a_k+A_k)-a_k^2\geq 0$, which allows us to neglect the last term on the right-hand side (which is then nonpositive).
Finally, since we have already assumed $a_k\geq 0$, requiring $\lambda_k(a_k+A_k)-a_k^2\geq 0$ corresponds to  
\[0\leq a_{k}\leq \frac{\lambda_k+\sqrt{\lambda_k^2+4A_k\lambda_k}}{2},\]
which yields the desired result.
\end{proof}

The convergence speed then follows from the same reasoning as for the proximal point method.

\begin{corollary}\label{corr:AHPE_bound}
Let $f \in \mathcal{F}_{0, \infty}$, $\{\lambda_i\}_{i\geq 0}$ be a sequence of nonnegative step sizes, and $\{x_i\}_{i\geq 0}$ be the corresponding sequence of iterates from Algorithm~\ref{alg:inexact-AccProx}. For all $k\in\mathbb{N}$, $k\geq 1$, it holds that
\[ f(x_k)-f_\star \leq \frac{2\|x_0-x_\star\|^2_2}{\left(\sum_{i=0}^{k-1}\sqrt{\lambda_{i}}\right)^2}.\]
\end{corollary}
\begin{proof}
Using the potential from Theorem~\ref{thm:pot_monteiro}:
\[\phi_k\defeq A_k(f(x_{k})-f(x_\star))+\frac{1}{2}\|z_{k}-x_\star\|^2_2\]
with $A_0=0$ as well as the chaining argument used in Corollary~\ref{cor:rate_prox}, we obtain \[ f(x_k)-f_\star \leq \frac{\|x_0-x_\star\|^2_2}{2A_k}.\]
The desired result then follows from
\[ A_{k+1}=A_k+a_k=A_k+\frac{\lambda_k+\sqrt{\lambda_k^2+4A_k\lambda_k}}{2}\geq A_k + \frac{\lambda_k}{2}+\sqrt{A_k\lambda_k}.\]
Hence, $A_{k}\geq \left(\sqrt{A_{k-1}}+\tfrac12\sqrt{\lambda_{k-1}}\right)^2\geq\tfrac14\left(\sum_{i=0}^{k-1}\sqrt{\lambda_{i}}\right)^2$.
\end{proof}

\section{Exploiting Strong Convexity}
In this section, we provide refined convergence results when the function to be minimized is $\mu$-strongly convex (all results from previous sections can be recovered by setting $\mu=0$). The algebra is slightly more technical, but the message and techniques are the same. While the proofs in the previous section can be seen as particular cases of the proofs presented below, we detail both versions separately to alleviate the algebraic barrier as much as possible.

\subsubsection{Proximal point algorithm under strong convexity}
We begin by refining the results on the proximal point algorithm. The same modification to the potential function is used to incorporate acceleration in the sequel.
\begin{theorem}\label{th:lyap-strong}
Let $f$ be a closed, $\mu$-strongly convex and proper function.  For any $k\in\mathbb{N}$, $A_k,\lambda_k\geq0$, any $x_k$, and $A_{k+1}= A_k(1+\lambda_k\mu)+\lambda_k$, it holds that
\begin{equation*}
\begin{aligned}
A_{k+1}(f(x_{k+1})-f(x_\star))&+\frac{1+\mu A_{k+1}}{2}\|x_{k+1}-x_\star\|^2_2\\&\leq A_k(f(x_{k})-f(x_\star))+\frac{1+\mu A_{k}}{2}\|x_{k}-x_\star\|^2_2,
\end{aligned}
\end{equation*}
with $x_{k+1}=\mathrm{prox}_{\lambda_k f}(x_k)$.
\end{theorem}
\begin{proof}
We perform a weighted sum of the following valid inequalities, which originate from our assumptions.
\begin{itemize}
    \item Strong convexity between $x_{k+1}$ and $x_\star$ with weight $A_{k+1}-A_k$:
    \[ f(x_\star)\geq f(x_{k+1})+\langle  g_f(x_{k+1});x_\star-x_{k+1}\rangle+\frac\mu2 \|x_\star-x_{k+1}\|^2_2,\]
    with some $g_f(x_{k+1})\in\partial f(x_{k+1})$ which we further use below.
    \item Convexity between $x_{k+1}$ and $x_k$ with weight $A_k$:
    \[f(x_k)\geq f(x_{k+1})+\langle  g_f(x_{k+1});x_k-x_{k+1}\rangle.\]
\end{itemize}
By performing the weighted sum of these two inequalities, we obtain the following valid inequality:
\begin{equation*}
    \begin{aligned}
    0\geq & (A_{k+1}-A_k) \bigg[f(x_{k+1})-f(x_\star)+\langle  g_f(x_{k+1});x_\star-x_{k+1}\rangle\\
    & \quad +\frac\mu2 \|x_\star-x_{k+1}\|^2_2\bigg]\\
    &+A_k \left[f(x_{k+1})-f(x_k)+\langle  g_f(x_{k+1});x_k-x_{k+1}\rangle \right].
    \end{aligned}
\end{equation*}
By matching the expressions term by term and after substituting $x_{k+1}=x_k-\lambda_k  g_f(x_{k+1})$ and $A_{k+1}=A_k(1+\lambda_k\mu)+\lambda_k$, one can check that the previous inequality can be rewritten as
\BEAS
&&A_{k+1} (f(x_{k+1})-f(x_\star))+\frac{1+\mu A_{k+1}}{2} \|x_{k+1}-x_\star\|^2_2\\
&\leq& A_k(f(x_k)-f(x_\star))+\frac{1+\mu A_{k}}{2}\|x_k-x_\star\|^2_2 \\ 
&&\quad - \lambda_k\frac{  A_k (2+\lambda_k  \mu )+\lambda_k }2\| g_f(x_{k+1})\|^2_2.
\EEAS
By neglecting the last term on the right-hand side (which is nonpositive), we reach the desired statement.
\end{proof}

To obtain the convergence speed guaranteed by the previous potential, we have to characterize the growth rate of $A_k$ again, observing that 
\[ A_{k+1}\geq A_k (1+\lambda_k\mu)= \frac{A_k}{1-\tfrac{\lambda_k\mu}{1+\lambda_k\mu}}.\]
The following corollary contains our final worst-case guarantee for the proximal point algorithm, which can be converted to its iteration complexity (details below).
\begin{corollary}\label{corr:PPA_strcvx}
Let $f$ be a closed, $\mu$-strongly convex, and proper function with $\mu\geq 0$, $\{\lambda_i\}_{i\geq 0}$ be a sequence of nonnegative step sizes, and $\{x_i\}_{i\geq 0}$ be the corresponding sequence of iterates from the proximal point algorithm. For all $k\in\mathbb{N}$, $k\geq 1$, it holds that
\[ f(x_k)-f_\star \leq \frac{\mu\|x_0-x_\star\|^2_2}{2(\Pi_{i=0}^{k-1}(1+\lambda_i\mu)-1)}.\]
\end{corollary}
\begin{proof}
Note that the recurrence for $A_k$ provided in Theorem~\ref{th:lyap-strong} has a simple solution
$A_k=([\Pi_{i=0}^{k-1}(1+\lambda_i\mu)]-1)/{\mu }$. By combining this with
\[ f(x_k)-f_\star\leq \frac{\|x_0-x_\star\|^2_2}{2A_k},\]
as provided by Theorem~\ref{th:lyap-strong}, and $A_0=0$, we reach the desired statement.
\end{proof}

As a particular case, note that we recover the case $\mu = 0$ from the previous corollary since $A_k\rightarrow\sum_{i=0}^{k-1}\lambda_i$ when $\mu$ goes to zero.

For arriving to the iteration complexity for obtaining an approximate solution $x_k$ satisfying $f(x_k)-f_\star\leq \epsilon$ with constant step sizes $\lambda_i=\lambda$, we use the following sufficient condition due to Corollary~\ref{corr:PPA_strcvx}:
\[\frac{\mu\|x_0-x_\star\|^2_2}{2(1+\lambda\mu)^k}=\left(1-\frac{\lambda\mu}{1+\lambda\mu}\right)^k\frac{\mu\|x_0-x_\star\|^2_2}{2} \leq \epsilon \,\Rightarrow\, f(x_k)-f_\star\leq \epsilon.\]
A few algebraic manipulations and taking logarithms allows obtaining the following equivalent sufficient condition
\[ k\geq \frac{\log\left(\frac{2\epsilon}{\mu\|x_0-x_\star\|^2_2}\right)}{\log\left(1-\frac{\lambda\mu}{1+\lambda\mu}\right)} \,\Rightarrow\, f(x_k)-f_\star\leq \epsilon.\]
Finally, using the bound $\log\left(1-\frac1{x}\right)\leq -\frac1x$ (for all $x\in(1,\infty)$) we arrive to
\[ k\geq \frac{1+\lambda\mu}{\lambda\mu}\log\left(\frac{\mu\|x_0-x_\star\|^2_2}{2\epsilon}\right) \,\Rightarrow\, f(x_k)-f_\star\leq \epsilon.\] We conclude that the accuracy $\epsilon$ is therefore achieved in $O\left(\tfrac{1+\lambda\mu}{\lambda\mu}\log\tfrac1\epsilon\right)$ iterations of the proximal point algorithm when the step size $\lambda_k=\lambda$ is kept constant. This contrasts with $O\left(\tfrac{1}{\lambda\epsilon}\right)$ in the non-strongly convex case.
\subsubsection{Proximal acceleration and inexactness under strong convexity}
To accelerate convergence while exploiting strong convexity, we upgrade Algorithm~\ref{alg:inexact-AccProx} to Algorithm~\ref{alg:inexact-AccProx-strconvex}, whose analysis follows the same lines as before. For simplicity, the algorithm is optimized for $\delta=\sqrt{1+\lambda_k\mu}$; it can be slightly improved by exploiting the case $0\leq\delta<\sqrt{1+\lambda_k\mu}$. This method is a simplified version of the A-HPE method of~\citep[Algorithm 5.1]{barre2021note}.
\begin{algorithm}[!ht]
  \caption{An inexact accelerated proximal point method}
  \label{alg:inexact-AccProx-strconvex}
  \begin{algorithmic}[1]
    \REQUIRE
      A ($\mu$-strongly) convex function $f$ and an initial point $x_0$.
    \STATE \textbf{Initialize} $z_0=x_0$ and $A_0=0$.
    \FOR{$k=0,\ldots$}
      \STATE Pick $A_{k+1}=A_k+ \tfrac {\lambda_k+2 A_k \lambda_k  \mu+\sqrt{4 A_k^2 \lambda_k  \mu  (\lambda_k  \mu +1)+4 A_k \lambda_k  (\lambda_k  \mu +1)+\lambda_k ^2}}{2}$
      \STATE $y_{k}=  x_k + \tfrac{(A_{k+1}-A_k) (A_k \mu +1)}{A_{k+1}+2 \mu A_k A_{k+1}-\mu  A_k^2} (z_k-x_k)$
      \STATE $x_{k+1}\approx_{\delta} \mathrm{prox}_{\lambda_k f}(y_k)$\hspace{.1cm}(see Eq.\eqref{eq:error_crit}, for some $\delta\in\, [0,\,\sqrt{1+\lambda_k\mu}]$)
      \STATE $z_{k+1} = z_k +\mu \tfrac{ A_{k+1}-A_k}{1+\mu A_{k+1}} (x_{k+1}-z_k)- \tfrac{A_{k+1}-A_k}{1+\mu A_{k+1}}  g_f(x_{k+1})$
    \ENDFOR
    \ENSURE Approximate solution $x_{k+1}$.
  \end{algorithmic}
\end{algorithm}

\begin{theorem} Let $f$ be a closed, $\mu$-strongly convex, and proper function. For any $k\in\mathbb{N}$ and $A_k,\lambda_k\geq0$, the iterates of Algorithm~\ref{alg:inexact-AccProx-strconvex} satisfy
\begin{equation*}
\begin{aligned}
A_{k+1}(f(x_{k+1})-f(x_\star))&+\frac{1+\mu A_{k+1}}{2}\|z_{k+1}-x_\star\|^2_2\\&\leq A_k(f(x_{k})-f(x_\star))+\frac{1+\mu A_{k}}{2}\|z_{k}-x_\star\|^2_2.
\end{aligned}
\end{equation*}
\end{theorem}
\begin{proof}
We perform a weighted sum of the following valid inequalities, which originate from our assumptions.
\begin{itemize}
    \item Strong convexity between $x_{k+1}$ and $x_\star$ with weight $A_{k+1}-A_k$:
    \[ f(x_\star)\geq f(x_{k+1})+\langle  g_f(x_{k+1});x_\star-x_{k+1}\rangle+\frac\mu2 \|x_\star-x_{k+1}\|^2_2,\]
    with some $g_f(x_{k+1})\in\partial f(x_{k+1})$, where this particular subgradient is used repetitively below.
    \item Strong convexity between $x_{k+1}$ and $x_k$ with weight $A_k$
    \[f(x_k)\geq f(x_{k+1})+\langle  g_f(x_{k+1});x_k-x_{k+1}\rangle+\frac\mu2\|x_k-x_{k+1}\|^2_2.\]
    \item Error magnitude with weight $\tfrac{A_{k+1}+2\mu A_k A_{k+1}-\mu A_k^2}{2\lambda_k (1+\mu  A_{k+1} )}$:
    \[ \|e_k\|^2_2\leq (1+\lambda_k\mu)\|x_{k+1}-y_k\|^2_2,\]
    which is valid for all $\delta\in[0,\,\sqrt{1+\lambda_k\mu}]$ in~\eqref{eq:error_crit}.
\end{itemize}
By performing a weighted sum of these three inequalities, with their respective weights, we obtain the following valid inequality:
\begin{equation*}
    \begin{aligned}
    0\geq &(A_{k+1}-A_k) \bigg[f(x_{k+1})-f(x_\star)+\langle  g_f(x_{k+1});x_\star-x_{k+1}\rangle\\
    & \quad\quad\quad\quad\quad\quad +\frac\mu2 \|x_\star-x_{k+1}\|^2_2\bigg]\\
    &+A_k \left[f(x_{k+1})-f(x_k)+\langle  g_f(x_{k+1});x_k-x_{k+1}\rangle+\frac\mu2\|x_k-x_{k+1}\|^2_2 \right]\\
    &+\frac{A_{k+1}+2\mu A_k A_{k+1}-\mu A_k^2}{2\lambda_k (1+\mu  A_{k+1} )} [\|e_k\|^2_2-(1+\lambda_k\mu)\|x_{k+1}-y_k\|^2_2].
    \end{aligned}
\end{equation*}
By matching the expressions term by term and by substituting the expressions for $y_k$, $x_{k+1}=y_k-\lambda_k  g_f(x_{k+1})+e_k$, and $z_{k+1}$, one can check that the previous inequality can be rewritten as (we advise against substituting $A_{k+1}$ at this stage):
\BEAS
&&A_{k+1}(f(x_{k+1})-f(x_\star))+\frac{1+\mu A_{k+1}}{2}\|z_{k+1}-x_\star\|^2_2\\
&\leq & A_k(f(x_{k})-f(x_\star))+\frac{1+\mu A_{k}}{2}\|z_{k}-x_\star\|^2_2\\
&&-\frac{(A_{k+1}+2\mu A_k A_{k+1}-\mu A_k^2)\lambda_k-(A_{k+1}-A_k)^2}{1  +\mu A_{k+1}} \, \frac12\| g_f(x_{k+1})\|^2_2\\
&&-\frac{A_k  (A_{k+1}-A_k) (1+\mu A_k)}{A_{k+1}+2\mu A_k A_{k+1}-\mu A_{k}^2}\, \frac{\mu}{2}\|x_k-z_k\|^2_2.
\EEAS
The conclusion follows from $A_{k+1}\geq A_k$ which allows us to discard the last term (which is then nonpositive). Positivity of the first residual term can be enforced by choosing $A_{k+1}$ such that
\[ (A_{k+1}+2\mu A_k A_{k+1}-\mu A_k^2)\lambda_k-(A_{k+1}-A_k)^2\geq 0.\]
The desired result is achieved by specifically choosing the largest root of the second-order polynomial in $A_{k+1}$, such that $A_{k+1}\geq A_k$.
\end{proof}

In contrast with the previous proximal point algorithm, this accelerated version requires 
\[
O\left(\sqrt{\tfrac{1+\lambda\mu}{\lambda\mu}}\log\tfrac1\epsilon\right)
\]
inexact proximal iterations to reach $f(x_k)-f(x_\star)\leq \epsilon$ when using a constant step size $\lambda_k=\lambda$. This follows from characterizing the growth rate of the sequence~$\{A_k\}_k$: 
\begin{equation}\label{eq:ak_acc_ppa}
A_{k+1}\geq A_k (1+\lambda_k \mu)+A_k\sqrt{\lambda_k\mu(1+\lambda_k\mu)}=\frac{A_k}{1-\sqrt{\frac{\lambda_k\mu}{1+\lambda_k\mu}}}.
\end{equation}
\begin{corollary}\label{corr:AHPE_bound_strcvx}
Let $f \in \mathcal{F}_{\mu, \infty}$ with $\mu \geq 0$, $\{\lambda_i\}_{i\geq 0}$ be a sequence of nonnegative step sizes, and $\{x_i\}_{i\geq 0}$ be the corresponding sequence of iterates from Algorithm~\ref{alg:inexact-AccProx-strconvex}. For all $k\in\mathbb{N}$, $k\geq 1$, it holds that
\[ f(x_k)-f_\star \leq \Pi_{i=1}^{k-1}\left(1-\sqrt{\tfrac{\lambda_i\mu}{1+\lambda_i\mu}} \right) \frac{\|x_0-x_\star\|^2_2}{2\lambda_0}.\]
\end{corollary}
\begin{proof}
The proof follows from the same arguments as before; that is,
\[ f(x_k)-f_\star\leq \frac{\|x_0-x_\star\|^2_2}{2A_k},\] and \[A_k\geq \frac{\lambda_0}{\Pi_{i=1}^{k-1}\left(1-\sqrt{\tfrac{\lambda_i\mu}{1+\lambda_i\mu}} \right)},\] where we used $A_0=0$ and $A_1=\lambda_0$. We then proceed with~\eqref{eq:ak_acc_ppa}.
\end{proof}

Before continuing to the next section, we note that combining Corollary~\ref{corr:AHPE_bound} with Corollary~\ref{corr:AHPE_bound_strcvx} shows that
\[ f(x_k)-f_\star \leq \min\left\{ \frac{\Pi_{i=1}^{k-1}\left(1-\sqrt{\tfrac{\lambda_i\mu}{1+\lambda_i\mu}} \right)}{\lambda_0},\frac{4}{\left(\sum_{i=0}^{k-1}\sqrt{\lambda_i}\right)^2}\right\}\frac{\|x_0-x_\star\|^2_2}{2}.\]

\section{Application: Catalyst Acceleration}

In what follows, we illustrate how to use proximal methods as meta-algorithms to improve the convergence of simple gradient-based first-order methods. The idea consists of using a base first-order method, such as gradient descent, to obtain approximations to the proximal point subproblems, within an accelerated proximal point method. This idea can be extended by embedding any algorithm that can solve the proximal subproblem.

There exist many notions of \emph{approximate solutions} to the proximal subproblems, giving rise to different types of guarantees together with slightly different methods. In particular, we required the approximate solution to have a small gradient. Other notions of approximate solutions are used, among others, in~\citep{guler1992new,schmidt2011convergence,villa2013accelerated}. Depending on the target application or on the target algorithm for solving inner problems, the \emph{natural} notion of an approximate solution to the proximal subproblem might change. A fairly general framework was developed by~\citet{monteiro2013accelerated} (where the error is controlled via a primal-dual gap on the proximal subproblem).

\subsection{Catalyst acceleration}\label{s:catalyst}
A popular application of the inexact accelerated proximal gradient is Catalyst acceleration \citep{lin2015universal}. For readability purposes, we do not present the general Catalyst framework but rather a simple instance. Stochastic versions of this acceleration procedure have also been developed, and we briefly summarize them in Section~\ref{s:stoch_cata}. The idea is again to use a base first-order method to approximate the proximal subproblem up to the required accuracy. For now, we assume that we want to minimize an $L$-smooth convex function~$f$, i.e.,
\[
    \min_{x\in\mathbb{R}^d} f(x),
\]
The corresponding proximal subproblem has the form
\BEQ
    \text{prox}_{\lambda f}(y) \defeq \argmin_x \left(f(x) + \frac{1}{2\lambda}\|x-y\|^2_2\right), \label{eq:prox_subproblem}
\EEQ
and it is therefore the minimization of an $(L+1/\lambda)$-smooth and $1/\lambda$-strongly convex function.  To solve such a problem, one can use a first-order method to approximate its solution. 

\subsubsection{Preliminaries}

In what follows, we consider using a method $\mathcal{M}$ to solve the proximal subproblem~\eqref{eq:prox_subproblem}. We assume that this method is guaranteed to converge linearly on any smooth strongly convex problem with minimizer $w_\star$, and more precisely that:
\begin{equation}
    \label{eq:lin_conv_req}
    \|w_k-w_\star\|_2\leq C_{\mathcal{M}}(1-\tau_{\mathcal{M}})^k \|w_0-w_\star\|_2
\end{equation}
(where $\{w_i\}_i$ are the iterates of $\mathcal{M}$) for some constant $C_{\mathcal{M}}\geq 0$ and some $0<\tau_{\mathcal{M}}\leq 1$. Note that we consider linear convergence in terms of $\|w_k-w_\star\|_2$ for convenience; other notions can be used, such as convergence in function values.

We distinguish the sequences $\{x_k\}_k$, $\{y_k\}_k$, and $\{z_k\}_k$, which are the iterates of the inexact accelerated proximal point algorithm (Algorithm~\ref{alg:inexact-AccProx}, or~\ref{alg:inexact-AccProx-strconvex}), and the sequence of iterates $\{w^{(k)}_i\}_i$, which are the iterates of $\mathcal{M}$, used to approximate $\text{prox}_{\lambda f}(y_k)$ in step 6 of Algorithm \ref{alg:inexact-AccProx} (or step 5 of Algorithm~\ref{alg:inexact-AccProx-strconvex}). We also use the \emph{warm-start strategy} $w_0^{(k)}=y_k$.

We can thus apply Algorithm~\ref{alg:inexact-AccProx} to minimize $f$ while (approximately) solving the proximal subproblems with ${\mathcal{M}}$. We first define four iteration counters:
\begin{enumerate}
    \item $N_{\text{outer}}$, the number of iterations of the inexact accelerated proximal point method (Algorithm~\ref{alg:inexact-AccProx}, or~\ref{alg:inexact-AccProx-strconvex}), which serves as the ``outer loop'' for the overall acceleration scheme. That is, the output of the overall method is $x_{N_{\text{outer}}}$ in the notation of Algorithm \ref{alg:inexact-AccProx} (or Algorithm~\ref{alg:inexact-AccProx-strconvex});
    \item $N_{\text{inner}}(k)$, the number of iterations needed by method $\mathcal{M}$ to approximately solve the proximal subproblem at iteration $k$ of the outer loop of the inexact proximal point method. That is, the number of iterations of $\mathcal{M}$ for approximating $\mathrm{prox}_{\lambda f}(y_k)$ to the target accuracy, when the initial iterate of $\mathcal{M}$ is set to $w_0^{(k)}=y_k$;
    \item $N_{\text{useless}}$, the number of iterations performed by $\mathcal{M}$ that did not result in an additional iteration of the inexact accelerated proximal point method. That is, if the user has a limited budget in terms of a total number of iterations for $\mathcal{M}$, it is likely that the last few iterations of $\mathcal{M}$ do not allow completing an iteration of the ``outer loop''. Thus, $N_{\text{useless}}<N_{\text{inner}}(N_{\text{outer}})$, i.e., the number of useless iterations of $\mathcal{M}$ is smaller than the number of iterations that would have lead to an additional outer iteration.
    \item $N_{\text{total}}$, the total number of iterations of method $\mathcal{M}$: \[N_{\text{total}}=N_{\text{useless}}+\sum_{k=0}^{N_{\text{outer}}-1}N_{\text{inner}}(k).\] Again, $N_{\text{useless}}$ is the number of iterations of $\mathcal{M}$ that did not allow an additional outer iteration to complete.
\end{enumerate}

\subsubsection{Overall complexity}

As we detail in the sequel, assuming that $\mathcal{M}$ satisfies~\eqref{eq:lin_conv_req}, the overall complexity of the combination of methods is guaranteed to be 
\[
f(x_{N_{\text{outer}}})-f_\star=O(N_{\text{total}}^{-2}),
\]
where $x_{N_{\text{outer}}}$ is the iterate produced after $N_{\text{outer}}$ iterations of the inexact accelerated proximal point method or equivalently, the iterate produced after a total number of iterations $N_{\text{total}}$  of method $\mathcal{M}$. More precisely, $x_{N_{\text{outer}}}$ is guaranteed to satisfy
\[ f(x_{N_{\text{outer}}})-f(x_\star)\leq \frac{2\|x_0-x_\star\|^2_2}{\lambda N_{\text{outer}}^2} \leq \frac{2\|x_0-x_\star\|^2_2}{\lambda \lfloor B_{\mathcal{M},\lambda}^{-1} N_{\text{total}}\rfloor^2},\]
(the first inequality follows from Corollary~\ref{corr:AHPE_bound} and the second one from the analysis below) where we used $N_{\text{inner}}(k)\leq B_{\mathcal{M},\lambda}$ for all $k\geq 0$ and hence $\lfloor\frac{N_{\text{total}}}{B_{\mathcal{M},\lambda}} \rfloor\leq  N_{\text{outer}}$, where the constant $B_{\mathcal{M},\lambda}$ depends solely on the choice of $\lambda$ and on properties of $\mathcal{M}$. This $B_{\mathcal{M},\lambda}$ represents the computational burden of approximately solving one proximal subproblem with $\mathcal{M}$, and it satisfies
\[
    B_{\mathcal{M},\lambda}\leq \frac{\log ( C_{\mathcal{M}}(\lambda L+2))}{\tau_{\mathcal{M}}}+1.
\]
We provide a few simple examples based on gradient methods for smooth strongly convex minimization. For all these methods, the embedding within the inexact proximal framework yields 
\BEQ\label{eq:cata-rate}
    N_{\text{total}} = O\left(B_{\mathcal{M},\lambda}\sqrt{\tfrac{L\|x_0-x_\star\|^2_2}{\epsilon}}\right)
\EEQ 
iteration complexity in terms of the total number of calls to $\mathcal{M}$ to find a point that satisfies $f(x_{N_{\text{outer}}})-f(x_\star)\leq \epsilon$. We can make this bound a bit more explicit depending on the choice of $\mathcal{M}$.

\begin{itemize}
    \item Let $\mathcal{M}$ be a regular gradient method with step size $1/(L+1/\lambda)$ that we use to solve the proximal subproblem. The method is known to converge linearly with $C_{\mathcal{M}}=1$ and $\tau_{\mathcal{M}}=\tfrac{1}{1+\lambda L}$ (the inverse condition ratio for the proximal subproblem), and it produces the accelerated rate in~\eqref{eq:cata-rate}. Note that directly applying the gradient method to the problem of minimizing $f$ yields a much worse iteration complexity: $O\left(\tfrac{L\|x_0-x_\star\|^2_2}{\epsilon}\right)$.
    \item Let $\mathcal{M}$ be a gradient method including an exact line-search. It is guaranteed to converge linearly with $C_{\mathcal{M}}=\lambda L + 1$ (the condition ratio of the proximal subproblem) and $\tau_{\mathcal{M}}=\tfrac{2}{2+\lambda L}$. The iteration complexity of applying this steepest descent scheme directly to $f$ is similar to that of vanilla gradient descent. One can also choose $\lambda$ to avoid having an excessively large $B_{\mathcal{M},\lambda}$; for example, $\lambda=1/L$.
    \item Let $\mathcal{M}$ be an accelerated gradient method specifically tailored for smooth strongly convex optimization, such as Nesterov's method with constant momentum; see Algorithm~\ref{alg:FGM_2_strcvx_constmomentum}. It is guaranteed to converge linearly with $C_{\mathcal{M}}=\lambda L + 1$ and $\tau_{\mathcal{M}}=\sqrt{\tfrac{1}{1+\lambda L}}$. Although there is no working guarantee for this method on the original minimization problem, if $f$ is not strongly convex, it can still be used to minimize $f$ through the inexact proximal point framework, as proximal subproblems are strongly convex.
\end{itemize}
To conclude, inexact accelerated proximal schemes produce accelerated rates for vanilla optimization methods that converge linearly for smooth strongly convex minimization. The idea of embedding a simple first-order method within an inexact accelerated scheme can be applied to a large array of settings, including to obtain acceleration in strongly convex problems or for stochastic minimization (see below). However, one should note that practical tuning of the corresponding numerical schemes (and particularly of the step size parameters) critically affects the overall performance, as discussed in, e.g.,~\citep{lin2017catalyst}. This makes effective implementation somewhat tricky. The analysis of non-convex settings is beyond the scope of this section, but examples of such results can be found in, e.g.,~\citep{paquette2018catalyst}.

\subsection{Detailed Complexity Analysis}
Recall that function value accuracies, e.g., in Corollary~\ref{corr:AHPE_bound} are expressed in terms of outer loop iterations. Therefore, to complete the analysis, we need to answer the following question: given a total budget of $N_{\text{total}}$ inner iterations of method $\mathcal{M}$, how many iterations of Algorithm~\ref{alg:inexact-AccProx}, $N_{\text{outer}}$, will we perform in the ideal strategy (in other words, what is  $B_{\mathcal{M},\lambda}$)? To answer this question, we start by analyzing the computational cost of solving a single proximal subproblem through~$\mathcal{M}$.

\paragraph{Computational cost of inner problems.}
Let
\[
    \Phi_k(x)\defeq f(x)+\frac{1}{2\lambda}\|x-y_k\|^2_2
\]
be the objective of the proximal subproblem that we aim to solve at iteration $k$ (line 6 of Algorithm~\ref{alg:inexact-AccProx}) centered at $y_k$. By construction, $\Phi_k(x)$  is $(L+1/\lambda)$-smooth and $1/\lambda$-strongly convex. Also denote by $w_0=y_k$ our (warm-started) initial iterate and by $w_0,w_1,\hdots,w_{N_{\text{inner}}(k)}$ the iterates of $\mathcal{M}$ used to solve $\min_x \Phi_k(x)$ (note that we drop the superscript $(k)$ for readability, avoiding the heavier notation $w_0^{(k)},w_1^{(k)},\hdots,w_{N_{\text{inner}}(k)}^{(k)}$). We also denote $w_\star(\Phi_k)\defeq\argmin_x\, \Phi_k(x)$. 

We need to compute an upper bound on the number of iterations $N_{\text{inner}}(k)$ required to satisfy the error criterion~\eqref{eq:error_crit}:
\begin{equation}\label{eq:cata_req}
    \|e_{N_{\text{inner}}(k)}\|_2= \lambda \|\nabla \Phi_k(w_{N_{\text{inner}}(k)})\|_2 \leq \|w_{N_{\text{inner}}(k)}-w_0\|_2,
\end{equation}
where we denote by $N_{\text{inner}}(k)=\inf\{i\, :\, \|\nabla \Phi_k(w_i)\|_2 \leq 1/\lambda \|w_i-w_0\|_2\}$ the index of the first iteration such that~\eqref{eq:cata_req} is satisfied: this is precisely the quantity we want to upper bound. We start with the following observations:
\begin{itemize}
    \item By $(L+1/\lambda)$-smoothness of $\Phi_k$, we have
\begin{equation}\label{eq:cata_req1}
 \| \nabla \Phi_k(w_i)\|_2 \leq (L+1/\lambda) \|w_i-w_\star(\Phi_k)\|_2,
\end{equation}
where $w_\star(\Phi_k)$ is the minimizer of $\Phi_k$.
\item The triangle inequality applied to $\|w_0-w_\star(\Phi_k)\|_2$ implies 
\begin{equation}\label{eq:cata_req2}
\|w_0-w_\star(\Phi_k)\|_2-\|w_i-w_\star(\Phi_k)\|_2 \leq \|w_0-w_i\|_2.
\end{equation}
\end{itemize}
Hence,~\eqref{eq:cata_req} is satisfied if the right-hand side of~\eqref{eq:cata_req1} is smaller than the left-hand side of~\eqref{eq:cata_req2} divided by $\lambda$. Thus, for any $i$ for which we can prove
\[(L+1/\lambda) \|w_i-w_\star(\Phi_k)\|_2\leq 1/\lambda (\|w_0-w_\star(\Phi_k)\|_2-\|w_\star(\Phi_k)-w_i\|_2), \]
we obtain $N_{\text{inner}}(k)\leq i$. Rephrasing this inequality leads to
\[ \|w_i-w_\star(\Phi_k)\|_2\leq \frac{1}{\lambda L+2}\|w_0-w_\star(\Phi_k)\|_2. \]
Therefore, by assumption on $\mathcal{M}$,~\eqref{eq:cata_req} is guaranteed to hold as soon as
\[ C_{\mathcal{M}}(1-\tau_{\mathcal{M}})^i\leq \frac{1}{\lambda L+2},\]
and thus~\eqref{eq:cata_req} holds for any $i$ that satisfies
\[ i \geq \left\lceil\frac{\log\left(C_{\mathcal{M}}(\lambda L+2)\right)}{ \log\left(1/(1-\tau_{\mathcal{M}})\right) }\right\rceil.\]
We conclude that~\eqref{eq:cata_req} is satisfied before this number of iterations is achieved; hence,
\[ N_{\text{inner}}(k)\leq \left\lceil \frac{\log\left(C_{\mathcal{M}}(\lambda L+2)\right)}{ \log\left(1/(1-\tau_{\mathcal{M}})\right) } \right\rceil\leq \frac{\log\left(C_{\mathcal{M}}(\lambda L+2)\right)}{ \log\left(1/(1-\tau_{\mathcal{M}})\right) } + 1.\]
Given that the right-hand side does not depend on $k$, we use the notation
\[ B_{\mathcal{M},\lambda}\defeq\frac{\log\left(C_{\mathcal{M}}(\lambda L+2)\right)}{ \log\left(1/(1-\tau_{\mathcal{M}})\right) } + 1\]
as our upper bound on the iteration cost of solving the proximal subproblem via $\mathcal{M}$.

\paragraph{Global complexity bound.} We have shown that the number of iterations in the inner loop is bounded above by a constant that depends on the specific choice of the regularization parameter and on the method $\mathcal{M}$. In other words, $N_{\text{inner}}(k)\leq B_{\mathcal{M},\lambda}$. Denoting by $N_{\text{total}}$ the total number of calls to the gradient of $f$, by $N_{\text{outer}}$ the number of iterations performed by Algorithm~\ref{alg:inexact-AccProx}, and by $N_{\text{useless}}$ the number of iterations of $\mathcal{M}$ that did not result in an additional outer iteration (see discussions in Section~\ref{s:catalyst} ``Preliminaries''), we conclude that
\[N_{\text{total}}=N_{\text{useless}}+\sum_{k=0}^{N_{\text{outer}}-1}N_{\text{inner}}(k)< (N_{\text{outer}}+1) B_{\mathcal{M},\lambda}.\]
Hence, $N_{\text{outer}} \geq  \lfloor B_{\mathcal{M},\lambda}^{-1} N_{\text{total}}\rfloor$ since $N_{\text{useless}}<B_{\mathcal{M},\lambda}$ (the number of useless iterations is smaller than the number of iterations that would lead to an additional outer iteration).
The conclusion follows from Corollary~\ref{corr:AHPE_bound}:
\[ f(x_{N_\text{outer}})-f(x_\star)\leq \frac{2\|x_0-x_\star\|^2_2}{\lambda N_{\text{outer}}^2}\leq \frac{2\|x_0-x_\star\|^2_2}{\lambda \lfloor B_{\mathcal{M},\lambda}^{-1} N_{\text{total}}\rfloor^2}.\]
That is, given a target accuracy $\epsilon$, the iteration complexity written in terms of the total number of approximate proximal minimizations in Algorithm~\ref{alg:inexact-AccProx} is $O(\sqrt{\tfrac{\|x_0-x_\star\|^2_2}{\lambda \epsilon}})$, and the total iteration complexity when solving the problem using $\mathcal{M}$ in the inner loops is simply the same bound multiplied by the cost of solving a single proximal subproblem, namely $O(B_{\mathcal{M},\lambda}\sqrt{\tfrac{\|x_0-x_\star\|^2_2}{\lambda \epsilon}})$.

\subsection{Catalyst for Strongly Convex Problems} 
The previous analysis holds for the convex (but not necessarily strongly convex) case. The iteration complexity of solving inner problem remains valid in the strongly convex case, and the expression for $B_{\mathcal{M},\lambda}$ can only be improved slightly---by taking into account the better strong convexity parameter $\mu+1/\lambda$ and the possibly larger acceptable error magnitude with the factor $\sqrt{1+\lambda\mu}$ in Algorithm~\ref{alg:inexact-AccProx-strconvex}. Therefore, the total number of iterations of Algorithm~\ref{alg:inexact-AccProx-strconvex} embedded with $\mathcal{M}$ remains bounded in a similar fashion, and the overall error decreases as $\left(1-\sqrt{\tfrac{\lambda\mu}{1+\lambda\mu}}\right)^{\lfloor N_{\text{outer}}/B_{\mathcal{M},\lambda}\rfloor}$, and the iteration complexity is therefore of order 
\begin{equation}\label{eq:cata_comp_strcvx}
O\left(B_{\mathcal{M},\lambda}\sqrt{\frac{1+\lambda\mu}{\lambda\mu}}\log\frac1\epsilon \right).
\end{equation}
It is thus natural to choose the value of $\lambda$ by optimizing the overall iteration complexity of Algorithm~\ref{alg:inexact-AccProx-strconvex} combined with $\mathcal{M}$. One way to proceed is by optimizing
\[ 
\left.{\sqrt{\frac{1+\lambda  \mu}{\lambda  \mu }}}\middle/{\tau_{\mathcal{M}} }\right.,
\]
essentially neglecting the factor $\log (C_{\mathcal{M}}(\lambda L+2))$ in the complexity estimate~\eqref{eq:cata_comp_strcvx}. Here are a few examples:
\begin{itemize}
    \item Gradient method with suboptimal tuning (e.g., when using backtracking or line-search techniques): $\tau_{\mathcal{M}}=\tfrac{\mu\lambda+1}{L\lambda+1}$. Optimizing the ratio leads to the choice $\lambda=\tfrac{1}{L-2\mu}$, and the ratio is equal to $2\sqrt{\tfrac{L}{\mu}-1}$. Assuming $C_{\mathcal{M}}=1$ (which is the case for the standard step size $1/L$), the overall iteration complexity is then $O\left(\sqrt{\tfrac{L}{\mu}}\log\tfrac1\epsilon \right),$ where we neglected the factor $\log (2\tfrac{1-\mu/L}{1-2\mu/L})\approx \log 2$ when $L/\mu$ is large enough.
    \item Gradient method with optimal tuning: $\tau_{\mathcal{M}}=\tfrac{2(\mu\lambda+1)}{(L\lambda+\mu\lambda+2)}$. The resulting choice is $\lambda=\tfrac{2}{L-3\mu}$ and the ratio is ${\sqrt{2}}\sqrt{\tfrac{L}{\mu }-1}$, thereby arriving at the same $O\left(\sqrt{\tfrac{L}{\mu}}\log\tfrac1\epsilon \right)$.
\end{itemize}

\subsection{Catalyst for Randomized/Stochastic Methods}\label{s:stoch_cata} Similar results hold for stochastic methods, assuming the convergence of $\mathcal{M}$ in expectation instead of~\eqref{eq:lin_conv_req}, such as in the form $\mathbb{E}\|w_k-w_\star\|_2\leq C_{\mathcal{M}}(1-\tau_{\mathcal{M}})^k \|w_0-w_\star\|_2$. Overall, the idea remains the same:
\begin{enumerate}
    \item Use the inexact accelerated proximal point algorithm (Algorithm~\ref{alg:inexact-AccProx} or~\ref{alg:inexact-AccProx-strconvex}) as if $\mathcal{M}$ were deterministic.
    \item Use the stochastic method $\mathcal{M}$ to obtain points that satisfy the accuracy requirement.
\end{enumerate}

Dealing with the computational burden of solving the inner problem is a bit more technical, but the overall analysis remains similar. One can bound the expected number of iterations needed to solve the inner problem $\mathbb{E}[N_{\text{inner}}(k)]$ by some constant $B_{\mathcal{M},\lambda}^{(\text{stoch})}$ of the form (details below)
    \[ B_{\mathcal{M},\lambda}^{(\text{stoch})}\defeq \frac{\log\left(C_{\mathcal{M}}(\lambda L+2)\right)}{ \log\left(1/(1-\tau_{\mathcal{M}})\right) } + 2,\]
    which is simply $B_{\mathcal{M},\lambda}^{(\text{stoch})}=B_{\mathcal{M},\lambda}+1$. A simple argument for obtaining this bound uses Markov's inequality as follows:
\begin{equation*}
\begin{aligned}
    \mathbb{P}(N_{\text{inner}}(k)>i)&\leq \mathbb{P}\left(\|w_i-w_\star(\Phi_k)\|_2>\tfrac{1}{\lambda L+2}\|w_0-w_\star(\Phi_k)\|_2\right)\\ &{\leq}\frac{\mathbb{E}[\|w_i-w_\star(\Phi_k)\|_2]}{\tfrac{1}{\lambda L+2}\|w_0-w_\star(\Phi_k)\|_2} \quad (\text{Markov})\\ &\leq\frac{C(1-\tau)^i\|w_0-w_\star(\Phi_k)\|_2}{\tfrac1{\lambda L+2}\|w_0-w_\star(\Phi_k)\|_2}=\frac{C(1-\tau)^i}{\tfrac1{\lambda L+2}}.
\end{aligned}
\end{equation*}
We then use a refined version of this bound: $\mathbb{P}(N_{\text{inner}}(k)>i)\leq \min\left\{1,(\lambda L+2){C(1-\tau)^i}\right\}$, and in order to bound $\mathbb{E}[N_{\text{inner}}(k)]$, we proceed with
\BEAS
\begin{aligned}
 \mathbb{E}[N_{\text{inner}}(k)]&=\sum_{t=1}^\infty \mathbb{P}(N_{\text{inner}}(k)\geq t)\\
 &\leq \int^{N_0}_0 1 dt+C(\lambda L+2)\int^\infty_{N_0}(1-\tau)^t dt,
\end{aligned}
\EEAS
where $N_0$ is such that $1=C(\lambda L+2)(1-\tau)^{N_0}$. Direct computation yields $\mathbb{E}[N_{\text{inner}}(k)]\leq B_{\mathcal{M},\lambda}^{(\text{stoch})}\defeq N_0+1$.

The overall expected iteration complexity is that of the inexact accelerated proximal point method multiplied by the expected computational burden of solving the proximal subproblems $B_{\mathcal{M},\lambda}^{(\text{stoch})}$. That is, the expected iteration complexity becomes 
\[
O\left(B_{\mathcal{M},\lambda}^{(\text{stoch})}\sqrt{\frac{\|x_0-x_\star\|^2_2}{\lambda \epsilon}}\right)
\]
in the smooth convex setting, and 
\[
O\left(B_{\mathcal{M},\lambda}^{(\text{stoch})}\sqrt{\frac{1+\lambda\mu}{\lambda\mu}}\log\frac1\epsilon \right)
\]
in the smooth strongly convex setting. The main argument of this section, namely the use of Markov's inequality, was adapted from~\citet[Appendix B.4]{lin2017catalyst} (merged with the arguments for the deterministic case above). Stochastic versions of Catalyst acceleration were also studied in~\citep{kulunchakov2019generic}.

\section{Notes and References}\label{s:notes_ref_prox}
In the optimization literature, the proximal operation is an essential algorithmic primitive at the heart of many practical optimization methods. Proximal point algorithms are also largely motivated by the fact that they offer a nice framework for obtaining ``meta'' (or high-level) algorithms. They naturally appear in augmented Lagrangian and splitting-based numerical schemes, among others.  We refer the reader to the excellent surveys in \citep{parikh2014proximal,ryu2016primer} for more details. 

\paragraph{Proximal point algorithms: accelerated and inexact variants.} Proximal point algorithms have a long history, dating back to the works of Moreau~(\citeyear{moreau1962proximite,moreau1965proximite}): they were introduced to the optimization community by Martinet~(\citeyear{martinet1970breve,martinet1972det}). Early interest in proximal methods was motivated by their connection to augmented Lagrangian techniques~\citep{rockafellar1973dual,rockafellar1976augmented,iusem1999augmented}; see also the helpful tutorial by~\citet{eckstein2013practical}). Among the many other successes and uses of proximal operations, one can cite the many \emph{splitting} techniques~\citep{lions1979,eckstein1989splitting}, for which there are sound surveys~\citep{boyd2011distributed,eckstein2012augmented,condat2019proximal}. In this context, inexact proximal operations had already been introduced by~\citet{rockafellar1976augmented} and were combined with acceleration much later by~\citet{guler1992new}---although not with a perfectly rigorous proof, which was later corrected in~\citep{salzo2012inexact,monteiro2013accelerated}.

\paragraph{Hybrid proximal extragradient (HPE) framework.} Whereas Catalyst acceleration is based on the idea of solving the proximal subproblem via a first-order method, the (related) hybrid proximal extragradient framework is also used together with a Newton scheme in~\citep{monteiro2013accelerated}. Furthermore, the accelerated hybrid proximal extragradient framework allows for an increasing sequence of step sizes, thereby leading to faster rates than those obtained via vanilla first-order methods. (That is, using an increasing sequence of $\{\lambda_i\}_i$, $(\sum_{i=1}^N\sqrt{\lambda_i})^2$ might grow much faster than $N^2$.)

The HPE framework was introduced by Solodov and Svaiter (\citeyear{solodov1999hybrid,solodov1999hybrid2,solodov2000error,solodov2001unified}) before it was embedded with acceleration techniques by~\citep{monteiro2013accelerated}.

\paragraph{Catalyst.} The variant presented in this section was chosen for simplicity of exposition; it is largely inspired by recent works on the topic in~\citep{lin2017catalyst,ivanova2019adaptive} along with~\citep{monteiro2013accelerated}. Efficient implementations of Catalyst can be found in the Cyanure package by~\citet{mairal2019cyanure}. In particular, most efficient practical implementations of Catalyst appear to rely on an \emph{absolute} inaccuracy criterion for the inexact proximal operation, instead of on \emph{relative} (or multiplicative) ones, as used in this section. In practice, the most convenient and efficient variants appear to be those that use a constant number of inner loop iterations to approximately solve the proximal subproblems. 

In this section, we chose the relative error model as we believe it allows for a slightly simpler exposition while relying on essentially the same techniques.
Catalyst was originally proposed by~\citet{lin2015universal} as a generic tool for reaching accelerated methods. Among others, it allowed for the acceleration of stochastic methods such as SVRG~\citep{johnson2013accelerating}, SAGA~\citep{defazio2014saga}, MISO~\citep{mairal2015incremental}, and Finito~\citep{defazio2014finito} before direct acceleration techniques had been developed for them~\citep{allen2017katyusha,zhou2018simple,zhou2019direct}.

\paragraph{Higher-order proximal subproblems.} Higher-order proximal subproblems of the form
\BEQ
\min_x\left\{ f(x)+\frac1{\lambda (p+1)}\lVert x-x_k\rVert^{p+1}_2\right\}
\EEQ
were used by~\citep{nesterov2020inexactAcc,nesterov2020inexact} as a new primitive for designing optimization schemes. These subproblems can also be solved approximately (via $p$th-order tensor methods~\citep{nesterov2019implementable}) while maintaining good convergence guarantees.

\paragraph{Optimized proximal point methods.} It is possible to develop optimized proximal methods in the spirit of optimized gradient methods. That is, given a computational budget---in the proximal setting, this consists of a number of iterations and a sequence of step sizes $\{\lambda_i\}_{0\leq i\leq N-1}$---one can choose algorithmic parameters to optimize the worst-case performance of a method of the type
\[ x_{k+1}=x_0-\sum_{i=1}^{k} \beta_i  g_f(x_i)-\lambda_{k} g_f(x_{k+1})\]
with respect to the $\beta_i$. The proximal equivalent of the optimized gradient method is G\"uler's second method~\citep[Section 6]{guler1992new}, which was obtained as an optimized proximal point method in~\citep{barre2020principled}. Alternatively, G\"uler's second method~\citep[Section 6]{guler1992new} can be obtained by applying the optimized gradient method (without its last iteration trick) to the Moreau envelope of the nonsmooth convex function $f$. More precisely, denoting by $\tilde{f}$ the Moreau envelope of $f$, one can apply the optimized gradient method without the last iteration trick to $\tilde{f}$ as $f(x_k)-f_\star=\tilde{f}(x_k)-\tilde{f}_\star-\tfrac1{2L}\|g_{\tilde{f}}(x_k)\|^2_2$, which corresponds precisely to the first term of the potential of the optimized gradient method (see Equation~\ref{eq:pot_ogm_preview}). In the more general setting of monotone inclusions, one can obtain alternate optimized proximal point methods for different criteria as in~\citep{kim2019accelerated,lieder2020convergence}. 

\paragraph{Proofs in this section.} The proofs of the potential inequalities in this section were obtained through the performance estimation methodology, introduced by~\citet{Dror14} and specialized to the study of inexact proximal operations by~\citet{barre2020principled}. More details can be found in Section~\ref{s:notes_ref_chapt_nest}, ``{On obtaining the proofs of this section}'' and in Appendix~\ref{a-WC_FO}. In particular, for reproducibility purposes, we provide code for symbolically verifying the algebraic reformulations of this section at \togglecodeurl{} together with those of Section~\ref{c-Nest}.

\chapter{Restart Schemes}\label{c-restart} 

In this section, we show that restart strategies can improve the performance of accelerated schemes when the objective function satisfies very generic H\"olderian error bounds (HEB) which generalize the notion of strong convexity, but only need to hold locally around the optimum. Restart schemes provide a convenient way to render standard first-order methods adaptive to the HEB parameters, and we will see that the cost of adaptation is only logarithmic.

\section{Introduction}\label{s:rest-intro}
\begin{figure}[ht]
\begin{center}
\begin{tabular}{ccc}
\raisebox{1.4\height}{\rotatebox{90}{$f-f_\star$}}{ }
\includegraphics[width=0.4\textwidth]{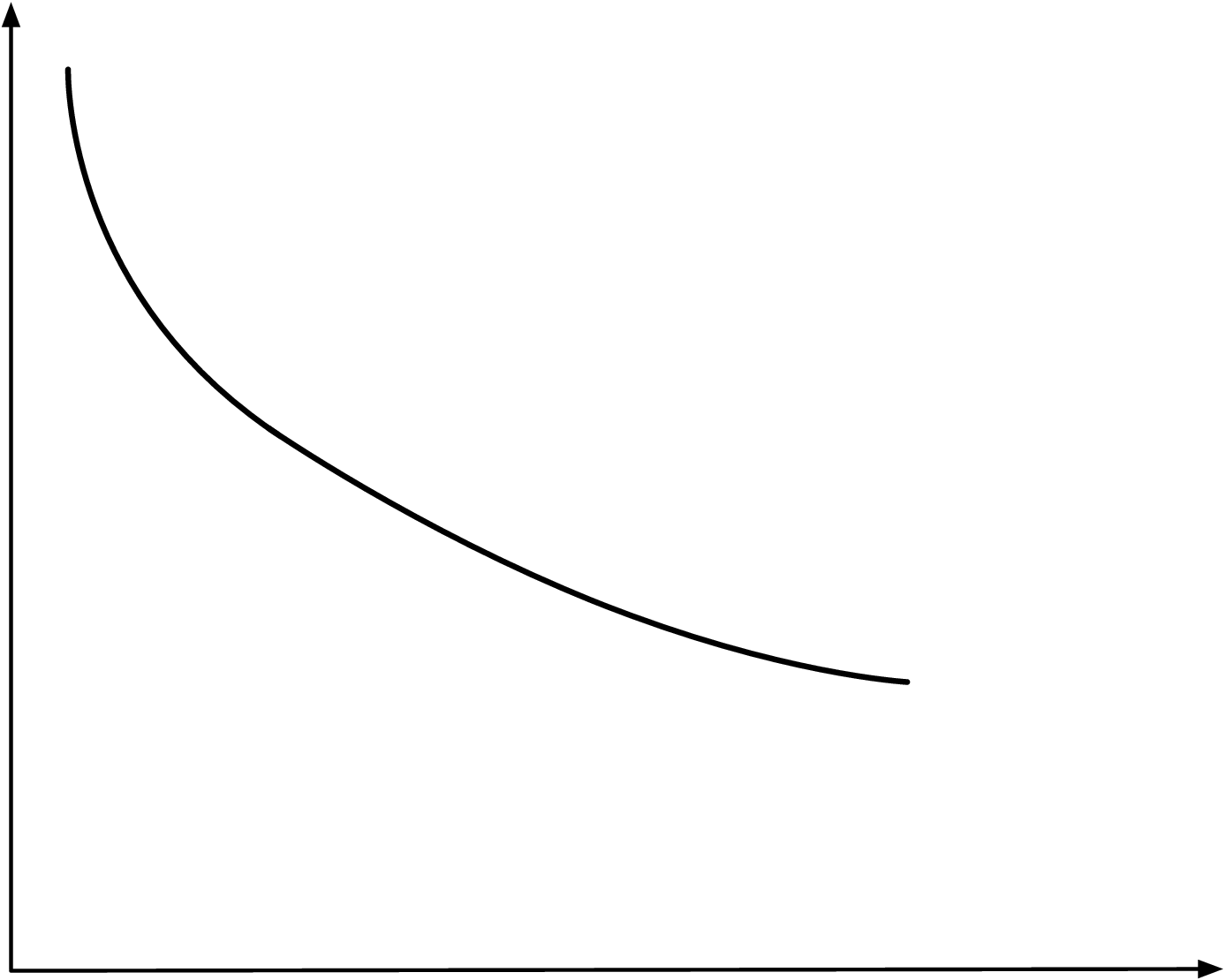}
&\qquad&
\raisebox{1.4\height}{\rotatebox{90}{$f-f_\star$}}{ }
\includegraphics[width=0.4\textwidth]{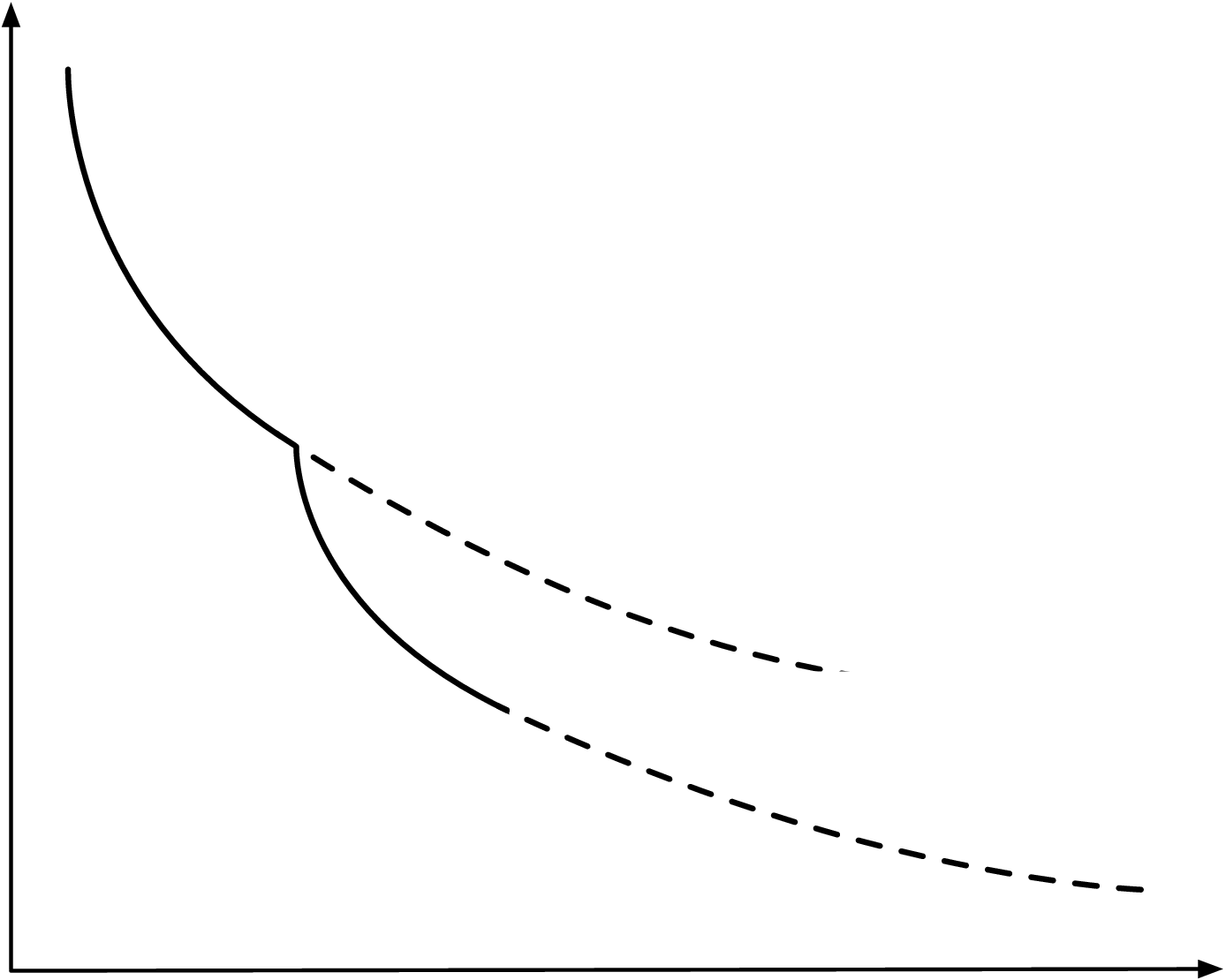}\\ 
{\footnotesize\quad  Iterations}
&\qquad&
{\footnotesize\quad Iterations}
\end{tabular}
\end{center}
\caption{{\em Left:} Sublinear convergence plot without restart. {\em Right:} Sublinear convergence plot with restart. \label{fig:rest-convex}}
\end{figure}

First-order methods typically exhibit a sublinear convergence, whose rate varies with gradient smoothness. The polynomial upper complexity bounds are typically convex functions of the number of iterations, so first-order methods converge faster in the beginning, then convergence tails off as iterations progress. This suggests that periodically restarting first-order methods, i.e., simply running more ``early'' iterations, could accelerate their convergence. We illustrate this concept in Figure~\ref{fig:rest-convex}.

Beyond this graphical argument, all accelerated methods have memory and look back at least one step to compute the next iterate. They iteratively form a model for the function around the optimum, and restarting allows this model to be periodically refreshed, thereby discarding outdated information as the algorithm converges towards the optimum. 

While the benefits of restart are immediately apparent in Figure~\ref{fig:rest-convex}, restart schemes raise several important questions: How many iterations should we run between restarts? What is the best complexity bound we can hope for using a restart scheme? What regularity properties of the problem drive the performance of restart schemes? Fortunately, all these questions have an explicit answer that stems from a simple and intuitive argument. We will see that restart schemes are also adaptive to unknown regularity constants and often reach near optimal convergence rates without observing these parameters. 

We begin by illustrating this adaptivity on the problem of minimizing a strongly convex function using the fixed step gradient method.

\subsection{The Strongly Convex Case}
We illustrate the main argument of this section when minimizing a strongly convex function using fixed step gradient descent. Suppose we seek to solve the minimization problem
\BEQ\label{eq:rest-prob}
\min_{x\in\reals^d}~f(x) 
\EEQ
Suppose that the gradient of $f$ is Lipschitz continuous with constant $L$ with respect to the Euclidean norm;
\BEQ\label{eq:rest-lip-grad}
\|\nabla f(y)-\nabla f(x)\|_2 \leq L \|y-x\|_2, \quad \mbox{for all $x,y\in\reals^d$}.
\EEQ

We can use the fixed step gradient method to solve problem~\eqref{eq:rest-prob}, as in Algorithm~\ref{alg:rest-grad} below.

\begin{algorithm}[ht]
  \caption{Gradient Method}
  \label{alg:rest-grad}
  \begin{algorithmic}[1]
    \REQUIRE
      A smooth convex function $f$ and an initial point $x_0$.
    \FOR{$k=0,\ldots$}
      \STATE $x_{k+1}=x_{k} - \frac{1}{L} \nabla f(x_k) $
    \ENDFOR
    \ENSURE An approximate solution $x_{k+1}$.
  \end{algorithmic}
\end{algorithm}

The smoothness assumption in~\eqref{eq:rest-lip-grad} ensures the complexity bound 
\BEQ\label{eq:rest-grad-bnd}
f(x_k)-f_\star \leq \frac{2L\|x_0-x_\star\|_2^2}{k+4}
\EEQ
after $k$ iterations (see Section~\ref{c-Nest} for a complete discussion).

Assume now that $f$ is also strongly convex with parameter $\mu$, with respect to the Euclidean norm. Strong convexity means that $f$ satisfies
\BEQ\label{eq:rest-strong-cvx}
\frac{\mu}{2} \|x-x_\star\|_2^2 \leq  f(x)-f_\star,
\EEQ
where $x_\star$ is an optimal solution to problem~\eqref{eq:rest-prob}, and $f_\star$ is the corresponding optimal objective value. Denote by $\mathcal{A}(x_0,k)$ the output of $k$ iterations of Algorithm~\ref{alg:rest-grad} started at $x_0$, and suppose that we periodically restart the gradient method according to the following scheme.

\begin{algorithm}
  \caption{Restart scheme}
  \label{alg:rest-scheme}
  \begin{algorithmic}[1]
    \REQUIRE
      A smooth convex function $f$, an initial point $x_0$ and an inner optimization algorithm $\mathcal{A}(x,k)$.
    \FOR{$i=0,\ldots,N-1$}
      \STATE Obtain $x_{i+1}$ by running $k_i$ iterations of the gradient method, starting at $x_i$, \ie
      \[
      x_{i+1}=\mathcal{A}(x_i,k_i)
      \]
    \ENDFOR
    \ENSURE An approximate solution $x_N$.
  \end{algorithmic}
\end{algorithm}

Combining the strong convexity bound in~\eqref{eq:rest-strong-cvx} with the complexity bound in~\eqref{eq:rest-grad-bnd} yields
\BEQ\label{eq:rest-chaining}
f(x_{i+1})-f_\star \leq \frac{2L\|x_i-x_\star\|_2^2}{k+4} \leq \frac{4L}{\mu(k+4)} (f(x_i)-f_\star)
\EEQ
after an iteration of the restart scheme in Algorithm~\ref{alg:rest-scheme} in which we run~$k$ (inner) iterations of the gradient method in Algorithm~\ref{alg:rest-grad}. This means that if we set 
\[
k_i= k = \left\lceil \frac{8L}{\mu} \right\rceil,
\]
then 
\[
f(x_{N})-f_\star \leq \left(\frac{1}{2}\right)^N (f(x_0)-f_\star)
\]
after $N$ iterations of the restart scheme in Algorithm~\ref{alg:rest-scheme}. Therefore, when running a total of $T=Nk$ gradient steps, we can rewrite the complexity bound in terms of the total number of gradient oracle calls (or inner iterations) as
\BEQ\label{eq:rest-lin-bnd}
f(x_{T})-f_\star \leq \left(\frac{1}{2^\frac{\mu}{8L}}\right)^T (f(x_0)-f_\star),
\EEQ
which proves linear convergence in the strongly convex case. 

Of course, the basic gradient method with fixed step size in Algorithm~\ref{alg:rest-grad} has no memory, so ``restarting'' it has no impact on the number of iterations or numerical performance. Invoking the restart scheme in Algorithm~\ref{alg:rest-scheme} simply allows us to produce a better complexity bound in the strongly convex case. Without information about the strong convexity parameter (since restart has no impact on the basic gradient method), whereas the classical bound yields sublinear convergence, while the restart method converges linearly.

Crucially here, the argument in~\eqref{eq:rest-chaining} can be significantly generalized to improve the convergence rate of several types of first-order methods. In fact, as we will see below, a local bound on the growth rate of the function akin to strong convexity holds almost generically, albeit with a different exponent than in~\eqref{eq:rest-strong-cvx}.

\subsection{Restart Strategies}
Empirical performance of restart schemes was studied at length in \citep{Beck12} and various restart strategies were explored to improve convergence of basic gradient methods by exploiting regularity properties of the objective function. \citep{Nest13} for example runs a bounded number of iterations between restarts to obtain linear convergence in the strongly convex case, while \citep{ODon15} obtain excellent empirical performance by restarting an accelerated method whenever convergence fails to be monotonic (accelerated methods typically exhibit oscillating convergence near the optimum). Below, we will describe the performance of a simple grid search on the restart strategy, attaining optimal performance while using a very limited number of grid points.

\section{H\"olderian Error Bounds}
We now recall several results related to subanalytic functions and H\"olderian error bounds of the form
\BEQ\label{eq:rest-loja}\tag{HEB}
    \frac{\mu}{r} d(x,X_\star)^r \leq  f(x)-f_\star,\quad\mbox{for all $x\in K,$}
\EEQ
for some $\mu,r>0$, where $d(x,X_\star)$ is the distance to the optimal set. We refer the reader to, e.g., \citep{Bolt07} for a more complete discussion. These results produce bounds akin to local versions of strong convexity, with various exponents, and they are known to hold under very generic conditions. In general of course, these values are neither observed nor known a priori, but as detailed below, restart schemes can be made adaptive to $\mu$ and $r$ and reach optimal convergence rates without any prior information.

\subsection{H\"olderian Error Bound and Smoothness} \label{sec:rest-sharp-smooth}
Let $f$ be a smooth convex function on $\reals^d$. Smoothness ensures that
\[
f(x) \leq f_\star+ \frac{L}{2}\|x-y\|^2_2,
\]
for any $x \in \reals^d$ and $y \in X_\star$. By setting $y$ to be the projection of $x$ on $X_\star$, this yields the following {\em upper bound} on suboptimality:
\BEQ\label{eq:rest-lb}
    f(x) - f_\star \leq \frac{L}{2} d(x,X_\star)^2.
\EEQ
Now, assume that $f$ satisfies the H\"olderian error bound~\eqref{eq:rest-loja} on a set $K$ with parameters $(r,\mu)$. Combining~\eqref{eq:rest-lb} and~\eqref{eq:rest-loja} leads to
\[
    \frac{2\mu}{rL} \leq d(x,X_\star)^{2-r},
\]
for every $x \in K$. This means that $2 \leq r$ by taking $x$ close enough to~$X_\star$. We will allow the gradient smoothness exponent of 2 to vary in later results, where we assume the gradient to be H\"older smooth, but we first detail the smooth case for simplicity. In what follows, we use the following notations:
\BEQ\label{eq:rest-kappa_tau}
\kappa \triangleq L/\mu^\frac{2}{r} \qquad \mbox{and} \qquad  \tau \triangleq 1-\frac{2}{r},
\EEQ
to define generalized condition numbers for the function $f$. Note that if $r=2$, then $\kappa$ matches the classical condition number of the function.

\subsection{Subanalytic Functions}
Subanalytic functions form a very broad class of functions for which we can demonstrate the H\"olderian error bounds as in~\eqref{eq:rest-loja}, akin to strong convexity. We recall some key definitions and refer the reader to, e.g., \citep{Bolt07} for a more complete discussion. 

\begin{definition}[Subanalyticity]
(i) A subset $A\subset\reals^d$ is called {\em semianalytic} if each point of $\reals^d$ admits a neighborhood $V$ for which $A \cap V$ assumes the following form
\[
\bigcup\limits_{i=1}^{p}\bigcap\limits_{j=1}^{q}\{x\in V: f_{ij}(x)=0, g_{ij}(x) >0\},
\]
where $f_{ij},g_{ij}:V\rightarrow \reals$ are real analytic functions for $1\leq i \leq p$, $1\leq j \leq q$.\\
(ii) A subset $A\subset\reals^d$ is called {\em subanalytic} if each point of $\reals^d$ admits a neighborhood $V$ such that
\[
A \cap V =\{x\in \reals^d: (x,y) \in B\}
\]
where $B$ is a bounded semianalytic subset of $\reals^d \times \reals^m$.\\
(iii) A function $f: \reals^d \rightarrow \reals \cup \{+\infty\}$ is called {\em subanalytic} if its graph is a subanalytic subset of $\reals^d \times \reals$.
\end{definition}

The class of subanalytic functions is, of course, very large, but the definition above suffers from one key shortcoming since the image and preimage of a subanalytic function are not generally subanalytic. To remedy this stability issue, we can define a notion of global subanalyticity. We first define the function $\beta_n$ with 
\[
\beta_d(x)\triangleq \left(\frac{x_1}{1+x_1^2},\ldots,\frac{x_d}{1+x_d^2}\right),
\]
and we have the following definition.

\begin{definition}[Global subanalyticity]
(i) A subset $A$ of $\reals^d$ is called {\em globally subanalytic} if its image under $\beta_d$ is a subanalytic subset of $\reals^d$.\\
(ii) A function $f:\reals^d \rightarrow \reals \cup \{+\infty\}$ is called {\em globally subanalytic} if its graph is a globally subanalytic subset of $\reals^d \times \reals$.
\end{definition}

We now recall the {\L}ojasiewicz factorization lemma, which gives us local growth bounds on the graph of a function around its minimum.

\begin{theorem}[{\L}ojasiewicz factorization lemma]\label{th:rest-loja-fact}
Let $K\subset\reals^d$ be a compact set and $g,h:K\rightarrow \reals$ two continuous globally subanalytic functions. If 
\[
h^{-1}(0) \subset g^{-1}(0),
\]
then
\BEQ
\frac{\mu}{r} |g(x)|^r \leq  |h(x)|,\quad\mbox{for all $x\in K,$}
\EEQ
for some $\mu,r>0$.
\end{theorem}

\begin{figure}[ht]
\begin{center}
\includegraphics[width=\textwidth]{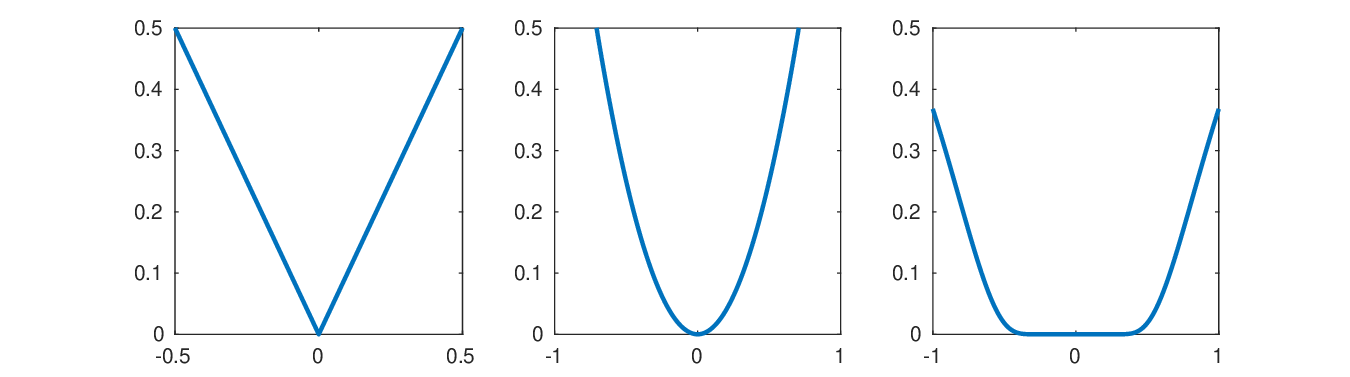}
\end{center}
\caption{{\em Left and center:} The functions $|x|$ and $x^2$ satisfy a growth condition around zero. {\em Right:} The function $\exp(-1/x^2)$ does not. \label{fig:rest-loja}}
\end{figure}

In an optimization context on a compact set $K\subset\reals^d$, we can set $h(x)=f(x)-f_\star$ and $g(x)=d(x,X_\star)$, the Euclidean distance from $x$ to the set $X_\star$, where $X_\star$ is the set of optimal solutions. In this case, we have $h^{-1}(0) \subset g^{-1}(0)$, and we can show that $g$ is globally subanalytic if $X_\star$ is globally subanalytic and $f$ is continuous and globally subanalytic. Theorem~\ref{th:rest-loja-fact} provides the following H\"olderian error bound,
\BEQ\label{eq:rest-loja-bis}
\frac{\mu}{r} d(x,X_\star)^r \leq  f(x)-f_\star,\quad\mbox{for all $x\in K,$}
\EEQ
for some $\mu,r>0$. Here, Theorem~\ref{th:rest-loja-fact} produces a bound on the growth rate of the function around the optimum, generalizing the strong convexity bound in~\eqref{eq:rest-strong-cvx}. We illustrate this in Figure~\ref{fig:rest-loja}. Overall, since continuity and subanalyticity are very weak conditions, Theorem~\ref{th:rest-loja-fact} shows that the H\"olderian error bound in~\eqref{eq:rest-loja} holds almost generically.

\section{Optimal Restart Schemes}
We now discuss how the H\"olderian error bounds detailed above can be exploited using restart schemes. Generic exponents beyond strong convexity, require restart schemes with a varying number of inner iterations (versus a constant one in the strongly convex case) and we we study here the cost of finding the best such scheme. Suppose again that we seek to solve the following unconstrained minimization problem:
\BEQ\label{eq:rest-prob-r}
\min_{x\in\reals^d}~f(x) 
\EEQ
where the gradient of $f$ is Lipschitz continuous with constant $L$ with respect to the Euclidean norm. The optimal method in~\eqref{cor:pot_fgm} detailed as Algorithm~\ref{alg:FGM_1} produces a point $x_k$ that satisfies 
\BEQ\label{eq:rest-smooth-bnd}
f(x_k)-f_\star \leq \frac{4L}{k^2}\|x_0-x_\star\|_2^2
\EEQ
after $k$ iterations.

Assuming that the function $f$ satisfies the H\"olderian error bound \eqref{eq:rest-loja}, we can use a chaining argument similar to that in~\eqref{eq:rest-chaining} to demonstrate improved convergence rates. While a constant number of inner iterations (between restarts) is optimal in the strongly convex case, the optimal restart scheme for $r>2$ involves a geometrically increasing number of inner iterations \citep{Nemic85,roulet2017sharpness}.

\begin{theorem}[Restart complexity]\label{th:rest-sched}
Let $f$ be a smooth convex function satisfying~\eqref{eq:rest-lip-grad} with parameter $L$ and~\eqref{eq:rest-loja} with parameters $(r,\mu)$ on a set $K$. Assume that we are given $x_0 \in \reals^d$ such that $\{x|\: f(x) \leq f(x_0) \} \subset K$. Run the restart scheme in Algorithm~\ref{alg:rest-scheme} from $x_0$ with iteration the schedule $k_i=C^\star_{\kappa,\tau}e^{\tau i}$, for $i=1,\ldots,R$,
where 
	\BEQ\label{def:C}
	C^\star_{\kappa,\tau} \triangleq e^{1-\tau}(c\kappa)^\frac{1}{2}(f(x_0)-f_\star)^{-\frac{\tau}{2}},
	\EEQ
with $\kappa$ and $\tau$ defined in~\eqref{eq:rest-kappa_tau} and $c = 4e^{2/e}$. The precision reached at the last point $\hat{x}$ is bounded by,
	\BEQ\label{eq:rest-conv-acc-bnd}
	f(\hat{x}) -f_\star \leq \frac{f(x_0)-f_\star}{\left(\tau e^{-1}(f(x_0)-f_\star)^{\frac{\tau}{2}}  (c\kappa)^{-\frac{1}{2}}N +1\right)^\frac{2}{\tau}} = O\left(N^{-\frac{2}{\tau}}\right),
	\EEQ
when $\tau >0$, where $N=\sum_{i=1}^R k_i$ is the total number of inner iterations.
\end{theorem} 

In the strongly convex case, i.e., when $\tau=0$, the bound above becomes 
\[
f(\hat{x}) - f_\star \leq \exp\left(-2e^{-1}(c\kappa)^{-\frac{1}{2}} N\right)(f(x_0)-f_\star) = O\left(\exp(-\kappa^{-\frac{1}{2}} N)\right)
\]
and we recover the classical linear convergence bound for Algorithm~\ref{alg:FGM_1_strconvex} in the strongly convex case. On the other hand, when $0<\tau<1$, bound \eqref{eq:rest-conv-acc-bnd} reveals a {\em faster convergence rate than accelerated gradient methods on non-strongly convex functions} (i.e., when $r>2$). The closer $r$ is to 2, the tighter the upper and lower bounds induced by smoothness and sharpness are, yielding a better model for the function and faster convergence. This property matches the lower bounds for optimizing smooth sharp functions~\citep{Nemic85} up to a constant factor. Moreover, setting $k_i = C^\star_{\kappa,\tau}e^{\tau i}$ yields continuous bounds on the precision, i.e.,~when $\tau \rightarrow 0$, bound \eqref{eq:rest-conv-acc-bnd} converges to the linear bound, which shows that for values of $\tau$ near zero, constant restart schemes are almost optimal. 

\section{Robustness and Adaptivity}
The previous restart schedules depend on the sharpness parameters $(r,\mu)$ in~\eqref{eq:rest-loja}. In general, of course, these values are neither observed nor known a priori. Making the restart scheme adaptive is thus crucial for practical performance. Fortunately, a simple logarithmic grid search on these parameters is enough to guarantee nearly optimal performance. In other words, as shown in \citep{roulet2017sharpness}, the complexity bound in~\eqref{eq:rest-conv-acc-bnd} is somewhat robust to misspecification of the inner iteration schedule $k_i$.

\subsection{Grid Search Complexity}
We can test several restart schemes in Algorithm~\ref{alg:rest-scheme}, each with a given number of inner iterations $N$ to perform a log-scale grid search on the values of $\tau$ and $\kappa$ in~\eqref{eq:rest-kappa_tau}. We see below that running $(\log_2 N)^2$ restart schemes suffices to achieve nearly optimal performance. We define these schemes as
\BEQ\label{algo:rest-grid}
\left\{\BA{l}
\mathcal{S}_{p,0}:\mbox{Restart Algorithm~\ref{alg:FGM_1} with $k_i = C_p$}, \\
\mathcal{S}_{p,q}:\mbox{Restart Algorithm~\ref{alg:FGM_1} with $k_i = C_p e^{\tau_q i}  $,} \\
\EA
\right.
\EEQ
where $C_p = 2^p$ and $\tau_q = 2^{-q}$. We stop these schemes when the total number of inner algorithm iterations exceeds $N$, i.e., at the smallest $R$ such that $\sum_{i=1}^R k_i \geq N$. The size of the grid search in $C_p$ is naturally bounded since as we cannot restart the algorithm after more than $N$ total inner iterations, so $p\in [1,\ldots,\lfloor \log_2 N \rfloor]$. Also, when $\tau$ is smaller than $1/N$, a constant schedule performs as well as the optimal, geometrically increasing schedule, which crucially means we can also choose $q \in [0,\ldots,\lceil \log_2 N \rceil ]$ and limits the cost of the grid search to $\log_2^2 N$. We have the following complexity bounds.

\begin{theorem}[Adaptive restart complexity]\label{th:rest-adap}
Let $f$ be a smooth convex function satisfying~\eqref{eq:rest-lip-grad} with parameter $L$ and \eqref{eq:rest-loja} with parameters $(r,\mu)$ on a set $K$. Assume that we are given $x_0 \in \reals^d$ such that $\{x |\: f(x) \leq f(x_0) \} \subset K$, and let $N$ be a given number of iterations. Run the schemes $\mathcal{S}_{p,q}$, defined in~\eqref{algo:rest-grid}, for $p\in [1,\ldots,\lfloor \log_2 N \rfloor]$ and $q \in [0,\ldots,\lceil \log_2 N \rceil]$, stopping each time after $N$ total inner algorithm iterations, i.e.,~for $R$ such that $\sum_{i=1}^R k_i \geq N$. Assume $N$ is large enough, such that $N\geq 2C^\star_{\kappa,\tau}$, and if $\frac{1}{N} >\tau >0$, $C^\star_{\kappa,\tau}>1$.\\
(i) If $\tau = 0$, there exists $p \in [1,\ldots,\lfloor \log_2 N \rfloor ]$ such that scheme $\mathcal{S}_{p,0}$ achieves a precision given by
	\[
	f(\hat{x}) -f_\star \leq \exp\left(-e^{-1}(c\kappa)^{-\frac{1}{2}} N\right)(f(x_0)-f_\star).
	\]
(ii) If $\tau > 0$, there exist $p \in [1,\ldots,\lfloor \log_2 N \rfloor ]$ and $q \in [1,\ldots,\lceil \log_2 N \rceil]$ such that scheme $\mathcal{S}_{p,q}$ achieves a precision given by
	\[
	f(\hat{x}) -f_\star \leq  \frac{f(x_0)-f_\star}{\left(\tau e^{-1}(c\kappa)^{-\frac{1}{2}} (f(x_0)-f_\star)^{\frac{\tau}{2}}  (N-1)/4 +1\right)^\frac{2}{\tau}}.
	\]
\end{theorem}

Overall, running the logarithmic grid search has a complexity that is $(\log_2 N )^2$ times higher than running $N$ iterations using the optimal scheme where we know the parameters in~\eqref{eq:rest-loja}, while the convergence rate is slowed down by roughly a factor four.

\section{Extensions}
We now discuss several extensions of the results above.

\subsubsection{H\"older Smooth Gradient} The results above can be extended somewhat directly to more general notions of regularity. In particular, if we assume that there exist $s \in [1,2]$ and $L>0$ on a set $J \subset \reals^d$, i.e.,
\BEQ\label{eq:rest-hold-grad}
\|\nabla f(x) -\nabla f(y)\|_2 \leq L \|x-y\|_2^{s-1}, \quad \mbox{for all $x,y \in J$}.
\EEQ
so that the gradient is H\"older smooth. Without further assumptions on~$f$, the optimal rate of convergence for this class of functions is bounded as $O(1/N^{\rho})$, where $N$ is the total number of iterations and
\BEQ\label{def:rest-q}
\rho = 3s/2-1,
\EEQ 
which gives $\rho =2$ for smooth functions and $\rho = 1/2$ for non-smooth functions. The universal fast gradient method \citep{Nest15} achieves this rate. It requires both a target accuracy $\epsilon$ and a starting point $x_0$ as inputs, and it outputs a point $x \triangleq \mathcal{U}(x_0,\epsilon,t)$ such that
\BEQ\label{eq:rest-algobound_gen}
f(x) - f_\star \leq  \frac{\epsilon}{2} + \frac{cL^\frac{2}{s} d(x_0,X_\star )^2}{\epsilon^{\frac{2}{s}}t^{\frac{2\rho}{s}}}  \frac{\epsilon}{2},
\EEQ
after $t$ iterations, where $c$ is a constant ($c = 2^\frac{4s-2}{s}$). We can extend the definition of $\kappa$ and $\tau$ in~\eqref{eq:rest-kappa_tau} to the case where the gradient is H\"older smooth, with
\BEQ\label{eq:rest-kappa-tau-Hold}
\kappa \triangleq \frac{L^\frac{2}{s}}{\mu^\frac{2}{r}} \qquad \mbox{and} \qquad  \tau \triangleq 1-\frac{s}{r}.
\EEQ
We see that $\tau$ acts as an analytical condition number that measures the tightness of upper and lower bound models. The key difference from the smooth case described above is that here we need to schedule {\em both} the target accuracy $\epsilon_i$ used by the algorithm {\em and} the number of iterations $k_i$ made at the $i\textsuperscript{th}$ run of the algorithm. Our scheme is described in Algorithm~\ref{alg:rest-scheduled_gen}.

\begin{algorithm}[ht]
	\caption{Universal scheduled restarts for convex minimization \label{alg:rest-scheduled_gen}}
	\begin{algorithmic}[1]	
		\REQUIRE $x_0\in\reals^d$, $\epsilon_0 \geq f(x_0)-f_\star$, $\gamma \geq 0$, a sequence $k_i$ and an inner algorithm $\mathcal{U}(x,\epsilon,k)$. 
		\FOR {$i=1,\ldots, R$}
		\STATE $\epsilon_i  =  e^{-\gamma} \epsilon_{i-1}$
		\STATE $x_i  =  \mathcal{U}(x_{i-1},\epsilon_i,k_i) $
		\ENDFOR
		\ENSURE An approximate solution $x_R$.
	\end{algorithmic}
\end{algorithm}

We choose a sequence $k_i$ that ensures
\[
f(x_i) - f_\star \leq \epsilon_i,
\]
for the geometrically decreasing sequence $\epsilon_i$. A grid search on the restart scheme still works in this case, but it requires knowledge of both $s$ and~$\tau$.

\begin{theorem}\label{th:rest-sched_gen}
Let $f$ be a convex function satisfying \eqref{eq:rest-hold-grad} with parameters $(s,L)$ on a set $J$ and \eqref{eq:rest-loja} with parameters $(r,\mu)$ on a set $K$. Given $x_0 \in \reals^d$ assume that $\{x | f(x) \leq f(x_0) \} \subset J\cap K$. Run the restart scheme in Algorithm~\ref{alg:rest-scheduled_gen} from $x_0$ for a given $\epsilon_0 \geq f(x_0)-f_\star$ with  
	\[ 
	\gamma =\rho, \qquad k_i = C^\star_{\kappa,\tau,\rho}e^{\tau i},\quad \mbox{where}\quad  C^\star_{\kappa,\tau,\rho} \triangleq e^{1-\tau}(c\kappa)^\frac{s}{2\rho}\epsilon_0^{-\frac{\tau}{\rho}} ,
	\]
and where $\rho$ is defined in \eqref{def:rest-q}, $\kappa$ and $\tau$ are defined in~\eqref{eq:rest-kappa-tau-Hold}, and $c = 8e^{2/e}$ here. The precision reached at the last point $x_R$ is given by
	\begin{align*}
	f(x_R) - f_\star ~\leq~ \exp\left(-\rho e^{-1} (c\kappa)^{-\frac{s}{2\rho}} N\right)\epsilon_0 ~=~ O\left(\exp(-\kappa^{-\frac{s}{2\rho}} N)\right),
	\end{align*}
when $\tau =0$, whereas when $\tau >0$,
	\begin{align*}
	f(x_R)-f_\star ~\leq~ \frac{\epsilon_0}
	{
		\left(\tau e^{-1} (c\kappa)^{-\frac{s}{2\rho}}\epsilon_0^{\frac{\tau}{\rho}}N
		+1
		\right)^{-\frac{\rho}{\tau}}} ~=~ O\left(\kappa^\frac{s}{2\tau}N^{-\frac{\rho}{\tau}}\right),
	\end{align*}
where $N=\sum_{i=1}^R k_i$ is the total number of iterations.
\end{theorem}

\subsubsection{Relative Smoothness} 
We can also extend the inequality defining condition~\eqref{eq:rest-loja} by replacing the distance to the optimal set by a more general Bregman divergence. Suppose $h(x)$ is a 1-strongly convex function with respect to the Euclidean norm. The Bregman divergence $D_h(x,y)$ is defined as
\BEQ\label{def:rest-bregman}
    D_h(x,y)\triangleq h(x)-h(y) - \langle  \nabla f(y) ; (x-y) \rangle,
\EEQ
and we say that a function $f$ is $L$-smooth with respect to $h$ on $\reals^d$ if 
\BEQ\label{eq:rest-L-smooth}
    f(y) \leq f(x) +\langle  \nabla f(x);(y-x)\rangle + L D_h(y,x),\quad \mbox{for all $x,y\in\reals^d$}.
\EEQ
We can then extend the H\"olderian error bound to the Bregman setting as follows. In an optimization context, on a compact set $K\subset\reals^d$, we can set $h(x)=f(x)-f_\star$ and $g(x)=D(x,X_\star)=\inf_{y\in X_\star}D_h(x,y)$, where $X_\star$ is the set of optimal solutions. In this case, we have $h^{-1}(0) \subset g^{-1}(0)$, and we can show that $g$ is globally subanalytic if $X_\star$ is subanalytic and if $f$ is continuous and globally subanalytic. Theorem~\ref{th:rest-loja-fact} shows that
\BEQ\label{eq:rest-loja-bregman}\tag{HEB-B}
    \frac{\mu}{r} D(x,X_\star)^r \leq  f(x)-f_\star,\quad\mbox{for all $x\in K,$}
\EEQ
for some $\mu,r>0$. This allows us to use the restart scheme complexity results above to accelerate proximal gradient methods.

\section{Calculus Rules}
In general, the exponent $r$ and the factor $\mu$ in the bounds~\eqref{eq:rest-loja} and~\eqref{eq:rest-strong_wolfe_primal_gap} are not observed and are difficult to estimate. Nevertheless, due to the robustness result in Theorem~\ref{th:rest-adap}, searching for the best restart scheme only introduces a $\log$ factor in the overall algorithm complexity. There are, however, a number of scenarios where we can produce much more precise estimates of $r$ and $\mu$ and hence both obtain refined a priori complexity bounds and reduce the cost of the grid search in~\eqref{algo:rest-grid}. 

In particular, \citep{Li18} provides ``calculus rules'' for the HEB exponent for a number of elementary operations using a related type of error bound known as the Kurdyka-\L ojasiewicz inequality; see~\citep[Theorem~5]{Bolt07} for more details on the relationship between these two notions. The results focus on the Kurdyka-\L ojasiewicz exponent $\alpha$, defined as follows.

\begin{definition}\label{def:rest-KL-exp}
A proper closed convex function has a Kurdyka-\L ojasie\-wicz (KL) exponent $\alpha$ if and only if for any point $\bar x \in \mathrm{dom}\partial f$ there is a neighborhood $\mathcal{V}$ of $\bar x$, a constant $\nu>0$, and a constant $c>0$ such that 
\BEQ\label{eq:KL}
D(\partial f(x),0) \geq c(f(x)-f(\bar x))^\alpha 
\EEQ
whenever $x \in \mathcal{V}$ and $f(\bar x) \leq f(x) \leq f(\bar x) + \nu$.
\end{definition}

In particular, \citep[Theorem\,3.3]{Bolt07} shows that~\eqref{eq:rest-loja} implies~\eqref{eq:KL} with exponent $\alpha=1-1/r$. The other way also holds with $r=1/(1-\alpha)$, but the constant is degraded; see~\citep[Section 3.1]{Bolt07}. Very briefly, the following calculus rules apply to the exponent $\alpha$.

\begin{itemize}
    \item If $f(x) = \min_i f_i(x)$ and each $f_i$ has the KL exponent $\alpha_i$, then $f$ has the KL exponent $\alpha=\max_i \alpha_i$ \citep[Corollary\,3.1]{Li18}.
    \item Let $f (x) = g \circ F (x)$, where g is a proper closed function and $F$ is a continuously differentiable mapping. Suppose in addition that g is a KL function with exponent $\alpha$ and that the Jacobian $JF(x)$ is a surjective mapping at some $\bar x \in \mathrm{dom} \partial f$. Then $f$ has the KL property at $\bar x$ with exponent $\alpha$ \citep[Theorem\,3.2]{Li18}.
    \item If $f(x) = \sum_i f_i(x_i)$ and each $f_i$ is continuous and has the KL exponent $\alpha_i$, then $f$ has KL exponent $\alpha=\max_i \alpha_i$ \citep[Corollary\,3.3]{Li18}.
    \item Let $f$ be a proper closed convex function with a KL exponent $\alpha \in [0,2/3]$. Suppose further that f is continuous on $\mathrm{dom} \partial f$. Fix $\lambda >0$ and consider
    \[
    F_\lambda(X) = \inf_y \left\{f(y)+ \frac{1}{2\lambda} \|x-y\|^2\right\}.
    \]
    Then $F_\lambda$ has the KL exponent $\alpha = \max \left\{\frac{1}{2} , \frac{\alpha}{2-2\alpha} \right\}$ \citep[Theorem\,3.4]{Li18}.
\end{itemize}

Note that a related notion of error bound in which the primal gap is replaced by the norm of the proximal step was studied in, e.g., \citep{Pang87,Luo92b,Tsen10,Zhou17}.

\section{Restarting Other First-Order Methods}
The restart argument can be readily extended to other optimization methods provided their complexity bound directly depends on some measure of distance to optimality. This is the case for instance for the Frank-Wolfe method, as detailed in~\citep{kerdreux2019restarting}. Suppose that we seek to solve the following constrained optimization problem
\BEQ\label{eq:rest-fw-prob}
\min_{x \in \mathcal{C}}  f(x).
\EEQ
The distance to optimality is now measured in terms of the strong Wolfe gap, defined as follows.

\begin{definition}[Strong Wolfe gap]\label{def:strong_wolfe_gap}
Let $f$ be a smooth convex function, $\mathcal{C}$ a polytope, and $x \in \mathcal{C}$ be arbitrary. Then the \emph{strong Wolfe gap $w(x)$ over $\mathcal{C}$} is defined as
\BEQ\label{eq:strong_wolfe_gap}
w(x) \triangleq\min_{S\in\mathcal{S}_x} ~ \max_{y \in S, z \in \mathcal{C}} \langle f(x) ; (y-z) \rangle,
\EEQ
where $x \in \Co(S)$ and 
\[
S_x = \{ S\subset\Ext(\mathcal{C}), \mbox{ finite, $x$ proper combination of elements of $S$}\},
\] 
is the set of proper supports of $x$.
\end{definition}

The inequality that plays the role of the H\"olderian error bound in~\eqref{eq:rest-loja} for the strong Wolfe gap is then written as follows.

\begin{definition}[Strong Wolfe primal bound]\label{def:rest-strong_wolfe_primal_gap}
Let $K$ be a compact neighborhood of $X_\star$ in $\mathcal{C}$, where $X_\star$ is the set of solutions of the constrained optimization problem \eqref{eq:rest-fw-prob}. A function $f$ satisfies an $r$-strong Wolfe primal bound on $K$, if and only if there exists $r\geq 1$ and $\mu>0$ such that for all $x\in K$
\BEQ\label{eq:rest-strong_wolfe_primal_gap}
f(x) - f_\star \leq \mu w(x)^{r},
\EEQ
where $f_\star$ it the minimum of $f$.
\end{definition}

Notice that this inequality is an upper bound on the primal gap $f(x) - f_\star$, whereas the H\"olderian error bound in~\eqref{eq:rest-loja} provides a lower bound. This is because the strong Wolfe gap can be understood as a gradient norm, such that~\eqref{eq:rest-strong_wolfe_primal_gap} is a {\L}ojasiewicz inequality as in \citep{Bolt07}, instead of a direct consequence of the {\L}ojasiewicz factorization lemma as in~\eqref{eq:rest-loja} above.

The regularity of $f$ is measured using the \emph{away curvature} as in \citep{lacoste2015global}, with
\BEQ\label{eq:away_curvature}
C_f^A\triangleq\underset{\substack{x,s,v\in\mathcal{C}\\\eta\in [0,1] \\ y=x+\eta (s-v)}}{\text{sup }}{\frac{2}{\eta^2}\big( f(y)-f(x)-\eta\langle\nabla f(x),s-v\rangle \big)},
\EEQ
allowing us to bound the performance the Fractional Away-Step Frank-Wolfe Algorithm in~\citep{kerdreux2019restarting}, as follows.

\begin{theorem}\label{th:rest-convergence_restart_fw}
Let $f$ be a smooth convex function with away curvature $C_f^{A}$. Assume the strong Wolfe primal bound in~\eqref{eq:rest-strong_wolfe_primal_gap} holds for some $1 \leq r \leq 2$. Let $\gamma > 0$ and assume $x_0\in \mathcal{C}$ is such that $e^{-\gamma} w(x_0,\mathcal{S}_0) / 2 \leq C_f^{A}$. With $\gamma_k = \gamma$, the output of the Fractional Away-Step Frank-Wolfe Algorithm satisfies
\begin{equation}\label{eq:cv_rates}
  \left\{
  \begin{array}{ll}
    f(x_T) - f_\star \leq w_0\frac{1}{\Big( 1 + \tilde{T} C_{\gamma}^r \Big)^{\frac{1}{2-r}}}   &  \text{when } 1 \leq r < 2\\
    f(x_T) - f_\star \leq  w_0 \exp\left(-\frac{\gamma}{e^{2\gamma}}\frac{\tilde{T}}{8 C_f^{A} \mu }\right)    & \text{when } r=2~,
  \end{array}
  \right.
\end{equation}
after $T$ steps, with $w_0=w(x_0,\mathcal{S}_0)$, $\tilde{T}\triangleq T - (|\mathcal{S}_0| - |\mathcal{S}_T|)$, and
\BEQ\label{eq:fct_gamma_rate}
C_{\gamma}^r \triangleq \frac{e^{\gamma (2-r)}-1}{ 8 e^{2\gamma} C_f^{A} \mu w(x_0,\mathcal{S}_0)^{r-2}}.
\EEQ
\end{theorem}

This result is similar to that of Theorem~\ref{th:rest-sched_gen}, and it shows that restart yields linear complexity bounds when the exponent in the strong Wolfe primal bound in~\eqref{eq:rest-strong_wolfe_primal_gap} matches that in the curvature (i.e., $r=2$) and that it yields to improved linear rates when the exponent $r$ satisfies $1 \leq r < 2$. Crucially, the method here is fully adaptive to the error bound parameters, so no prior knowledge of these parameters is required to get the accelerated rates in Theorem~\ref{th:rest-convergence_restart_fw}, and no log-scale grid search is required.

\section{Application: Compressed Sensing}
In some applications such as compressed sensing, under some classical assumptions on the problem data, the exponent $r$ is equal to one and the constant $\mu$ can be directly computed from quantities controlling recovery performance. In such problems, a single parameter thus controls both signal recovery and computational performance.

Consider, for instance, a sparse recovery problem using the $\ell_1$ norm. Given a matrix $A\in\reals^{n\times p}$ and observations $b=Ax_\star$ on a signal $x_\star\in\reals^p$, recovery is performed by solving the $\ell_1$ minimization program
\BEQ\label{eq:rest-l1-recov}\tag{$\ell_1$ recovery}
\BA{ll}
\mbox{minimize} & \|x\|_1\\
\mbox{subject to} & Ax=b
\EA
\EEQ
in the variable $x\in\reals^p$. A number of conditions on $A$ have been derived to guarantee that~\eqref{eq:rest-l1-recov} recovers the true signal whenever it is sparse enough. Among these, the null space property (see \cite{Cohe06} and references therein) is defined as follows. 

\begin{definition}[Null space property]\label{def:nsp_prop} The matrix $A$ satisfies the Null Space Property (NSP) \emph{on support $S\subset \{ 1,p \}$} with constant $\alpha \geq 1$ if for any $z\in \Null(A)\setminus\{0\}$,
	\BEQ\label{def:nsp}\tag{NSP}
	\alpha \|z_S\|_1 < \|z_{S^c}\|_1.
	\EEQ
	The matrix $A$ satisfies the Null Space Property \emph{at order $s$} with constant $\alpha \geq 1$ if it satisfies it on every support $S$ of cardinality at most $s$.
\end{definition}

The null space property is a necessary and sufficient condition for the convex program \eqref{eq:rest-l1-recov} to recover all signals up to some sparsity threshold. We have, the following proposition directly linking the null space property and the H\"olderian error bound~\eqref{eq:rest-loja}.

\begin{proposition}\label{prop:recov}
	Given a coding matrix $A\in\reals^{n \times p}$ satisfying \eqref{def:nsp} at order $s$ with constant $\alpha \geq 1$, if the original signal $x_\star$ is $s$-sparse, then for any $x\in\reals^p$ satisfying $Ax=b$, $x\neq x_\star$, we have  
	\BEQ
	\|x\|_1-\|x_\star\|_1 > \frac{\alpha-1}{\alpha+1} \|x-x_\star\|_1.
	\EEQ
	This implies signal recovery, i.e. optimality of $x_\star$ for \eqref{eq:rest-l1-recov}\, and the H\"olderian error bound~\eqref{eq:rest-loja} with $\mu = \frac{\alpha-1}{\alpha+1}$.
\end{proposition}

\section{Notes and References}

The optimal complexity bounds and exponential restart schemes detailed here can be traced back to \citep{Nemic85}. Restart schemes were extensively benchmarked in the numerical toolbox TFOCS by \citep{Beck12}, with a particular focus on compressed sensing applications. The robustness result showing that a log scale grid search produces near optimal complexity bounds is due to \citep{roulet2017sharpness}.

Restart schemes based on the gradient norm as a termination criterion also reach nearly optimal complexity bounds and adapt to strong convexity \citep{Nest13} or HEB parameters \citep{Ito21}.

H\"olderian error bounds for analytic functions can be traced back to the work of \citet{Loja63}. They were extended to much broader classes of functions by \citep{Kurd98,Bolt07}. Several examples of problems in signal processing where this condition holds can be found in, e.g., \citep{Zhou15,Zhou17}. Calculus rules for the exponent are discussed in details in, e.g.,~\citep{Li18}.

Restarting is also helpful in the stochastic setting, with \citep{Davi19} showing recently that stochastic algorithms with geometric step decay converge linearly on functions satisfying H\"olderian error bounds. This validates a classical empirical acceleration trick, which is to restarts every few epochs after adjusting the step size (aka the learning rate in machine learning terminology).

\appendix
\chapter{Useful Inequalities}\label{a-inequalities}

In this appendix, we prove basic inequalities involving smooth strongly convex functions. Most of these inequalities are not used in our developments. Nevertheless, we believe they are useful for gaining intuition about smooth strongly convex of functions, as well as for comparisons with the literature.

Also note that these inequalities can be considered standard (see, e.g.,~\citep[Theorem 2.1.5]{Nest03a}.

\section{Smoothness and Strong Convexity in Euclidean spaces}\label{s:ineq_eucl}
In this section, we consider a Euclidean setting, where $\|x\|_2^2=\langle x;x\rangle$ and $\langle .;.\rangle:\mathbb{R}^d\times\mathbb{R}^d\rightarrow\mathbb{R}$ is a dot product.

The following theorem summarizes known inequalities that characterize the class of smooth convex functions. Note that these characterizations of $f\in\mathcal{F}_{0,L}$ are all equivalent assuming that $f\in\mathcal{F}_{0,\infty}$ since convexity is not implied by some of the points below. In particular, (i), (ii), (v), (vi), and (vii) do not encode the convexity of $f$ when taken on their own, whereas (iii) and (iv) encode both smoothness and convexity.
\begin{theorem}\label{thm:ineqs_smooth} Let $f:\mathbb{R}^d\rightarrow\mathbb{R}$ be a differentiable convex function. The following statements are equivalent for inclusion in $\mathcal{F}_{0,L}$.
\begin{itemize}
    \item[(i)] $\nabla f$ satisfies a Lipschitz condition: for all $x,y\in\mathbb{R}^d$,
    \[\| \nabla f(x)-\nabla f(y)\|_2\leq L\|x-y\|_2.\]
    \item[(ii)] $f$ is upper bounded by quadratic functions: for all $x,y\in\mathbb{R}^d$,
    \[ f(x)\leq f(y)+\langle \nabla f(y);x-y\rangle +\frac{L}{2}\|x-y\|_2^2.\]
    \item[(iii)] $f$ satisfies, for all $x,y\in\mathbb{R}^d$,
    \[f(x)\geq f(y)+\langle \nabla f(y);x-y\rangle+\frac{1}{2L}\|\nabla f(x)-\nabla f(y)\|_2^2.\]
    \item[(iv)] $\nabla f$ is cocoercive: for all $x,y\in\mathbb{R}^d$,
    \[ \langle \nabla f(x)-\nabla f(y);x-y\rangle\geq \frac1L\|\nabla f(x)-\nabla f(y)\|_2^2.\]
    \item[(v)] $\nabla f$ satisfies, for all $x,y\in\mathbb{R}^d$,
    \[ \langle \nabla f(x)-\nabla f(y);x-y\rangle\leq L\|x-y\|_2^2.\]
    \item[(vi)] $\frac{L}{2}\|x\|_2^2-f(x)$ is convex.
    \item[(vii)] $f$ satisfies, for all $\lambda\in[0,1]$,
    \[ f(\lambda x+(1-\lambda)y)\geq \lambda f(x)+(1-\lambda)f(y)-\lambda(1-\lambda)\frac{L}{2}\|x-y\|_2^2.\]
\end{itemize}
\end{theorem}
\begin{proof}
We start with (i)$\Rightarrow$(ii). We use the first-order expansion
\[ f(y)=f(x)+\int_{0}^1 \langle \nabla f(x+\tau(y-x));y-x\rangle d\tau.\]
The quadratic upper bound then follows from algebraic manipulations and from upper bounding the integral term:
\begin{equation*}
\begin{aligned}
f(y)&=f(x)+\langle \nabla f(x);y-x\rangle \\
& \quad + \int_0^1 \langle \nabla f(x+\tau(y-x))-\nabla f(x);y-x\rangle d\tau\\
&\leq f(x)+\langle \nabla f(x);y-x\rangle \\
& \quad +\int_0^1 \|\nabla f(x+\tau(y-x))-\nabla f(x)\|_2 \|y-x\|_2 d\tau\\
&\leq f(x)+\langle \nabla f(x);y-x\rangle+ L\|x-y\|_2^2\int_0^1 \tau d\tau\\
&= f(x)+\langle \nabla f(x);y-x\rangle+ \frac{L}{2}\|x-y\|_2^2.
\end{aligned}
\end{equation*}

We proceed with (ii)$\Rightarrow$(iii). The idea is to require the quadratic upper bound to be everywhere above the linear lower bound arising from the convexity of $f$. That is, for all $x,y,z\in\mathbb{R}^d$,
\[ f(y)+\langle \nabla f(y);z-y\rangle \leq f(z) \leq f(x)+\langle \nabla f(x);z-x\rangle +\frac{L}{2}\|x-z\|_2^2.\]
In other words, for all $z\in\mathbb{R}^d$, we must have
\begin{equation*}
\begin{aligned}
f(y)+\langle \nabla f(y);z-y\rangle \leq f(x)+\langle \nabla f(x);z-x\rangle +\frac{L}{2}\|x-z\|_2^2\\
\Leftrightarrow f(y)-f(x)+\langle \nabla f(y);z-y\rangle-\langle \nabla f(x);z-x\rangle-\frac{L}{2}\|x-z\|_2^2\leq 0\\
\Leftrightarrow f(y)-f(x)+\max_{z\in\mathbb{R}^d} \langle \nabla f(y);z-y\rangle-\langle \nabla f(x);z-x\rangle-\frac{L}{2}\|x-z\|_2^2\leq 0\\
\Leftrightarrow f(y)-f(x)+\langle \nabla f(y);x-y\rangle +\frac1{2L}\|\nabla f(x)-\nabla f(y)\|_2^2\leq 0,
\end{aligned}
\end{equation*}
where the last line follows from the explicit maximization on $z$. That is, we pick $z=x-\tfrac1L(\nabla f(x)-\nabla f(y))$ and reach the desired result after base algebraic manipulations.

We continue with  (iii)$\Rightarrow$(iv), which simply follows from adding
\begin{equation*}
\begin{aligned}
f(x)&\geq f(y)+\langle \nabla f(y);x-y\rangle +\frac1{2L}\|\nabla f(x)-\nabla f(y)\|_2^2\\
f(y)&\geq f(x)+\langle \nabla f(x);y-x\rangle +\frac1{2L}\|\nabla f(x)-\nabla f(y)\|_2^2.
\end{aligned}
\end{equation*}
To obtain (iv)$\Rightarrow$(i), one can use Cauchy-Schwartz:
\begin{align*}
\frac1L \|\nabla f(x)-\nabla f(y)\|_2^2 \leq& \langle \nabla f(x)-\nabla f(y);x-y\rangle\\
\leq&\|\nabla f(x)-\nabla f(y)\|_2\|x-y\|_2,
\end{align*}
which allows us to conclude that $\|\nabla f(x)-\nabla f(y)\|_2\leq L\|x-y\|_2$, thus reaching the final statement.

To obtain (ii)$\Rightarrow$(v), we simply add 
\begin{equation*}
\begin{aligned}
f(x)&\leq f(y)+\langle \nabla f(y);x-y\rangle +\frac{L}{2}\|x-y\|_2^2\\
f(y)&\leq f(x)+\langle \nabla f(x);y-x\rangle +\frac{L}{2}\|x-y\|_2^2
\end{aligned}
\end{equation*}
and reorganize the resulting inequality.

To obtain (v)$\Rightarrow$(ii), we again use a first-order expansion:
\[ f(y)=f(x)+\int_{0}^1 \langle \nabla f(x+\tau(y-x));y-x\rangle d\tau.\]
The quadratic upper bound then follows from algebraic manipulations and from upper bounding the integral term. (We use the intermediate variable $z_\tau=x+\tau(y-x)$ for convenience)
\begin{equation*}
\begin{aligned}
f(y)&=f(x)+\langle \nabla f(x);y-x\rangle \\
& \quad +\int_0^1 \langle \nabla f(x+\tau(y-x))-\nabla f(x);y-x\rangle d\tau\\
&=f(x)+\langle \nabla f(x);y-x\rangle+\int_0^1 \frac1{\tau}\langle \nabla f(z_\tau)-\nabla f(x);z_\tau-x\rangle d\tau\\
&\leq f(x)+\langle \nabla f(x);y-x\rangle+\int_0^1 \frac{L}{\tau} \|z_\tau-x\|_2^2 d\tau\\
&= f(x)+\langle \nabla f(x);y-x\rangle+ L\|x-y\|_2^2\int_0^1 \tau d\tau\\
&= f(x)+\langle \nabla f(x);y-x\rangle+ \frac{L}{2}\|x-y\|_2^2.
\end{aligned}
\end{equation*}

For the equivalence (vi)$\Leftrightarrow$(ii), simply define $h(x)=\frac{L}{2}\|x\|_2^2-f(x)$ (and hence $\nabla h(x)=Lx-\nabla f(x)$) and observe that for all $x,y\in\mathbb{R}^d$,
\[ h(x)\geq h(y)+\langle \nabla h(y);x-y\rangle\,\Leftrightarrow\,f(x)\leq f(y)+\langle \nabla f(y);x-y\rangle+\frac{L}{2}\|x-y\|_2^2,\]
which follows from base algebraic manipulations.

Finally, the equivalence (vi)$\Leftrightarrow$(vii) follows the same $h(x)=\frac{L}{2}\|x\|_2^2-f(x)$ (and hence $\nabla h(x)=Lx-\nabla f(x)$) and the observation that for all $x,y\in\mathbb{R}^d$ and $\lambda\in[0,1]$, we have
\begin{equation*}
\begin{aligned}
h(\lambda x+(1-\lambda)y)&\leq \lambda h(x)+(1-\lambda)h(y)\\&\Leftrightarrow \\
f(\lambda x+(1-\lambda)y)&\geq \lambda f(x)+(1-\lambda)f(y)-\lambda(1-\lambda)\frac{L}{2}\|x-y\|_2^2,
\end{aligned}
\end{equation*}
which follows from base algebraic manipulations.
\end{proof}

To obtain the corresponding inequalities in the strongly convex case, one can rely on Fenchel conjugation between smoothness and strong convexity; see, for example,~\citep[Proposition 12.6]{rockafellar2009variational}. The following inequalities are stated without proofs; they can be obtained either as direct consequences of the definitions or from Fenchel conjugation along with the statements of Theorem~\ref{thm:ineqs_smooth}.
\begin{theorem}\label{thm:ineqs_smooth_scvx} Let $f:\mathbb{R}^d\rightarrow\mathbb{R}$ be a closed convex proper function. The following statements are equivalent for inclusion in $\mathcal{F}_{\mu,L}$.
\begin{itemize}
    \item[(i)] $\nabla f$ satisfies a Lipschitz and an inverse Lipschitz condition: for all $x,y\in\mathbb{R}^d$,
    \[\mu\|x-y\|_2\leq \| \nabla f(x)-\nabla f(y)\|_2\leq L\|x-y\|_2.\]
    \item[(ii)] $f$ is lower and upper bounded by quadratic functions: for all $x,y\in\mathbb{R}^d$,
    \begin{align*}
    f(y)+&\langle \nabla f(y);x-y\rangle +\frac{\mu}{2}\|x-y\|_2^2 \\&\leq f(x)\leq f(y)+\langle \nabla f(y);x-y\rangle +\frac{L}{2}\|x-y\|_2^2.
    \end{align*}
    \item[(iii)] $f$ satisfies, for all $x,y\in\mathbb{R}^d$,
    \begin{equation*}
    \begin{aligned}
    f(y)+&\langle \nabla f(y);x-y\rangle+\frac{1}{2 L}\|\nabla f(x)-\nabla f(y)\|_2^2\\&\leq f(x)\leq\\ f(y)&+\langle \nabla f(y);x-y\rangle+\frac{1}{2\mu}\|\nabla f(x)-\nabla f(y)\|_2^2.
    \end{aligned}
    \end{equation*}
    \item[(iv)] $\nabla f$ satisfies, for all $x,y\in\mathbb{R}^d$,
    \begin{align*}
    \frac1L& \| \nabla f(x)-\nabla f(y)\|_2^2\\&\leq \langle \nabla f(x)-\nabla f(y);x-y\rangle
    \leq \frac1\mu \| \nabla f(x)-\nabla f(y)\|_2^2.
    \end{align*}
    \item[(v)] $\nabla f$ satisfies, for all $x,y\in\mathbb{R}^d$,
    \[ \mu \| x-y\|_2^2\leq \langle \nabla f(x)-\nabla f(y);x-y\rangle\leq L \| x-y\|_2^2.\]
    \item[(vi)] For all $\lambda\in[0,1]$,
    \begin{equation*}
    \begin{aligned}
    \lambda f(x)+&(1-\lambda)f(y)-\lambda(1-\lambda)\frac{L}{2}\|x-y\|_2^2\\&\leq f(\lambda x+(1-\lambda)y)\leq\\
    \lambda f(x)+&(1-\lambda)f(y)-\lambda(1-\lambda)\frac{\mu}{2}\|x-y\|_2^2.
    \end{aligned}
    \end{equation*}
    \item[(vii)] $f(x)-\frac{\mu}{2}\|x\|_2^2$ and $\frac{L}{2}\|x\|_2^2-f(x)$ are convex and $(L-\mu)$-smooth.
\end{itemize}
\end{theorem}
Finally, we mention that the existence of an inequality that allows us to encode both smoothness and strong convexity together. This inequality is also known as an \emph{interpolation} inequality~\citep{taylor2017smooth}, and it turns out to be particularly useful for proving worst-case guarantees.
\begin{theorem}\label{thm:interp_proof} Let $f:\mathbb{R}^d\rightarrow\mathbb{R}$ be a differentiable function. $f$ is $L$-smooth $\mu$-strongly convex if and only if
\begin{equation}\label{eq:interp_app}
\begin{aligned}
f(x)\geq f(y)+&\langle \nabla f(y);x-y\rangle+\frac{1}{2L}\|\nabla f(x)-\nabla f(y)\|_2^2\\
+&\frac{\mu}{2(1-\mu/L)}\|x-y-\frac1L(\nabla f(x)-\nabla f(y))\|_2^2.
\end{aligned}
\end{equation}
\end{theorem}
\begin{proof}
($f\in\mathcal{F}_{\mu,L}\Rightarrow$~\eqref{eq:interp_app}) The idea is to require the quadratic upper bound from smoothness to be everywhere above the quadratic lower bound arising from strong convexity. That is, for all $x,y,z\in\mathbb{R}^d$
\begin{align*}
f(y)+\langle \nabla f(y);z-y\rangle +\frac{\mu}{2}\|z-y\|_2^2\leq f(z) \leq& f(x)+\langle \nabla f(x);z-x\rangle \\
& +\frac{L}{2}\|x-z\|_2^2.
\end{align*}
In other words, for all $z\in\mathbb{R}^d$, we must have
\begin{equation*}
\begin{aligned}
f(y)+& \langle \nabla f(y);z-y\rangle +\frac{\mu}{2}\|z-y\|_2^2 \leq f(x)\\
& \quad\quad\quad +\langle \nabla f(x);z-x\rangle +\frac{L}{2}\|x-z\|_2^2\\
\Leftrightarrow & f(y)-f(x)+\langle \nabla f(y);z-y\rangle +\frac{\mu}{2}\|z-y\|_2^2-\langle \nabla f(x);z-x\rangle\\
& \quad\quad\quad -\frac{L}{2}\|x-z\|_2^2\leq 0\\
\Leftrightarrow & f(y)-f(x)+\max_{z\in\mathbb{R}^d} \bigg(\langle \nabla f(y);z-y\rangle +\frac{\mu}{2}\|z-y\|_2^2 \\
& \quad\quad\quad -\langle \nabla f(x);z-x\rangle-\frac{L}{2}\|x-z\|_2^2\bigg)\leq 0
\end{aligned}
\end{equation*}
explicit maximization over $z$. That is, picking $z=\frac{Lx-\mu y}{L-\mu}-\tfrac1{L-\mu}(\nabla f(x)-\nabla f(y))$ allows the desired inequality to be reached by base algebraic manipulations.

(\eqref{eq:interp_app}$\Rightarrow f\in\mathcal{F}_{\mu,L}$) $f\in\mathcal{F}_{0,L}$ is direct by observing that~\eqref{eq:interp_app} is stronger than Theorem~\ref{thm:ineqs_smooth}(iii); $f\in\mathcal{F}_{\mu,L}$ is then direct by reformulating~\eqref{eq:interp_app} as
\begin{equation*}
\begin{aligned}
f(x)\geq f(y)+&\langle \nabla f(y);x-y\rangle+\frac{\mu}{2}\|x-y\|_2^2\\
+&\frac{1}{2L(1-\mu/L)}\|\nabla f(x)-\nabla f(y)-\mu(x-y)\|_2^2,
\end{aligned}
\end{equation*}
which is stronger than $f(x)\geq f(y)+\langle \nabla f(y);x-y\rangle+\frac{\mu}{2}\|x-y\|_2^2$.
\end{proof}

\begin{remark}\label{rem:restrictedset}
It is crucial to recall that some of the inequalities above are only valid when $\dom f=\mathbb{R}^d$---in particular, this holds for Theorem~\ref{thm:ineqs_smooth}(iii \& iv), Theorem~\ref{thm:ineqs_smooth_scvx}(iii\&iv), and Theorem~\ref{thm:interp_proof}. We refer to~\citep{drori2018properties} for an illustration that some inequalities are not valid when restricted on some $\dom f\neq \mathbb{R}^d$. Most standard inequalities, however, do hold even in the case of restricted domains, as established in, e.g.,~\citep{Nest03a}. Some other inequalities, such as Theorem~\ref{thm:ineqs_smooth}(iv) and Theorem~\ref{thm:ineqs_smooth_scvx}(iv), do hold under the additional assumption of twice continuous differentiability(see, for example,~\citep{de2020worst}).
\end{remark}
\let\mysectionmark\sectionmark
\renewcommand\sectionmark[1]{}
\section{Smoothness for General Norms and Restricted Sets}
\let\sectionmark\mysectionmark
\sectionmark{Smoothness for General Norms and Restricted Sets}
\label{app:inequalities_nonEuclid}

In this section, we show that requiring a Lipschitz condition on $\nabla f$, on a convex set $C\subseteq \mathbb{R}^d$, implies a quadratic upper bound on $f$. That is, requiring that for all $x,y\in C$,
\[ \|\nabla f(x)-\nabla f(y)\|_*\leq L\|x-y\|,\]
where $\|.\|$ is some norm and $\|.\|_*$ is the corresponding dual norm, implies a quadratic upper bound $\forall x,y\in C$:
\[ f(x)\leq f(y)+\langle \nabla f(y);x-y\rangle+\frac{L}{2}\|x-y\|^2.\]
\begin{theorem} Let $f:\mathbb{R}^d\rightarrow \mathbb{R}\cup\{+\infty\}$ be continuously differentiable on some open convex set $C\subseteq \mathbb{R}^d$, and let it satisfy a Lipschitz condition
\[ \|\nabla f(x)-\nabla f(y)\|_*\leq L\|x-y\|,\]
for all $x,y\in C$. Then, it holds that
\[ f(x)\leq f(y)+\langle \nabla f(y);x-y\rangle+\frac{L}{2}\|x-y\|^2,\]
for all $x,y\in C$.
\end{theorem}
\begin{proof}
The desired result is obtained from a first-order expansion:
\[ f(y)=f(x)+\int_{0}^1 \langle \nabla f(x+\tau(y-x));y-x\rangle d\tau.\]
The quadratic upper bound then follows from algebraic manipulations and from upper bounding the integral term
\begin{equation*}
\begin{aligned}
f(y)&=f(x)+\langle \nabla f(x);y-x\rangle \\
& \quad +\int_0^1 \langle \nabla f(x+\tau(y-x))-\nabla f(x);y-x\rangle d\tau\\
&\leq f(x)+\langle \nabla f(x);y-x\rangle \\
& \quad +\int_0^1 \|\nabla f(x+\tau(y-x))-\nabla f(x)\|_* \|y-x\| d\tau\\
&\leq f(x)+\langle \nabla f(x);y-x\rangle+ L\|x-y\|^2\int_0^1 \tau d\tau\\
&= f(x)+\langle \nabla f(x);y-x\rangle+ \frac{L}{2}\|x-y\|^2.\qedhere
\end{aligned}
\end{equation*}
\end{proof}

\clearpage
\chapter{Variations on Nesterov Acceleration}\label{a-backandforth} 

\section{Relations between Acceleration Methods}

\subsection{Optimized Gradient Method: Forms I \& II}\label{s:eq_OGM}
In this short section, we show that Algorithm~\ref{alg:OGM_1} and Algorithm~\ref{alg:OGM_2} generate the same sequence $\{y_k\}_k$. A direct consequence of this statement is that the sequences $\{x_k\}_k$ also match, as in both cases they are generated from simple gradient steps on $\{y_k\}_k$.

For this purpose we show that Algorithm~\ref{alg:OGM_2} is a reformulation of Algorithm~\ref{alg:OGM_1}.
\begin{proposition} The sequence $\{y_k\}_k$ generated by Algorithm~\ref{alg:OGM_1} is equal to that generated by Algorithm~\ref{alg:OGM_2}.
\end{proposition}
\begin{proof}
We first observe that the sequences are initiated the same way in both formulations of the OGM. Furthermore, consider one iteration of the OGM in form I:
\[y_k=\left(1-\frac{1}{\theta_{k,N}}\right)x_k+\frac1{\theta_{k,N}}z_k.\]
Therefore, we clearly have $z_k=\theta_{k,N}y_k+(1-\theta_{k,N})x_k$. At the next iteration, we have
\begin{equation*}
\begin{aligned}
y_{k+1}=&\left(1-\frac{1}{\theta_{k+1,N}}\right)x_{k+1}+\frac1{\theta_{k+1,N}}\left(z_k-\frac{2\theta_{k,N}}{L}\nabla f(y_k)\right)\\
=&\left(1-\frac{1}{\theta_{k+1,N}}\right)x_{k+1}\\
& \quad +\frac1{\theta_{k+1,N}}\left(\theta_{k,N}y_k+(1-\theta_{k,N})x_k-\frac{2\theta_{k,N}}{L}\nabla f(y_k)\right),
\end{aligned}
\end{equation*}
where we substituted $z_k$ by its equivalent expression from the previous iteration. Now, by noting that $-\frac1L \nabla f(y_k)=x_{k+1}-y_k$, we reach
\begin{equation*}
\begin{aligned}
y_{k+1}&=\frac{\theta_{k+1,N}-1}{\theta_{k+1,N}}x_{k+1}+\frac1{\theta_{k+1,N}}\left((1-\theta_{k,N})x_k+2\theta_{k,N}x_{k+1}-\theta_{k,N}y_k\right)\\
&=x_{k+1}+\frac{\theta_{k,N}-1}{\theta_{k+1,N}}(x_{k+1}-x_k)+\frac{\theta_{k,N}}{\theta_{k+1,N}}(x_{k+1}-y_k),
\end{aligned}
\end{equation*}
where we reorganized the terms to achieve the same format as in Algorithm~\ref{alg:OGM_2}.
\end{proof}

\subsection{Nesterov's Method: Forms I, II, and III}\label{s:eq_FGM}\label{s:eq_FGM_3}
\begin{proposition} The two sequences $\{x_k\}_k$ and $\{y_k\}_k$ generated by Algorithm~\ref{alg:FGM_1} are equal to those generated by Algorithm~\ref{alg:FGM_2}.
\end{proposition}
\begin{proof} In order to prove the result, we use the identities $A_{k+1}=a_k^2$ as well as $A_k=\sum_{i=0}^{k-1} a_i$, and  $a_{k+1}^2=a_k^2+a_{k+1}$.

Given that the sequences $\{x_k\}_k$ are obtained from gradient steps on $y_k$ in both formulations, it is sufficient to prove that the sequences $\{y_k\}_k$ match. The equivalence is clear for $k=0$, as both methods generate $y_1=x_0-\frac{1}{L}\nabla f(x_0)$. For $k\geq 0$, from Algorithm~\ref{alg:FGM_1}, one can write iteration $k$ as 
\[ y_k = \frac{A_k}{A_{k+1}}x_k + \left(1-\frac{A_k}{A_{k+1}}\right)z_k,\]
and hence,
\begin{equation*}
\begin{aligned}
z_k&=\frac{A_{k+1}}{A_{k+1}-A_{k}} y_k+\left(1-\frac{A_{k+1}}{A_{k+1}-A_{k}}\right)x_k\\
&=a_k y_k+\left(1-a_k\right)x_k.
\end{aligned}
\end{equation*}
Substituting this expression in that for iteration $k+1$, we reach
\begin{equation*}
\begin{aligned}
 y_{k+1}=&\frac{A_{k+1}}{A_{k+2}}x_{k+1}+\frac{A_{k+2}-A_{k+1}}{A_{k+2}}\left(z_k-\frac{A_{k+1}-A_k}{L}\nabla f(y_k)\right)\\
        =&\frac{a_{k}^2}{a_{k+1}^2}x_{k+1}+\frac{1}{a_{k+1}}\left(a_k y_k+\left(1-a_k\right)x_k-\frac{a_k}{L}\nabla f(y_k)\right)\\
        =&\frac{a_{k}^2}{a_{k+1}^2}x_{k+1}+\frac{1}{a_{k+1}}\left(a_k x_{k+1}+\left(1-a_k\right)x_k\right)\\
        =&x_{k+1}+\frac{a_k-1}{a_{k+1}}(x_{k+1}-x_k),
\end{aligned}
\end{equation*}
where we substituted the expression for $z_k$ and used previous identities to reach the desired statement.
\end{proof}

\enlargethispage{\baselineskip}
The same relationship holds with Algorithm~\ref{alg:FGM_3}, as provided by the next proposition.
\begin{proposition} The three sequences $\{z_k\}_k$, $\{x_k\}_k$ and $\{y_k\}_k$ generated by Algorithm~\ref{alg:FGM_1} are equal to those generated by Algorithm~\ref{alg:FGM_3}.
\end{proposition}
\begin{proof}
Clearly, we have $x_0=z_0=y_0$ in both methods. Let us assume that the sequences match up to iteration $k$, that is, up to $y_{k-1}$, $x_k$, and $z_k$. Clearly, both $y_k$ and $z_{k+1}$ are computed in the same way in both methods. It remains to compare the update rules for $x_{k+1}$: in Algorithm~\ref{alg:FGM_3}, we have
\begin{equation*}
\begin{aligned}
x_{k+1}&=\frac{A_k}{A_{k+1}}x_k+\left(1-\frac{A_k}{A_{k+1}}\right)z_{k+1}\\
&=y_k-\left(1-\frac{A_k}{A_{k+1}}\right)\frac{A_{k+1}-A_k}{L}\nabla f(y_k),
\end{aligned}
\end{equation*}
where we used the update rule for $z_{k+1}$. Further simplifications, along with the identity $(A_{k+1}-A_k)^2=A_{k+1}$ allows us to arrive at
\begin{equation*}
\begin{aligned}
x_{k+1}&=y_k-\frac{(A_{k+1}-A_k)^2}{LA_{k+1}}\nabla f(y_k)\\
&=y_k-\frac{1}{L}\nabla f(y_k),
\end{aligned}
\end{equation*}
which is clearly the same update rule as that of Algorithm~\ref{alg:FGM_1}. Hence, all sequences match and the desired statement is proved. 
\end{proof}

\subsection{Nesterov's Accelerated Gradient Method (Strongly Convex Case): Forms I, II, and III}\label{s:eq_FGM_strcvx}
In this short section, we provide alternate, equivalent, formulations for Algorithm~\ref{alg:FGM_1_strconvex}.
\begin{algorithm}[!ht]
  \caption{Nesterov's method, form II}
  \label{alg:FGM_2_strconvex}
  \begin{algorithmic}[1]
    \REQUIRE $L$-smooth $\mu$-strongly convex function $f$ and initial point $x_0$.
    \STATE \textbf{Initialize} $z_0=x_0$; $\cond=\mu/L$, $A_0=0$, and $A_1=(1-\cond)^{-1}$.
    \FOR{$k=0,\ldots$}
      \STATE $A_{k+2}=\frac{2 A_{k+1}+1+\sqrt{4 A_{k+1}+4 \cond A_{k+1}^2 +1}}{2 \left(1-\cond\right)}$
      \STATE $x_{k+1}=y_{k}-\frac1L \nabla f(y_{k})$
      \STATE $y_{k+1}=x_{k+1}+\beta_k(x_{k+1}-x_k)$
      \STATE with $\beta_k=\frac{(A_{k+2}-A_{k+1}) \left(A_{k+1} \left(1-\cond\right)-A_{k}-1\right)}{A_{k+2} \left(2\cond A_{k+1}+1\right)-\cond A_{k+1}^2}$
    \ENDFOR
    \ENSURE Approximate solution $x_{N}$.
  \end{algorithmic}
\end{algorithm}

\begin{proposition} The two sequences $\{x_k\}_k$ and $\{y_k\}_k$ generated by Algorithm~\ref{alg:FGM_1_strconvex} are equal to those generated by Algorithm~\ref{alg:FGM_2_strconvex}.
\end{proposition}
\begin{proof}
Without loss of generality, we can consider that a third sequence $z_k$ is present in Algorithm~\ref{alg:FGM_2_strconvex} (although it is not computed).

Obviously, we have $x_0=z_0=y_0$ in both methods. Let us assume that the sequences match up to iteration $k$, that is, up to $y_{k}$, $x_k$, and $z_k$. Clearly, $x_{k+1}$ is computed in the same way in both methods as a gradient step from $y_k$, and it remains to compare the update rules for $y_{k+1}$. In Algorithm~\ref{alg:FGM_1_strconvex}, we have
\begin{align*}
y_{k+1}=& x_k + \left(\tau_k-\tau_{k+1}(\tau_k-1) (1-\cond  \delta_k)\right)(z_k-x_k)\\
& \quad -\frac{(\delta_k-1) \tau_{k+1}+1}{L}\nabla f(y_k),
\end{align*}
whereas in Algorithm~\ref{alg:FGM_1_strconvex}, we have
\[y_{k+1}=x_k+(\beta_k+1) \tau_k(z_k-x_k)-\frac{1+\beta_k}{L}\nabla f(y_k).\]
By noting that $\beta_k=\tau_{k+1}(\delta_k-1)$, we see that the coefficients in front of $\nabla f(y_k)$ match in both expressions. It remains to check that
\[ (\beta_k+1) \tau_k-\left(\tau_k-\tau_{k+1}(\tau_k-1) (1-\cond  \delta_k)\right)\]
is identically $0$ to reach the desired statement. By substituting $\beta_k=\tau_{k+1}(\delta_k-1)$, this expression reduces to
\[{\tau_{k+1}(\delta_k (\tau_k (1-\cond )+\cond )-1)},\]
and we have to verify that $(\delta_k (\tau_k (1-\cond )+\cond )-1)$ is zero.
Substituting and reworking this expression using the expressions for $\tau_k$, and $\delta_k$, we arrive at
\[ \frac{\tau_{k} \left((A_{k+1}-A_{k})^2-A_{k+1}-\cond  A_{k+1}^2\right)}{(A_{k+1}-A_{k}) (1+\cond  A_{k+1})}=0,\]
as we recognize that $(A_{k+1}-A_{k})^2-A_{k+1}-\cond  A_{k+1}^2=0$ (which is the expression we used to select $A_{k+1}$).
\end{proof}
\begin{algorithm}[!ht]
  \caption{Nesterov's method, form III}
  \label{alg:FGM_3_strconvex}
  \begin{algorithmic}[1]
    \REQUIRE $L$-smooth $\mu$-strongly convex function $f$ and initial point $x_0$.
    \STATE \textbf{Initialize} $z_0=x_0$ and $A_0=0$; $\cond=\mu/L$.
    \FOR{$k=0,\ldots$}
      \STATE $A_{k+1}=\frac{2 A_k+1+\sqrt{4 A_k+4 \cond A_k^2+1}}{2 \left(1-\cond\right)}$
      \STATE set $\tau_k=\frac{(A_{k+1}-A_k) (1+\cond A_k)}{A_{k+1}+2\cond A_k A_{k+1}-\cond A_k^2 }$ and $\delta_{k}=\frac{A_{k+1}-A_{k}}{1+\cond A_{k+1}}$
      \STATE $y_{k}=  x_k+\tau_k (z_k-x_k)$
      \STATE $z_{k+1}=(1-\cond\delta_k)z_k+ \cond\delta_k y_k- \frac{\delta_k}{L}\nabla f(y_{k})$
      \STATE $x_{k+1}=\frac{A_k}{A_{k+1}}x_k+(1-\frac{A_k}{A_{k+1}})z_{k+1}$
    \ENDFOR
    \ENSURE Approximate solution $x_{N}$.
  \end{algorithmic}
\end{algorithm}
\begin{proposition} The three sequences $\{z_k\}_k$, $\{x_k\}_k$, and $\{y_k\}_k$ generated by Algorithm~\ref{alg:FGM_1_strconvex} are equal to those generated by Algorithm~\ref{alg:FGM_3_strconvex}.
\end{proposition}
\begin{proof}
Clearly, we have $x_0=z_0=y_0$ in both methods. Let us assume that the sequences match up to iteration $k$, that is, up to $y_{k-1}$, $x_k$, and $z_k$. Since $y_k$ and $z_{k+1}$ are clearly computed in the same way in both methods, we only have to verify that the update rules for $x_{k+1}$ match. In other words, we have to verify that
\[\frac{A_k}{A_{k+1}}x_k+(1-\frac{A_k}{A_{k+1}})z_{k+1}=y_k-\frac1L\nabla f(y_k),\]
which, using the update rules for $z_{k+1}$ and $y_k$, amounts to verifying that
\[ -\frac{(A_{k+1}-A_{k})^2-A_{k+1}-\cond  A_{k+1}^2}{LA_{k+1}(1+\cond A_{k+1})} \nabla f(y_k)=0.\]
This statement is true since we recognize $(A_{k+1}-A_{k})^2-A_{k+1}-\cond  A_{k+1}^2=0$ as the expression used to select $A_{k+1}$.
\end{proof}

\section{Conjugate Gradient Method}\label{s:conj_methods}
Historically, Nesterov's accelerated gradient method~\citep{Nest83} was preceded by a few other methods with optimal worst-case convergence rates $O(N^{-2})$ for smooth convex minimization. However, the alternate schemes required the capability to optimize exactly over a few dimensions---plane-searches were used in~\citep{Book:NemirovskyYudin,nemirovski1983information} and line-searches were used in~\citep{nemirovski1982orth}; unfortunately these references are not available in English, and we refer to~\citep{narkiss2005sequential} for related discussions. 

In this vein, accelerated methods can be obtained through their links with conjugate gradients (Algorithm~\ref{alg:CG_1}), as a by-product of the worst-case analysis. In this section, we illustrate the absolute perfection of the connection between the OGM and conjugate gradients is absolutely perfect: an identical proof (achieving the lower bound) is valid for both methods.
\begin{algorithm}[!ht]
  \caption{Conjugate gradient method}
  \label{alg:CG_1}
  \begin{algorithmic}[1]
    \REQUIRE
     $L$-smooth convex function $f$, initial point $y_0$, and budget $N$.
    \FOR{$k=0,\ldots,N-1$}
      \STATE $y_{k+1}=\mathrm{argmin}_x\{ f(x):\, x\in y_0+\mathrm{span}\{\nabla f(y_0),\hdots,$ $\nabla f(y_k)\}\}$
    \ENDFOR
    \ENSURE Approximate solution $y_{N}$.
  \end{algorithmic}
\end{algorithm}
The conjugate gradient (CG) method for solving quadratic optimization problems is known to have an efficient form that does not require span-searches (which are in general too expensive to be of any practical interest); see, for example,~\citep{nocedal2006numerical}. Beyond quadratics, it is generally not possible to reformulate the CG method in an efficient way. However, it is possible to find other methods for which the same worst-case analysis applies, and it turns out that the OGM is one of them---see~\citep{drori2019efficient} for details. Similarly, by slightly weakening the analysis of the CG method, one can find other methods, such as Nesterov's accelerated gradient (see Remark~\ref{rem:nest_cg} below for more details).

More precisely, recall the previous definition for the sequence $\{\theta_{k,N}\}_k$, defined in~\eqref{eq:thetas}:
\begin{equation*}
\begin{aligned}
\theta_{k+1,N}= \left\{\begin{array}{ll}
        \frac{1+\sqrt{4\theta_{k,N}^2+1}}{2}  \, & \text{if } k\leq N-2  \\
        \frac{1+\sqrt{8\theta_{k,N}^2+1}}{2}  \, & \text{if } k=N-1.
      \end{array}\right.
\end{aligned}
\end{equation*}
As a result of the worst-case analysis presented below, all methods satisfying 
\begin{equation}\label{eq:super_LS}
\begin{aligned}
 \langle \nabla f(y_i); y_i-&\Bigg[ \left(1-\frac{1}{\theta_{i,N}}\right)\left(y_{i-1}-\tfrac1L \nabla f(y_{i-1})\right)\\
 & \quad +\frac1{\theta_{i,N}}\left(y_0-\tfrac2L \sum_{j=0}^{i-1}\theta_{j,N}\nabla f(y_j) \right) \Bigg]\rangle \leq 0
\end{aligned}
\end{equation}
achieve the optimal worst-case complexity of smooth convex minimization that is provided by Theorem~\ref{thm:smooth_LB}. On the one hand, the CG ensures that this inequality holds thanks to its span-searches (which ensure the orthogonality of successive search directions); that is,
\begin{equation*}
\begin{aligned}
\langle \nabla f(y_i); y_i-y_{i-1}+\frac{1}{\theta_{i,N}}(y_{i-1}-y_0)\rangle&=0\\
\langle \nabla f(y_i); \nabla f(y_0)\rangle &=0\\
&\vdots\\
\langle \nabla f(y_i); \nabla f(y_{i-1})\rangle&=0.
\end{aligned}
\end{equation*}
On the other hand, the OGM enforces this inequality by using
\[ y_i=\left(1-\frac{1}{\theta_{i,N}}\right)\left(y_{i-1}-\tfrac1L \nabla f(y_{i-1})\right)+\frac1{\theta_{i,N}}\left(y_0-\tfrac2L \sum_{j=0}^{i-1}\theta_{j,N}\nabla f(y_j) \right) .\]

\subsubsection{Optimized and Conjugate Gradient Methods: Worst-case Analyses}

The worst-case analysis below relies on the same potentials used for the optimized gradient method; see Theorem~\ref{thm:pot_OGM} and Lemma~\ref{thm:pot_OGM_final}.
\begin{theorem} Let $f$ be an $L$-smooth convex function, $N\in\mathbb{N}$ and some $x_\star\in\mathrm{argmin}_x\, f(x)$. The iterates of the conjugate gradient method (CG, Algorithm~\ref{alg:CG_1}) and of all methods whose iterates are compliant with~\eqref{eq:super_LS} satisfy
\[ f(y_N)-f(x_\star) \leq \frac{L\lVert y_0-x_\star\rVert_2^2}{2\theta_{N,N}^2},\]
for all $y_0\in\mathbb{R}^d$.
\end{theorem}

\bgroup
\addtolength{\jot}{-0.1em}
\begin{proof} 
The result is obtained from the same potential as that used for the OGM, obtained from further inequalities. That is, we first perform a weighted sum of the following inequalities.
\begin{itemize}   
    \item Smoothness and convexity of $f$ between $y_{k-1}$ and $y_k$ with weight $\lambda_1={2\theta_{k-1,N}^2}$:
    \begin{align*}
    0\geq& f(y_k)-f(y_{k-1})+\langle \nabla f(y_k); y_{k-1}-y_k\rangle \\
    & \quad +\frac1{2L}\lVert \nabla f(y_k)-\nabla f(y_{k-1})\rVert_2^2.
    \end{align*}
    \item Smoothness and convexity of $f$ between $x_\star$ and $y_k$ with weight $\lambda_2={2\theta_{k,N}}$:
    \[ 0\geq f(y_k)-f(x_\star)+\langle \nabla f(y_k); x_\star-y_k\rangle +\frac1{2L}\lVert \nabla f(y_k)\rVert_2^2. \]
    \item Search procedure to obtain $y_k$, with weight $\lambda_3=2 \theta_{k,N}^2$:
    \[ 0\geq \langle \nabla f(y_k); y_k-\Bigg[ \left(1-\frac{1}{\theta_{k,N}}\right)\left(y_{k-1}-\tfrac1L \nabla f(y_{k-1})\right)+\frac1{\theta_{k,N}}z_k \Bigg]\rangle ,\]
 where we used $z_k:=y_0-\tfrac2L \sum_{j=0}^{k-1}\theta_{j,N}\nabla f(y_j)$.
\end{itemize}
The weighted sum is a valid inequality:
\begin{equation*}
\begin{aligned}
0\geq& \lambda_1[f(y_k)-f(y_{k-1})+\langle \nabla f(y_k); y_{k-1}-y_k\rangle \\
& \quad \quad +\frac1{2L}\lVert \nabla f(y_k)-\nabla f(y_{k-1})\rVert_2^2]\\
&+\lambda_2[ f(y_k)-f(x_\star)+\langle \nabla f(y_k); x_\star-y_k\rangle +\frac1{2L}\lVert \nabla f(y_k)\rVert_2^2]\\
&+\lambda_3[\langle \nabla f(y_k); y_k-\Bigg[ \left(1-\frac{1}{\theta_{k,N}}\right)\left(y_{k-1}-\tfrac1L \nabla f(y_{k-1})\right)\\
& \quad \quad\quad\quad\quad\quad\quad\quad\quad +\frac1{\theta_{k,N}}z_k \Bigg]\rangle].
\end{aligned}
\end{equation*} Substituting $z_{k+1}$, the previous inequality can be reformulated exactly as
\begin{equation*}
\begin{aligned}
0\geq& 2\theta_{k,N}^2\left(f(y_{k})-f_\star-\frac{1}{2L}\lVert \nabla f(y_{k})\rVert_2^2\right)+\frac{L}{2}\|z_{k+1}-x_\star\|_2^2\\
&-2\theta_{k-1,N}^2\left(f(y_{k-1})-f_\star-\frac1{2L}\lVert \nabla f(y_{k-1})\rVert_2^2\right)-\frac{L}{2}\|z_{k}-x_\star\|_2^2\\
&+2 \left(\theta_{k-1,N}^2-\theta_{k,N}^2+\theta_{k,N}\right) \left(f(y_k)-f_\star+\frac1{2L}\|\nabla f(y_k)\|_2^2\right)\\
&+2 \left(\theta_{k-1,N}^2-\theta_{k,N}^2+\theta_{k,N}\right) \langle \nabla f(y_k); y_{k-1}-\tfrac{1}{L}\nabla f(y_{k-1})-y_k\rangle.
\end{aligned}
\end{equation*}
We reach the desired inequality by selecting $\theta_{k,N}$ that satisfies $\theta_{k,N}\geq\theta_{k-1,N}$ and \[\theta_{k-1,N}^2-\theta_{k,N}^2+\theta_{k,N}=0,\]
thereby reaching the same potential as in Theorem~\ref{thm:pot_OGM}.

To obtain the technical lemma that allows us to bound the final $f(y_N)-f_\star$, we follow the same steps with the following inequalities.
\begin{itemize}
    \item Smoothness and convexity of $f$ between $y_{k-1}$ and $y_k$ with weight $\lambda_1={2\theta_{N-1,N}^2}$:
    \begin{align*}
    0\geq& f(y_N)-f(y_{N-1})+\langle \nabla f(y_N); y_{N-1}-y_N\rangle \\
    & \quad +\frac1{2L}\lVert \nabla f(y_N)-\nabla f(y_{N-1})\rVert_2^2.
    \end{align*}
    \item Smoothness and convexity of $f$ between $x_\star$ and $y_k$ with weight $\lambda_2=\theta_{N,N}$:
    \[ 0\geq f(y_N)-f(x_\star)+\langle \nabla f(y_N); x_\star-y_N\rangle +\frac1{2L}\lVert \nabla f(y_N)\rVert_2^2. \]
    \item Search procedure to obtain $y_N$, with weight $\lambda_3=\theta_{N,N}^2$:
    \begin{align*}
    0\geq \langle \nabla f(y_N); y_N-\Bigg[& \left(1-\frac{1}{\theta_{N,N}}\right)\left(y_{N-1}-\tfrac1L \nabla f(y_{N-1})\right)\\
    & \quad +\frac1{\theta_{N,N}}z_N \Bigg]\rangle.
    \end{align*}
    \end{itemize}
The weighted sum can then be reformulated as:
\begin{equation*}
\begin{aligned}
0\geq &\theta_{N,N}^2 \left(f(y_N)-f_\star\right)+\frac{L}{2}\|z_{N}-\tfrac{\theta_{N,N}}{L}\nabla f(y_N)-x_\star\|_2^2\\
&-2\theta_{N-1,N}^2\left(f(y_{N-1})-f_\star-\frac1{2L}\lVert \nabla f(y_{N-1})\rVert_2^2\right)-\frac{L}{2}\|z_{N}-x_\star\|_2^2\\
&+\left(2 \theta_{N-1,N}^2-\theta_{N,N}^2+\theta_{N,N}\right)\left(f(y_N)-f_\star+\frac{1}{2 L}\|\nabla f(y_N)\|_2^2\right)\\
&+\left(2 \theta_{N-1,N}^2-\theta_{N,N}^2+\theta_{N,N}\right) \langle \nabla f(y_N);y_{N-1}-\tfrac{1}{L}\nabla f(y_{N-1})-y_{N}\rangle,
\end{aligned}
\end{equation*}
thus reaching the desired inequality, as in Lemma~\ref{thm:pot_OGM_final}, by selecting $\theta_{N,N}$ that satisfies $\theta_{N,N}\geq \theta_{N-1,N}$ and
\[2 \theta_{N-1,N}^2-\theta_{N,N}^2+\theta_{N,N}.\]
Hence, the potential argument from Corollary~\ref{cor:OGM_bound} applies as such, and we reach the desired conclusion. In other words, for all $k\in\{0,\hdots,N\}$, one can define
\[ \phi_k\defeq 2\theta_{k-1,N}^2\left(f(y_{k-1})-f_\star-\frac{1}{2L}\lVert \nabla f(y_{k-1})\rVert_2^2\right)+\frac{L}{2}\|z_{k}-x_\star\|_2^2\]
and \[ \phi_{N+1}\defeq \theta_{N,N}^2\left(f(y_N)-f_\star\right)+\frac{L}{2}\|z_{N}-\tfrac{\theta_{N,N}}{L}\nabla f(y_N)-x_\star\|_2^2\]
and reach the desired statement by chaining the inequalities:\[\theta_{N,N}^2(f(y_N)-f_\star)\leq \phi_{N+1}\leq\phi_{N}\leq\hdots\leq \phi_0=\frac{L}2\|y_0-x_\star\|_2^2.\qedhere\]
\end{proof}
\egroup

\begin{remark}\label{rem:nest_cg}It is possible to further exploit the conjugate gradient method to design practical accelerated methods in different settings, such as that of~\citet{Nest83}. This point of view has been exploited in~\citep{narkiss2005sequential,karimi2016unified,karimi2017single,diakonikolas2019conjugate}, among others. The link between the CG method and the OGM presented in this section is due to~\citet{drori2019efficient}, though with a different presentation that does not involve the potential function.
\end{remark}

\let\mysectionmark\sectionmark
\renewcommand\sectionmark[1]{}
\section{Acceleration Without Monotone Backtracking}
\let\sectionmark\mysectionmark
\sectionmark{Acceleration Without Monotone Backtracking}

\subsection{FISTA without Monotone Backtracking}\label{s:Fista_nonmonotone_BT}

In this section, we show how to incorporate backtracking strategies that may not satisfy $L_{k+1}\geq L_k$, which is important in practice. The developments are essentially the same; one possible trick is to incorporate all the knowledge about $L_k$ in $A_k$. That is, we use a rescaled shape for the potential function:
\[ \phi_k\defeq B_k (f(x_k)-f_\star)+\frac{1+\mu B_k}{2}\|z_k-x_\star\|_2^2,\]
where without the backtracking strategy, $B_k=\tfrac{A_k}{L}$. This seemingly cosmetic change allows $\phi_k$ to depend on $L_k$ solely via $B_k$, and it applies to both backtracking methods presented in Section~\ref{c-Nest} (Section~\ref{s:backtracking}). 

The idea used to obtain both methods below is that one can perform the same computations as in Algorithm~\ref{alg:FGM_1_strconvex}, replacing $A_k$ by $L_{k+1}B_k$ and $A_{k+1}$ by $L_{k+1}A_{k+1}$ at iteration $k$. Thus, as in previous versions, only the current approximate Lipschitz constant $L_{k+1}$ is used at iteration $k$: previous approximations were only used to compute $B_k$.

\begin{algorithm}[!ht]
  \caption{Strongly convex FISTA (general initialization of $L_{k+1}$)}
  \label{alg:fista_1_general_BT}
  \begin{algorithmic}[1]
    \REQUIRE An $L$-smooth (possibly $\mu$-strongly) convex function $f$, a convex function $h$ with proximal operator available, an initial point $x_0$, and an initial estimate $L_0>\mu$.
    \STATE \textbf{Initialize} $z_0=x_0$, $B_0=0$, and some $\alpha>1$.
    \FOR{$k=0,\ldots$}
        \STATE Pick $L_{k+1}\in[L_0,L_k]$.
        \LOOP 
            \STATE set $\cond_{k+1}=\mu/L_{k+1}$,
            \STATE $B_{k+1}=\frac{2 L_{k+1}B_k+1+\sqrt{4 L_{k+1}B_k+4\mu L_{k+1}B_k^2 +1}}{2 \left(L_{k+1}-\mu\right)}$
            \STATE set $\tau_k=\frac{(B_{k+1}-B_k) (1+\mu B_k )}{(B_{k+1}+2\mu B_k B_{k+1}-\mu B_k^2)}$ and $\delta_k=L_{k+1}\frac{B_{k+1}-B_{k}}{1+\mu B_{k+1}}$
            \STATE $y_{k}= x_k+\tau_k(z_k-x_k)$
            \STATE $x_{k+1}=\mathrm{prox}_{h/L_{k+1}}\left(y_{k}-\tfrac{1}{L_{k+1}}\nabla f(y_{k})\right)$
            \STATE $z_{k+1}=(1-\cond_{k+1}\delta_k)z_k+\cond_{k+1}\delta_k y_k+\delta_k \left(x_{k+1}-y_k\right)$
            \IF{\eqref{eq:backtracking_cond} holds}
                \STATE \textbf{break} \COMMENT{Iterates accepted; $k$ will be incremented.}
            \ELSE
                \STATE $L_{k+1}=\alpha L_{k+1}$ \COMMENT{Iterates  not accepted; compute new $L_{k+1}$.}
            \ENDIF 
        \ENDLOOP
    \ENDFOR
    \ENSURE Approximate solution $x_{k+1}$.
  \end{algorithmic}
\end{algorithm}

\enlargethispage{-\baselineskip}
The proof follows the same lines as used for FISTA (Algorithm~\ref{thm:FISTA_strcvx}). In this case, $f$ is assumed to be smooth and convex over $\mathbb{R}^d$ (i.e., it has full domain, $\dom f=\mathbb{R}^d$), and we are therefore allowed to evaluate gradients of $f$ outside of the domain of $h$.

\begin{theorem}\label{thm:FISTA_strcvx_general_BT} Let $f\in\mathcal{F}_{\mu,L}$ (with full domain, $\dom f=\mathbb{R}^d$), $h$ be a closed convex proper function, $x_\star \in\mathrm{argmin}_x\, \{F(x)\defeq f(x)+h(x)\}$, and $k\in\mathbb{N}$. For any $x_k,z_k\in\mathbb{R}^d$ and $B_k\geq 0$, the iterates of Algorithm~\ref{alg:fista_1_general_BT} that satisfy~\eqref{eq:backtracking_cond} also satisfy
\begin{equation*}
\begin{aligned}
B_{k+1}&(F(x_{k+1})-F_\star)+\frac{1+\mu B_{k+1}}{2}\|z_{k+1}-x_\star\|_2^2\\&\leq B_k (F(x_k)-F_\star)+\frac{1+\mu B_{k}}2 \|z_k-x_\star\|_2^2,
\end{aligned}
\end{equation*}
with $B_{k+1}=\frac{2 L_{k+1}B_k+1+\sqrt{4L_{k+1} B_k+4 \mu L_{k+1}B_k^2 +1}}{2 \left(L_{k+1}-\mu\right)}$.
\end{theorem}\begin{proof}
The proof consists of a weighted sum of the following inequalities.
\begin{itemize}
    \item Strong convexity of $f$ between $x_\star$ and $y_{k}$ with weight $\lambda_1=B_{k+1}-B_k$:
    \[ f_\star \geq f(y_{k})+\langle \nabla f(y_{k});x_\star-y_k\rangle+\frac{\mu}{2}\|x_\star-y_k\|_2^2.\]
    \item Strong convexity of $f$ between $x_k$ and $y_{k}$ with weight $\lambda_2=B_k$:
    \[ f(x_k) \geq f(y_{k})+\langle \nabla f(y_{k});x_k-y_{k}\rangle.\]
    \item Smoothness of $f$ between $y_{k}$ and $x_{k+1}$ (\emph{descent lemma}) with weight $\lambda_3=B_{k+1}$:
    \[ f(y_{k}) + \langle \nabla f(y_{k}); x_{k+1}-y_{k}\rangle +\frac{L_{k+1}}{2}\|x_{k+1}-y_{k}\|_2^2 \geq f(x_{k+1}).\]
    \item Convexity of $h$ between $x_\star$ and $x_{k+1}$ with weight $\lambda_4=B_{k+1}-$ $B_k$:
    \[ h(x_\star)\geq h(x_{k+1})+\langle g_h(x_{k+1});x_\star-x_{k+1}\rangle,\]
    with $g_h(x_{k+1})\in\partial h(x_{k+1})$ and $x_{k+1}=y_k-\tfrac1{L_{k+1}} (\nabla f(y_k)+g_h$ $(x_{k+1}))$.
    \item Convexity of $h$ between $x_k$ and $x_{k+1}$ with weight $\lambda_5=B_k$:
    \[ h(x_k)\geq h(x_{k+1})+\langle g_h(x_{k+1});x_k-x_{k+1}\rangle.\]
\end{itemize}
We obtain the following inequality:
\begin{equation*}
\begin{aligned}
0\geq &\lambda_1 [f(y_{k})-f_\star+\langle \nabla f(y_{k});x_\star-y_k\rangle+\frac{\mu}{2}\|x_\star-y_k\|_2^2]\\&+\lambda_2[f(y_{k})-f(x_k)+\langle \nabla f(y_{k});x_k-y_{k}\rangle]\\
&+\lambda_3[f(x_{k+1})-(f(y_{k}) + \langle \nabla f(y_{k}); x_{k+1}-y_{k}\rangle \\
& \quad +\frac{L_{k+1}}{2}\|x_{k+1}-y_{k}\|_2^2)]\\
&+\lambda_4[ h(x_{k+1})-h(x_\star)+\langle g_h(x_{k+1});x_\star-x_{k+1}\rangle]\\&+\lambda_5[ h(x_{k+1})-h(x_k)+\langle g_h(x_{k+1});x_k-x_{k+1}\rangle].
\end{aligned}
\end{equation*}
Substituting the $y_k$, $x_{k+1}$, and $z_{k+1}$ with
\begin{equation*}
\begin{aligned}
y_k&=x_k+\tau_k  (z_k-x_k)\\
x_{k+1}&=y_k-\frac{1}{L_{k+1}}(\nabla f(y_k)+g_h(x_{k+1})) \\
z_{k+1}&=(1-\cond_{k+1}\delta_k)z_k+\cond_{k+1}\delta_k y_k+\delta_k \left(x_{k+1}-y_k\right),
\end{aligned}
\end{equation*}
after some basic but tedious algebra, yields
\begin{equation*}
\begin{aligned}
B_{k+1}&(f(x_{k+1})+h(x_{k+1})-f(x_\star)-h(x_\star))+\frac{1+B_{k+1}\mu}{2}\|z_{k+1}-x_\star\|_2^2\\
\leq& B_k (f(x_k)+h(x_k)-f(x_\star)-h(x_\star))+\frac{1+B_k\mu}{2}\|z_k-x_\star\|_2^2\\
&+\frac{L_{k+1}(B_{k}-B_{k+1})^2-B_{k+1}-\mu  B_{k+1}^2}{1+\mu B_{k+1}}\\
& \quad\quad\times \frac1{2L_{k+1}}\|\nabla f(y_k)+g_h(x_{k+1})\|_2^2\\
&-\frac{B_{k}^2 (B_{k+1}-B_{k}) (1+\mu  B_{k}) (1+\mu  B_{k+1})}{\left( B_{k+1} +2 \mu B_{k}B_{k+1}-\mu B_{k}^2\right)^2}\frac{\mu}{2}\|x_k-z_k\|_2^2.
\end{aligned}
\end{equation*}
Then, choosing $B_{k+1}$ such that $B_{k+1}\geq B_k$ and
\[ L_{k+1}(B_{k}-B_{k+1})^2-B_{k+1}-\mu  B_{k+1}^2 = 0,\]
yields the desired result:
\begin{equation*}
\begin{aligned}
B_{k+1}&(f(x_{k+1})+h(x_{k+1})-f(x_\star)-h(x_\star))+\frac{1+B_{k+1}\mu}{2}\|z_{k+1}-x_\star\|_2^2\\
\leq& B_k (f(x_k)+h(x_k)-f(x_\star)-h(x_\star))+\frac{1+B_k\mu}{2}\|z_k-x_\star\|_2^2.\qedhere
\end{aligned}
\end{equation*}
\end{proof}

\enlargethispage{\baselineskip}
Finally, we obtain a complexity guarantee by adapting the potential argument~\eqref{eq:chained_pot} and by noting that $B_{k+1}$ is a decreasing function of $L_{k+1}$ (whose maximal value is $\alpha L$, assuming $L_0<L$; otherwise, its maximal value is $L_0$). The growth rate of $B_k$ in the smooth convex setting remains unchanged (see~\eqref{eq:conv_Ak_fgm}) since we have
\[B_{k+1}\geq \frac{\left(\frac{1}{2}+\sqrt{B_k L_{k+1}}\right)^2}{L_{k+1}}, \]
and hence, $\sqrt{B_{k+1}}\geq \tfrac{1}{2\sqrt{L_{k+1}}}+\sqrt{B_k}$. Therefore, $B_k\geq \left(\tfrac{k}{2\sqrt{\ell}}\right)^2$ with $\ell=\max\{L_0,\alpha L\}$ and $L_{k+1}\leq \ell$. As for the geometric rate, we similarly obtain
\[ B_{k+1}\geq B_k\frac{ \left(1+\sqrt{\frac{\mu }{L_{k+1}}}\right)}{1-\frac{\mu }{L_{k+1}}}=\frac{B_k}{1-\sqrt{\frac{\mu }{L_{k+1}}}},\]
and therefore, $B_{k+1}\geq (1-\sqrt{\tfrac{\mu}{\ell}})^{-1}B_k$.
\begin{corollary}\label{cor:FISTA_general_BT} Let $f\in\mathcal{F}_{\mu,L}(\mathbb{R}^d)$ (with full domain, $\dom f=\mathbb{R}^d$), $h$ be a closed convex proper function and $x_\star\in\mathrm{argmin}_x\, \{F(x)\defeq f(x)+h(x)\}$. For any $N\in\mathbb{N}$, $N\geq 1$, and $x_0\in\mathbb{R}^d$, the output of Algorithm~\ref{alg:fista_1_general_BT} satisfies
\[F(x_N)-F_\star\leq \min\left\{\frac{2}{N^2},\left(1-\sqrt{\frac{\mu}{\ell}}\right)^{N}\right\}\ell\|x_0-x_\star\|_2^2,\]
with $\ell=\max\{\alpha L,L_0\}$.
\end{corollary}
\begin{proof} We assume that $L>L_0$ since otherwise, $f\in\mathcal{F}_{\mu,L_0}$ and the proof directly follows from the case without backtracking. The chained potential argument~\eqref{eq:chained_pot} can be used as before. Using $B_0=0$, we reach
\[F(x_N)-F_\star\leq \frac{\|x_0-x_\star\|_2^2}{2 B_N}.\]
Our previous bounds on $B_N$ yields the desired result, using \[B_1=\frac{1}{L_{k+1}-\mu}\geq \frac{2\ell^{-1}}{1-\tfrac{\mu}{\ell}}=\frac{2\ell^{-1}}{\left(1-\sqrt{\tfrac{\mu}{\ell}}\right)\left(1+\sqrt{\tfrac{\mu}{\ell}}\right)}\geq \frac{\ell^{-1}}{1-\sqrt{\tfrac{\mu}{\ell}}},\]
and hence, $B_N\geq \ell^{-1} \left(1-\sqrt{\tfrac{\mu}{\ell}}\right)^{-N}$ as well as $B_k\geq \left(\tfrac{k}{2\sqrt{\ell}}\right)^2$.
\end{proof}

\subsection[Another Method without Monotone Backtracking]{Another Accelerated Method without Monotone \\ Backtracking}\label{s:Acc_nonmonotone_BT}
Just as for FISTA, we can perform the same cosmetic change to Algorithm~\ref{alg:MST_strconvex} for incorporating a non-monotonic estimations of the Lipschitz constant.  The proof is therefore essentially that of Algorithm~\ref{alg:MST_strconvex}.
\begin{algorithm}[!ht]
  \caption{A proximal accelerated gradient (general initialization of $L_{k+1}$)}
  \label{alg:MST_strconvex_general_BT}
  \begin{algorithmic}[1]
    \REQUIRE $h\in\mathcal{F}_{0,\infty}$ with proximal operator available, $f\in\mathcal{F}_{\mu,L}(\dom h)$, an initial point $x_0\in\dom h$, and an initial estimate $L_0>\mu$.
    \STATE \textbf{Initialize} $z_0=x_0$, $A_0=0$, and some $\alpha>1$.
    \FOR{$k=0,\ldots$}
        \STATE Pick $L_{k+1}\in [L_0,L_k]$.
        \LOOP
            \STATE Set $\cond_{k+1}=\mu/L_{k+1}$,
            \STATE $B_{k+1}=\frac{2 L_{k+1}B_k+1+\sqrt{4 L_{k+1}B_k+4\mu L_{k+1}B_k^2 +1}}{2 \left(L_{k+1}-\mu\right)}$
            \STATE Set $\tau_k=\frac{L_{k+1}(B_{k+1}-B_k) (1+\mu B_k )}{L_{k+1}(B_{k+1}+2\mu B_k B_{k+1}-\mu B_k^2)}$ and $\delta_k=L_{k+1}\frac{B_{k+1}-B_{k}}{1+\mu B_{k+1}}$
            \STATE $y_{k}=  x_k+\tau_k (z_k-x_k)$
            \STATE $z_{k+1}=\mathrm{prox}_{\delta_k h/L_{k+1}}\bigg((1-\cond_{k+1}\delta_k)z_k+ \cond_{k+1}\delta_k y_k- \frac{\delta_k}{L_{k+1}}\nabla$ $f(y_{k})\bigg)$
            \STATE $x_{k+1}=\frac{A_k}{A_{k+1}}x_k+(1-\frac{A_k}{A_{k+1}})z_{k+1}$
            \IF{\eqref{eq:backtracking_cond} holds}
                \STATE \textbf{break} \COMMENT{Iterates accepted; $k$ will be incremented.}
            \ELSE
                \STATE $L_{k+1}=\alpha L_{k+1}$ \COMMENT{Iterates  not accepted; compute new $L_{k+1}$.}
            \ENDIF 
        \ENDLOOP
    \ENDFOR
    \ENSURE An approximate solution $x_{k+1}$.
  \end{algorithmic}
\end{algorithm}

\begin{theorem}\label{thm:MST_strcvx_general_BT} Let $h\in\mathcal{F}_{0,\infty}$, $f\in\mathcal{F}_{\mu,L}(\dom h)$, $x_\star\in\mathrm{argmin}_x\, \{F(x)\defeq f(x)+h(x)\}$, and $k\in\mathbb{N}$. For any $x_{k},z_k\in\mathbb{R}^d$ and $B_k\geq0$, the iterates of Algorithm~\ref{alg:MST_strconvex_general_BT} that satisfy~\eqref{eq:backtracking_cond} also satisfy
\begin{equation*}
\begin{aligned}
B_{k+1}&(F(x_{k+1})-F_\star)+\frac{1+\mu B_{k+1}}{2}\|z_{k+1}-x_\star\|_2^2\\&\leq B_k (F(x_k)-F_\star)+\frac{1+\mu B_{k}}2 \|z_k-x_\star\|_2^2,
\end{aligned}
\end{equation*}
with $B_{k+1}=\frac{2 L_{k+1}B_k+1+\sqrt{4 L_{k+1}B_k+4\mu L_{k+1}B_k^2 +1}}{2 \left(L_{k+1}-\mu\right)}$.
\end{theorem}
\begin{proof}  {First, $\{z_k\}_k$ is in $\dom h$ by construction---it is the output of a proximal/projection step. Furthermore, we have $0\leq \tfrac{B_{k}}{B_{k+1}}\leq 1$ given that $B_{k+1}\geq B_k\geq 0$. A direct consequence is that since $z_0=x_0\in\dom h$, all subsequent $\{y_k\}_k$ and $\{x_k\}_k$ are also in $\dom h$ (as they are obtained from convex combinations of feasible points).}

The rest of the proof consists of a weighted sum of the following inequalities (which are valid due to the feasibility of the iterates).
\begin{itemize}
    \item Strong convexity of $f$ between $x_\star$ and $y_{k}$ with weight $\lambda_1=B_{k+1}-B_k$:
    \[ f(x_\star) \geq f(y_{k})+\langle \nabla f(y_{k});x_\star-y_k\rangle+\frac{\mu}{2}\|x_\star-y_k\|_2^2.\]
    \item Convexity of $f$ between $x_k$ and $y_{k}$ with weight $\lambda_2=B_k$: 
    \[ f(x_k) \geq f(y_{k})+\langle \nabla f(y_{k});x_k-y_{k}\rangle.\]
    \item Smoothness of $f$ between $y_{k}$ and $x_{k+1}$ (\emph{descent lemma}) with weight $\lambda_3=B_{k+1}$:
    \[ f(y_{k}) + \langle \nabla f(y_{k}); x_{k+1}-y_{k}\rangle +\frac{L_{k+1}}{2}\|x_{k+1}-y_{k}\|_2^2 \geq f(x_{k+1}).\]
    \item Convexity of $h$ between $x_\star$ and $z_{k+1}$ with weight $\lambda_4=B_{k+1}-B_k$:
    \[ h(x_\star)\geq h(z_{k+1})+\langle g_h(z_{k+1});x_\star-z_{k+1}\rangle,\]
    with $g_h(z_{k+1})\in\partial h(z_{k+1})$ and $z_{k+1}=(1-\cond\delta_k)z_k+ \cond\delta_k y_k- \frac{\delta_k}{L_{k+1}}(\nabla f(y_{k})+g_h(z_{k+1}))$.
    \item Convexity of $h$ between $x_k$ and $x_{k+1}$ with weight $\lambda_5=B_k$:
    \[ h(x_k)\geq h(x_{k+1})+\langle g_h(x_{k+1});x_k-x_{k+1}\rangle,\]
    with $g_h(x_{k+1})\in\partial h(x_{k+1})$.
    \item Convexity of $h$ between $z_{k+1}$ and $x_{k+1}$ with weight $\lambda_6=B_{k+1}-B_k$:
    \[ h(z_{k+1})\geq h(x_{k+1})+\langle g_h(x_{k+1});z_{k+1}-x_{k+1}\rangle.\]
\end{itemize}
We obtain the following inequality:
\begin{equation*}
\begin{aligned}
0\geq &\lambda_1 [f(y_{k})-f_\star+\langle \nabla f(y_{k});x_\star-y_k\rangle+\frac{\mu}{2}\|x_\star-y_k\|_2^2]\\
&+\lambda_2[f(y_{k})-f(x_k)+\langle \nabla f(y_{k});x_k-y_{k}\rangle]\\
&+\lambda_3[f(x_{k+1})-(f(y_{k}) + \langle \nabla f(y_{k}); x_{k+1}-y_{k}\rangle \\
& \quad +\frac{L_{k+1}}{2}\|x_{k+1}-y_{k}\|_2^2)]\\
&+\lambda_4[ h(z_{k+1})-h(x_\star)+\langle g_h(z_{k+1});x_\star-z_{k+1}\rangle]\\
&+\lambda_5[ h(x_{k+1})-h(x_k)+\langle g_h(x_{k+1});x_k-x_{k+1}\rangle]\\
&+\lambda_6[ h(x_{k+1})-h(z_{k+1})+\langle g_h(x_{k+1});z_{k+1}-x_{k+1}\rangle].
\end{aligned}
\end{equation*}
Substituting the $y_k$, $z_{k+1}$, and $x_{k+1}$ by
\begin{equation*}
\begin{aligned}
y_k&=x_k+\tau_k  (z_k-x_k)\\
z_{k+1}&= (1-\cond_{k+1} \delta_k )z_k + \cond_{k+1} \delta_k y_k - \frac{\delta_k}{L_{k+1}} (\nabla f(y_k) + g_h(z_{k+1}))\\
x_{k+1}&= \frac{B_k}{B_{k+1}} x_k + \left(1 - \frac{B_k}{B_{k+1}}\right) z_{k+1},
\end{aligned}
\end{equation*}
and algebra allows us to obtain the following reformulation:
\begin{equation*}
\begin{aligned}
B_{k+1}&(f(x_{k+1})+h(x_{k+1})-f(x_\star)-h(x_\star))+\frac{1+\mu B_{k+1}}{2}\|z_{k+1}-x_\star\|_2^2\\
\leq& B_k (f(x_k)+h(x_k)-f(x_\star)-h(x_\star))+\frac{1+\mu B_k}{2}\|z_k-x_\star\|_2^2\\
&+\frac{(B_k-B_{k+1})^2 \left(L_{k+1}(B_{k}-B_{k+1})^2- B_{k+1}-\mu  B_{k+1}^2\right)}{ B_{k+1} (1+\mu B_{k+1})^2} \\
& \quad \quad \quad\times \frac{1}{2}\|\nabla f(y_k)+g_h(z_{k+1})\|_2^2\\
&-\frac{B_k^2 (B_{k+1}-B_k) (1+ \mu B_k) (1+\mu B_{k+1})}{\left(B_{k+1}+2 \mu  B_k B_{k+1}-\mu  B_k^2\right)^2} \frac{\mu}{2}\|x_k-z_k\|_2^2.
\end{aligned}
\end{equation*}
The desired inequality follows from selecting $B_{k+1}$ such that $B_{k+1}\geq B_k$ and
\[ L_{k+1}(B_{k}-B_{k+1})^2-B_{k+1}-\mu  B_{k+1}^2 = 0,\]
thereby yielding
\begin{equation*}
\begin{aligned}
B_{k+1}&(f(x_{k+1})+h(x_{k+1})-f(x_\star)-h(x_\star))+\frac{1+\mu B_{k+1}}{2}\|z_{k+1}-x_\star\|_2^2\\
\leq& B_k (f(x_k)+h(x_k)-f(x_\star)-h(x_\star))+\frac{1+B_k\mu}{2}\|z_k-x_\star\|_2^2.\qedhere
\end{aligned}
\end{equation*}
\end{proof}

The final corollary follows from the same arguments as those used for Corollary~\ref{cor:FISTA_general_BT}. It provides the final bound for Algorithm~\ref{alg:MST_strconvex_general_BT}.

\begin{corollary} Let $h\in\mathcal{F}_{0,\infty}$, $f\in\mathcal{F}_{\mu,L}(\dom h)$, and $x_\star\in\mathrm{argmin}_x\, \{F(x)$ $\defeq f(x)+h(x)\}$. For any $N\in\mathbb{N}$, $N\geq 1$, and $x_0\in\mathbb{R}^d$, the output of Algorithm~\ref{alg:MST_strconvex_general_BT} satisfies
\[F(x_N)-F_\star\leq \min\left\{\frac{2}{N^2},\left(1-\sqrt{\frac{\mu}{\ell}}\right)^{-N}\right\}\ell\|x_0-x_\star\|_2^2,\]
with $\ell=\max\{\alpha L,L_0\}$.
\end{corollary}
\begin{proof} The proof follows the same arguments as those for Corollary~\ref{cor:FISTA_general_BT}, using the potential from Theorem~\ref{thm:MST_strcvx_general_BT} and the fact that the output of the algorithm satisfies~\eqref{eq:backtracking_cond}.
\end{proof}

\clearpage
\chapter{On Worst-case Analyses for First-order Methods}\label{a-WC_FO}

\section{Principled Approaches to Worst-case Analyses}

In this section, we show that obtaining convergence rates and proofs can be framed as finding feasible points to certain convex problems. More precisely, all convergence guarantees from Section~\ref{c-Nest} and Section~\ref{c-prox} can be obtained as feasible points to certain linear matrix inequalities (LMI). As we see in what follows, this approach can be seen as a \emph{principled} approach to worst-case analysis of first-order methods: the approach fails only when no such guarantees can be found. The purpose of this section is to provide complete examples of the LMIs for a few cases of interest: analyses of gradient and accelerated gradient methods, as well as pointers to the relevant literature. We provide a full derivation for the base case, and leave advanced ones as exercises for the reader. Notebooks for obtaining the corresponding LMIs are provided in Section~\ref{s:ref}.

\enlargethispage{\baselineskip}
The elements of this section are largely inspired by the presentation of Taylor and Bach~(\citeyear{taylor19bach}) with elements borrowed from the presentation of Taylor, Hendrickx and Glineur~(\citeyear{taylor2017smooth}), which is itself largely inspired by that of Drori and Teboulle~(\citeyear{Dror14}). The arguments are also similar to the line of work by Lessard, Recht and Packard~(\citeyear{lessard2016analysis}) and follow-up works, see, e.g.,~\citep{fazlyab2018analysis,hu2017dissipativity}. The latter line of works is similar in spirit to the former, but framed in control-theoretic terms, via so-called \emph{integral quadratic constraints}, popularized by~\citet{megretski1997system}.

These techniques are analogous and mostly differs in their presentation styles. Roughly speaking, they can be seen as \emph{dual} to each others. That is, whereas the \emph{performance estimation} viewpoint stems from the problem of computing worst-case scenarios and approaches worst-case guarantees as feasible point to the corresponding dual problems, the \emph{integral quadratic constraint} approach directly starts from the problem of performing linear combination of inequalities, which is exactly the dual problem to that of computing worst-case scenarios. Depending on the background of the researchers involved in a work on one of those topics, things might therefore be named in different ways. We insist on the fact that those are really two facets of the same coin with only subtle differences in terms of presentations.

We choose to take the performance estimation viewpoint as using the definition of a ``worst-case'' allows to carefully select the most appropriate set of inequalities to be used. Informally, this advantageous construction allows certifying the approach to provide meaningful worst-case guarantees: either the approach provides a satisfying worst-case guarantee, or there exists a non-satisfying counterexample, invalidating the existence of any satisfying guarantee of the desired form. 

Further discussions and a more thorough list of references are provided in Section~\ref{s:ref}. Readability in mind, the presentation focuses on some examples of interest rather than on a general framework. We refer to~\citep{Dror14,taylor2017exact,taylor2017smooth} for more details.

\section{Worst-case Analysis as Optimization/Feasibility Problems}

In this section, we provide examples illustrating the type of problems that can be used for obtaining worst-case guarantees. The base idea underlying the technique is that worst-case scenarios are by definition solutions to certain optimization problems. In the context of first-order convex optimization methods, those worst-case scenarios correspond to solutions to linear semidefinite programs (SDP), which are convex; see, e.g.,~\citep{vandenberghe1999applications}. It nicely follows from this theory that any worst-case guarantee (i.e., any upper bound on a worst-case performance) can be formulated as a feasible point to the dual problem to that of finding worst-case scenarios. Equivalently, those dual solutions correspond to appropriate weighted sums of inequalities, whose weights correspond to the values of the dual variables. Proofs from Section~\ref{c-Nest} and Section~\ref{c-prox} correspond to such dual certificates.

Those statements are made more precise in the next sections. We begin by providing a few examples of LMIs that can be used for designing worst-case guarantees.

\paragraph{Preview: worst-case guarantees via LMIs.} Perhaps the most basic LMI that can be presented for obtaining worst-case guarantees concerns gradient descent and its convergence in terms of distance to an optimal point. We present it for simplicity, as the corresponding LMI only involves very few variables. This LMI  has also relatively simple solutions. As our target here is to present the approach, we let finding their solutions as exercises. We present the LMIs in their most \emph{raw} forms, even without a few direct simplifications.

Note that those LMIs always involve $n(n-1)$ ``dual'' variables (the precise meaning of \emph{dual} becomes clear in the sequel), where $n$ is the number of points at which the type of guarantee under consideration requires using or specifying a function or gradient evaluation (either in the algorithm or for computing the value of the guarantee). In the following example, we need two dual variables because the guarantee only requires using two gradients of $f$, namely $\nabla f(x_k)$ (for expressing a gradient step $x_{k+1}=x_k-\gamma_k \nabla f(x_k)$) and $\nabla f(x_\star)$ (for expressing optimality of $x_\star$ as $\nabla f(x_\star)=0$).

\begin{theorem}\label{thm:dist_gd}
Let $\tau\geq 0$ and $\gamma_k\in\mathbb{R}$. The inequality \begin{equation}\label{eq:dist_gd}
    \|x_{k+1}-x_\star\|^2_2\leq {\tau} \|x_k-x_\star\|^2_2
\end{equation} holds for all $d\in\mathbb{N}$, all $f\in\mathcal{F}_{\mu,L}(\mathbb{R}^d)$, all $x_k,x_{k+1},x_\star\in\mathbb{R}^d$ (such that $x_{k+1}=x_k-\gamma_k \nabla f (x_k)$ and $\nabla f(x_\star)=0$), if and only if
\begin{equation}\label{eq:gd_dist_LMI}
    \exists {\lambda_1,\lambda_2}\geq 0:\left\{\begin{array}{ll}
    \lambda_1=\lambda_2\\
    0\preceq \begin{bmatrix}
								\tau -1+\frac{\mu  L (\lambda_1+\lambda_2)}{2 (L-\mu )} & \gamma_k-\frac{L \lambda_1+\mu  \lambda_2}{2 (L- \mu) } \\
 \gamma_k -\frac{L \lambda_1+\mu  \lambda_2}{2(L-\mu) } & -\gamma_k^2+\frac{\lambda_1+\lambda_2}{2 (L-\mu )}  \\
								\end{bmatrix}.
    \end{array}\right.
\end{equation}
\end{theorem}
We emphasize that the message underlying Theorem~\ref{thm:dist_gd} is that verifying a worst-case convergence guarantee of the form~\eqref{eq:dist_gd} boils down to verifying the feasibility of a certain convex problem. It is relatively straightforward to convert a feasible point of~\eqref{eq:gd_dist_LMI} to  a proof that only consists of a weighted linear combination of inequalities, see, e.g.,~\citep[Theorem 3.1]{taylor2018exact}. The corresponding weights are the values of the multipliers (that is, in Theorem~\ref{thm:dist_gd}, the weights are $\lambda_1$ and $\lambda_2$) as showcased in Section~\ref{c-Nest} and Section~\ref{c-prox}.

As we see in Section~\ref{s:gd_obtain_lmi}, changing the Lyapunov, or potential, function to be verified also changes the LMI to be solved. The desired LMI can be obtained following a principled approach presented in the sequel. In particular, the following result is slightly more complicated and corresponds to verifying the potential provided by Theorem~\ref{thm:pot_GM}. One should note that those LMIs can be solved numerically, providing nice guides for choosing appropriate analytical weights. Symbolic computations and computer algebra software might also help.

The following LMI relies on $6$ \emph{dual variables} $\lambda_1,\hdots,\lambda_6$ as it involves gradients and/or function values of $f(\cdot)$ at three points: $x_k$, $x_{k+1}$, and $x_\star$, thereby fixing $n=3$ and hence $n(n-1)=6$ dual variables.
\begin{theorem}\label{thm:lyap_gd}
Let $A_{k+1},A_k\geq 0$ and $\gamma_k\in\mathbb{R}$. The inequality \[A_{k+1}(f(x_{k+1})-f_\star)+\tfrac{L}{2}\|x_{k+1}-x_\star\|_2^2\leq A_k(f(x_k)-f_\star)+\tfrac{L}{2}\|x_k-x_\star\|_2^2\] holds for all $d\in\mathbb{N}$, all $f\in\mathcal{F}_{L}(\mathbb{R}^d)$, all $x_k,x_{k+1},x_\star\in\mathbb{R}^d$ (such that $x_{k+1}=x_k-\gamma_k \nabla f (x_k)$ and $\nabla f(x_\star)=0$) if and only if
\begin{equation*}
\begin{aligned}
    \exists \lambda_1,&\lambda_2,\hdots,\lambda_6\geq 0:\\
    &\left\{\begin{array}{l}
     0=A_k+\lambda_1+\lambda_2-\lambda_4-\lambda_6 \\
     0=-A_{k+1}-\lambda_2+\lambda_3+\lambda_4-\lambda_5 \\
    0\preceq \begin{bmatrix}
	0 & \star & \star \\
 \frac{1}{2} (\gamma_k L-\lambda_1) & \frac{\lambda_1+\lambda_2+\lambda_4+\lambda_6-\gamma_k^2 L^2-2 \gamma_k L \lambda_2}{2 L} &\star \\
 -\frac{\lambda_3}{2} & \frac{1}{2} \left(\gamma_k (\lambda_3+\lambda_4)-\frac{\lambda_2+\lambda_4}{L}\right) & \frac{\lambda_2+\lambda_3+\lambda_4+\lambda_5}{2 L} \\
 \end{bmatrix},
    \end{array}\right.
\end{aligned}
\end{equation*}
(where $\star$'s denote symmetric elements in the matrix).
\end{theorem}
\begin{remark} The LMIs of this section are put in their ``raw'' forms, for simplicity of the presentation (which does not focus on solving those LMIs analytically. Of course, a few simplifications are relatively direct: for instance, any feasible point will have $\lambda_1=\gamma_kL$ and $\lambda_3=0$, as the corresponding matrix could not be positive semidefinite otherwise.
\end{remark}
As we discuss in the sequel (see Remark~\ref{rem:weaker}), it is also relatively straightforward to obtain weaker versions of those LMIs which are then only sufficient for obtaining valid worst-case guarantees. Those simplified LMIs might be simpler to solve analytically, and might therefore be advantageous in certain contexts. Brief discussions and pointers for this topic are provided in Remark~\ref{rem:weaker} and Section~\ref{s:ref}.

A strongly convex version of Theorem~\ref{thm:lyap_gd} is provided in Theorem~\ref{thm:lyap_gd_strcvx}. It is slightly more algrebaic in its vanilla form, but allows recovering the results of Theorem~\ref{thm:pot_GM_strcvx} as a feasible point. Analyses of accelerated methods can be obtained in a similar way, as illustrated by the following LMI. The latter uses on $12$ \emph{dual variables} $\lambda_1,\hdots,\lambda_{12}$, as it relies on evaluating gradients and/or function values of $f(\cdot)$ at four points: 
$y_k$, $x_k$, $x_{k+1}$, and $x_\star$, so $n=4$ and hence $n(n-1)=12$. Although this LMI might appear as a bit of a brutal approach to worst-case analysis, one might observe that many of elements of the LMI can be set to zero due to the structure of the problem.
\begin{theorem}\label{thm:lyap_fgm}
Let $A_{k+1},A_k\geq 0$ and $\alpha_k,\gamma_k,\tau_k\in\mathbb{R}$, and consider the iteration
\begin{equation}\label{eq:acc_app}
\begin{aligned}
y_{k}&=x_k+\tau_k (z_k-x_k)\\
x_{k+1}&=y_{k}-\alpha_k \nabla f(y_{k})\\
z_{k+1}&=z_k-\gamma_k \nabla f(y_{k}).
\end{aligned}
\end{equation}
The inequality \[A_{k+1}(f(x_{k+1})-f_\star)+\tfrac{L}{2}\|z_{k+1}-x_\star\|_2^2\leq A_k(f(x_k)-f_\star)+\tfrac{L}{2}\|z_k-x_\star\|_2^2\] holds for all $d\in\mathbb{N}$, all $f\in\mathcal{F}_{L}(\mathbb{R}^d)$, and all $x_k,x_{k+1},z_k,z_{k+1},x_\star\in\mathbb{R}^d$ (such that $x_{k+1},z_{k+1}$ are generated by~\eqref{eq:acc_app} and $\nabla f(x_\star)=0$) if and only if
\begin{equation*}
\begin{aligned}
    \exists \lambda_1,&\lambda_2,\hdots,\lambda_{12}\geq 0:\\
    &\left\{\begin{array}{l}
\begin{aligned}
    0&=A_k+\lambda_{1}+\lambda_{2}-\lambda_{4}-\lambda_{6}-\lambda_{8}+\lambda_{11}\\
    0&=-A_{k+1}-\lambda_{2}+\lambda_{3}+\lambda_{4}-\lambda_{5}-\lambda_{9}+\lambda_{12}\\
    0&=\lambda_{7}+\lambda_{8}+\lambda_{9}-\lambda_{10}-\lambda_{11}-\lambda_{12}\\
    0&\preceq \begin{bmatrix}
    0 & 0 & S_{1,3} & S_{1,4} & S_{1,5}\\
    0 & 0 & S_{2,3} & S_{2,4} & S_{2,5}\\
    S_{1,3} & S_{2,3} & S_{3,3} & S_{3,4} & S_{3,5}\\
    S_{1,4} & S_{2,4} & S_{3,4} & S_{4,4} & S_{4,5}\\
    S_{1,5} & S_{2,5} & S_{3,5} & S_{4,5} & S_{5,5}
 \end{bmatrix},
\end{aligned}
    \end{array}\right.
\end{aligned}
\end{equation*}
with
\begin{equation*}
    \begin{aligned}
    S_{1,3}&=\tfrac{1}{2} (\lambda_{7} (\tau_k-1)+\lambda_{8} \tau_k),\\S_{1,4}&=-\tfrac{1}{2} (\lambda_{1}+\tau_k (\lambda_{2}+\lambda_{11})),\,\,
    S_{1,5}=\tfrac{1}{2} (\lambda_{3} (\tau_k-1)+\lambda_{4} \tau_k),\\
    S_{2,3}&=\tfrac{1}{2} (\gamma_k L-\tau_k (\lambda_{7}+\lambda_{8})),\\
    S_{2,4}&=\tfrac{1}{2} \tau_k (\lambda_{2}+\lambda_{11}),\,\,
    S_{2,5}=-\tfrac{1}{2} \tau_k (\lambda_{3}+\lambda_{4}),\\
    S_{3,3}&=\frac{\lambda_{7}+\lambda_{8}+\lambda_{9}+\lambda_{10}+\lambda_{11}+\lambda_{12}-\gamma_k^2 L^2-2 \alpha_k L \lambda_{9}}{2 L},\\
    S_{3,4}&=-\frac{\alpha_k L \lambda_{2}+\lambda_{8}+\lambda_{11}}{2 L},\,\,
    S_{3,5}=\tfrac{1}{2} \left(\alpha_k (\lambda_{3}+\lambda_{4}+\lambda_{12})-\frac{\lambda_{9}+\lambda_{12}}{L}\right),\\
    S_{4,4}&=\frac{\lambda_{1}+\lambda_{2}+\lambda_{4}+\lambda_{6}+\lambda_{8}+\lambda_{11}}{2 L},\,\,
    S_{4,5}=-\frac{\lambda_{2}+\lambda_{4}}{2 L},\\
    S_{5,5}&=\frac{\lambda_{2}+\lambda_{3}+\lambda_{4}+\lambda_{5}+\lambda_{9}+\lambda_{12}}{2 L}.
    \end{aligned}
\end{equation*}
\end{theorem}

\section{Analysis of Gradient Descent via Linear Matrix Inequalities}\label{s:gd_obtain_lmi}

In this section, we detail the approach to obtain LMIs such as those of Theorem~\ref{thm:dist_gd}, Theorem~\ref{thm:lyap_gd} and Theorem~\ref{thm:lyap_fgm}. We provide full details for gradient descent. The same technique is presented in a more expeditious way for its accelerated versions afterwards.

\subsection{Linear Convergence of Gradient Descent}
We consider gradient descent for minimizing smooth strongly convex functions. For exposition purposes, we investigate a type of one-iteration worst-case convergence guarantee in terms of the distance to the optimum (see Theorem~\ref{thm:dist_gd}) for gradient descent, of the form:
\begin{equation}\label{eq:WC_grad_dist}
    \|x_{k+1}-x_\star\|_2^2 \leq \tau \|x_k-x_\star\|_2^2
\end{equation}
which are valid for all $d\in\mathbb{N}$, $x_k,x_{k+1},x_\star\in\mathbb{R}^d$ and all $f\in\mathcal{F}_{\mu,L}(\mathbb{R}^d)$ ($L$-smooth $\mu$-strongly convex function) when $x_{k+1}=x_k-\gamma_k \nabla f(x_k)$ (gradient descent) and $\nabla f(x_\star)=0$ ($x_\star$ is optimal for $f$). In this context, we denote by $\tau_\star$ (we omit the dependence on $\gamma_k$, $\mu$, and $L$ for convenience) the smallest value $\tau$ for which~\eqref{eq:WC_grad_dist} is valid. By definition, this value can be formulated as the solution to an optimization problem looking for worst-case scenarios:
\begin{equation}\label{eq:def_pep}
    \begin{aligned}
    \tau_\star \defeq \max_{\substack{d,f\\x_k,x_{k+1},x_\star}} &\frac{\|x_{k+1}-x_\star\|_2^2}{\|x_k-x_\star\|_2^2}\\ \text{s.t.}\,\,& d\in\mathbb{N},f\in\mathcal{F}_{\mu,L}(\mathbb{R}^d)\\&x_k,x_{k+1},x_\star\in\mathbb{R}^d\\
    &x_{k+1}=x_k-\gamma_k \nabla f(x_k)\\
    &\nabla f(x_\star)=0.
    \end{aligned}
\end{equation}
As it is, this problem does not look quite practical. However, it actually admits an equivalent formulation as a linear semidefinite program. As a first step for reaching this formulation, the previous problem can be formulated in an equivalent \emph{sampled} manner. That is, we sample $f$ at the points where the first-order information is explicitly used: 
\begin{equation}
    \begin{aligned}
    \tau_\star = \max_{\substack{d\\f_k,f_\star\\g_k,,g_\star\\x_k,x_{k+1},x_\star}} &\frac{\|x_{k+1}-x_\star\|_2^2}{\|x_k-x_\star\|_2^2}\\ \text{s.t.}\,\,& d\in\mathbb{N},f_k,f_\star\in\mathbb{R}\\&x_k,x_{k+1},x_\star,g_k,g_\star\in\mathbb{R}^d\\
    &\exists f\in\mathcal{F}_{\mu,L}:\,\left\{\begin{array}{l}
     f_k=f(x_k) \text{ and } g_k=\nabla f(x_k)  \\
     f_\star=f(x_\star) \text{ and } g_\star=\nabla f(x_\star)
    \end{array}\right.\\
    &g_\star=0\\
    &x_{k+1}=x_k-\gamma_k g_k,
    \end{aligned}
\end{equation}
and $f$ is now represented in terms of its samples at $x_\star$ and $x_k$.

A second stage in this reformulation consists of replacing the existence of a certain $f\in\mathcal{F}_{\mu,L}$ interpolating (or extending) the samples by an equivalent explicit condition provided by the following theorem.
\begin{theorem}[$\mathcal{F}_{\mu,L}$-interpolation, Theorem 4 in~\citep{taylor2017smooth}]\label{thm:interp} Let $L>\mu\geq0$, $I$ be an index set and $S=\{(x_i,g_i,f_i)\}_{i\in I}\subseteq \mathbb{R}^d\times\mathbb{R}^d\times \mathbb{R}$ be a set of triplets. There exists $f\in\mathcal{F}_{\mu,L}$ satisfying $f(x_i)=f_i$ and $g_i\in\partial f(x_i)$ for all $i\in I$ if and only if
\begin{equation}\label{eq:interp_thm}
\begin{aligned}
f_i\geq f_j&+\langle g_j;x_i-x_j\rangle+\frac{1}{2L}\|g_i-g_j\|_2^2\\&+\frac{\mu}{2(1-\mu/L)}\|x_i-x_j-\tfrac{1}{L}(g_i-g_j)\|_2^2
\end{aligned}
\end{equation}
holds for all $i,j\in I$.
\end{theorem}
Theorem~\ref{thm:interp} conveniently allows replacing the existence constraint by a set of quadratic inequalities, reaching:
\begin{equation}
    \begin{aligned}
    \tau_\star = \max_{\substack{d\\f_k,f_\star\\g_k,x_k,x_\star}} &\frac{\|x_k-\gamma_k g_k-x_\star\|_2^2}{\|x_k-x_\star\|_2^2}\\ \text{s.t.}\,\,& d\in\mathbb{N},f_k,f_\star\in\mathbb{R}\\&x_k,x_\star,g_k\in\mathbb{R}^d\\
    &f_\star\geq f_k+\langle{g_k};{x_\star-x_k}\rangle+\frac{1}{2L}\|{g_k}\|_2^2\\& \quad\quad\quad\quad+\frac{\mu}{2(1-\mu/L)}\|{x_k-\frac{1}{L}g_k-x_\star}\|_2^2\\
	&f_k\geq f_\star+\frac{1}{2L}\|{g_k}\|_2^2\\& \quad\quad\quad\quad+\frac{\mu}{2(1-\mu/L)}\|{x_k-\frac{1}{L}g_k-x_\star}\|_2^2,
    \end{aligned}
\end{equation}
where we also substituted $x_{k+1}$ and $g_\star$ by their respective expressions. Finally, we arrive to a first (convex) semidefinite reformulation of the problem via new variables: $G\succeq 0$ and $F$ defined as
			\begin{align*}
			G \defeq \begin{bmatrix}
			\|x_k-x_\star\|_2^2 & \langle g_k,x_k-x_\star\rangle\\
			\langle g_k, x_k-x_\star\rangle & \| g_k\|_2^2
			\end{bmatrix},\quad 	F \defeq 			f_k-f_\star.
			\end{align*}
The problem turns out to be linear in $G$ and $F$:
\begin{equation}\label{eq:SDP_denom}
    \begin{aligned}
    \tau_\star = \max_{\substack{G,F}}\,\, &\frac{G_{1,1}+\gamma_k^2 G_{2,2}-2\gamma_k G_{1,2}}{G_{1,1}}\\ \text{s.t.}\,\,& F\in\mathbb{R},G\in\mathbb{S}^2\\
	&G\succeq 0\\
    &F + \tfrac{L\mu}{2(L-\mu)} G_{1,1}+\tfrac{1}{2(L-\mu)}G_{2,2}-\tfrac{L}{L-\mu}G_{1,2}\leq 0\\
	&-F + \tfrac{L\mu}{2(L-\mu)} G_{1,1}+\tfrac{1}{2(L-\mu)}G_{2,2}-\tfrac{\mu}{L-\mu}G_{1,2}\leq 0.
    \end{aligned}
\end{equation}
Finally, a simple homogeneity argument (for any feasible $(G,F)$ to~\eqref{eq:SDP_denom}, the pair $(\tilde{G},\tilde{F})\defeq (G/G_{1,1},F/G_{1,1})$ is also feasible with the same objective value, with $\tilde{G}_{1,1}=1$ so we can assume without loss of generality that $G_{1,1}=1$ without changing the optimal value of the problem---note that it is relatively straightforward to establish that the optimal solution satisfies $G_{1,1}\neq 0$) allows arriving to the equivalent:
\begin{equation}
    \begin{aligned}
    \tau_\star = \max_{\substack{G,F}}\,\, &G_{1,1}+\gamma_k^2 G_{2,2}-2\gamma_k G_{1,2}\\ \text{s.t.}\,\,& F\in\mathbb{R},G\in\mathbb{S}^2\\
	&G\succeq 0\\
    &F + \tfrac{L\mu}{2(L-\mu)} G_{1,1}+\tfrac{1}{2(L-\mu)}G_{2,2}-\tfrac{L}{L-\mu}G_{1,2}\leq 0\\
	&-F + \tfrac{L\mu}{2(L-\mu)} G_{1,1}+\tfrac{1}{2(L-\mu)}G_{2,2}-\tfrac{\mu}{L-\mu}G_{1,2}\leq 0\\
	&G_{1,1}=1.
    \end{aligned}
\end{equation}
For arriving to the desired LMI, it remains to dualize the problem. That is, we perform the following primal-dual associations:
\begin{equation*}
\begin{aligned}
    &F + \tfrac{L\mu}{2(L-\mu)} G_{1,1}+\tfrac{1}{2(L-\mu)}G_{2,2}-\tfrac{L}{L-\mu}G_{1,2}\leq 0 &&:\lambda_1,\\
	&-F + \tfrac{L\mu}{2(L-\mu)} G_{1,1}+\tfrac{1}{2(L-\mu)}G_{2,2}-\tfrac{\mu}{L-\mu}G_{1,2}\leq 0 &&:\lambda_2,\\
	&G_{1,1}=1&&:\tau.
\end{aligned}
\end{equation*}
Standard Lagrangian duality allows arriving to
\begin{equation}\label{eq:dual_GD_dist}
\begin{aligned}
\tau_\star=\min_{\lambda_1,\lambda_2,\tau\geq 0} &\tau\\
\text{s.t. }&\lambda_1=\lambda_2\\
    &0\preceq \begin{bmatrix}
							\tau-1+\frac{\mu L (\lambda_1+\lambda_2)}{2(L-\mu )}  & {\gamma_k}-\frac{\lambda_1 L+\lambda_2 \mu}{2 (L-\mu )} \\
								{\gamma_k}-\frac{\lambda_1 L+\lambda_2 \mu}{2 (L-\mu )} & -{\gamma_k^2}+\frac{\lambda_1+\lambda_2}{2(L-\mu)} \\
			    \end{bmatrix},
\end{aligned}
\end{equation}
where we used the fact there is no duality gap, as one can show the existence of a Slater point~\citep{boyd2004convexopt}. One such Slater point can be obtained by applying gradient descent on the function $f(x)=\tfrac{1}{2}x^\top\mathrm{diag}(L,\mu)x$ (i.e., $d=2$) with $x_k=(1,1)$. A formal statement is provided in~\citep[Theorem 6]{taylor2017smooth}. 

Theorem~\ref{thm:dist_gd} is now a direct consequence of the dual reformulation~\eqref{eq:dual_GD_dist}, as provided by the following proof.

\begin{proof}[Proof of Theorem~\ref{thm:dist_gd}]
(Sufficiency, $\Leftarrow$) If there exists a feasible point $(\tau,\lambda_1,\lambda_2)$ for~\eqref{eq:gd_dist_LMI}, weak duality implies that it is an upper bound on $\tau_\star$ by construction.

(Necessity, $\Rightarrow$) For any $\tau$ such that there exists no $\lambda_1,\lambda_2\geq 0$ for which $(\tau,\lambda_1,\lambda_2)$ is feasible for~\eqref{eq:gd_dist_LMI}, it follows that $\tau\leq \tau_\star$, and strong duality implies that there exists a problem instance ($f\in\mathcal{F}_{\mu,L}$, $d\in\mathbb{N}$, and $x_k\in\mathbb{R}^d$) on which $\|x_{k+1}-x_\star\|_2^2=\tau_\star\|x_k-x_\star\|_2^2\geq \tau \|x_k-x_\star\|_2^2$.
\end{proof}
\begin{remark}
Following similar lines as those of this section, one can verify other types of inequalities, beyond~\eqref{eq:dist_gd}, simply by changing the objective in~\eqref{eq:def_pep}. This allows obtaining the statement from Theorem~\ref{thm:lyap_gd} and Theorem~\ref{thm:lyap_fgm}.
\end{remark}
\begin{remark}
Finding analytical solutions to such LMIs (parametrized by the algorithm and problem parameters) might be challenging. For gradient descent, the solution is provided in e.g.,~\citep[Section 4.4]{lessard2016analysis} and~\citep[Theorem 3.1]{taylor2018exact}. For more complicated cases, one can rely on numerical inspiration for finding analytical solutions (or upper bounds on it).
\end{remark}
\begin{remark}\label{rem:weaker}
It is possible to obtain ``weaker'' LMIs based on other sets of inequalities (which are necessary but not sufficient for interpolation). Those LMIs are then only sufficient for finding worst-case guarantees. Those alternate LMIs might enjoy simpler analytical solutions, but this comes at the cost of loosing a priori tightness guarantees. This is in general not a problem if the worst-case guarantee is satisfying, but the subtle consequence is that those LMIs might then fail to provide a satisfying guarantee even when there exists one.
\end{remark}

\subsection{Potential Function for Gradient Descent} 

For formulating the LMI for verifying potential functions as those of Theorem~\ref{thm:pot_GM} and Theorem~\ref{thm:pot_GM_strcvx}, one essentially has to follow the same steps as in the previous section. The strongly convex version is a bit heavy and is provided below. In short, verifying that 
\[ \phi_k\defeq A_k(f(x_k)-f_\star)+\tfrac{L+\mu A_k}{2}\|x_k-x_\star\|^2_2\]
is a potential function, that is, $\phi_{k+1}\leq\phi_k$ (for all $f\in\mathcal{F}_{\mu,L}$, $d\in\mathbb{N}$, and $x_k\in\mathbb{R}^d$), amount to verify that
\begin{equation*}
\begin{aligned}
0\geq \max\left\{\phi_{k+1}-\phi_k:\right.&\, d\in\mathbb{N},\,f\in\mathcal{F}_{\mu,L},\,x_k,x_{k+1},x_\star\in\mathbb{R}^d,\\& \left.x_{k+1}=x_k-\gamma_k\nabla f(x_k),\text{ and } \nabla f(x_\star)=0\right\},
\end{aligned}
\end{equation*}
where the maximum is taken over $d$, $f$, $x_k$, $x_{k+1}$ and $x_\star$. This problem can be reformulated as in Section~\ref{s:gd_obtain_lmi} using the same technique with more samples. More precisely, this formulation requires sampling the function $f$ at three points (instead of two): $x_\star$, $x_k$, and $x_{k+1}$, and hence $6$ dual variables are required (because $6$ inequalities of the form~\eqref{eq:interp_thm} are used for describing the sampled version of the function $f$). The formal statement is provided by the following theorem, without a proof.

\begin{theorem}\label{thm:lyap_gd_strcvx}
Let $A_{k+1},A_k\geq 0$ and $\gamma_k\in\mathbb{R}$. The inequality 
\begin{equation*}
\begin{aligned}
A_{k+1}(f(x_{k+1})-f_\star)&+\tfrac{L+\mu A_{k+1}}{2}\|x_{k+1}-x_\star\|_2^2\\&\leq A_k(f(x_k)-f_\star)+\tfrac{L+\mu A_k}{2}\|x_k-x_\star\|_2^2
\end{aligned}
\end{equation*}
holds for all $d\in\mathbb{N}$, all $f\in\mathcal{F}_{\mu,L}(\mathbb{R}^d)$, all $x_k,x_{k+1},x_\star\in\mathbb{R}^d$ (such that $x_{k+1}=x_k-\gamma_k\nabla f (x_k    )$ and $\nabla f(x_\star)=0$) if and only if
\begin{equation*}
    \exists \lambda_1,\lambda_2,\hdots,\lambda_6\geq 0:\,\,\left\{\begin{array}{l}
     0=A_k+\lambda_1+\lambda_2-\lambda_4-\lambda_6\\
     0=-A_{k+1}-\lambda_2+\lambda_3+\lambda_4-\lambda_5\\
    0\preceq \begin{bmatrix}
	 S_{1,1} & S_{1,2} & S_{1,3} \\
 S_{1,2} & S_{2,2} & S_{2,3} \\
 S_{1,3} & S_{2,3} & S_{3,3} \\
 \end{bmatrix},
    \end{array}\right.
\end{equation*}
with 
\begin{equation*}
\begin{aligned}
S_{1,1}&=\tfrac{1}{2} \mu  \left(A_{k}-A_{k+1}+\tfrac{L (\lambda_1+\lambda_3+\lambda_5+\lambda_6)}{L-\mu }\right)\\
S_{1,2}&=-\tfrac{\gamma_k (\mu  A_{k+1} (\mu -L)+L (\mu  (\lambda_3+\lambda_5+1)-L))+\lambda_6 \mu +\lambda_1 L}{2 (L-\mu )}\\
S_{1,3}&=-\tfrac{\lambda_5 \mu +\lambda_3 L}{2 (L- \mu) }\\
S_{2,2}&=\tfrac{\gamma_k^2 (\mu  A_{k+1} (\mu -L)+L (\mu  (\lambda_2+\lambda_3+\lambda_4+\lambda_5+1)-L))-2 \gamma_k (\lambda_4 \mu +\lambda_2 L)+\lambda_1+\lambda_2+\lambda_4+\lambda_6}{2 (L-\mu )}\\
S_{2,3}&=\tfrac{\gamma_k (\mu  (\lambda_2+\lambda_5)+L (\lambda_3+\lambda_4))-\lambda_2-\lambda_4}{2 (L-\mu )}\\
S_{3,3}&=\tfrac{\lambda_2+\lambda_3+\lambda_4+\lambda_5}{2 (L- \mu) }.
\end{aligned}
\end{equation*}
\end{theorem}

Note again that a notebook is provided in Section~\ref{s:ref} for obtaining and verifying this LMI formulation via symbolic computations.
 \section{Accelerated Gradient Descent via Linear Matrix Inequalities}
We provide the main ideas for formulating the LMI for verifying potential functions as those of Theorem~\ref{thm:pot_GM_strcvx} and Theorem~\ref{thm:Nest_first_strcvx}. In short, verifying that
\[ \phi_k\defeq A_k(f(x_k)-f_\star)+\tfrac{L+\mu A_k}{2}\|z_k-x_\star\|_2^2\]
is a potential function, that is, $\phi_{k+1}\leq\phi_k$ (for all $f\in\mathcal{F}_{\mu,L}$, $d\in\mathbb{N}$, and $x_k,z_k,x_\star\in\mathbb{R}^d$, $\nabla f(x_\star)=0$), amounts to verify that
\begin{equation*}
\begin{aligned}
0\geq \max\left\{\phi_{k+1}-\phi_k:\right.&\, d\in\mathbb{N},\,f\in\mathcal{F}_{\mu,L},\,z_k,x_k,x_\star\in\mathbb{R}^d,\,\nabla f(x_\star)=0,\\& \left.\text{and }y_k,x_{k+1},z_{k+1}\in\mathbb{R}^d \text{ generated by~\eqref{eq:FGM_1_strcvx}}\right\},
\end{aligned}
\end{equation*}
where the maximum is taken over $d$, $f$, the iterates, as well as $x_\star$. This problem can be cast as a SDP using the same ideas as in Section~\ref{s:gd_obtain_lmi} with more samples, again. More precisely, this formulation requires sampling the function~$f$ at four points: $x_\star$, $x_k$, $x_{k+1}$, and $y_k$. The case $\mu=0$ is covered by Theorem~\ref{thm:lyap_fgm}.

\section{Notes and References}\label{s:ref}

\paragraph{General frameworks.} The whole idea of using semidefinite programming for analyzing first-order methods dates back to~\citep{Dror14} (more details and examples in~\citep{drori2014contributions,drori2016optimal}). The principled approach to worst-case analysis using performance estimation problems with interpolation/extension arguments was proposed in in~\citep{taylor2017smooth}, and generalized to more problem setups in~\citep{taylor2017exact}. The integral quadratic approach to first-order methods was proposed in~\citep{lessard2016analysis}, specifically for studying linearly converging methods (focus on strong convexity and related notions). Those two related methodologies were then further extended and linked in different setup~\citep{hu2017dissipativity,hu2017unified,taylor2018lyapunov,fazlyab2018analysis,taylor19bach,lieder2020convergence,aybat2020robust,hu2021analysis,Ryu20,dragomir2019optimal}. Among those developments, some works performed analyses via ``weaker'' LMIs, based on other sets of inequalities which are necessary but not sufficient for interpolation; see, e.g.,~\citep{ryupark2021}. The advantage of this approach is that it is often simpler to obtain analytical solutions to some of those LMIs, at the cost of loosing tightness guarantees (which might not be a problem when the guarantee is satisfying). This is in general the case for IQC-based works. In those cases, non-tightness is usually coupled with the search for a Lyapunov function. In general, it is possible to simultaneously look for a tight guarantee and a Lyapunov/potential function, see e.g.,~\citep{taylor2018lyapunov,taylor19bach}. A simplified approach to performance estimation problems was implemented in the performance estimation toolbox~\citep[PESTO]{pesto2017}.

\paragraph{Designing methods using semidefinite programming.} The optimized gradient method (OGM) was apparently the first method obtained by optimizing its worst-case using SDPs/LMIs. It was obtained as a solution to a convex optimization problem by~\citet{Dror14}, which was later solved analytically by~\citet{kim2016optimized}. The same method was obtained through an analogy with the conjugate gradient method~\citep{drori2019efficient}, which might serve as a strategy for designing method in various setups. Optimized methods can be developed for other criteria and setups as well. As an example, optimized methods for gradient norms $\|\nabla f(x_N)\|_2^2$ are studied in~\citet{kim2018optimizing,kim2018generalizing}, in the smooth convex setting. See also Section~\ref{s:ogm_sc} and Section~\ref{s:tmm}; in particular, the \emph{Triple Momentum Method} (TMM)~\citep{van2017fastest} was designed as a time-independent optimized gradient method, through Lyapunov arguments (and IQCs). See also~\citep{lessard2020direct,zhou2020boosting,gramlich2020convex,drori2021exact} for different ways of recovering the TMM. Optimized methods were also developed in other setups, such as fixed-point iteration~\citep{lieder2020convergence} and monotone inclusions~\citep{kim2019accelerated} (which turned out to be a particular case of~\citep{lieder2020convergence}).

\paragraph{Specific methods.} The SDP/LMI approaches were taken further for studying first-order methods in a few different contexts. It was originally used for studying gradient-type methods (see, e.g.,~\citep{Dror14,drori2014contributions,lessard2016analysis,taylor2017smooth}) and accelerated/fast gradient-type methods (see, e.g.,~\citep{Dror14,drori2014contributions,lessard2016analysis,taylor2017smooth,taylor2017exact,hu2017dissipativity,van2017fastest,cyrus2018robust,safavi2018explicit,aybat2020robust}) for convex minimization. It was used later for analyzing, among others, nonsmooth setups~\citep{drori2016optimal,drori2019efficient}, stochastic~\citep{hu2017unified,hu2018dissipativity,taylor19bach,hu2021analysis}, coordinate-descent~\citep{shi2017better,taylor19bach}, nonconvex setups~\citep{abbaszadehpeivasti2021exact,abbaszadehpeivasti2021rate}, proximal methods~\citep{taylor2017exact,kim2018another,kim2018optimizing,barre2020principled}, splitting methods~\citep{ryu2020finding,ryu2020operator,taylor2018exact}, monotone inclusions and variational inequalities~\citep{ryu2020operator,gu2019optimal,gu2020tight,zhang2021unified}, fixed-point iterations~\citep{lieder2020convergence}, and distributed/decentralized optimization~\citep{sundararajan2020analysis,colla2021automated}.

\paragraph{Obtaining and solving the LMIs.} For solving the LMIs, standard numerical semidefinite optimization packages can be used, see, e.g.,~\citep{Yalmip,Sedumi,Mosek,SdpT3}. For obtaining and verifying analytical solutions, symbolic computing might also be a great asset. For the purpose of reproducibility, we provide notebooks for obtaining the LMI formulations of this section symbolically, and for solving them numerically, at \togglecodeurl.

\begin{acknowledgements}
The authors would like to warmly thank Rapha\"el Berthier, Mathieu Barr\'e, Aymeric Dieuleveut, Fabian Pedregosa and Baptiste Goujaud for comments on early versions of this manuscript; for spotting a few typos; and for discussions and developments related to Section~\ref{c-Cheb}, Section~\ref{c-Nest}, and Section~\ref{c-prox}. We are also greatly indebted to Lena\"ic Chizat, Laurent Condat, Jelena Diakonikolas, Alexander Gasnikov, 
Shuvomoy Das Gupta, Pontus Giselsson, Crist\'obal Guzm\'an, Julien Mairal, and Ir\`ene Waldspurger for spotting a few typos and inconsistencies in the first version of the manuscript.

We further want to thank Francis Bach, S\'ebastien Bubeck, Radu-Alexandru Dragomir, Yoel Drori, Hadrien Hendrikx, Reza Babanezhad, Claude Brezinski, Pavel Dvurechensky, Herv\'e Le Ferrand, Georges Lan, Adam Ouorou, Michela Redivo-Zaglia, Simon Lacoste-Julien, Vincent Roulet, and Ernest Ryu for fruitful discussions and pointers, which largely simplified the writing and revision process of this manuscript.

AA is also extremely grateful to the French ministry of education and \'ecole Etienne Marcel for keeping school mostly open during the 2020-2021 pandemic. 

AA is at the D\'epartement d’informatique de l’ENS, \'Ecole normale sup\'erieure, UMR CNRS 8548, PSL Research University, 75005 Paris, France and INRIA. AA would like to acknowledge support from the ML and Optimisation joint research initiative with the funds AXA pour la Recherche and Kamet Ventures, a Google focused award, as well as funding from the French government under the management of the Agence Nationale de la Recherche as part of the  ``Investissements d’avenir'' program, reference ANR-19-P3IA-0001 (PRAIRIE 3IA Institute). AT is at INRIA and the D\'epartement d’informatique de l’ENS, \'Ecole normale sup\'erieure, CNRS, PSL Research University, 75005 Paris, France. AT acknowledges support from the European Research Council (ERC grant SEQUOIA 724063).
\end{acknowledgements}  

\backmatter

\printbibliography

\end{document}